\documentclass[openright,12pt,twoside,makeidx]{book}
\usepackage{a4}
\usepackage[T1]{fontenc}
\usepackage[all]{xy}
\usepackage{verbatim}
\usepackage{stmaryrd}
\usepackage{amsfonts}
\usepackage{amsmath,amssymb,amsthm}
\usepackage{enumerate}
\usepackage{xspace}
\usepackage{ulem}
\usepackage{euscript}
\usepackage{epsfig}
\usepackage{fancyhdr}
\usepackage{makeidx}
\usepackage{indentfirst}
\newtheoremstyle{note}
{20pt}
{20pt}
{\itshape}
{}
{\bfseries}
{.}
{0.5em}
{}
\theoremstyle{note}
\newtheorem{theo}{Th{\'e}or{\`e}me}[section]
\newtheorem{defi}[theo]{D{\'e}finition}

\newtheorem{lem}[theo]{Lemme}
\newtheorem{lemf}[theo]{Lemme fondamental}
\newtheorem{exem}[theo]{Exemple}
\newtheorem{theor}{Th{\'e}or{\`e}me}
\newtheorem{lemfo}[theor]{Lemme fondamental}

\newtheorem{rem}[theo]{Remarque}
\newtheorem{propri}[theo]{Propri{\'e}t{\'e}}
\newtheorem{fait}[theo]{Fait}
\newtheorem{rems}[theo]{Remarques}
\newtheorem{coro}[theo]{Corollaire}

\newtheorem{notat}[theo]{Notation}
\newtheorem{pro}[theo]{Proposition}
\newtheorem{propo}[theor]{Proposition}
\newtheorem{theosub}{Th{\'e}or{\`e}me}[subsection]
\newtheorem{defisub}[theosub]{D{\'e}finition}
\newtheorem{lemsub}[theosub]{Lemme}
\newtheorem{exemsub}[theosub]{Exemple}
\newtheorem{remsub}[theosub]{Remarque}
\newtheorem{faitsub}[theosub]{Fait}

\newtheorem{corosub}[theosub]{Corollaire}

\newtheorem{prosub}[theosub]{Proposition}
\input cyracc.def
\DeclareFontFamily{U}{russian}{}
\DeclareFontShape{U}{russian}{m}{n}
{ <5><6> wncyr5
<7><8><9> wncyr7
<10><10.95><12><14.4><17.28><20.74><24.88> wncyr10 }{}
\DeclareSymbolFont{Russian}{U}{russian}{m}{n}
\DeclareSymbolFontAlphabet{\mathcyr}{Russian}
\makeatletter
\let\@math@cyr\mathcyr
\makeindex
\begin{document}
\setcounter{secnumdepth}{3}
\pagenumbering{Roman}
\title{Descente de torseurs, gerbes et points rationnels}
\author{St{\'e}phane Zahnd}
\maketitle
\newpage
\thispagestyle{empty}
$\ $
\newpage
\thispagestyle{empty}
\tableofcontents
\vspace{80mm}
\begin{flushleft}
$\ $ \\
\normalsize{\textbf{Mots-cl{\'e}s.} Points rationnels, (bi-)torseurs, champs, gerbes, cohomologie non-ab{\'e}lienne, obstruction de Brauer-Manin.}
\end{flushleft}
\thispagestyle{empty}
\chapter*{Remerciements}
\addcontentsline{toc}{chapter}{Remerciements}
\markboth{REMERCIEMENTS}{REMERCIEMENTS}
\thispagestyle{empty}
Les mots sont faibles pour exprimer ma reconnaissance {\`a} Jean-Claude DOUAI pour tout ce qu'il m'a apport{\'e} durant l'{\'e}laboration de cette th{\`e}se. Sa passion et sa soif de d{\'e}couvrir, peut-{\^e}tre plus encore que sa disponibilit{\'e} de tous les instants et son impressionnante culture math{\'e}matique, ont transform{\'e} toutes nos discussions en v{\'e}ritables moments de bonheur. A ses c{\^o}t{\'e}s, la d{\'e}couverte des gerbes et de la cohomologie non-ab{\'e}lienne m'est apparue comme un fantastique voyage.\vspace{2mm}

J'ai bien s{\^u}r {\'e}t{\'e} tr{\`e}s honor{\'e} que Jean GIRAUD ait accept{\'e} la lourde t{\^a}che de rapporteur, et je tiens ici {\`a} le remercier d'avoir {\'e}t{\'e} aussi pr{\'e}cis et consciencieux dans ses observations. Le pr{\'e}sent texte porte la trace de son magnifique ouvrage sur la cohomologie non-ab{\'e}lienne, et de ses commentaires tr{\`e}s riches et d{\'e}taill{\'e}s sur mon travail.\vspace{2mm}

Le temps que David HARARI a eu la gentillesse de bien vouloir me consacrer, les indications et explications pr{\'e}cises qu'il m'a fournies sur les questions arithm{\'e}tiques en g{\'e}n{\'e}ral et sur l'obstruction de Brauer-Manin en particulier, m'ont permis d'{\'e}norm{\'e}ment progresser dans la compr{\'e}hension de ces probl{\`e}mes. Qu'il trouve ici l'expression de ma profonde gratitude.\vspace{2mm}

Je tiens {\`a} remercier chaleureusement Pierre D{\`E}BES, pour les conseils et encouragements qu'il m'a prodigu{\'e}s tout au long de ce travail. J'ai beaucoup appris sur les rev{\^e}tements (entre autres) {\`a} son contact, et nos discussions ont toujours {\'e}t{\'e}, sur de nombreux sujets, tr{\`e}s {\'e}clairantes.\vspace{2mm}

Je remercie Michel EMSALEM de m'avoir consacr{\'e} tant de temps, et d'avoir accept{\'e} de se pencher avec moi sur des probl{\`e}mes aussi effrayants (ah! le $gr$-champ des auto-{\'e}quivalences de la gerbe des $G$-torseurs\ldots). Manipuler les cat{\'e}gories fibr{\'e}es et les champs me para{\^i}t plus facile gr{\^a}ce {\`a} ses indications.\vspace{2mm}

Le cours de DEA de G{\'e}om{\'e}trie Alg{\'e}brique de Dimitri MARKUSHEVICH m'a donn{\'e} le go{\^u}t de cette discipline et du travail bien fait, et c'est pourquoi je souhaite ici lui dire avec reconnaissance: $\mathcyr{beskon}\acute{\mathcyr{e}}\mathcyr{ch}\mathcyr{no}$ $\mathcyr{spas}$$\acute{\mathcyr{i}}$$\mathcyr{b}$$\mathcyr{o}$.\vspace{2mm}

La th{\'e}orie des champs occupe une place importante dans ce travail. Je tiens donc {\`a} remercier chaleureusement Laurent MORET-BAILLY de m'avoir fait l'honneur d'accepter de faire partie de ce jury. Ses travaux dans ce domaine apportent {\`a} de nombreux endroits de la pr{\'e}sente th{\`e}se un {\'e}clairage particuli{\`e}rement int{\'e}ressant.\vspace{8mm}

Je tiens {\`a} remercier ma famille de son soutien inconditionnel durant ces trois ann{\'e}es. Evidemment, je m'adresse tout d'abord {\`a} Sophie, pour avoir support{\'e} sans sourciller mes envies de solitude, et pour avoir accept{\'e} que je me plonge de longues soir{\'e}es dans la lecture des \textquotedblleft{grands classiques}\textquotedblright.\vspace{2mm}

Je remercie Valerio VASSALLO, qui m'a enseign{\'e} les fondements de la G{\'e}om{\'e}trie Alg{\'e}brique, et dont les id{\'e}es m'ont permis d'aiguiser mon intuition et mon d{\'e}sir d'aller plus loin dans ce domaine.\vspace{2mm}

Enfin, je ne veux {\`a} aucun prix oublier mes amis: vos encouragements (en particulier ces derni{\`e}res semaines), et les conversations que nous avons eues ensemble, que ce soit devant un tableau noir, autour d'un caf{\'e} ou sur une aire d'autoroute quelque part en Allemagne, m'ont aid{\'e} {\`a} conserver intacte ma motivation. S{\'e}verine (la mise en page de cette th{\`e}se porte ta signature!), mes \textquotedblleft{compagnons autrichiens}\textquotedblright\ Yann et Salah, mes camarades de jeu en cohomologie {\'e}tale Ben et Diallo, c'est un plaisir de vous associer {\`a} ces remerciements.
\newpage
\thispagestyle{empty}
$\ $
\vspace{110mm}
\begin{verse}
	Mais, si vous en croyez tout le monde savant,\\L'esprit doit sur le corps prendre le pas devant,\\Et notre plus grand soin, notre premi{\`e}re instance,\\Doit {\^e}tre {\`a} le nourrir du suc de la science.
\end{verse}
\vspace{3mm}
$\ $
\hspace{50mm}
Moli{\`e}re, \textit{Les femmes savantes, (Acte II, Sc{\`e}ne VII).}
\newpage
\thispagestyle{empty}
\chapter*{Introduction}
\addcontentsline{toc}{chapter}{Introduction}
\markboth{INTRODUCTION}{INTRODUCTION}
\thispagestyle{empty}
La question qui motive ce travail est la suivante:
\begin{quotation}
	\textit{Soient $k$ un corps et $X$ un $k$-sch{\'e}ma. $X$ poss{\`e}de-t-il des points $k$-rationnels?}
\end{quotation}

Il semble {\'e}videmment sans espoir de r{\'e}pondre enti{\`e}rement {\`a} ce probl{\`e}me, puisque cela {\'e}quivaudrait {\`a} prouver d'un seul coup tous les {\'e}nonc{\'e}s du type \textquotedblleft{Th{\'e}or{\`e}me de Fermat}\textquotedblright\ possibles et imaginables. L'objectif plus raisonnable que nous nous fixons ici est le suivant:
\begin{quotation}
	\textit{Soient $k$ un corps de caract{\'e}ristique nulle et $X$ un $k$-sch{\'e}ma. D{\'e}terminer des obstructions cohomologiques {\'e}tales {\`a} l'existence de points $k$-rationnels sur $X$.}
\end{quotation}

Pour arriver {\`a} nos fins, nous utilisons les gerbes, introduites par Grothendieck et Giraud, mais fort peu utilis{\'e}es depuis (sauf peut-{\^e}tre par les Physiciens, qui se servent des gerbes ab{\'e}liennes en th{\'e}orie des cordes, \textit{cf.} \cite{Hi2}), ce qui nous semble {\^e}tre une flagrante injustice. En effet, les gerbes apparaissent naturellement dans des probl{\`e}mes aussi nombreux que vari{\'e}s, et dont nous donnons maintenant quelques exemples, en commen\c{c}ant par celui qui est au c{\oe}ur de cet expos{\'e}:
\vspace{1mm}
\newline

\uline{\textbf{Probl{\`e}me central}}\textbf{:} soient $k$ un corps de caract{\'e}ristique nulle, $X$ un $k$-sch{\'e}ma g{\'e}om{\'e}triquement connexe, et $G$ un $k$-groupe alg{\'e}brique lin{\'e}aire. Soient encore:
	\[\bar{P}\longrightarrow \bar{X}=X\otimes_{k}\bar{k}
\]
un $\bar{G}_{X}$-torseur, o{\`u}:\label{barGX}
	\[G_{X}=G\times_{\textup{Spec}\;k}X\ \ \textup{et}\ \ \bar{G}_{X}=G_{X}\times_{X}\bar{X}.
\]
et $\bar{k}$ d{\'e}signe une cl{\^o}ture alg{\'e}brique de $k$ fix{\'e}e {\`a} l'avance. On suppose que pour tout $\sigma\in\textup{Gal}\left(\bar{k}/k\right)$, il existe un isomorphisme de $\bar{G}_{X}$-torseurs:
	\[\varphi_{\sigma}:\;^{\sigma}\bar{P}\longrightarrow \bar{P}
\]
o{\`u} $^{\sigma}\bar{P}$ d{\'e}signe le sch{\'e}ma obtenu par pullback {\`a} partir de $\widetilde{\sigma}$ (l'automorphisme de $\textup{Spec}\;\bar{k}$ induit par $\sigma$); autrement dit, le carr{\'e} ci-dessous est cart{\'e}sien:
	\[\xymatrix@C=15pt@R=15pt{^{\sigma}\bar{P}\ar[dd]\ar[rr]^{\varphi_{\sigma}}&&\bar{P}\ar[dd]\\&\Box \\\textup{Spec}\;\bar{k}\ar[rr]_{\widetilde{\sigma}}&&\textup{Spec}\;\bar{k}}
\]

Employant une terminologie propre {\`a} la th{\'e}orie des rev{\^e}tements, nous dirons que le torseur $\bar{P}\rightarrow \bar{X}$ est de \textit{corps des modules $k$}.\index{corps des modules}

Le probl{\`e}me est alors de savoir s'il existe un $G_{X}$-torseur $Q\rightarrow X$ tel que $\bar{Q}\approx\bar{P}$, ce dernier isomorphisme vivant dans la cat{\'e}gorie des $\bar{G}_{X}$-torseurs sur $\bar{X}$. Nous appellerons \textit{mod{\`e}le} de $\bar{P}$ sur $X$ un tel torseur $Q\rightarrow X$. Si $\bar{P}$ poss{\`e}de un mod{\`e}le sur $X$, nous dirons qu'il est \textit{d{\'e}fini sur $k$}. \index{d{\'e}fini sur $k$}

Supposons par exemple que $G$ soit un $k$-groupe ab{\'e}lien de type multiplicatif, et que $X$ soit un $k$-sch{\'e}ma propre (\textit{e.g.} une $k$-vari{\'e}t{\'e} projective). Alors la suite exacte {\`a} 5 termes:
\begin{equation}
	0\longrightarrow H^{1}\left(k,G\right)\longrightarrow H^{1}\left(X,G_{X}\right)\stackrel{u}{\longrightarrow} H^{1}\left(\bar{X},\bar{G}_{X}\right)^{\Gamma}\stackrel{\delta^{1}}{\longrightarrow}H^{2}\left(k,G\right)\stackrel{}{\longrightarrow}H^{2}\left(X,G_{X}\right)
\end{equation}
(o{\`u} $\Gamma$ d{\'e}signe le groupe de Galois absolu de $k$) d{\'e}duite de la suite spectrale de Leray:
	\[E_{2}^{p,q}=H^{p}\left(k,R^{q}\pi_{\ast}G_{X}\right)\Longrightarrow H^{p+q}\left(X,G_{X}\right)=E^{p+q}
\]
montre que l'obstruction {\`a} ce que $\bar{P}$ (qui repr{\'e}sente une classe dans $H^{1}\left(\bar{X},\bar{G}_{X}\right)^{\Gamma}$) soit d{\'e}fini sur $k$ (\textit{i.e}: $\left[\bar{P}\right]\in \textup{im}\;u$) est mesur{\'e}e par une classe vivant dans $H^{2}\left(k,G\right)$. Une telle classe est une (classe d'{\'e}quivalence de) gerbe sur le site {\'e}tale de $k$.

Par ailleurs, sous les m{\^e}mes hypoth{\`e}ses, on a une interpr{\'e}tation en termes de \textit{type}\index{type d'un torseur} de ce probl{\`e}me. On dispose en effet de la suite exacte {\`a} 5 termes introduite par Colliot-Th{\'e}l{\`e}ne et Sansuc dans \cite{CTS}:
\begin{equation}
	H^{1}\left(k,G\right)\longrightarrow H^{1}\left(X,G_{X}\right)\stackrel{\chi}{\longrightarrow} \textup{Hom}_{\Gamma}\left(\widehat{G},\textup{Pic}\;\bar{X}\right)\stackrel{\partial}{\longrightarrow}H^{2}\left(k,G\right)\stackrel{}{\longrightarrow}H^{2}\left(X,G_{X}\right)
\end{equation}
celle-ci {\'e}tant d{\'e}duite de la suite spectrale:
	\[E_{2}^{p,q}=\textup{Ext}^{p}_{\Gamma}\left(\widehat{G},H^{q}\left(\bar{X},\mathbb{G}_{m,\bar{X}}\right)\right)\Longrightarrow H^{p+q}\left(X,G\right)=E^{p+q}
\]

Les suites exactes (1) et (2) sont canoniquement isomorphes (\textit{cf.} \cite{HS}). On appelle type d'un $G_{X}$-torseur $P\rightarrow X$ l'image de $\left[P\right]$ par le morphisme $\chi$ de la suite (2). Plus g{\'e}n{\'e}ralement, on peut associer {\`a} tout $\bar{G}_{X}$-torseur $\bar{P}\rightarrow\bar{X}$ un type not{\'e} $\tau_{\bar{P}}$, appartenant \textit{a priori} {\`a} $\textup{Hom}\left(\widehat{G},\textup{Pic}\;\bar{X}\right)$; si de plus $\bar{P}$ est de corps des modules $k$, alors $\tau_{\bar{P}}$ est $\Gamma$-{\'e}quivariant (\textit{cf.} \cite{HS} 3.7). Dans ce cas, l'existence d'un mod{\`e}le pour le $\bar{G}_{X}$-torseur $\bar{P}\rightarrow\bar{X}$ est {\'e}quivalente {\`a} l'existence d'un $G_{X}$-torseur sur $X$ de type $\tau_{\bar{P}}$.
\vspace{1mm}
\newline

Un probl{\`e}me tr{\`e}s similaire est celui de l'existence de mod{\`e}les pour un rev{\^e}tement galoisien:
\vspace{1mm}
\newline

\uline{\textquotedblleft{\bfseries{Rev{\^e}tements alg{\'e}briques: corps des modules contre corps de d{\'e}finition}}\textquotedblright}\textbf{:}\footnote{Titre honteusement plagi{\'e} sur \cite{DD2}.} soient $X$ une vari{\'e}t{\'e} alg{\'e}brique d{\'e}finie sur un corps $K$, et soit $K^{sep}$ une cl{\^o}ture s{\'e}parable de $K$. Soient $G$ un groupe fini, et $\bar{f}:\bar{Y}\rightarrow\bar{X}$ un $G$-rev{\^e}tement. On suppose que $\bar{f}$ est isomorphe {\`a} tous ses conjugu{\'e}s par l'action de $\textup{Gal}\left(K^{sep}/K\right)$.

Le $G$-rev{\^e}tement $\bar{f}$ est-il d{\'e}fini sur $K$? 
\vspace{1mm}
\newline

Une situation l{\'e}g{\`e}rement diff{\'e}rente des deux pr{\'e}c{\'e}dentes est li{\'e}e {\`a} ce que nous conviendrons d'appeler la
\vspace{1mm}
\newline

\uline{\textbf{Conjecture de Grothendieck sur les groupes de Brauer}}\textbf{:}\index{conjecture de Grothendieck sur les groupes de Brauer@Conjecture de Grothendieck sur les groupes de Brauer} soit $X$ un sch{\'e}ma. Notons $\textup{Br}_{Az}X$\index{groupe!de Brauer-Azumaya}\label{Braz} le groupe des classes d'{\'e}quivalence d'alg{\`e}bres d'Azumaya sur $X$ (\textit{cf.} \cite{G4}), et notons $\textup{Br}\;X$\index{groupe!de Brauer!cohomologique}\label{BrX} la partie de torsion du second groupe de cohomologie {\'e}tale $H^{2}\left(X,\mathbb{G}_{m}\right)$ (lorsque $X$ est r{\'e}gulier, on a: $\textup{Br}\;X=H^{2}\left(X,\mathbb{G}_{m}\right)$ d'apr{\`e}s \cite{G4} II.1.4). On a toujours un morphisme injectif de groupes:
	\[\Delta:\textup{Br}_{Az}X\longrightarrow \textup{Br}\;X
\]

La surjectivit{\'e} de ce morphisme n'est pas connue en g{\'e}n{\'e}ral, mais elle l'est dans le cas o{\`u} $X$ est:\vspace{2mm}
\begin{itemize}
\item le spectre d'un corps;\vspace{3mm}
\item un sch{\'e}ma affine, ou la r{\'e}union de deux sch{\'e}mas affines ayant une intersection affine (Gabber \cite{Ga});\vspace{3mm}
\item une vari{\'e}t{\'e} ab{\'e}lienne (Hoobler \cite{Ho1});\vspace{3mm}
\item une vari{\'e}t{\'e} torique lisse (Demeyer et Ford \cite{DF});\vspace{3mm}
\item une surface alg{\'e}brique s{\'e}par{\'e}e g{\'e}om{\'e}triquement normale (Schröer \cite{Sc});
\end{itemize}
$\ $

Le probl{\`e}me de d{\'e}terminer si $\Delta$ est surjective ou non est {\'e}videmment li{\'e} aux gerbes. Consid{\'e}rons en effet une classe $\left[c\right]\in \textup{Br}\;X$; comme $\left[c\right]$ est de torsion, il existe un entier $n$ tel que $\left[c\right]$ soit repr{\'e}sentable par une classe $\left[c'\right]\in H^{2}\left(X,\mu_{n}\right)$. Un repr{\'e}sentant de $\left[c'\right]$ est une $\mu_{n}$-gerbe sur le site {\'e}tale de $X$ (au passage, $\mathcal{G}$ est en particulier un champ de Deligne-Mumford), et $\mathcal{G}$ appartient {\`a} l'image de $\Delta$ si et seulement si $\mathcal{G}$ est isomorphe {\`a} un champ quotient (\textit{cf.} \cite{EHKV} 3.6).
\vspace{1mm}
\newline

Ce sont les m{\^e}mes types de consid{\'e}rations qui interviennent lorsque l'on {\'e}tudie la
\vspace{1mm}
\newline

\uline{\textbf{Relation de domination de Springer}}\textbf{:}\index{relation de domination de Springer} soient $K$ un corps, $K^{sep}$ une cl{\^o}ture s{\'e}parable de $K$, $G$ un $K$-groupe alg{\'e}brique, et $H$ un sous-$K$-groupe alg{\'e}brique de $G$. Il existe une relation \cite{Sp}:
	\[H^{1}\left(\textup{Gal}\left(K^{sep}/K\right),G\right) \multimap H^{1}\left(\textup{Gal}\left(K^{sep}/K\right);G,H\right)
\]
entre l'ensemble des classes d'isomorphie de $G$-torseurs ({\`a} droite) sur $K$, et celui des $K$-espaces homog{\`e}nes sous l'action ({\`a} droite) de $G$ avec isotropie $H$. Soit $V$ un tel espace homog{\`e}ne. Une fois encore, l'obstruction {\`a} ce que $V$ soit domin{\'e} par un $G$-torseur, \textit{i.e.} l'obstruction {\`a} ce que $\left[V\right]$ appartienne {\`a} l'image de la relation est mesur{\'e}e par une gerbe sur $k$, localement li{\'e}e par $H$.
\vspace{1mm}
\newline

Le dernier exemple que nous donnons se distingue des pr{\'e}c{\'e}dents par le fait que c'est purement un probl{\`e}me de G{\'e}om{\'e}trie Alg{\'e}brique, et non un probl{\`e}me d'Arithm{\'e}tique. Ceci {\'e}tant dit, m{\^e}me si les contextes sont diff{\'e}rents, les gerbes sont encore pr{\'e}sentes:
\vspace{1mm}
\newline

\uline{\textbf{Classe de Chern d'un fibr{\'e} en droites}}\textbf{:} soient $X$ une vari{\'e}t{\'e} alg{\'e}brique complexe projective lisse. On note $X^{an}$ la vari{\'e}t{\'e} analytique canoniquement associ{\'e}e {\`a} $X$ (\textit{cf.} \cite{Fu}). De la suite exponentielle, on d{\'e}duit la suite de cohomologie:
	\[\xymatrix{H^{1}\left(X^{an},\mathcal{O}_{X^{an}}\right) \ar[r]^{exp} & H^{1}\left(X^{an},\mathcal{O}^{\ast}_{X^{an}}\right) \ar[r]^{\ \ \delta^{1}} & H^{2}\left(X^{an},\mathbb{Z}\right)}
\]

D'apr{\`e}s le th{\'e}or{\`e}me de Serre (\textit{cf.} \cite{Har} p.440), qui {\'e}tablit l'{\'e}quivalence entre la cat{\'e}gorie des faisceaux coh{\'e}rents sur $X$ et celle des faisceaux analytiques coh{\'e}rents sur $X^{an}$, on en d{\'e}duit la suite exacte:
	\[\xymatrix{H^{1}\left(X,\mathcal{O}_{X}\right) \ar[r]^{\ \ \ \eta} & \textup{Pic}\;X \ar[r]^{\delta^{1}\ \ \ } & H^{2}\left(X^{an},\mathbb{Z}\right)}
\]

L'image d'un fibr{\'e} en droites $\mathcal{L}$ sur $X$ (\textit{i.e.} d'un repr{\'e}sentant d'une classe de $\textup{Pic}\;X$) par $\delta^{1}$ est appel{\'e}e la \textit{classe de Chern} de $\mathcal{L}$, et elle est not{\'e}e $c_{1}\left(\mathcal{L}\right)$.\footnote{Ce n'est pas la seule mani{\`e}re de d{\'e}finir la classe de Chern d'un fibr{\'e} en droites. C'est en tout cas celle pr{\'e}sent{\'e}e dans \cite{Fu}, 19.3.1; pour un autre point de vue, nous renvoyons {\`a} l'appendice A de \cite{Har}, et pour un retour aux sources {\`a} \cite{G2}.} On peut {\'e}videmment voir cette classe $c_{1}\left(\mathcal{L}\right)$ comme une gerbe, sur le site analytique de $X^{an}$. Sa non-nullit{\'e} est une obstruction {\`a} ce que $\mathcal{L}$ appartienne {\`a} l'image du morphisme $\eta$.
\vspace{2mm}
\newline

Apr{\`e}s ces exemples qui illustrent la diversit{\'e} des contextes dans lesquels les gerbes interviennent, voici le plan que nous suivrons dans ces notes:
\vspace{2mm}
\newline

\uline{\scshape{Chapitre 1}}: nous y introduisons les outils et le langage n{\'e}cessaires dans toute la suite. On commence par rappeler la d{\'e}finition de site en g{\'e}n{\'e}ral, en ayant {\`a} l'esprit le fait que seule la notion de site {\'e}tale d'un sch{\'e}ma nous sera vraiment utile. On rappelle ensuite les notions de cat{\'e}gorie fibr{\'e}e, pr{\'e}champ et champ afin de pouvoir pr{\'e}senter les actrices principales de ce travail: les gerbes. Cependant, plut{\^o}t que de travailler sur un site g{\'e}n{\'e}ral (comme dans \cite{Gi1} et \cite{Gi2}), nous {\'e}tudierons plus particuli{\`e}rement les propri{\'e}t{\'e}s des gerbes sur le site {\'e}tale d'un sch{\'e}ma, ce qui rendra peut-{\^e}tre plus facile leur interpr{\'e}tation. Enfin nous rappellerons les notions de lien et de cohomologie {\`a} valeurs dans un lien, en donnant des exemples de situations o{\`u} le $H^{2}$ est \textquotedblleft{facilement}\textquotedblright\ calculable.
\vspace{2mm}
\newline

\uline{\scshape{Chapitre 2}}: nous nous attaquons au probl{\`e}me central par son versant le plus facile: le cas ab{\'e}lien. Plus pr{\'e}cis{\'e}ment, on s'int{\'e}resse {\`a} la situation suivante: on consid{\`e}re un corps de caract{\'e}ristique nulle $k$, dont on fixe une cl{\^o}ture alg{\'e}brique $\bar{k}$, un $k$-sch{\'e}ma $X$ g{\'e}om{\'e}triquement irr{\'e}ductible, et un $k$-groupe alg{\'e}brique ab{\'e}lien $G$. On peut alors, en ajoutant de peu contraignantes hypoth{\`e}ses, d{\'e}terminer une obstruction cohomologique ab{\'e}lienne {\`a} l'existence d'un point $k$-rationnel sur $X$. Explicitement, on suppose satisfaite la condition $\bar{G}_{X}\left(\bar{X}\right)=\bar{G}\left(\bar{k}\right)$ (\textit{e.g.} $X$ projective et $G=\mathbb{G}_{m}$). On a alors la suite exacte {\`a} 5 termes:
	\[H^{1}\left(k,G\right)\longrightarrow H^{1}\left(X,G_{X}\right)\stackrel{u}{\longrightarrow} H^{1}\left(\bar{X},\bar{G}_{X}\right)^{\Gamma}\stackrel{\delta^{1}}{\longrightarrow}H^{2}\left(k,G\right)\stackrel{v}{\longrightarrow}H^{2}\left(X,G_{X}\right)
\]

L'existence d'un point $k$-rationnel sur $X$ entra{\^i}ne l'existence d'une r{\'e}traction de l'application $v$, ce qui force le cobord $\delta^{1}$ {\`a} {\^e}tre nul, cette nullit{\'e} entra{\^i}nant la surjectivit{\'e} de l'application $u$. On en d{\'e}duit donc une premi{\`e}re obstruction:
\begin{theor}[Obstruction ab{\'e}lienne {\`a} l'existence d'un point rationnel]Soient $k$ un corps de caract{\'e}ristique nulle, $X$ un $k$-sch{\'e}ma, quasi-compact et quasi-s{\'e}par{\'e}, et $G$ un $k$-groupe alg{\'e}brique ab{\'e}lien tels que $\bar{G}_{X}\left(\bar{X}\right)=\bar{G}\left(\bar{k}\right)$.

Si $X\left(k\right)\neq\emptyset$, alors dans la suite exacte pr{\'e}c{\'e}dente, le morphisme $u$ est surjectif. Autrement dit, tout $\bar{G}_{X}$-torseur sur $\bar{X}$ de corps des modules $k$ est d{\'e}fini sur $k$.
\end{theor}

De cet {\'e}nonc{\'e}, on extrait les cons{\'e}quences suivantes:
\begin{propo} Le corps $k$, le sch{\'e}ma $X$ et le groupe $G$ {\'e}tant comme indiqu{\'e}s dans l'{\'e}nonc{\'e} pr{\'e}c{\'e}dent, on a la suite longue de cohomologie (toujours d{\'e}duite de la suite spectrale de Leray\index{suite spectrale!de Leray}):
\vspace{2mm}
\begin{flushleft}
$0\longrightarrow H^{1}\left(k,G\right)\longrightarrow H^{1}\left(X,G_{X}\right)\stackrel{u}{\longrightarrow} H^{1}\left(\bar{X},\bar{G}_{X}\right)^{\Gamma}\stackrel{\delta^{1}}{\longrightarrow}H^{2}\left(k,G\right)$
\end{flushleft}
\vspace{2mm}
\begin{flushright}
$\stackrel{v}{\longrightarrow}H^{2}\left(X,G_{X}\right)^{tr}\longrightarrow H^{1}\left(k,H^{1}\left(\bar{X},\bar{G}_{X}\right)\right)\stackrel{\delta^{2}}{\longrightarrow}H^{3}\left(k,G\right)$
\end{flushright}
\vspace{1mm}
$\ $
\newline
o{\`u} $H^{2}\left(X,G_{X}\right)^{tr}=\ker\left\{H^{2}\left(X,G_{X}\right)\rightarrow H^{2}\left(\bar{X},\bar{G}_{X}\right)\right\}$. Supposons que $X$ poss{\`e}de un point $k$-rationnel. Alors:
\begin{enumerate}[(i)]
\item les suites:
	\[0\longrightarrow H^{1}\left(k,G\right)\longrightarrow H^{1}\left(X,G_{X}\right)\stackrel{u}{\longrightarrow} H^{1}\left(\bar{X},\bar{G}_{X}\right)^{\Gamma}\longrightarrow 0
\]
et
	\[0\longrightarrow H^{2}\left(k,G\right)\stackrel{v}{\longrightarrow} H^{2}\left(X,G_{X}\right)^{tr}\longrightarrow H^{1}\left(k,H^{1}\left(\bar{X},\bar{G}_{X}\right)\right)\longrightarrow0
\]
sont exactes;
\item le morphisme $H^3\left(k,G\right)\longrightarrow H^3\left(X,G_{X}\right)$ est injectif (c'est l'edge: $E_{2}^{3,0}\rightarrow E^3$).
\end{enumerate}

Dans le cas particulier o{\`u} $G=\mathbb{G}_{m,k}$, on a alors:
\begin{enumerate}[($\mathbb{G}_{m}$-i)]
\item tout fibr{\'e} en droites sur $\bar{X}$ de corps des modules $k$ est d{\'e}fini sur $k$;
\item la suite :
	\[0\longrightarrow \textup{Br}\;k\longrightarrow \textup{Br}^{tr}X\longrightarrow H^{1}\left(k,\textup{Pic}\;\bar{X}\right)\longrightarrow 0
\]
est exacte;
\item le morphisme $H^{3}\left(k,\mathbb{G}_{m}\right)\longrightarrow H^{3}\left(X,\mathbb{G}_{m,X}\right)$ est injectif.
\end{enumerate}
\end{propo}

Signalons que ce r{\'e}sultat est connu de longue date, puisqu'{\`a} peu de choses pr{\`e}s, l'{\'e}nonc{\'e} de cette proposition est celui du lemme 6.3 de \cite{Sa}. Notons aussi que la condition $\bar{G}_{X}\left(\bar{X}\right)=\bar{G}\left(\bar{k}\right)$ est remplie, lorsque $G$ est de type multiplicatif, si $X$ est tel que $\bar{k}\left[X\right]^{\ast}=\bar{k}^{\ast}$. Les $k$-vari{\'e}t{\'e}s satisfaisant cette derni{\`e}re condition constituent d'ailleurs, pour citer Skorobogatov \textquotedblleft{une classe raisonnable de vari{\'e}t{\'e}s pour lesquels les m{\'e}thodes de descente fonctionnent bien}\textquotedblright\ (\textit{cf.} \cite{Sk} p.407). Cependant, toutes les vari{\'e}t{\'e}s ne satisfont pas cette condition, et il est tr{\`e}s instructif de regarder ce qui arrive sur un exemple (d{\^u} {\`a} J.-L. Colliot-Th{\'e}l{\`e}ne et O. Gabber) de vari{\'e}t{\'e} pour laquelle l'obstruction ab{\'e}lienne ne tient pas. Plus pr{\'e}cis{\'e}ment, nous {\'e}tudions une vari{\'e}t{\'e} $X$, telle que $\bar{k}\left[X\right]^{\ast}=\bar{k}^{\ast}\oplus\mathbb{Z}$, poss{\'e}dant un point $k$-rationnel, mais telle que le morphisme:
	\[\textup{Pic}\;X\longrightarrow \left(\textup{Pic}\;\bar{X}\right)^{\Gamma}
\]
n'est pas surjectif.

Il est {\'e}galement int{\'e}ressant de remarquer que l'absence de point rationnel n'est pas une obstruction {\`a} la descente des torseurs en g{\'e}n{\'e}ral\footnote{\textit{I.e.} la r{\'e}ciproque du th{\'e}or{\`e}me 1 est fausse en g{\'e}n{\'e}ral.}. Concr{\`e}tement, d'apr{\`e}s \cite{CTS}, la descente des $\bar{G}_{X}$-torseurs ($G$ {\'e}tant ab{\'e}lien) sur $\bar{X}$ est possible d{\`e}s que $X$ poss{\`e}de des points ad{\'e}liques d'un certain type (ce qui est plus faible que de demander l'existence de points rationnels). Nous verrons que les vari{\'e}t{\'e}s dont l'obstruction de Brauer-Manin est nulle ont justement des points ad{\'e}liques de ce type.
\vspace{2mm}
\newline

\uline{\scshape{Chapitre 3}}: l'objet de ce chapitre est de s'inspirer du cas ab{\'e}lien pour obtenir une obstruction non-ab{\'e}lienne {\`a} l'existence de point rationnel. Bien entendu, il n'est plus question d'utiliser des suites spectrales, et il ne subsiste de la suite exacte {\`a} 5 termes (1) que la suite exacte au sens des ensembles point{\'e}s:
	\[0\longrightarrow H^{1}\left(k,\pi_{\ast}G_{X}\right)\longrightarrow H^{1}\left(X,G_{X}\right)\stackrel{u}{\longrightarrow} H^{1}\left(\bar{X},\bar{G}_{X}\right)^{\Gamma}
\]

Pour esp{\'e}rer prolonger cette suite et obtenir un analogue de la suite {\`a} 5 termes dans le cas non-ab{\'e}lien, on d{\'e}finit, pour tout $\bar{P}\rightarrow\bar{X}$ repr{\'e}sentant une classe dans $H^{1}\left(\bar{X},\bar{G}_{X}\right)^{\Gamma}$ sa gerbe des mod{\`e}les, not{\'e}e $D\left(\bar{P}\right)$. C'est une gerbe sur $k$, dont la neutralit{\'e} est n{\'e}cessaire et suffisante pour que $\bar{P}$ soit d{\'e}fini sur $k$ (ou encore pour que $\bar{P}$ ait un mod{\`e}le sur $X$). La difficult{\'e} est {\'e}videmment que cette gerbe ne repr{\'e}sente en g{\'e}n{\'e}ral pas une classe de \textquotedblleft{$H^{2}\left(k,G\right)$}\textquotedblright, ni m{\^e}me de $H^{2}\left(k,\textup{lien}\;G'\right)$, o{\`u} $G'$ serait une $k$-forme de $G$. Le lemme suivant, qui est trivial mais d'une importance capitale, explique en partie l'origine de cette difficult{\'e}:
\begin{lemfo} Soient $S$ un sch{\'e}ma, $G_{S}$ un sch{\'e}ma en groupes sur $S$, et $P$ un $G_{S}$-torseur sur $S$. Alors $\textup{ad}_{G_{S}}\left(P\right)$ est une $S$-forme int{\'e}rieure de $G_{S}$; autrement dit:
	\[\textup{ad}_{G_{S}}\left(P\right)\textit{ repr{\'e}sente une classe de } H^{1}\left(S,\textup{Int}\;G_{S}\right).
\]

En particulier, si $G_{S}$ est ab{\'e}lien, alors: 
	\[\textup{ad}_{G_{S}}\left(P\right)\approx G_{S}
\]
\end{lemfo}

Malgr{\'e} cela, nous montrons que l'existence d'un point $k$-rationnel $x$ sur $X$ entra{\^i}ne justement la neutralit{\'e} de $D\left(\bar{P}\right)$, lorsque l'on fait l'hypoth{\`e}se suppl{\'e}mentaire:\vspace{2mm}

\uline{\textbf{Condition} $\left(\star_{\bar{P}}\right)$}: \textit{En notant $G'=\textup{ad}_{G_{\bar{X}}}\left(\bar{P}\right)$ et $\bar{x}$ un point g{\'e}om{\'e}trique associ{\'e} au point $k$-rationnel $x$, on a:}
	\[H^{0}\left(\bar{X},G'\right)=H^{0}\left(\bar{k},G'_{\bar{x}}\right)
\]\vspace{2mm}

qui est l'analogue non-ab{\'e}lien de la condition:
	\[\bar{G}_{X}\left(\bar{X}\right)=\bar{G}\left(\bar{k}\right)
\]

Sous cette hypoth{\`e}se, on a le r{\'e}sultat suivant:
\begin{theor} Soient $X$ un $k$-sch{\'e}ma, et $G$ un $k$-groupe lin{\'e}aire, et $\bar{P}\rightarrow \bar{X}$ un $G_{\bar{X}}$-torseur de corps des modules $k$. On suppose que $X$ poss{\`e}de un point $k$-rationnel $x$, et on suppose satisfaite la condition suivante:

\textbf{Condition $\left(\star_{\bar{P}}\right)$:} En notant $G'=\textup{ad}_{G_{\bar{X}}}\left(\bar{P}\right)$ et $\bar{x}$ un point g{\'e}om{\'e}trique associ{\'e} au point $k$-rationnel $x$, on a:
	\[H^{0}\left(\bar{X},G'\right)=H^{0}\left(\bar{k},G'_{\bar{x}}\right)
\]

Alors $\bar{P}\rightarrow\bar{X}$ est d{\'e}fini sur $k$.
\end{theor}

D'autre part, de nombreux obstacles se pr{\'e}sentent lorsque l'on cherche {\`a} obtenir des informations pr{\'e}cises sur le lien de la gerbe des mod{\`e}les d'un torseur. Cependant, dans le cas o{\`u} le groupe $G$ est fini, on retrouve l'{\'e}nonc{\'e} suivant, {\`a} rapprocher de celui de Harari et Skorobogatov \cite{HS}, th{\'e}or{\`e}me 2.5 et section 3.1:
\begin{theor} Soient $k$ un corps de caract{\'e}ristique nulle, $X$ un $k$-sch{\'e}ma g{\'e}om{\'e}triquement connexe, et $G$ un $k$-groupe fini. Si $X\left(k\right)\neq\emptyset$, alors tout $\bar{G}_{X}$-torseur sur $\bar{X}$ de corps des modules $k$ est d{\'e}fini sur $k$.
\end{theor}

Enfin, nous essayons d'expliquer ce qui emp{\^e}che d'obtenir une description pr{\'e}cise du lien de $D\left(\bar{P}\right)$, et ce qui emp{\^e}che donc \textit{a priori} de raffiner les {\'e}nonc{\'e}s pr{\'e}c{\'e}dents. Nous {\'e}tudierons aussi ce qui se passe si l'on n'impose plus la condition corps des modules sur les torseurs que l'on cherche {\`a} descendre. Pr{\'e}cis{\'e}ment, on peut toujours associer {\`a} un torseur $\bar{P}$ sur $\bar{X}$ repr{\'e}sentant une classe dans $H^{1}\left(\bar{X},\bar{G}_{X}\right)$ sa gerbe des mod{\`e}les, mais celle-ci n'est plus n{\'e}cessairement une $k$-gerbe. De fait, c'est une $k_{\bar{P}}$-gerbe, o{\`u} $k_{\bar{P}}$ d{\'e}signe le corps des modules\footnote{C'est la plus petite extension {\'e}tale $L$ de $k$ telle que:
	\[^{\sigma}\bar{P}\approx\bar{P},\ \forall\ \sigma\in\textup{Gal}\left(\bar{k}/L\right).
\]
} de $\bar{P}$. Le torseur $\bar{P}$ repr{\'e}sente un point $\xi_{\bar{P}}$ du $k$-champ $\pi_{\ast}\textup{Tors}\left(X,G_{X}\right)$ obtenu {\`a} partir de la gerbe des $G_{X}$-torseurs sur $X$ par image directe \textit{via} le morphisme structural $\pi:X\rightarrow\textup{Spec}\;k$. Ce champ joue le r{\^o}le de champ des modules grossiers pour les $\bar{G}_{X}$-torseurs sur $\bar{X}$, et on a une interpr{\'e}tation en termes de gerbe r{\'e}siduelle pour $D\left(\bar{P}\right)$. En utilisant la terminologie de Giraud, on peut {\'e}galement interpr{\'e}ter $D\left(\bar{P}\right)$ comme une section au-dessus de l'ouvert $\left(\textup{Spec}\;k_{\bar{P}}\rightarrow \textup\;k\right)$ du faisceau des sous-gerbes maximales du $k$-champ $\pi_{\ast}\textup{Tors}\left(X,G_{X}\right)$. Nous essayons d'expliquer le lien entre ces deux fa\c{c}ons d'aborder les choses. Il nous semble que le diagramme commutatif ci-dessous r{\'e}sume assez bien la situation:
	\[\xymatrix@C=25pt@R=20pt{D\left(\bar{P}\right)\ar@{^{(}->}[rr]^{i\ \ \ \ \ } \ar[dd]_{\pi} && \pi_{\ast}\textup{Tors}\left(X,G_{X}\right)\ar[dd]^{\left[\bullet\right]}\\\\\textup{Spec}\;k_{\bar{P}}\ar[rr]_{\left[\bar{P}\right]}&&R^{1}\pi_{\ast}G_{X}}
\]

Dans ce diagramme, on a not{\'e}:\vspace{2mm}
\begin{itemize}
\item $\pi:D\left(\bar{P}\right)\longrightarrow\textup{Spec}\;k_{\bar{P}}$ le morphisme structural de la gerbe des mod{\`e}les de $\bar{P}$;\vspace{2mm}
\item $i:D\left(\bar{P}\right)\longrightarrow \pi_{\ast}\textup{Tors}\left(X,G_{X}\right)$ d{\'e}signe le monomorphisme canonique (rendu explicite dans le chapitre 3);\vspace{2mm}
\item $\left[\bullet\right]: \pi_{\ast}\textup{Tors}\left(X,G_{X}\right)\longrightarrow R^{1}\pi_{\ast}G_{X}$ est d{\'e}fini en envoyant une section du champ $\pi_{\ast}\textup{Tors}\left(X,G_{X}\right)$, \textit{i.e.} un $G_{X_{L}}$-torseur $T_{L}\rightarrow X_{L}$ sur sa classe d'isomorphie $\left[T_{L}\right]$;\vspace{2mm}
\item enfin, le choix de $\bar{P}$ donne naissance {\`a} un point de $R^{1}\pi_{\ast}G_{X}$ {\`a} valeurs dans $\textup{Spec}\;k_{\bar{P}}$; c'est ce point que nous avons not{\'e} $\left[\bar{P}\right]$.
\end{itemize}
\vspace{2mm}

\uline{\scshape{Chapitre 4}}: les espaces homog{\`e}nes sur un corps sont un exemple de vari{\'e}t{\'e}s pour lesquelles l'existence de points rationnels a {\'e}t{\'e} et est encore particuli{\`e}rement {\'e}tudi{\'e}; un probl{\`e}me source de nombreuses activit{\'e}s est celui de la validit{\'e} du principe de Hasse. Plus pr{\'e}cis{\'e}ment, il serait int{\'e}ressant de savoir si l'obstruction de Brauer-Manin est la seule pour les espaces homog{\`e}nes sous l'action d'un groupe lin{\'e}aire $G$ avec isotropie $H$. D'apr{\`e}s Sansuc \cite{Sa}, on sait d{\'e}j{\`a} que c'est le cas pour les torseurs sous des groupes connexes (\textit{i.e.} $G$ connexe et $H=\left\{1\right\}$), et d'apr{\`e}s Borovoï \cite{Bo2} c'est aussi le cas pour des espaces homog{\`e}nes sous des groupes connexes avec isotropie connexe, ou pour des espaces homog{\`e}nes sous des groupes simplement connexes avec isotropie ab{\'e}lienne finie.

Pour fixer les id{\'e}es, consid{\'e}rons un corps de nombres $k$, et $H$ un sous-groupe de $SL_{n}$. L'obstruction {\`a} ce qu'un $k$-espace homog{\`e}ne $V$ sous $SL_{n}$ avec isotropie $H$ poss{\`e}de un point $k$-rationnel est mesur{\'e}e par une gerbe $\mathcal{G}$ sur $k$, localement li{\'e}e par $H$; la neutralit{\'e} de $\mathcal{G}$ (\textit{i.e.} l'existence d'un point $k$-rationnel sur $\mathcal{G}$) est {\'e}quivalente {\`a} l'existence d'un point $k$-rationnel sur $V$. Or la m{\^e}me gerbe $\mathcal{G}$ peut correspondre {\`a} plusieurs espaces homog{\`e}nes. D'o{\`u} l'id{\'e}e, dans un souci d'{\'e}conomie, de travailler directement avec les gerbes plut{\^o}t qu'avec les espaces homog{\`e}nes.

Dans ce chapitre, fruit d'un travail en commun avec Jean-Claude Douai et Michel Emsalem, on commence donc par d{\'e}finir l'obstruction de Brauer-Manin d'une gerbe\footnote{Pour des raisons techniques, on doit consid{\'e}rer des gerbes qui sont des champs de Deligne-Mumford, ce qui n'est pas g{\^e}nant dans nos applications, puisque le cas int{\'e}ressant est justement celui o{\`u} l'isotropie est finie.}. On montre que pour tout $k$-espace homog{\`e}ne $V$ sous $SL_{n}$ (ou n'importe quel autre groupe semi-simple simplement connexe) avec isotropie finie, on a:
	\[m_{\mathcal{H}}\left(V\right)=m_{\mathcal{H}}\left(\mathcal{G}_{V}\right)
\]
o{\`u} $\mathcal{G}_{V}$ est la gerbe des trivialisations de $V$ (pour faire le lien avec les travaux de Springer, c'est aussi l'obstruction {\`a} ce que $V$ soit domin{\'e} par un $k$-torseur sous $SL_{n}$; mais un tel torseur est toujours trivial). On obtient, gr{\^a}ce {\`a} cette obstruction, une nouvelle interpr{\'e}tation du th{\'e}or{\`e}me de Tate-Poitou (quand $H$ est ab{\'e}lien fini), et un \textquotedblleft{demi-th{\'e}or{\`e}me de Tate-Poitou}\textquotedblright\ lorsque $H$ est non-ab{\'e}lien fini.
\vspace{5mm}

\newpage

\thispagestyle{empty}
\chapter*{Notations}
\addcontentsline{toc}{chapter}{Notations}
\markboth{NOTATIONS}{NOTATIONS}
\thispagestyle{empty}
$\ $
\vspace{6mm}
\newline

Etant donn{\'e} un corps $k$, nous appellerons \textit{$k$-groupe alg{\'e}brique} un $k$-groupe alg{\'e}brique lin{\'e}aire, et \textit{$k$-groupe alg{\'e}brique r{\'e}ductif} un $k$-groupe alg{\'e}brique r{\'e}ductif et connexe.\vspace{10mm}

Pour tout corps $k$ et tout $k$-sch{\'e}ma $Y$, nous noterons $Y\left(k\right)$ l'ensemble des points $k$-rationnels de $Y$. Pour toute extension $L$ de $k$, $Y\left(L\right)$ d{\'e}signera l'ensemble des points de $Y$ {\`a} valeurs dans $\textup{Spec}\;L$.\vspace{10mm}

Si $Y$ est un sch{\'e}ma et $G$ un sch{\'e}ma en groupes (\textit{resp.} un sch{\'e}ma en groupes ab{\'e}liens) sur $Y$, nous noterons $H^{0}\left(Y,G\right)$ et $H^{1}\left(Y,G\right)$ (\textit{resp.} $H^{i}\left(Y,G\right)$, $i\geq0$) les ensembles (\textit{resp.} les groupes) de cohomologie {\'e}tale $H^{0}_{\acute{e}t}\left(Y,G\right)$ et $H^{1}_{\acute{e}t}\left(Y,G\right)$ (\textit{resp.} $H^{i}_{\acute{e}t}\left(Y,G\right)$, $i\geq0$), o{\`u} $G$ est identifi{\'e} au faisceau de groupes qu'il repr{\'e}sente sur le site {\'e}tale de $Y$.\vspace{10mm}

Soient $k$ un corps de caract{\'e}ristique nulle, $\bar{k}$ une cl{\^o}ture alg{\'e}brique de $k$, $\Gamma=\textup{Gal}\left(\bar{k}/k\right)$, $\pi:X\rightarrow\textup{Spec}\;k$ un $k$-sch{\'e}ma, et $G$ un $k$-groupe alg{\'e}brique. Soit encore $\bar{P}\rightarrow\bar{X}$ un $\bar{G}_{X}$-torseur ($\bar{G}_{X}=\left(G\times_{\textup{Spec}\;k}\textup{Spec}\;\bar{k}\right)\times_{\textup{Spec}\;\bar{k}}\bar{X}$). Nous dirons que ce torseur est de \textit{corps des modules $k$} lorsque $\bar{P}\rightarrow\bar{X}$ repr{\'e}sente une classe de $H^{1}\left(\bar{X},\bar{G}_{X}\right)^{\Gamma}$ (notons que cette appellation repr{\'e}sente un abus par rapport {\`a} la d{\'e}finition usuelle de corps des modules dans la th{\'e}orie des rev{\^e}tements; avec celle-ci en effet, la condition $\left[\bar{P}\right]\in H^{1}\left(\bar{X},\bar{G}_{X}\right)^{\Gamma}$ assure seulement que le corps des modules de $\bar{P}\rightarrow\bar{X}$ est \textit{inclus} dans $k$).\vspace{10mm}

Lorsque $\mathcal{C}$ est une cat{\'e}gorie, nous noterons $\textup{Ob}\left(\mathcal{C}\right)$ la classe de ses objets. Etant donn{\'e}s deux objets $A$ et $B$ de $\mathcal{C}$, nous noterons $\textup{Hom}_{\mathcal{C}}\left(A,B\right)$ la classe des morphismes (ou des fl{\`e}ches) de $\mathcal{C}$ entre $A$ et $B$. Si de plus $\mathcal{C}$ est un groupoïde, $\textup{Isom}_{\mathcal{C}}\left(A,B\right)$ d{\'e}signera la classe des isomorphismes de $\mathcal{C}$ entre $A$ et $B$. Nous conviendrons que la cat{\'e}gorie vide est un groupoïde. Enfin, nous noterons $\mathfrak{Ens}$ (\textit{resp.} $\mathfrak{Gr}$, \textit{resp.} $\mathfrak{Ab}$, \textit{resp.} $FAGR\left(Y\right)$, \textit{resp.} $FAGRAB\left(Y\right)$) la cat{\'e}gorie des ensembles (\textit{resp.} des groupes, \textit{resp.} des groupes ab{\'e}liens, \textit{resp.} des faisceaux de groupes sur le site {\'e}tale de $Y$, \textit{resp.} des faisceaux de groupes ab{\'e}liens sur le site {\'e}tale de $Y$).
\newpage
\thispagestyle{empty}
$\ $
\vspace{5mm}
\newline
$\ $

\pagenumbering{arabic}
\begin{chapter}{Champs et gerbes}
\thispagestyle{empty}
Dans ce chapitre, nous rappelons les notions de site et de topos, qui fournissent une g{\'e}n{\'e}ralisation de la notion d'espace topologique. Cette g{\'e}n{\'e}ralisation permet comprendre pourquoi les probl{\`e}mes {\'e}voqu{\'e}s dans l'introduction (descente de torseurs ou de rev{\^e}tements, banalisation d'une alg{\`e}bre d'Azumaya,\ldots) sont de m{\^e}me nature, dans le sens o{\`u} on peut tous les interpr{\'e}ter comme des probl{\`e}mes de recollement, moyennant le choix d'un site idoine.

Notons tout de suite que nous nous r{\'e}duirons tr{\`e}s vite en pratique au site {\'e}tale d'un sch{\'e}ma (qui sera d'ailleurs souvent le site {\'e}tale d'un corps). Cette restriction est motiv{\'e}e d'une part par le fait que notre probl{\`e}me central appara{\^i}t naturellement comme un probl{\`e}me de descente sur le site {\'e}tale d'un corps, et d'autre part parce que la manipulation des champs et des gerbes est assez d{\'e}licate et lourde sur des sites g{\'e}n{\'e}raux.

On dispose sur les sites (et sur les topoï, qui sont des sites particuliers) des m{\^e}mes objets et outils que sur les espaces topologiques. En particulier, on a une notion de faisceau, et la deuxi{\`e}me section est consacr{\'e}e {\`a} l'{\'e}tude de faisceaux particuliers: les torseurs (sur le site {\'e}tale d'un sch{\'e}ma), qui jouent un des premiers r{\^o}les dans ce travail. Les exemples de torseurs sont nombreux \textquotedblleft{dans la nature}\textquotedblright, aussi bien en G{\'e}om{\'e}trie Alg{\'e}brique (\textit{e.g}: les fibr{\'e}s vectoriels de rang $n$ sur une vari{\'e}t{\'e} alg{\'e}brique $X$ \textquotedblleft{sont}\textquotedblright\ les $GL_{n,X}$-torseurs\footnote{Ce n'est pas tout-{\`a} fait vrai: pour {\^e}tre pr{\'e}cis, le m{\^e}me ensemble $H^{1}\left(X,GL_{n,X}\right)$ peut {\^e}tre interpr{\'e}t{\'e} soit comme l'ensemble des classes d'isomorphie de $GL_{n,X}$-torseurs sur $X$, soit comme l'ensemble des classes d'isomorphie de fibr{\'e}s vectoriels de rang $n$ sur $X$.} sur le site de Zariski de $X$) qu'en Arithm{\'e}tique (\textit{e.g}: les alg{\`e}bres simples centrales d'indice $n$ sur un corps $K$ \textquotedblleft{sont}\textquotedblright\ les $PGL_{n}$-torseurs sur le site {\'e}tale de $K$). Apr{\`e}s avoir rappel{\'e} quelques-unes de leurs propri{\'e}t{\'e}s, nous constaterons avec d{\'e}pit le manque de structure de l'ensemble des classes d'isomorphie de $G$-torseurs, lorsque $G$ n'est pas ab{\'e}lien. Ce manque de structure peut toutefois {\^e}tre partiellement combl{\'e} en sym{\'e}trisant la situation, \textit{via} les bitorseurs. En anticipant un peu, disons que les bitorseurs sont particuli{\`e}rement bien adapt{\'e}s {\`a} notre probl{\`e}me central, puisque, dans un sens que nous pr{\'e}ciserons dans le chapitre III, ils permettent de \textquotedblleft{ne pas perdre d'informations}\textquotedblright.

Munis de ces outils, on peut enfin donner la notion de gerbe, dont nous verrons qu'elle est tr{\`e}s fortement li{\'e}e {\`a} celles de torseur et de bitorseur. Il est impossible de parler de gerbes sans {\'e}voquer les champs, ce qui justifie les rappels sur les cat{\'e}gories fibr{\'e}es, les pr{\'e}champs\ldots\ Pour avoir un premier aper\c{c}u de ce que peuvent {\^e}tre les gerbes, nous d{\'e}crivons ensuite celles qui sont associ{\'e}es aux divers probl{\`e}mes {\'e}voqu{\'e}s dans l'introduction. Enfin, nous rappelons les d{\'e}finitions de lien et de 2-cohomologie {\`a} valeurs dans un lien. Dans leur {\'e}crasante majorit{\'e}, les propri{\'e}t{\'e}s concernant les gerbes que nous rappelons ici sont issues de \cite{Gi1} et \cite{Gi2}, mais elles sont peut-{\^e}tre rendues un peu plus explicites par notre choix de sites particuliers.
\begin{section}{Sites et topoï}
Comme nous venons de le signaler, les sites (et les topoï) g{\'e}n{\'e}ralisent les espaces topologiques\footnote{Afin d'illustrer cette affirmation, citons Grothendieck (\cite{G2}, p.301): \textquotedblleft{$\left[\ldots\right]$la notion de topos, d{\'e}riv{\'e} naturel du point de vue faisceautique en Topologie, constitue un {\'e}largissement substantiel de la notion d'espace topologique, englobant un grand nombre de situations qui autrefois n'{\'e}taient pas consid{\'e}r{\'e}es comme relevant de l'intuition topologique. Le trait caract{\'e}ristique de telles situations est qu'on y dispose d'une notion de \textquotedblleft{localisation}\textquotedblright, notion qui est formalis{\'e}e pr{\'e}cis{\'e}ment par la notion de site et, en derni{\`e}re analyse, par celle de topos$\left[\ldots\right]$}\textquotedblright.}. Grossi{\`e}rement, un site est une cat{\'e}gorie munie d'une topologie, et plus pr{\'e}cis{\'e}ment:
\begin{defi} Soit $\mathcal{C}$ une cat{\'e}gorie. Pour tout objet $U$ de $\mathcal{C}$, on se donne des familles de morphismes, et on note $\textup{Cov}\left(\mathcal{C}\right)$ la r{\'e}union de ces familles. On dit que $\textup{Cov}\left(\mathcal{C}\right)$ est une \textbf{topologie de Grothendieck}\index{topologie de Grothendieck} si les conditions suivantes sont satisfaites:\begin{enumerate}[(TG 1)]
\item Si $U\stackrel{\phi}{\rightarrow}V$ est un isomorphisme, alors $\phi\in{\textup{Cov}\left(\mathcal{C}\right)}$.
\item(Caract{\`e}re local). Si $\left(T_{\alpha}\xrightarrow{\phi_{\alpha}}T\right)_{\alpha}$ est une famille dans $\textup{Cov}\left(\mathcal{C}\right)$, et si pour tout indice $\alpha$, $\left(U_{\alpha\beta}\xrightarrow{\phi_{\alpha\beta}}T_{\alpha}\right)_{\beta}$ est une famille de $\textup{Cov}\left(\mathcal{C}\right)$, alors la famille 
	\[\left(U_{\alpha\beta}\xrightarrow{\phi_{\alpha}\circ\phi_{\alpha\beta}}T\right)_{\alpha\beta}
\]
obtenue \textquotedblleft{par composition}\textquotedblright, appartient {\`a} $\textup{Cov}\left(\mathcal{C}\right)$.
\item(Stabilit{\'e} par changement de base). Si $\left(T_{\alpha}\xrightarrow{\phi_{\alpha}}T\right)_{\alpha}$ est une famille dans $\textup{Cov}\left(\mathcal{C}\right)$, et si $V$ est un objet de $\mathcal{C}$ sur\footnote{\textit{I.e}: il existe un morphisme de domaine $V$ et de codomaine $T$.} $T$, alors les produits fibr{\'e}s $T_{\alpha}\times_{T}V$ existent dans $\mathcal{C}$, et la famille 
	\[\left(\left(\phi_{\alpha}\right)_{V}:T_{\alpha}\times_{T}V\rightarrow{V}\right)_{\alpha}
\]
appartient {\`a} $\textup{Cov}\left(\mathcal{C}\right)$.
\end{enumerate}
On appelle \textbf{site}\index{site} une cat{\'e}gorie munie d'une topologie de Grothendieck.
\end{defi}
\begin{exem}\textup{Le site $\EuScript{O}uv_{X}$ associ{\'e} {\`a} un espace topologique $X$: on note $\EuScript{O}uv$ la cat{\'e}gorie dont les objets sont les ouverts de $X$, et dont les morphismes sont les inclusions. On choisit pour tout ouvert $U$ de $X$ la famille de toutes les familles d'ouverts $\left(V_{i}\hookrightarrow U\right)_{i\in I}$ dont la r{\'e}union est {\'e}gale {\`a} $U$. On note $\textup{Cov}\left(\EuScript{O}uv\right)$ la r{\'e}union de toutes ces familles; alors $\textup{Cov}\left(\EuScript{O}uv\right)$ est une topologie de Groethendieck sur $\EuScript{O}uv$; en effet, les trois propri{\'e}t{\'e}s de la d{\'e}finition se traduisent dans cet exemple de la fa\c{c}on suivante:
\begin{enumerate}[(TG 1)]
\item Le singleton $\left\{U\xrightarrow{id}U\right\}$ constitue un recouvrement ouvert de $U$;
\item \textquotedblleft{Un recouvrement ouvert d'un recouvrement ouvert est un recouvrement ouvert}\textquotedblright;
\item Si $\left(U_{\alpha}\right)_{\alpha}$ est un recouvrement ouvert d'un ouvert $U$, et si $V$ est un ouvert inclus dans $U$, alors\footnote{Dans cette cat{\'e}gorie, le produit fibr{\'e} n'est autre que l'intersection.} $\left(U_{\alpha}\cap{V}\right)_{\alpha}$ est un recouvrement ouvert de $V$.
\end{enumerate}
Par suite $\EuScript{O}uv_{X}=\left(\EuScript{O}uv,\textup{Cov}\left(\EuScript{O}uv\right)\right)$ est un site. En particulier, si on prend pour espace topologique $X$ un sch{\'e}ma, muni de la topologie de Zariski, on obtient le \textit{site de Zariski}\index{site!de Zariski} de $X$, not{\'e} $X_{\textup{Zar}}$.}
\end{exem}
\begin{exem}\textup{Le site $\left(Sch\right)$: c'est la cat{\'e}gorie des sch{\'e}mas, munie de la topologie de Grothendieck pour laquelle une famille $\left(S_{\alpha}\xrightarrow{\phi_{\alpha}}S\right)_{\alpha}$ est couvrante (\textit{i.e.} est dans $\textup{Cov}\left(Sch\right)$) si et seulement si elle est surjective. Il est de nouveau imm{\'e}diat (moyennant l'existence du produit fibr{\'e} dans la cat{\'e}gorie des sch{\'e}mas, qui elle, n'est pas imm{\'e}diate, \textit{cf.} \cite{EGA1} 3.2.1\ldots) que les propri{\'e}t{\'e}s (TG 1), (TG 2) et (TG 3) sont v{\'e}rifi{\'e}es.}
\end{exem}
\begin{exem}\textup{Le site $S_{\acute{e}t}$\label{Setsite} des sch{\'e}mas {\'e}tales sur un sch{\'e}ma $S$: c'est la cat{\'e}gorie $\left(Et/S\right)$ des $S$-sch{\'e}mas {\'e}tales ($S$ {\'e}tant un sch{\'e}ma quelconque), munie de la topologie de Grothendieck pour laquelle une famille de morphismes {\'e}tales est couvrante si et seulement si elle est surjective. Les propri{\'e}t{\'e}s (TG 1), (TG 2) et (TG 3) sont v{\'e}rifi{\'e}es, car l'identit{\'e} est un morphisme {\'e}tale, et le caract{\`e}re {\'e}tale est stable par composition et par changement de base (\textit{cf.} \cite{Mi} proposition I.3.3). On appelle ce site le \textit{site {\'e}tale}\index{site!etale d'un sch{\'e}ma@{\'e}tale d'un sch{\'e}ma} de $S$, et on le note $S_{\acute{e}t}$.}

\textup{Dans le cas particulier o{\`u} $S$ est le spectre d'un corps $k$, on obtient le site {\'e}tale de $k$: un objet de $\left(\textup{Spec}\:k\right)_{\acute{e}t}$ est une $k$-alg{\`e}bre {\'e}tale, c'est-{\`a}-dire une $k$-alg{\`e}bre isomorphe {\`a} un produit d'extensions s{\'e}parables finies de $k$.}
\end{exem}
\begin{defi} Soient $\mathcal{S}=\left(\mathcal{C},\textup{Cov}\left(\mathcal{C}\right)\right)$ et $\mathcal{S}'=\left(\mathcal{C}',\textup{Cov}\left(\mathcal{C'}\right)\right)$ deux sites. On appelle \textbf{morphisme de sites}\index{morphisme!de sites} la donn{\'e}e d'un foncteur $F:\mathcal{C}\rightarrow\mathcal{C}'$ satisfaisant les propri{\'e}t{\'e}s suivantes:
\begin{enumerate}[(i)]
\item Si $\left(U_{\alpha}\xrightarrow{\phi_{\alpha}}U\right)_{\alpha}$ est une famille de $\textup{Cov}\left(\mathcal{C}\right)$, alors $\left(F\left(U_{\alpha}\right)\xrightarrow{F\left(\phi_{\alpha}\right)}F\left(U\right)\right)_{\alpha}$ est une famille de $\textup{Cov}\left(\mathcal{C}'\right)$;
\item Si $\left(U_{\alpha}\xrightarrow{\phi_{\alpha}}U\right)_{\alpha}$ est une famille de $\textup{Cov}\left(\mathcal{C}\right)$ et si $V\rightarrow U$ est un morphisme (dans la cat{\'e}gorie $\mathcal{C}$), alors le morphisme canonique:
	\[F\left(U_{\alpha}\times_{U}V\right)\longrightarrow F\left(U_{\alpha}\right)\times_{F\left(U\right)}F\left(V\right)
\]
est un isomorphisme pour tout indice $\alpha$.
\end{enumerate}
\end{defi}
\begin{exem}\textup{Soient $X$ et $Y$ deux espaces topologiques. Une application continue $f:X\rightarrow Y$ induit {\'e}videmment un morphisme de sites:}
	\[F:\EuScript{O}uv_{Y}\longrightarrow \EuScript{O}uv_{X}
\]
\textup{le foncteur $F$ associant {\`a} un ouvert $U\subset Y$ (\textit{resp.} {\`a} l'inclusion $V\hookrightarrow U\subset Y$) l'ouvert $f^{-1}\left(U\right)\subset X$ (\textit{resp.} l'inclusion $f^{-1}\left(V\right)\hookrightarrow f^{-1}\left(U\right)\subset X$).}
\end{exem}
\begin{exem}\textup{Soit $f:Y\rightarrow S$ un morphisme de sch{\'e}mas. On peut lui associer un morphisme de sites (\textit{cf}. \cite{SGA4-7} 1.4):}
	\[f_{\acute{e}t}:S_{\acute{e}t}\longrightarrow Y_{\acute{e}t}
\]
\textup{d{\'e}fini en associant {\`a} un ouvert {\'e}tale $\left(S'\rightarrow S\right)$ de $S$ l'ouvert {\'e}tale de $Y$ obtenu par changement de base:}
	\[Y\times_{S}S'\longrightarrow Y
\]
\end{exem}
\begin{defi} Soit $\mathcal{S}=\left(\mathcal{C},\textup{Cov}\left(\mathcal{C}\right)\right)$ un site. On appelle \textbf{pr{\'e}faisceau}\index{pr{\'e}faisceau sur un site} d'ensembles (\textit{resp.} de groupes,\ldots) sur $\mathcal{S}$ un foncteur:
	\[\mathcal{P}:\mathcal{C}^{0}\longrightarrow \mathfrak{Ens}\textup{  (\textit{resp.} $\mathfrak{Grp}$,\ldots)}
\]

On appelle \textbf{faisceau}\index{faisceau!sur un site} d'ensembles (\textit{resp.} de groupes,\ldots) sur $\mathcal{S}$ un pr{\'e}faisceau satisfaisant la condition suivante: pour toute famille $\left(U_{\alpha}\rightarrow U\right)_{\alpha}$ de $\textup{Cov}\left(\mathcal{C}\right)$, le diagramme:
	\[\xymatrix{\mathcal{F}\left(U\right)\ar[r]&\prod_{\alpha}\mathcal{F}\left(U_{\alpha}\right)\ar@<0.5ex>[r]\ar@<-0.5ex>[r]&\prod_{\alpha,\beta}\mathcal{F}\left(U_{\alpha}\times_{U}U_{\beta}\right)}
\]
est exact.
\end{defi}
\begin{exem}\textup{Dans la situation usuelle o{\`u} $X$ est un espace topologique, un faisceau d'ensembles au sens usuel sur $X$ est un faisceau d'ensembles sur le site $\EuScript{O}uv_{X}$ (\textit{cf.} exemple 1.1.2).}
\end{exem}

Nous donnons maintenant quelques exemples de faisceaux que nous utiliserons dans les applications. Il s'agit de faisceaux sur le site {\'e}tale d'un sch{\'e}ma $S$; pour plus de d{\'e}tails concernant les propri{\'e}t{\'e}s de ces faisceaux, nous renvoyons {\`a} \cite{Mi} ou \cite{Tam}.
\begin{exem}\textup{Le faisceau $\mathbb{G}_{a,S}$:\label{Gadd} c'est le faisceau dont le groupe des sections au-dessus d'un ouvert {\'e}tale $\left(S'\rightarrow S\right)$, not{\'e} simplement $\mathbb{G}_{a,S}\left(S'\right)$ est:}
	\[\mathbb{G}_{a,S}\left(S'\right)=\Gamma\left(S',\mathcal{O}_{S'}\right)=\mathcal{O}_{S'}\left(S'\right)
\]

\textup{Le faisceau $\mathbb{G}_{a,S}$ est appel{\'e} le \textit{groupe additif}\index{groupe!additif d'un sch{\'e}ma} de $S$.}
\end{exem}
\begin{exem}\textup{Le faisceau $\mathbb{G}_{m,S}$\label{Gmult} est le faisceau dont le groupe des sections au-dessus d'un ouvert {\'e}tale $\left(S'\rightarrow S\right)$, not{\'e} $\mathbb{G}_{m,S}\left(S'\right)$ est cette fois constitu{\'e} des fonctions inversibles d{\'e}finies globalement sur $S'$:}
	\[\mathbb{G}_{m,S}\left(S'\right)=\Gamma\left(S',\mathcal{O}^{\ast}_{S'}\right)=\mathcal{O}^{\ast}_{S'}\left(S'\right)
\]

\textup{Le faisceau $\mathbb{G}_{m,S}$ est appel{\'e} le \textit{groupe multiplicatif}\index{groupe!multiplicatif d'un sch{\'e}ma} de $S$.}
\end{exem}
\begin{exem}\textup{Le faisceau $\mu_{n,S}$\label{mun} est le faisceau des racines $n$-i{\`e}mes de l'unit{\'e}\index{faisceau!des racines $n$-i{\`e}mes de l'unit{\'e}}. Les sections de celui-ci au-dessus d'un ouvert {\'e}tale $\left(S'\rightarrow S\right)$ sont donn{\'e}es par:}
	\[\mu_{n,S}\left(S'\right)=\left\{f\in \mathcal{O}^{\ast}_{S'}\left(S'\right)/ f^{n}=1 \right\}
\]
\end{exem}
\begin{exem}\textup{On peut associer {\`a} tout $S$-sch{\'e}ma $Y$ un faisceau, not{\'e} temporairement $\mathcal{F}_{Y}$, de la mani{\`e}re suivante: pour tout ouvert {\'e}tale $\left(S'\rightarrow S\right)$, on prend comme sections de $\mathcal{F}_{Y}$:}
	\[\mathcal{F}_{Y}\left(S'\right)=\textup{Hom}_{S}\left(S',Y\right)
\]
\textup{Un faisceau $\mathcal{F}$ sur $S_{\acute{e}t}$ est dit \textit{repr{\'e}sentable}\index{faisceau!repr{\'e}sentable} s'il existe un $S$-sch{\'e}ma $Y$ tel que $\mathcal{F}=\mathcal{F}_{Y}$. Les trois pr{\'e}c{\'e}dents exemples de faisceaux sont justement repr{\'e}sentables:\vspace{2mm}
\begin{itemize}
\item le faisceau $\mathbb{G}_{a,S}$ est repr{\'e}sent{\'e} par le sch{\'e}ma:
	\[\textup{Spec}\;\mathbb{Z}\left[T\right]\times_{\textup{Spec}\;\mathbb{Z}}S
\]
\item le faisceau $\mathbb{G}_{m,S}$ est repr{\'e}sent{\'e} par le sch{\'e}ma:
	\[\textup{Spec}\;\mathbb{Z}\left[T,T^{-1}\right]\times_{\textup{Spec}\;\mathbb{Z}}S
\]
\item et enfin, le faisceau $\mu_{n,S}$ est repr{\'e}sent{\'e} par le sch{\'e}ma:
	\[\textup{Spec}\;\left(\frac{\mathbb{Z}\left[T\right]}{\left(T^{n}-1\right)}\right)\times_{\textup{Spec}\;\mathbb{Z}}S
\]
\end{itemize}}
\vspace{1mm}

\textup{Dans nos applications, nous consid{\`e}rerons un $k$-groupe alg{\'e}brique $G$, et un $k$-sch{\'e}ma $X$. Nous noterons $G_{X}$ le produit fibr{\'e} de $X$ et de $G$ au-dessus de $\textup{Spec}\;k$:}
	\[\xymatrix{G_{X}\ar[r]\ar[d]&G\ar[d]\\X\ar[r]_{\pi\ \ }&\textup{Spec}\;k}
\]
\textup{D'apr{\`e}s ce qui pr{\'e}c{\`e}de, le $k$-groupe alg{\'e}brique $G$ (\textit{resp.} le $X$-sch{\'e}ma en groupes $G_{X}$) repr{\'e}sente un faisceau sur le site {\'e}tale de $k$ (\textit{resp.} de $X$). Nous identifierons souvent $G$ et $G_{X}$ avec les faisceaux qu'ils repr{\'e}sentent.}
\end{exem}

On dispose {\'e}videmment d'une notion de morphisme de (pr{\'e})faisceaux au-dessus d'un site $\mathcal{S}$ donn{\'e}. Il est non moins {\'e}vident que les faisceaux d'ensembles sur $\mathcal{S}$ et les morphismes entre iceux constituent une cat{\'e}gorie, que l'on note $\widetilde{\mathcal{S}}$.

Nous concluons cette premi{\`e}re section en introduisant la notion de topos \cite{SGA4-4}:
\begin{defi} On appelle \textbf{topos}\index{topos} une cat{\'e}gorie $\mathcal{T}$ telle qu'il existe un site $\mathcal{S}$ telle que $\mathcal{T}$ soit {\'e}quivalente {\`a} la cat{\'e}gorie $\widetilde{\mathcal{S}}$ des faisceaux d'ensembles sur $\mathcal{S}$.
\end{defi}

Le principal topos auquel nous nous int{\'e}resserons dans ce travail est le \textit{topos {\'e}tale}\index{topos!etale@{\'e}tale d'un sch{\'e}ma} d'un sch{\'e}ma $S$, not{\'e} $\widetilde{S}_{\acute{e}t}$\label{Settopos}: c'est le topos des faisceaux d'ensembles sur le site {\'e}tale de $S$. Nous renvoyons {\`a} \cite{SGA4-7} pour une {\'e}tude approfondie des topoï {\'e}tales des sch{\'e}mas. Faisons juste une remarque \textquotedblleft{fonctorielle}\textquotedblright: soit:
	\[f:Y\longrightarrow S
\]
un morphisme de sch{\'e}mas. On a vu dans l'exemple 1.1.7 que $f$ induit un morphisme de sites:
	\[f_{\acute{e}t}:Y_{\acute{e}t}\longrightarrow S_{\acute{e}t}
\]
correspondant au foncteur:
	\[f^{-1}:\left(S'\rightarrow S\right)\rightsquigarrow \left(Y\times_{S}S'\rightarrow Y\right)
\]
qui associe {\`a} un ouvert {\'e}tale de $S$ l'ouvert {\'e}tale de $Y$ obtenu par changement de base. Alors $f$ induit {\'e}galement un morphisme entre les topos {\'e}tales de $Y$ et de $S$:
	\[\widetilde{f}_{\acute{e}t}:\widetilde{Y}_{\acute{e}t}\longrightarrow \widetilde{S}_{\acute{e}t}
\]
Le morphisme $\widetilde{f}_{\acute{e}t}$ associe {\`a} tout faisceau d'ensembles sur le site {\'e}tale de $Y$ le faisceau sur le site {\'e}tale de $S$ obtenu par image directe gr{\^a}ce {\`a} $f$.
\begin{notat}\textup{Dans nos applications, le sch{\'e}ma de base $S$ sera souvent $X$ ou $\textup{Spec}\;k$. Pour all{\'e}ger un peu la terminologie, nous dirons simplement, lorsque le contexte est clair, \textquotedblleft{faisceau sur $X$ (\textit{resp.} sur $k$)}\textquotedblright\ au lieu de \textquotedblleft{faisceau d'ensembles sur le site {\'e}tale de $X$ (\textit{resp.} de $\textup{Spec}\;k$)}\textquotedblright. D'ailleurs, plus g{\'e}n{\'e}ralement nous parlerons de pr{\'e}champ, champ, gerbe\ldots\ sur un sch{\'e}ma $S$ (\textit{resp.} sur un corps $k$) pour d{\'e}signer un pr{\'e}champ, un champ, une gerbe\ldots\ sur le site {\'e}tale de $S$ (\textit{resp.} sur le site {\'e}tale de $\textup{Spec}\;k$).}
\end{notat}
\end{section}
\begin{section}{Torseurs}
Nous nous limitons aux torseurs sur le site {\'e}tale d'un sch{\'e}ma et nous renvoyons au chapitre III de \cite{Gi2} pour l'{\'e}tude des torseurs sur un site g{\'e}n{\'e}ral.
\begin{defi}Soient $S$ un sch{\'e}ma et $G_{S}$ un sch{\'e}ma en groupes sur $S$. Un \textbf{$G_{S}$-torseur}\index{torseur} {\`a} droite sur $S$ est un faisceau $P$ sur $S$, muni d'une action {\`a} droite du faisceau de groupes $G_{S}$, tel que:

il existe un recouvrement {\'e}tale $\left(S_{i}\rightarrow S\right)_{i\in I}$ tel que l'ensemble $P\left(S_{i}\right)$ soit principal homog{\`e}ne sous l'action du groupe $G_{S}\left(S_{i}\right)$, pour tout $i\in I$.

Un \textbf{morphisme}\index{morphisme!de torseurs} $\varphi:P\rightarrow P'$ de $G_{S}$-torseurs est un morphisme de faisceaux $G_{S}$-{\'e}quivariant. Etant donn{\'e} un $G_{S}$-torseur $P$, on note:\label{adG}
	\[\textup{ad}_{G_{S}}\left(P\right)
\]
le faisceau des automorphismes du $G_{S}$-torseur $P$.

Enfin on note:
	\[\textup{EHP}\left(G_{S}/S\right)
\]
l'ensemble des classes d'isomorphie de $G_{S}$-torseurs sur $S$.
\end{defi}
\begin{rem}\textup{Rappelons que d'apr{\`e}s nos conventions et notations:}

\textup{\begin{itemize}
\item $P\left(S_{i}\right)$ (\textit{resp.} $G_{S}\left(S_{i}\right)$) d{\'e}signe l'ensemble (\textit{resp.} le groupe) des sections du faisceau $P$ (\textit{resp.} $G_{S}$) au-dessus de l'ouvert {\'e}tale $\left(S_{i}\rightarrow S\right)$;\vspace{2mm}
\item le $S$-sch{\'e}ma en groupes $G_{S}$ est identifi{\'e} au faisceau de groupes qu'il repr{\'e}sente;\vspace{2mm}
\item la locution \textquotedblleft{$P$ est un faisceau sur $S$}\textquotedblright\ est mise pour \textquotedblleft{$P$ est un faisceau d'ensembles sur le site {\'e}tale de $S$}\textquotedblright.
\end{itemize}}
\end{rem}
\begin{rem}\textup{On peut {\'e}videmment d{\'e}finir de la m{\^e}me mani{\`e}re les $G_{S}$-torseurs {\`a} gauche. Notons que par la suite, nous appellerons simplement $G_{S}$-torseur un $G_{S}$-torseur {\`a} droite. Les $G_{S}$-torseurs sur $S$ et leurs morphismes constituent une cat{\'e}gorie, que l'on note $Tors\left(S,G_{S}\right)$\label{Torscat}.}
\end{rem}

Apr{\`e}s cette d{\'e}finition \textquotedblleft{faisceautique}\textquotedblright\ de torseur, notons que l'on dispose, moyennant des hypoth{\`e}ses tr{\`e}s raisonnables sur le sch{\'e}ma $G_{S}$, d'une interpr{\'e}tation g{\'e}om{\'e}trique de ces objets. Plus pr{\'e}cis{\'e}ment:
\begin{pro} Soit $S$ un sch{\'e}ma, et $G_{S}$ un sch{\'e}ma en groupes affine sur $S$. Alors tout $G_{S}$-torseur sur $S$ est repr{\'e}sentable par un sch{\'e}ma.
\end{pro}

\uline{\textsc{Preuve}}: c'est le (a) du th{\'e}or{\`e}me III.4.3 de \cite{Mi}.
\begin{flushright}
$\Box$
\end{flushright}

En particulier, dans la situation de notre probl{\`e}me central (\textit{cf.} page XI), les $G_{X}$-torseurs sont repr{\'e}sentables. En effet, du fait que $G$ est un $k$-groupe alg{\'e}brique lin{\'e}aire, le morphisme structural:
	\[G\longrightarrow \textup{Spec}\;k
\]
est affine. Il s'ensuit que le morphisme:
	\[G_{X}\longrightarrow X
\]
l'est aussi, le caract{\`e}re affine {\'e}tant stable par changement de base quelconque (\textit{cf.} \cite{EGA1} I.9.1.16).
\vspace{1mm}
\newline

Donnons maintenant quelques exemples de torseurs:
\begin{exem}\textup{Etant donn{\'e}s un sch{\'e}ma $S$ et un $S$-sch{\'e}ma en groupes $G_{S}$, on appelle \textit{$G_{S}$-torseur trivial}\index{torseur!trivial} et on note $G_{S,d}$ le $G_{S}$-torseur obtenu en faisant op{\'e}rer $G_{S}$ sur lui-m{\^e}me {\`a} droite par translations. Par d{\'e}finition m{\^e}me, tout $G_{S}$-torseur est isomorphe {\`a} $G_{S,d}$, localement pour la topologie {\'e}tale sur $S$. En continuant d'enfoncer des portes ouvertes, il s'ensuit que deux $G_{S}$-torseurs sont toujours localement isomorphes pour la topologie {\'e}tale.}
\end{exem}
\begin{exem}\textup{Avec les m{\^e}mes notations que ci-dessus, on peut associer naturellement {\`a} tout $G_{S}$-torseur ({\`a} droite) $P$ un $\textup{ad}_{G_{S}}\left(P\right)$-torseur {\`a} gauche, le faisceau des automorphismes ($G_{S}$-{\'e}quivariants) de $P$ op{\'e}rant {\`a} gauche sur $P$ de mani{\`e}re {\'e}vidente. Nous reviendrons sur cette remarque dans la section concernant les bitorseurs.}
\end{exem}
\begin{exem}\textup{Les classes d'isomorphie de $\mathbb{G}_{m,S}$-torseurs (\textit{resp.} de $GL_{n,S}$-torseurs) sur $S$ coïncident avec les classes d'isomorphie de fibr{\'e}s en droites (\textit{resp.} de fibr{\'e}s vectoriels de rang $n$) sur $S$.}
\end{exem}
\begin{exem}\textup{Soient $n$ un entier $\geq2$ et $k$ un corps quelconque. On peut associer au groupe $\mathbb{Z}/n\mathbb{Z}$ un faisceau de groupes sur la droite projective $\mathbb{P}^{1}_{k}$ (faisceau constant). Un $\mathbb{Z}/n\mathbb{Z}$-torseur sur $\mathbb{P}^{1}_{k}$ est alors un rev{\^e}tement {\'e}tale de la droite projective de groupe $\mathbb{Z}/n\mathbb{Z}$.}
\end{exem}

En vue de nos applications, il est {\'e}videmment indispensable de disposer d'un moyen de calculer pratiquement l'ensemble $\textup{EHP}\left(G_{S}/S\right)$. Ce moyen est fourni par l'{\'e}nonc{\'e} ci-dessous, qui est un cas particulier du corollaire III.4.7 de \cite{Mi}:
\begin{pro} Soient $k$ un corps de caract{\'e}ristique nulle, $G$ un $k$-groupe alg{\'e}brique lin{\'e}aire, et $X$ un $k$-sch{\'e}ma. Alors les ensembles $\textup{EHP}\left(G_{X}/X\right)$ et $H^{1}_{\acute{e}t}\left(X,G_{X}\right)$ sont en bijection.
\end{pro}

L'exemple ci-dessous est une application de cette proposition aux alg{\`e}bres simples centrales et aux vari{\'e}t{\'e}s de Severi-Brauer.
\begin{exem}\textup{Soit $S$ un sch{\'e}ma. Rappelons que l'on appelle alg{\`e}bre d'Azumaya sur $S$ un faisceau de $\mathcal{O}_S$-alg{\`e}bres $\mathcal{A}$, pour lequel il existe un recouvrement {\'e}tale $\left(S_i\rightarrow S\right)_{i\in I}$ de $S$ tel que pour tout $i\in I$, il existe un entier $n_i$ tel que:
	\[\mathcal{A}\otimes_{\mathcal{O}_S}\mathcal{O}_{S_i}\approx{M_{n_i}\left(\mathcal{O}_{S_i}\right)}
\]
Il revient au m{\^e}me (d'apr{\`e}s le thm. 5.1 p.57 de \cite{G4}) de dire que $\mathcal{A}$ est une $\mathcal{O}_{S}$-alg{\`e}bre localement libre telle que:
\begin{enumerate}[(i)]
\item pour tout point $s$ de $S$, $\mathcal{A}\otimes_{\mathcal{O}_{S,s}}k\left(s\right)$ est une alg{\`e}bre simple centrale;
\item le morphisme naturel $\mathcal{A}\otimes_{\mathcal{O}_S}\mathcal{A}^{opp}\rightarrow{\textup{End}_{\mathcal{O}_{S}}\left(\mathcal{A}\right)}$ est un isomorphisme de $\mathcal{O}_{S}$-alg{\`e}bres.
\end{enumerate}}

\textup{Nous appellerons alg{\`e}bre d'Azumaya\index{alg{\`e}bre!d'Azumaya} d'indice\index{indice!d'une alg{\`e}bre d'Azumaya} $n$ une alg{\`e}bre d'Azumaya pour lesquels tous les $n_i$ de la d{\'e}finition ci-dessus sont {\'e}gaux {\`a} $n$. Qu'une alg{\`e}bre d'Azumaya d'indice $n$ sur un sch{\'e}ma $S$ soit un $PGL_{n,S}$-torseur sur $S$ d{\'e}coule de l'{\'e}nonc{\'e} suivant, qui g{\'e}n{\'e}ralise aux sch{\'e}mas le lemme de Skolem-Noether:
\begin{theo}[Auslander-Goldman, \cite{G4}] Soit $\mathcal{A}$ une alg{\`e}bre d'Azumaya sur $S$, $u$ un automorphisme de $\mathcal{A}$. Alors, localement pour la topologie {\'e}tale, $u$ est int{\'e}rieur, \textit{i.e.} de la forme:
	\[u\left(s\right)=asa^{-1}
\]
o{\`u} $a$ est une section inversible de $\mathcal{A}$, d{\'e}termin{\'e}e d'ailleurs de fa\c{c}on unique modulo multiplication par une section de $\mathcal{O}_{S}^{\ast}$. De mani{\`e}re {\'e}quivalente, le sch{\'e}ma des automorphismes de l'alg{\`e}bre associative $M_n\left(\mathcal{O}_{S}\right)$ est canoniquement isomorphe au groupe projectif $PGL_{n,S}$.
\end{theo}
}

\textup{Dans le cas particulier o{\`u} $S=\textup{Spec}\;k$, on retrouve la d{\'e}finition d'alg{\`e}bre simple centrale: une alg{\`e}bre d'Azumaya d'indice $n$ sur $k$ est une alg{\`e}bre simple centrale\index{alg{\`e}bre!simple centrale}\index{indice!d'une alg{\`e}bre simple centrale} d'indice $n$, soit une $k$-forme de l'alg{\`e}bre de matrices $M_{n}$. D'apr{\`e}s la proposition 1.2.8, l'ensemble des classes d'isomorphie de $k$-alg{\`e}bres simples centrales est en bijection avec l'ensemble $H^{1}_{\acute{e}t}\left(k,PGL_{n}\right)$. Remarquons alors que les multiples facettes de $PGL_{n}$ fournissent des interpr{\'e}tations diff{\'e}rentes de cet ensemble. Explicitement, comme:}\index{vari{\'e}t{\'e}!de Severi-Brauer}
	\[PGL_{n}=\textup{Aut}\;M_{n}=\textup{Int}\;SL_{n}=\textup{Aut}\;\mathbb{P}^{n-1}
\]
\textup{on en d{\'e}duit une correspondance entre les $k$-formes de $M_{n}$, celles de $\mathbb{P}^{n-1}$ et les $k$-formes int{\'e}rieures de $SL_{n}$:}
	\[\xymatrix{\left\{\textup{$k$-alg{\`e}bres simples centrales d'indice }n\right\}\ar@{<->}[d]\\\left\{\textup{$k$-formes int{\'e}rieures de $SL_{n}$}\right\}\ar@{<->}[d]\\\left\{\textup{Vari{\'e}t{\'e}s de Severi-Brauer de dimension }n-1\right\}}
\]
\end{exem}
\begin{exem}\textup{Suivant la remarque 2.1 de \cite{SGA4-7}, lorsque $G$ est un groupe commutatif ordinaire (\textit{i.e.} $G$ est un faisceau constant), on {\'e}crira simplement $H^{1}\left(X,G\right)$ au lieu de $H^{1}\left(X,G_{X}\right)$. Si par exemple $G$ est ab{\'e}lien fini, alors $H^{1}\left(X,G\right)$ est l'ensemble des classes d'isomorphie de $G$-rev{\^e}tements galoisiens, pour lesquels les probl{\`e}mes de descente ont {\'e}t{\'e} largement {\'e}tudi{\'e}s dans \cite{DD1} et \cite{DD2}. Lorsque $X$ est connexe et muni d'un point g{\'e}om{\'e}trique $x$, on a l'isomorphisme canonique:}
	\[H^{1}\left(X,G\right)=\textup{Hom}\left(\Pi_{1}\left(X,x\right),G\right)
\]
\end{exem}
\begin{exem}\textup{Pour faire le lien avec la terminologie de Springer (\cite{Sp}), disons encore que si $G$ est un $k$-groupe alg{\'e}brique, les $G$-torseurs sur $k$ sont {\'e}videmment des $k$-espaces homog{\`e}nes avec isotropie triviale.}
\end{exem}

Le r{\'e}sultat suivant, qui concerne les faisceaux d'automorphismes des torseurs, est assez {\'e}vident mais d'une importance capitale pour comprendre la diff{\'e}rence entre les chapitres 2 et 3 (\textit{i.e.} entre les situations ab{\'e}lienne et non-ab{\'e}lienne):
\begin{lemf}\label{lemf} Soient $S$ un sch{\'e}ma, $G_{S}$ un sch{\'e}ma en groupes sur $S$, et $P$ un $G_{S}$-torseur sur $S$. Alors $\textup{ad}_{G_{S}}\left(P\right)$ est une $S$-forme int{\'e}rieure de $G_{S}$; autrement dit:
	\[\textup{ad}_{G_{S}}\left(P\right)\textit{ repr{\'e}sente une classe de } H^{1}\left(S,\textup{Int}\;G_{S}\right).
\]

En particulier, si $G_{S}$ est ab{\'e}lien, alors: 
	\[\textup{ad}_{G_{S}}\left(P\right)\approx G_{S}
\]
\end{lemf}

\uline{\textsc{Preuve}}: par d{\'e}finition, il existe un recouvrement {\'e}tale $\left(S_{i}\rightarrow S\right)$ trivialisant $P$. On choisit alors une section $p_{i}$ de $P_{\left|S_{i}\right.}$, pour tout $i\in I$.

Pour tout $\left(i,j\right)\in I\times I$, on note $S_{ij}$ le sch{\'e}ma $S_{i}\times_{S}S_{j}$. Du fait que l'action de ${G_{S}}_{\left|S_{ij}\right.}$ sur $P_{\left|S_{ij}\right.}$ est simplement transitive\footnote{Cet abus de langage signifie en fait que le groupe $G_{S}\left(S'\right)$ op{\`e}re simplement transitivement ({\`a} droite, d'apr{\`e}s nos conventions sur les torseurs) sur l'ensemble $P\left(S'\right)$, pour tout $S_{ij}$-sch{\'e}ma {\'e}tale $S'$.}, il existe $\gamma_{ij}\in {G_{S}}\left(S_{ij}\right)$ tel que:
	\[{p_{j}}_{\left|S_{ij}\right.}={p_{i}}_{\left|S_{ij}\right.}.\gamma_{ij}
\]
(la famille $\left(\gamma_{ij}\right)_{i,j}$ est justement un 1-cocycle repr{\'e}sentant la classe de $P$ dans $H^{1}\left(S,G_{S}\right)$). Afin de ne pas surcharger les notations, nous abandonnons {\`a} partir de maintenant les indices \textquotedblleft{$_{\left|S_{ij}\right.}$}\textquotedblright\ . Nous r{\'e}{\'e}crivons donc l'{\'e}galit{\'e} pr{\'e}c{\'e}dente:
	\[p_{j}=p_{i}.\gamma_{ij}\ \ \ \ \left(E1\right)
\]

Soit maintenant $f$ un automorphisme de $P$. Comme $f$ est $G_{S}$-{\'e}quivariante:
	\[\exists\;g_{i}\in G_{S}\left(S_{i}\right)\ \textup{tel que}:\ f\left(p_{i}\right)=p_{i}.g_{i},\ \forall\;i\in I.
\]

Pour tout $\left(i,j\right)\in I\times I$, on obtient, en utilisant la relation $\left(E1\right)$:
	\[f\left(p_{j}\right)=p_{j}.g_{j}=p_{i}.\gamma_{ij}\;g_{j}\ \ \ \ \left(E2\right)
\]

D'autre part en utilisant la relation $\left(E2\right)$ et la $G_{S}$-{\'e}quivariance de $f$, on a:
	\[f\left(p_{j}\right)=f\left(p_{i}.\gamma_{ij}\right)=f\left(p_{i}\right).\gamma_{ij}=p_{i}.g_{i}\;\gamma_{ij}\ \ \ \left(E3\right)
\]

En comparant les relations $\left(E2\right)$ et $\left(E3\right)$ et en utilisant une nouvelle fois la simple transitivit{\'e} de l'action, on obtient:
	\[g_{i}\;\gamma_{ij}=\gamma_{ij}\;g_{j}
\]
soit finalement:
	\[g_{i}=\gamma_{ij}\;g_{j}\;\gamma_{ij}^{-1}=\textup{conj}\left(\gamma_{ij}\right)\left(g_{j}\right)
\]

Notons tout de suite qu'il est imm{\'e}diat que la famille $\left(\textup{conj}\left(\gamma_{ij}\right)\right)_{i,j}$ est un 1-cocycle {\`a} valeurs dans $\textup{Int}\;G_{S}$, puisque $\left(\gamma_{ij}\right)_{i,j}$ est un 1-cocycle {\`a} valeurs dans $G_{S}$. 

De plus, la relation ci-dessus traduit le fait que le faisceau $\textup{ad}_{G_{S}}\left(P\right)$ est la donn{\'e}e des faisceaux locaux ${G_{S}}_{\left|S_{i}\right.}$ recoll{\'e}s au-dessus des $S_{ij}$ par les automorphismes int{\'e}rieurs $\textup{conj}\left(\gamma_{ij}\right)$. Autrement dit, c'est une $S$-forme int{\'e}rieure de $G_{S}$.
\begin{flushright}
$\Box$
\end{flushright}

Dans la suite, nous aurons besoin de conna{\^i}tre l'influence d'un morphisme de sch{\'e}mas (plus pr{\'e}cis{\'e}ment l'influence du morphisme structural $\pi:X\rightarrow \textup{Spec}\;k$) sur les torseurs; c'est l'int{\'e}r{\^e}t de l'{\'e}nonc{\'e} suivant:
\begin{lem} Soient $k$ un corps de caract{\'e}ristique nulle, $X$ un $k$-sch{\'e}ma, $\pi:X\rightarrow \textup{Spec}\;k$ le morphisme structural, $G$ un $k$-groupe alg{\'e}brique lin{\'e}aire, et $P$ (\textit{resp.} $Q$) un $G_{X}$-torseur (\textit{resp.} un $G$-torseur) sur $X$ (\textit{resp.} sur $\textup{Spec}\;k$). Alors:
\begin{enumerate}[(i)]
\item le faisceau image directe $\pi_{\ast}P$ est un $\pi_{\ast}G_{X}$-pseudo-torseur sur $\textup{Spec}\;k$;
\item si $\bar{G}_{X}\left(\bar{X}\right)=\bar{G}\left(\bar{k}\right)$, alors $\pi_{\ast}P$ est un $G$-pseudo-torseur sur $\textup{Spec}\;k$;
\item le faisceau image inverse $\pi^{\ast}Q$ est un $\pi^{\ast}G$-torseur sur $X$.
\end{enumerate}
\end{lem}

\uline{\textsc{Preuve}}: pour le point (i), il suffit de prouver que pour toute extension {\'e}tale $L$ de $k$, l'ensemble $\pi_{\ast}P\left(\textup{Spec}\;L\rightarrow\textup{Spec}\;k\right)$ est vide ou principal homog{\`e}ne sous l'action du groupe $\pi_{\ast}G_{X}\left(\textup{Spec}\;L\rightarrow\textup{Spec}\;k\right)$. Ceci provient du fait que l'ensemble:
	\[\pi_{\ast}P\left(\textup{Spec}\;L\rightarrow\textup{Spec}\;k\right)=P\left(X_{L}\rightarrow X\right)=\textup{Hom}_{X}\left(X_{L},P\right)
\]
est vide ou principal homog{\`e}ne sous le groupe:
	\[\pi_{\ast}G_{X}\left(\textup{Spec}\;L\rightarrow\textup{Spec}\;k\right)=G_{X}\left(X_{L}\rightarrow X\right)=\textup{Hom}_{X}\left(X_{L},G_{X}\right)
\]
car $P$ est un $G_{X}$-torseur sur $X$.

Pour {\'e}tablir le point (ii), il suffit de comparer les faisceaux (sur $\left(\textup{Spec}\;k\right)_{\acute{e}t}$) $\pi_{\ast}G_{X}$ et $G$; pour ce faire, il suffit d'{\'e}tudier leurs fibres $\left(\pi_{\ast}G_{X}\right)_{\bar{k}}$ et $G_{\bar{k}}$. Par d{\'e}finition:
	\[\left(\pi_{\ast}G_{X}\right)_{\bar{k}}=\lim_{\xrightarrow[L]{}}\left[G_{X}\left(X_{L}\rightarrow X\right)\right]=\textup{Hom}_{X}\left(\bar{X},G_{X}\right)
\]
la limite directe {\'e}tant prise sur les extensions {\'e}tales $L$ de $k$. De la m{\^e}me fa\c{c}on:
	\[\left(G\right)_{\bar{k}}=\lim_{\xrightarrow[L]{}}\left[G\left(\textup{Spec}\;L\rightarrow \textup{Spec}\;k\right)\right]=\textup{Hom}_{k}\left(\textup{Spec}\;\bar{k},G\right)
\]

Par cons{\'e}quent, les faisceaux $\pi_{\ast}G_{X}$ et $G$ sont isomorphes si et seulement si:
	\[\textup{Hom}_{X}\left(\bar{X},G_{X}\right)=\textup{Hom}_{k}\left(\textup{Spec}\;\bar{k},G\right)
\]

Cette derni{\`e}re condition est {\'e}quivalente {\`a}:
	\[\left[\textup{Hom}_{\bar{X}}\left(\bar{X},\bar{G}_{X}\right)=\right]\ \ \bar{G}_{X}\left(\bar{X}\right)=\bar{G}\left(\bar{k}\right)\ \ \left[=\textup{Hom}_{\bar{k}}\left(\textup{Spec}\;\bar{k},\bar{G}\right)\right]
\]

En effet, les ensembles $\textup{Hom}_{X}\left(\bar{X},G_{X}\right)$ et $\textup{Hom}_{\bar{X}}\left(\bar{X},\bar{G}_{X}\right)$ sont naturellement en bijection; {\`a} tout {\'e}l{\'e}ment $\varphi\in \textup{Hom}_{\bar{X}}\left(\bar{X},\bar{G}_{X}\right)$ on associe l'{\'e}l{\'e}ment $p_{G}\circ\varphi$ de $\textup{Hom}_{X}\left(\bar{X},G_{X}\right)$:
	\[\xymatrix@C=50pt@R=40pt{\bar{G}_{X} \ar[r]^{p_{G}} \ar[d] &G_{X} \ar[d] \\ \bar{X} \ar@/^20pt/[u]^{\varphi} \ar[r] \ar@{-->}[ur]^{p_{G}\circ\varphi} & X}
\]

R{\'e}ciproquement, on peut associer {\`a} tout $f\in \textup{Hom}_{X}\left(\bar{X},G_{X}\right)$ le morphisme (dont l'existence est assur{\'e}e par la propri{\'e}t{\'e} universelle du produit fibr{\'e}) $\varphi\in \textup{Hom}_{\bar{X}}\left(\bar{X},\bar{G}_{X}\right)$ du diagramme ci-dessous:
	\[\xymatrix@C=25pt@R=20pt{\bar{X} \ar[drrr]^{f} \ar[dddr]_{id} \ar@{-->}[dr]^{\varphi} \\ & \bar{G}_{X} \ar[rr]^{p_{G}} \ar[dd] && G_{X} \ar[dd] \\ \\ & \bar{X} \ar[rr] && X}
\]
\pagebreak

Le point (iii) est {\'e}vident, {\`a} partir du moment o{\`u} l'on a remarqu{\'e} que comme $G$ est lin{\'e}aire, le $G$-torseur $Q$ est repr{\'e}sentable par un sch{\'e}ma, que l'on note encore $Q$, et le faisceau $\pi^{\ast}Q$ est alors repr{\'e}sent{\'e} par le $X$-sch{\'e}ma $\pi^{\ast}Q$ d'apr{\`e}s la proposition II.3.1.3 de \cite{Tam}.
\begin{flushright}
$\Box$
\end{flushright}
\begin{rem}\textup{Nous aurons l'occasion de revenir en d{\'e}tail sur l'utilit{\'e} de la condition $\bar{G}_{X}\left(\bar{X}\right)=\bar{G}\left(\bar{k}\right)$ dans les probl{\`e}mes de descente.}
\end{rem}
\begin{rem}\textup{Avec les notations du lemme pr{\'e}c{\'e}dent, il n'est en g{\'e}n{\'e}ral pas vrai que le faisceau image inverse $\pi^{\ast}G$ est repr{\'e}sent{\'e} par le sch{\'e}ma en groupes $G_{X}$. C'est cependant le cas si $G$ est un groupe fini, d'apr{\`e}s la remarque 3.1.(d) page 69 de \cite{Mi}.}
\end{rem}
\end{section}
\begin{section}{Bitorseurs}
Lorsque l'on veut {\'e}tudier les $G_{S}$-torseurs sur un sch{\'e}ma $S$, il est essentiel de conna{\^i}tre leurs faisceaux d'automorphismes. Dans le cas ab{\'e}lien, d'apr{\`e}s le lemme 1.2.13, il n'y a rien {\`a} faire. En revanche dans le cas g{\'e}n{\'e}ral, d'apr{\`e}s le m{\^e}me lemme, ces faisceaux sont des formes int{\'e}rieures de $G_{S}$, et il convient donc de tenir compte de cette information suppl{\'e}mentaire. C'est ce qui motive l'introduction des bitorseurs\footnote{C'est {\'e}galement ce qui justifie le fait qu'il est d{\'e}sesp{\'e}r{\'e} d'arriver {\`a} une g{\'e}n{\'e}ralisation vraiment \textquotedblleft{id{\'e}ale}\textquotedblright\ de la suite {\`a} 5 termes dans le cas non-ab{\'e}lien avec les \textquotedblleft{$H^{1}$}\textquotedblright\ usuels; puisque par exemple l'ensemble $H^{1}_{\acute{e}t}\left(S,G_{S}\right)$ correspond certes aux classes d'isomorphie de $G_S$-torseurs sur $S$, mais oublie d{\'e}lib{\'e}r{\'e}ment la structure {\`a} gauche de ces objets.}. Comme d'habitude, nous nous restreignons ici aux bitorseurs sur un sch{\'e}ma, et nous renvoyons {\`a} \cite{B1} pour la d{\'e}finition de bitorseur sur un topos en g{\'e}n{\'e}ral.
\begin{defi} Soient $S$ un sch{\'e}ma, $G_{S}$ et $H_{S}$ deux sch{\'e}mas en groupes sur $S$. On appelle \textbf{$\left(H_{S},G_{S}\right)$-bitorseur\index{bitorseur} sur $S$} un faisceau d'ensembles $B$ sur le site {\'e}tale de $S$ muni d'une action {\`a} gauche (\textit{resp.} {\`a} droite) de $H_{S}$ (\textit{resp.} de $G_{S}$), ces actions commutant entre elles, tel que $B$ soit un $H_{S}$-torseur {\`a} gauche et un $G_{S}$-torseur {\`a} droite.

Lorsque $G_{S}=H_{S}$, nous dirons simplement $G_{S}$-bitorseur pour d{\'e}signer un $\left(G_{S},G_{S}\right)$-bitorseur.

Soient $B_{1}$ et $B_{2}$ deux $\left(H_{S},G_{S}\right)$-bitorseurs sur $S$. Un \textbf{morphisme}\index{morphisme!de bitorseurs} de $\left(H_{S},G_{S}\right)$-bitorseurs est morphisme de faisceaux d'ensembles $\varphi:B_{1}\rightarrow B_{2}$ qui est $\left(H_{S},G_{S}\right)$-{\'e}quivariant.
\end{defi}

Les $\left(H_{S},G_{S}\right)$-bitorseurs (\textit{resp.} les $G_{S}$-bitorseurs) sur $S$ et leurs morphismes constituent une cat{\'e}gorie, not{\'e}e:\label{BitorsSHG}
	\[Bitors\left(S;H_{S},G_{S}\right)\ \ \ \left(resp.\ Bitors\left(S;G_{S}\right)\right)
\]
\begin{exem}
\textup{Le \textit{$G_{S}$-bitorseur trivial}\index{bitorseur!trivial}, not{\'e} $G_{S,bitriv}$, est obtenu en faisant op{\'e}rer $G_{S}$ sur lui-m{\^e}me par translations {\`a} gauche et {\`a} droite.}
\end{exem}
\begin{exem}\textup{Soit $u$ un automorphisme de $G_{S}$. On lui associe le $G_{S}$-bitorseur not{\'e} $\left(G_{S},u\right)$, d{\'e}fini en faisant op{\'e}rer $G_{S}$ sur lui-m{\^e}me:}
\vspace{2mm}
\textup{\begin{itemize}
\item {\`a} droite par translations;\vspace{2mm}
\item {\`a} gauche en posant: $g\bullet_{gauche}g'=u\left(g\right).g'$
\end{itemize}}
\end{exem}
\begin{exem} \textup{A tout $G_{S}$-bitorseur, on associe naturellement un $G_{S}$-torseur\index{torseur!associ{\'e} {\`a} un bitorseur} en \textquotedblleft{oubliant}\textquotedblright\ l'action de $G_{S}$ {\`a} gauche. R{\'e}ciproquement, {\`a} tout $G_{S}$-torseur $P$, on associe un $\left(\textup{ad}_{G_{S}}\left(P\right),G_{S}\right)$-bitorseur\index{bitorseur!associ{\'e} {\`a} un torseur} en faisant op{\'e}rer {\`a} gauche $\textup{ad}_{G_{S}}\left(P\right)$ sur $P$ de fa\c{c}on {\'e}vidente.}
\end{exem}
\begin{exem} \textup{A tout $\left(H_{S},G_{S}\right)$-bitorseur $B$ on peut associer un $\left(G_{S},H_{S}\right)$-bitorseur not{\'e} $B^{opp}$ et appel{\'e} \textit{bitorseur oppos{\'e} {\`a} $B$}\index{bitorseur!oppos{\'e}}, la nouvelle action {\`a} gauche (resp. {\`a} droite) de $G_{S}$ (resp. de $H_{S}$) {\'e}tant d{\'e}finie en posant:
	\[g\boxdot b=b\bullet g^{-1}\ \ \ \left(resp.\ b\boxcircle h=h\ast b\right),\ \forall\;g\in G_{S},\ \forall\;h\in H_{S},\ \forall\;b\in B;
\]
o{\`u} $\bullet$ (resp. $\ast$) d{\'e}signe l'action {\`a} droite (resp. {\`a} gauche) de $G_{S}$ (resp. de $H_{S}$) donn{\'e}e par la structure de $\left(H_{S},G_{S}\right)$-bitorseur de $B$.}
\end{exem}

\begin{rem} \textup{Pour classifier les bitorseurs, on doit introduire la cohomologie {\`a} valeurs dans les modules crois{\'e}s. On appelle \textit{module crois{\'e}} ({\`a} gauche)\index{module crois{\'e}} sur $S$ la donn{\'e}e d'un morphisme de sch{\'e}mas en groupes sur $S$:
	\[\varphi:G_{S}\longrightarrow H_{S}
\]
et d'une action {\`a} gauche de $H_{S}$ sur $G_{S}$ tels que:}
\vspace{1mm}
\textup{\begin{enumerate}[(i)]
\item $\varphi\left(^{h}g\right)=h.\varphi\left(g\right).h^{-1},\ \forall\;g\in G_{S},\ \forall\;h\in H_{S};$
\item $^{\varphi\left(g\right)}g'=g.g'.g^{-1},\ \forall\;g,g'\in G_{S};$
\end{enumerate}}
\vspace{2mm}

\textup{Par exemple, pour tout sch{\'e}ma en groupes $G_{S}$ sur $S$, le morphisme {\'e}vident:
	\[conj:G_{S}\longrightarrow \textup{Aut}\;G_{S}
\]
et l'action naturelle de $\textup{Aut}\;G_{S}$ sur $G_{S}$ donnent lieu {\`a} un module crois{\'e} sur $S$. De plus, dans \cite{B1} ou \cite{B2}, Breen d{\'e}finit la cohomologie {\`a} valeurs dans ce module crois{\'e} $\left(G_{S}\rightarrow \textup{Aut}\;G_{S}\right)$. Sans rentrer dans les d{\'e}tails de cette construction\footnote{Que l'on peut d'ailleurs rapprocher de la cohomologie {\`a} valeurs dans un syst{\`e}me de coefficients introduite dans \cite{Do1}.}, signalons juste que:}

\textup{\begin{itemize}
\item $H^{0}\left(S,G_{S}\rightarrow \textup{Aut}\;G_{S}\right)$ correspond aux classes d'isomorphie de $G_{S}$-bitorseurs sur $S$, et\vspace{2mm}
\item $H^{1}\left(S,G_{S}\rightarrow \textup{Aut}\;G_{S}\right)$ correspond aux classes d'{\'e}quivalence de gerbes sur $S$ localement li{\'e}es par $G_{S}$. En anticipant un peu, cet ensemble diff{\`e}re du $H^{2}$ de Giraud, puisque deux gerbes {\'e}quivalentes (et non lien $G_{S}$-{\'e}quivalentes) donnent la m{\^e}me classe dans ce $H^{1}$, et pas n{\'e}cessairement dans le $H^{2}$ de Giraud.
\end{itemize}}
\end{rem}

La suite exacte (\textit{cf.} \cite{B1}) suivante contient l'essentiel des propri{\'e}t{\'e}s des bitorseurs dont nous aurons besoin dans la suite de notre propos:
\vspace{2mm}
\begin{flushleft}
$\xymatrix@C=12pt{&H^{0}\left(\textup{Int}\;G_{S}\right)\ar[dr]\\H^{0}\left(G_{S}\right)\ar[ur] \ar[rr]^{conj.} && H^{0}\left(\textup{Aut}\;G_{S}\right)\ar[r]^{b\ \ \ \ } & H^{0}\left(G_{S}\rightarrow \textup{Aut}\;G_{S}\right)\ar@{-}[r]^{\ \ \ \ \ \ \ \ \ \ AP}& \cdots}$
\end{flushleft}
\begin{flushright}
$\xymatrix@C=12pt{&&H^{1}\left(\textup{Int}\;G_{S}\right)\ar[dr]\\\cdots\ar[r]&H^{1}\left(G_{S}\right)\ar[ur] \ar[rr]^{\textup{ad}_{G_{S}}\ } && H^{1}\left(\textup{Aut}\;G_{S}\right)\ar[dr]^{\lambda}\\&&&&H^{1}\left(\textup{Out}\;G_{S}\right)}$
\end{flushright}
\vspace{1mm}
o{\`u} l'on a not{\'e} $H^{i}\left(\bullet\right)$ l'ensemble de cohomologie $H^{i}\left(S,\bullet\right)$ et:\vspace{2mm}
\vspace{1mm}
\begin{itemize}
\item $conj.:H^{0}\left(G_{S}\right)\longrightarrow H^{0}\left(\textup{Aut}\;G_{S}\right)$ est le morphisme {\'e}vident;\vspace{2mm}
\item pour tout $u\in H^{0}\left(\textup{Aut}\;G_{S}\right)$, $b\left(u\right)$ est le bitorseur $\left(G_{S},u\right)$ de l'exemple 1.3.3;\vspace{2mm}
\item $AP$ est le foncteur \textquotedblleft{amn{\'e}sie partielle}\textquotedblright\ qui associe {\`a} tout $G_{S}$-bitorseur $B$ le $G_{S}$-torseur {\`a} droite sous-jacent;\vspace{2mm}
\item pour tout $G_{S}$-torseur $P$, $\textup{ad}_{G_{S}}\left(P\right)$ est comme d'habitude le faisceau des $G_{S}$- automorphismes de $P$ (dans la preuve du lemme 1.2.13, nous avons explicitement d{\'e}crit ce morphisme);\vspace{2mm}
\item enfin l'application:
	\[\lambda:H^{1}\left(\textup{Aut}\;G_{S}\right)\longrightarrow H^{1}\left(\textup{Out}\;G_{S}\right)
\]
est celle qui associe {\`a} une forme de $G_{S}$ le lien qu'elle repr{\'e}sente; remarquons l'image par $\lambda$ d'une forme $G'_{S}$ de $G_{S}$ est {\'e}gale {\`a} la classe privil{\'e}gi{\'e}e de $H^{1}\left(\textup{Out}\;G_{S}\right)$ si et seulement si $G'_{S}$ est une $S$-forme \textit{int{\'e}rieure} de $G_{S}$.
\end{itemize}
\vspace{2mm}

En particulier, si $G_{S}$ est un sch{\'e}ma en groupes ab{\'e}liens, on a la suite exacte:
	\[0\longrightarrow H^{0}\left(\textup{Aut}\;G_{S}\right)\longrightarrow H^{0}\left(G_{S}\rightarrow \textup{Aut}\;G_{S}\right)\longrightarrow H^{1}\left(G_{S}\right)\longrightarrow 0
\]

Dans ce cas, les $G_{S}$-bitorseurs constituent une extension des $G_{S}$-torseurs par $\textup{Aut}\ G_{S}$.

Pour clore cette section, notons que l'on dispose sur les bitorseurs une loi de composition partiellement d{\'e}finie, gr{\^a}ce au produit contract{\'e}. Plus explicitement:
\begin{defi} Soient $S$ un sch{\'e}ma, et $G_{S}$, $H_{S}$ et $L_{S}$ des sch{\'e}mas en groupes sur $S$. Soient encore $B$ (resp. $C$) un $\left(H_{S},G_{S}\right)$-bitorseur (resp. un $\left(G_{S},L_{S}\right)$-bitorseur). On appelle \textbf{produit contract{\'e} de $B$ et de $C$}\index{produit contract{\'e}}, et on note:
	\[B\wedge^{G_{S}}C
\]
le faisceau $B\times C$ quotient{\'e} par la relation:
	\[\left(b.g,c\right)=\left(b,g.c\right)
\]

D'apr{\`e}s \cite{B1}, l'action {\`a} gauche (resp. {\`a} droite) de $H_{S}$ (resp. de $L_{S}$) fait de $B\wedge^{G_{S}}C$ un $\left(H_{S},L_{S}\right)$-bitorseur.
\end{defi}

Le bitorseur trivial est l'{\'e}l{\'e}ment neutre pour ce produit, et le produit contract{\'e} d'un bitorseur $B$ avec son oppos{\'e} $B^{opp}$ est isomorphe au bitorseur trivial.
\end{section}
\begin{section}{Pr{\'e}champs et champs}
Il est bien connu (par exemple gr{\^a}ce {\`a} \cite{LMB}) que les champs alg{\'e}briques fournissent une g{\'e}n{\'e}ralisation de la notion de sch{\'e}ma. Ils apparaissent naturellement dans des probl{\`e}mes de modules, o{\`u} le foncteur que l'on {\'e}tudie n'est pas repr{\'e}sentable par un sch{\'e}ma, du fait de la pr{\'e}sence d'automorphismes (voir par exemple l'{\'e}tude des courbes de genre $g\geq2$ sur un sch{\'e}ma $S$ \cite{DM}). En clair, on ne peut obtenir un espace des modules fin qu'en \textquotedblleft{rigidifiant}\textquotedblright\ la situation, par l'interm{\'e}diaire de structures suppl{\'e}mentaires (pour {\^e}tre un peu plus concret, cela se fait \textit{via} les structures de niveau pour les vari{\'e}t{\'e}s ab{\'e}liennes). De la m{\^e}me mani{\`e}re, il est raisonnable de voir les champs en g{\'e}n{\'e}ral (sur le site {\'e}tale d'un sch{\'e}ma $S$) comme une g{\'e}n{\'e}ralisation des faisceaux d'ensembles sur $S$. 

La d{\'e}finition de champ sur un sch{\'e}ma $S$ n{\'e}cessite plusieurs {\'e}tapes. Premi{\`e}rement, on appelle \textit{cat{\'e}gorie fibr{\'e}e sur $S$}\index{cat{\'e}gorie fibr{\'e}e} la donn{\'e}e:
\begin{enumerate}[(i)]
\item pour tout ouvert {\'e}tale $\left(S_{1}\rightarrow S\right)$ d'une cat{\'e}gorie, not{\'e}e $\mathcal{G}\left(S_{1}\right)$ et appel{\'e}e \textit{cat{\'e}gorie fibre} de $\mathcal{G}$ au-dessus de l'ouvert {\'e}tale $\left(S_{1}\rightarrow S\right)$;
\item pour toute inclusion\footnote{Pour faire le parall{\`e}le avec la d{\'e}finition usuelle de pr{\'e}faisceau, nous appelons ici \textquotedblleft{inclusion}\textquotedblright\ un morphisme entre ouverts {\'e}tales de $S$; cet abus est justifi{\'e} par le fait qu'une cat{\'e}gorie fibr{\'e}e est grossi{\`e}rement un pr{\'e}faisceau en cat{\'e}gories, comme il est d'ailleurs indiqu{\'e} dans \cite{B3}. Un tel morphisme n'est cependant pas un monomorphisme en g{\'e}n{\'e}ral.}:
	\[\xymatrix@C=5pt{S_{2}\ar[rrrr]\ar[rrd]&&&&S_{1}\ar[lld]\\&&S}
\]
entre ouverts {\'e}tales de $S$, d'un foncteur (restriction {\`a} $S_{2}$):
	\[\xymatrix@R=5pt@C=20pt{\rho_{S_{1}S_{2}}:&\mathcal{G}\left(S_{1}\right)\ar[r]&\mathcal{G}\left(S_{2}\right)\\& g_{1} \ar@{|->}[r]&{g_{1}}_{\left|S_{2}\right.}}
\]
\item pour toute double inclusion:
	\[\xymatrix@C=5pt{S_{3}\ar[rrr]\ar[rrrd]&&&S_{2}\ar[d]\ar[rrr]&&&S_{1}\ar[llld]\\&&&S}
\]
d'une transformation naturelle:
	\[\tau_{S_{3}S_{2}S_{1}}:\rho_{S_{1}S_{3}}\Longrightarrow \rho_{S_{2}S_{3}}\circ\rho_{S_{1}S_{2}}
\]
de telle sorte que, pour toute triple inclusion:
	\[\xymatrix@C=5pt{S_{4}\ar[rr]\ar[rrrd]&&S_{3}\ar[rr]\ar[rd]&&S_{2}\ar[ld]\ar[rr]&&S_{1}\ar[llld]\\&&&S}
\]
les transformations naturelles obtenues par composition:
	\[\rho_{S_{1}S_{4}}\Longrightarrow \rho_{S_{3}S_{4}}\circ \rho_{S_{1}S_{3}}\Longrightarrow \rho_{S_{3}S_{4}}\circ\left(\rho_{S_{2}S_{3}}\circ \rho_{S_{1}S_{2}}\right)
\]
et
	\[\rho_{S_{1}S_{4}}\Longrightarrow \rho_{S_{2}S_{4}}\circ \rho_{S_{1}S_{2}}\Longrightarrow \left(\rho_{S_{3}S_{4}}\circ \rho_{S_{2}S_{3}}\right)\circ\rho_{S_{1}S_{2}}
\]
coïncident.
\end{enumerate}

$\mathcal{G}$ est appel{\'e}e une \textit{cat{\'e}gorie fibr{\'e}e en groupoïdes} sur $S$, si la cat{\'e}gorie fibre $\mathcal{G}\left(S'\right)$ au-dessus de tout ouvert {\'e}tale $\left(S'\rightarrow S\right)$ est un groupoïde\footnote{Un groupoïde est une cat{\'e}gorie o{\`u} toute fl{\`e}che est inversible.}. Le foncteur naturel:
	\[\mathcal{G}\longrightarrow S_{\acute{e}t}
\]
est appel{\'e} \textit{projection} ou \textit{morphisme structural}\index{morphisme!structural} de $\mathcal{G}$. Commen\c{c}ons par un exemple naïf:
\begin{exem}\textup{Soit $\pi:X\rightarrow \textup{Spec}\;k$ un $k$-sch{\'e}ma. On peut associer {\`a} $X$ une cat{\'e}gorie fibr{\'e}e sur $k$,\footnote{\textit{I.e.} une cat{\'e}gorie fibr{\'e}e sur le site {\'e}tale de $\textup{Spec}\;k$.} not{\'e}e provisoirement $\mathcal{G}_{X}$, en prenant pour toute extension {\'e}tale $L$ de $k$ le groupoïde discret dont l'ensemble d'objets est:
	\[\mathcal{G}_{X}\left(L\right)=\textup{Hom}_{\textup{Spec}\;k}\left(\textup{Spec}\;L,X\right)
\]
Autrement dit, la cat{\'e}gorie fibre de $\mathcal{G}_{X}$ au-dessus de $L$ a pour objets les points {\`a} valeurs dans $\textup{Spec}\;L$ du $k$-sch{\'e}ma $X$, et les seuls morphismes sont les identit{\'e}s des objets (en particulier $\mathcal{G}_{X}$ est une cat{\'e}gorie fibr{\'e}e en groupoïdes!).}
\end{exem}

Nous donnerons un peu plus loin des exemples plus int{\'e}ressants (heureusement!) de cat{\'e}gories fibr{\'e}es. La totalit{\'e} des cat{\'e}gories fibr{\'e}es que nous consid{\`e}rerons seront des cat{\'e}gories fibr{\'e}es en groupoïdes.

Un \textit{morphisme $F:\mathcal{G}\rightarrow\mathcal{G'}$}\index{morphisme!de cat{\'e}gories fibr{\'e}es} entre deux cat{\'e}gories fibr{\'e}es sur $S$ est la donn{\'e}e pour tout ouvert {\'e}tale $\left(S_{1}\rightarrow S\right)$ d'un foncteur $F_{S_{1}}:\mathcal{G}\left(S_{1}\right)\rightarrow\mathcal{G}'\left(S_{1}\right)$ naturellement compatible avec les restrictions, dans le sens o{\`u} pour toute inclusion $\left(S_{2}\rightarrow S_{1}\right)$ le carr{\'e} suivant:
\[\xymatrix@C=50pt@R=50pt{\mathcal{G}\left(S_{1}\right)\ar[d]_{\rho^{\mathcal{G}}_{S_{1}S_{2}}}\ar[r]^{F_{S_{1}}}&\mathcal{G}'\left(S_{1}\right)\ar[d]^{\rho_{S_{1}S_{2}}^{\mathcal{G'}}}\\\mathcal{G}\left(S_{2}\right)\ar[r]_{F_{S_{2}}}&\mathcal{G}'\left(S_{2}\right)}
\]
commute, {\`a} un isomorphisme de foncteurs pr{\`e}s; \textit{i.e.} il existe une transformation naturelle:
	\[\theta_{S_{1}S_{2}}:\rho_{S_{1}S_{2}}^{\mathcal{G}'}\circ F_{S_{1}}\Longrightarrow F_{S_{2}}\circ \rho_{S_{1}S_{2}}^{\mathcal{G}}
\]
cette transformation naturelle satisfaisant {\`a} sont tour des conditions de compatibilit{\'e} avec les transformations naturelles \textquotedblleft{$\tau_{S_{3}S_{2}S_{1}}$}\textquotedblright. Un tel foncteur est aussi appel{\'e} \textit{foncteur cart{\'e}sien}\index{foncteur cart{\'e}sien} entre les cat{\'e}gories $\mathcal{G}$ et $\mathcal{G}'$.

Un \textit{pr{\'e}champ de groupoïdes sur $S$}\index{pr{\'e}champ} (ou simplement pr{\'e}champ sur $S$) est une cat{\'e}gorie fibr{\'e}e en groupoïdes $\mathcal{G}$ sur $S$ o{\`u} les isomorphismes se recollent: explicitement, cela signifie qu'{\'e}tant donn{\'e}s un ouvert {\'e}tale $\left(S'\rightarrow S\right)$, $x$ et $y$ deux objets de $\mathcal{G}\left(S'\right)$, $\left(S_{i}\rightarrow S'\right)_{i\in I}$ un recouvrement {\'e}tale de $S'$, la suite d'ensembles:
	\[\xymatrix{\textup{Isom}_{\mathcal{G}\left(S'\right)}\left(x,y\right)\ar[r]& \prod_{i\in I}\textup{Isom}_{\mathcal{G}\left(S_{i}\right)}\left(x_{\left|S_{i}\right.},y_{\left|S_{i}\right.}\right)\ar@<0.5ex>[r]\ar@<-0.5ex>[r]& \prod_{i,j\in I}\textup{Isom}_{\mathcal{G}\left(S_{ij}\right)}\left(x_{\left|S_{ij}\right.},y_{\left|S_{ij}\right.}\right)}
\]
est exacte (on a not{\'e} $S_{ij}=S_{i}\times_{S'}S_{j}$). Plus succinctement, il revient au m{\^e}me de dire que $\textup{Isom}\left(x,y\right)$ est un faisceau sur le site {\'e}tale de $S'$.

Un \textit{champ de groupoïdes sur $S$}\index{champ} (ou juste champ sur $S$, ou encore $S$-champ) est un pr{\'e}champ de groupoïdes $\mathcal{G}$ sur $S$ o{\`u} toute donn{\'e}e de descente sur les objets est effective\index{donn{\'e}e!de descente!effective}, ce qui signifie ceci: soient:
\vspace{1mm}
\begin{itemize}
\item $\left(S'\rightarrow S\right)$ un ouvert {\'e}tale de $S$;\vspace{2mm}
\item $\left(S_{i}\rightarrow S'\right)_{i\in I}$ un recouvrement {\'e}tale de $S'$;\vspace{2mm}
\item $\left(g_{i}\right)_{i\in I}$ une famille de sections de $\mathcal{G}$, plus pr{\'e}cis{\'e}ment:
	\[g_{i}\in Ob\left(\mathcal{G}\left(S_{i}\right)\right),\ \forall\ i\in I;
\]\vspace{2mm}
\item pour tout couple $\left(i,j\right)\in I\times I$, un isomorphisme:
	\[\varphi_{ij}:{g_{j}}_{\left|S_{ij}\right.}\longrightarrow {g_{i}}_{\left|S_{ij}\right.}
\]
de telle sorte que:
	\[\varphi_{ik}=\varphi_{ij}\circ \varphi_{jk},\ \forall\;\left(i,j,k\right)\in I\times I\times I
\]
cette derni{\`e}re {\'e}galit{\'e} ayant lieu dans $\textup{Isom}\left(\mathcal{G}\left(S_{ijk}\right)\right)$.
\end{itemize}
\vspace{2mm}

Alors il existe (effectivement) un objet $g'\in \mathcal{G}\left(S'\right)$ et des isomorphismes:
	\[\eta_{i}:g'_{\left|S_{i}\right.}\longrightarrow g_{i}
\]
compatibles avec les isomorphismes de recollement, dans le sens o{\`u}:
	\[\varphi_{ij}\circ\eta_{j}=\eta_{i}\;\ \textup{sur }S_{ij},\ \forall\;\left(i,j\right)\in I\times I.
\]

Grossi{\`e}rement, un champ est donc une cat{\'e}gorie fibr{\'e}e dans laquelle les morphismes \textit{et} les objets se recollent. Nous en donnons maintenant quelques exemples.
\begin{exem}[Le champ $\textup{Tors}\left(k,G\right)$\label{ChTorskG} des $G$-torseurs sur un corps $k$]\index{champ!des torseurs} \textup{Soient $k$ un corps et $G$ un $k$-groupe alg{\'e}brique. Pour tout ouvert {\'e}tale $\left(\textup{Spec}\;L\rightarrow \textup{Spec}\;k\right)$, les $G_{L}$-torseurs sur le site {\'e}tale de $L$ sont les objets de la cat{\'e}gorie $Tors\left(L,G_{L}\right)$; cette cat{\'e}gorie est d'ailleurs un groupoïde, puisque tout morphisme de torseurs est un isomorphisme, du fait de la simple transitivit{\'e} de l'action. En outre, si $\left(\textup{Spec}\;L'\rightarrow \textup{Spec}\;k\right)$ est un ouvert {\'e}tale de $k$ inclus dans le premier, \textit{i.e.} si le diagramme ci-dessous est commutatif:}
	\[\xymatrix{\textup{Spec}\;L'\ar[rr]\ar[rd]&&\textup{Spec}\;L\ar[ld]\\&\textup{Spec}\;k}
\]
\textup{alors on a un foncteur:}
	\[\rho_{LL'}:\ Tors\left(L,G_{L}\right)\longrightarrow Tors\left(L',G_{L'}\right)
\]
\textup{qui associe {\`a} un $G_{L}$-torseur $P\rightarrow \textup{Spec}\;L$ le torseur $P_{L'}\rightarrow\textup{Spec}\;L'$ obtenu par changement de base; le carr{\'e} suivant est donc cart{\'e}sien:}
	\[\xymatrix@C=15pt@R=15pt{P_{L'}\ar[dd]\ar[rr]&&P\ar[dd]\\&\Box \\\textup{Spec}\;L'\ar[rr]&&\textup{Spec}\;L}
\]
\textup{Par ailleurs, le foncteur $\rho_{LL'}$ associe {\`a} un morphisme $f:P\rightarrow Q$ de $G_{L}$-torseurs sur $L$ le morphisme $f_{L'}:P_{L'}\rightarrow Q_{L'}$ obtenu lui aussi par changement de base; le diagramme ci-dessous est donc commutatif:}
	\[\xymatrix@C=10pt@R=15pt{&&Q_{L'}\ar[rrrr]\ar[dddl]
	&&&&Q\ar[dddl]\\P_{L'}\ar[rru]^{f_{L'}} \ar[rrrr]|{\ \ \ } \ar[ddr]&&&&P\ar[ddr]\ar[urr]^{f}\\\\&\textup{Spec}\;L'\ar[rrrr]&&&&\textup{Spec}\;L}
\]

\textup{La collection des groupoïdes $Tors\left(L,G_{L}\right)$ (pour $L$ parcourant les extensions {\'e}tales de $k$) et les foncteurs $\rho_{LL'}$ constituent une cat{\'e}gorie fibr{\'e}e en groupoïdes sur $k$. Cette cat{\'e}gorie fibr{\'e}e est not{\'e}e $\textup{Tors}\left(k,G\right)$: par d{\'e}finition, le groupoïde fibre $\textup{Tors}\left(k,G\right)\left(L\right)$ de cette cat{\'e}gorie fibr{\'e}e au-dessus d'un ouvert {\'e}tale $\left(\textup{Spec}\;L\rightarrow \textup{Spec}\;k\right)$ est justement le groupoïde $Tors\left(L,G_{L}\right)$. La cat{\'e}gorie fibr{\'e}e $\textup{Tors}\left(k,G\right)$ est un pr{\'e}champ puisque les morphismes se recollent, et c'est un champ car les torseurs (qui sont des faisceaux) se recollent (toujours par la th{\'e}orie classique de la descente).}

\textup{Si l'on consid{\`e}re le cas particulier o{\`u} $G=PGL_{n}$, on obtient le champ $\textup{Tors}\left(k,PGL_{n}\right)$, et suivant que l'on interpr{\`e}te $PGL_{n}$ comme le groupe d'automorphismes de $M_{n}$ ou de $\mathbb{P}^{n-1}$, on obtient le $k$-champ $\textup{Asc}\left(k,n\right)$\label{Asckn} des $k$-alg{\`e}bres simples centrales d'indice $n$, ou celui des vari{\'e}t{\'e}s des $k$-vari{\'e}t{\'e}s de Severi-Brauer\index{vari{\'e}t{\'e}!de Severi-Brauer} de dimension $n-1$, not{\'e} $\textup{SB}\left(k,n-1\right)$\label{SBkn-1}. En particulier, les champs $\textup{Asc}\left(k,n\right)$ et $\textup{SB}\left(k,n-1\right)$ sont isomorphes, et l'isomorphisme de champs n'est autre que le foncteur qui permet d'associer {\`a} une alg{\`e}bre simple centrale $A$ la vari{\'e}t{\'e} de Severi-Brauer $X_{A}$ correspondante, ce foncteur {\'e}tant d{\'e}crit dans le chapitre 5 de \cite{Jah}, ou encore dans \cite{J} III.3.5.}
\end{exem}
\begin{exem}[Le champ $\textup{Tors}\left(S,G_{S}\right)$\label{TorsSG} des $G_{S}$-torseurs sur un sch{\'e}ma $S$] \textup{Soient $S$ un sch{\'e}ma et $G_{S}$ un $S$-sch{\'e}ma en groupes. On d{\'e}finit de la m{\^e}me mani{\`e}re que pr{\'e}c{\'e}demment le champ $\textup{Tors}\left(S,G_{S}\right)$ des $G_{S}$-torseurs sur $S$. Le groupoïde fibre de ce champ au-dessus d'un ouvert {\'e}tale $\left(S'\rightarrow S\right)$ a pour objets les $G_{S'}$-torseurs\footnote{On a not{\'e} $G_{S'}$ le sch{\'e}ma en groupes $G_{S}\times_{S}S'$.} sur $S'$, et pour fl{\`e}ches les isomorphismes de $G_{S'}$-torseurs sur $S'$. Lorsque l'on prend $G_{S}=\mathbb{G}_{m,S}$ (\textit{resp.} $G_{S}=GL_{n,S}$, \textit{resp.} $G_{S}=PGL_{n,S}$), on obtient le champ $\textup{LBun}\left(S\right)$ des fibr{\'e}s en droites sur $S$ (\textit{resp.} le champ $\textup{VBun}\left(n,S\right)$ des fibr{\'e}s vectoriels de rang $n$ sur $S$, \textit{resp.} le champ $\textup{Az}\left(n,S\right)$ des alg{\`e}bres d'Azumaya d'indice $n$ sur $S$).}
\end{exem}
\begin{exem}[Champ associ{\'e} {\`a} un sch{\'e}ma]\index{champ!associ{\'e} {\`a} un sch{\'e}ma} \textup{Soient $S$ un sch{\'e}ma et $Y$ un $S$-sch{\'e}ma. On peut associer un $S$-champ {\`a} $Y$, que l'on note encore $Y$; le groupoïde fibre de ce champ $Y$ au-dessus d'un ouvert {\'e}tale $\left(S'\rightarrow S\right)$ le groupoïde discret\footnote{\textit{I.e.} les seules fl{\`e}ches sont les identit{\'e}s.} dont l'ensemble d'objets est: $\textup{Hom}_{S}\left(S',Y\right)$. Un $S$-champ $\mathfrak{X}$ est dit \textit{repr{\'e}sentable} s'il existe un $S$-sch{\'e}ma $Y$ tel que le $S$-champ associ{\'e} {\`a} $Y$ par le proc{\'e}d{\'e} d{\'e}crit ci-dessus et le $S$-champ $\mathfrak{X}$ sont isomorphes. La pr{\'e}sence d'automorphismes non-triviaux est donc clairement un obstacle {\`a} ce qu'un $S$-champ soit repr{\'e}sentable par un sch{\'e}ma.}

\textup{Cette correspondance entre sch{\'e}mas et champs nous permet d'interpr{\'e}ter le foncteur projection d'un champ comme un morphisme de champs. Consid{\'e}rons par exemple le champ des $G$-torseurs sur un corps $k$. Le foncteur projection relatif {\`a} ce champ:}
	\[p:\textup{Tors}\left(k,G\right)\longrightarrow \left(\textup{Spec}\;k\right)_{\acute{e}t}
\]
\textup{est celui qui associe {\`a} un torseur l'extension de $k$ au-dessus de laquelle il est d{\'e}fini. Or, d'apr{\`e}s ce qui pr{\'e}c{\`e}de, on peut associer au $k$-sch{\'e}ma $\textup{Spec}\;k$ (!) un $k$-champ. En adoptant ce point de vue, le foncteur $p$ devient un morphisme de champs.}

\textup{Enfin, la construction du $S$-champ associ{\'e} {\`a} un $S$-sch{\'e}ma donne lieu {\`a} un foncteur pleinement fid{\`e}le:
	\[\left(\textup{Sch/S}\right)\longrightarrow \left(\textup{Champ/S}\right)
\]
de la cat{\'e}gorie des $S$-sch{\'e}mas dans celle des $S$-champs, qui se factorise d'ailleurs par la cat{\'e}gorie des $S$-champs alg{\'e}briques (\textit{cf.} le chapitre 4 de \cite{LMB}).}
\end{exem}
\begin{exem}[Le champ $\mathcal{M}_{g}$ des courbes stables de genre g ($g\geq{2}$)]\textup{Soit $S$ un sch{\'e}ma; rappelons que l'on appelle courbe stable (suivant \cite{DM}) de genre $g$ sur $S$ un morphisme 
	\[\pi:C\longrightarrow{S}
\]
propre, plat, dont les fibres g{\'e}om{\'e}triques sont des sch{\'e}mas $C_{S}$ de dimension 1, r{\'e}duits, connexes, et tels que:\vspace{2mm}
\begin{itemize}
\item $C_{S}$ n'a que des points doubles ordinaires;\vspace{2mm}
\item si $E$ est une composante rationnelle non-singuli{\`e}re de $C_{S}$, alors $E$ rencontre les autres composantes de $C_{S}$ en au moins trois points.
\end{itemize}}

\textup{On note $\left(\textup{SchEt}\right)$ le site dont la cat{\'e}gorie sous-jacente est celle des sch{\'e}mas, munie de la topologie {\'e}tale. Pour tout sch{\'e}ma $S$, on d{\'e}finit un groupoïde que l'on note $\mathcal{M}_{g,S}$ en posant:\vspace{2mm}
\begin{itemize}
\item [\textit{\uline{Objets de $\mathcal{M}_{g,S}$}:}] courbes stables de genre g sur $S$;\vspace{2mm}
\item [\textit{\uline{Isomorphismes de $\mathcal{M}_{g,S}$}:}] $S$-isomorphismes de sch{\'e}mas.\vspace{3mm}
\end{itemize}}

\textup{La collection des $\mathcal{M}_{g,S}$ constitue une cat{\'e}gorie fibr{\'e}e $\mathcal{M}_{g}$ sur $\left(\textup{SchEt}\right)$, et $\mathcal{M}_{g}$ est un champ de groupoïdes sur $\left(\textup{SchEt}\right)$ (\textit{cf.} \cite{Do5}). En outre, ce champ est isomorphe au quotient de $H_{g}$ par $PGL\left(5g-6\right)$, o{\`u} $H_{g}$ est le sous-sch{\'e}ma de $\textup{Hilb}_{5g-6}^{P_{g}}$ des courbes stables tricanoniquement plong{\'e}es dans $\mathbb{P}^{5g-6}$, o{\`u} $P_{g}\left(n\right)=\left(6n-1\right)\left(g-1\right)$ est le polyn{\^o}me de Hilbert (\textit{cf.} \cite{Do5} p.127-128).}
\end{exem}
\begin{exem}[Champ associ{\'e} {\`a} un pr{\'e}champ]\index{champ!associ{\'e} {\`a} un pr{\'e}champ} \textup{Soit $S$ un sch{\'e}ma. Pour construire le faisceau (sur le site {\'e}tale de $S$) associ{\'e} {\`a} un pr{\'e}faisceau sur $S$, il suffit de rendre locale la d{\'e}finition de section, de telle sorte que le recollement des sections devient possible. De la m{\^e}me fa\c{c}on, on construit le $S$-champ associ{\'e} {\`a} un $S$-pr{\'e}champ en rendant locale la d{\'e}finition d'objet, de mani{\`e}re {\`a} forcer l'effectivit{\'e} des donn{\'e}es de descente. Ainsi, on obtient un foncteur: 
	\[\left(\textup{Pr{\'e}champ/S}\right)\longrightarrow \left(\textup{Champ/S}\right)
\]
et tout $S$-champ appartient localement {\`a} l'image essentielle de celui-ci \cite{B3}. Plus pr{\'e}cis{\'e}ment, si $\mathcal{P}$ est un pr{\'e}champ sur $S$, et si $\left(S'\rightarrow S\right)$ est un ouvert {\'e}tale, une section du champ $\mathcal{P}^{+}$ associ{\'e} {\`a} $\mathcal{P}$ au-dessus de cet ouvert est la donn{\'e}e:}
\textup{\begin{enumerate}[(i)]
\item d'un recouvrement {\'e}tale $\left(S'_{i}\rightarrow S'\right)_{i\in I}$;
\item d'une famille de sections $\left(p_{i}\right)_{i\in I}$ du pr{\'e}champ $\mathcal{P}$ relativement {\`a} ce recouvrement: pr{\'e}cis{\'e}ment, on demande que:
	\[p_{i}\in\textup{Ob}\left(\mathcal{P}\left(S'_{i}\right)\right),\ \forall\;i\in I;
\]
\item (\textit{donn{\'e}e de recollement})\index{donn{\'e}e!de recollement} d'un isomorphisme de $\mathcal{P}\left(S'_{i}\times_{S'}S'_{j}\right)$:
	\[\varphi_{ij}:{p_{j}}_{\left|S'_{i}\times_{S'}S'_{j}\right.}\longrightarrow {p_{i}}_{\left|S'_{i}\times_{S'}S'_{j}\right.},\ \forall\;\left(i,j\right)\in I\times I;
\]
\item (\textit{donn{\'e}e de descente})\index{donn{\'e}e!de descente} les isomorphismes $\varphi_{ij}$ satisfaisant la condition de 1-cocycle\footnote{C'est une mani{\`e}re rapide de dire que:
	\[{\varphi_{ij}}_{\left|S'_{i}\times_{S'}S'_{j}\times_{S'}S'_{k}\right.}\circ{\varphi_{jk}}_{\left|S'_{i}\times_{S'}S'_{j}\times_{S'}S'_{k}\right.}={\varphi_{ik}}_{\left|S'_{i}\times_{S'}S'_{j}\times_{S'}S'_{k}\right.},\ \forall\;\left(i,j,k\right)\in I\times I\times I,
\]
cette {\'e}galit{\'e} ayant lieu dans $\textup{Isom}\left(\mathcal{P}\left(S'_{i}\times_{S'}S'_{j}\times_{S'}S'_{k}\right)\right)$.}:
	\[\varphi_{ij}\circ\varphi_{jk}=\varphi_{ik},\ \forall\;\left(i,j,k\right)\in I\times I\times I.
\]
\end{enumerate}}

\textup{De fait, il existe une mani{\`e}re plus intrins{\`e}que de d{\'e}finir le champ associ{\'e} {\`a} un pr{\'e}champ (nous renvoyons {\`a} la d{\'e}finition II.2.1.1 de \cite{Gi2} pour celle-ci.)}
\end{exem}
\end{section}
\begin{section}{Gerbes} 
\begin{defi}Soit $S$ un sch{\'e}ma. On appelle \textbf{gerbe sur $S$}\index{gerbe} (ou $S$-gerbe) un champ de groupoïdes $\mathcal{G}$ sur le site {\'e}tale de $S$, localement non-vide et localement connexe, soit:\vspace{2mm}
\begin{itemize}
\item ($\mathcal{G}$ est localement non-vide). Il existe un recouvrement {\'e}tale $\left(S_{i}\rightarrow S\right)_{i\in I}$ tel que:
	\[\textup{Ob}\left(\mathcal{G}\left(S_{i}\right)\right)\neq\emptyset,\ \forall\ i\in I.
\]
\item ($\mathcal{G}$ est localement connexe). Soient $\left(S'\rightarrow S\right)$ un ouvert {\'e}tale, et soient $x$ et $y$ deux objets de $\mathcal{G}\left(S'\right)$. Alors il existe un recouvrement {\'e}tale $\left(S'_{j}\rightarrow S'\right)_{j\in J}$ tel que:
	\[x_{\left|S'_{j}\right.}\approx y_{\left|S'_{j}\right.}
\]
cet isomorphisme vivant dans le groupoïde $\mathcal{G}\left(S'_{j}\right)$, $\forall\ j\in J$.\vspace{3mm}
\end{itemize}

Une $S$-gerbe est dite \textbf{neutre}\index{gerbe!neutre} si elle a une section au-dessus de $S$, \textit{i.e.} si le groupoïde fibre $\mathcal{G}\left(S\right)$ est non-vide.

Un \textbf{morphisme de gerbes}\index{morphisme!de gerbes} sur $S$ est un morphisme de champs dont la source et le but sont des $S$-gerbes.
\end{defi}
\begin{exem}[Gerbes de torseurs]\index{gerbe!de torseurs} \textup{Soient $S$ un sch{\'e}ma, et $G_{S}$ un $S$-sch{\'e}ma en groupes. Le champ $\textup{Tors}\left(S,G_{S}\right)$ de l'exemple 1.4.3 est une gerbe. En effet:}

\textup{\begin{itemize}
\item le champ $\textup{Tors}\left(S,G_{S}\right)$ est localement non-vide, puisqu'il l'est m{\^e}me globalement: le $G_{S}$-torseur trivial $S\times_{S}G_{S}\rightarrow S$ est en effet une section de ce champ au-dessus de $S$;\vspace{2mm}
\item le champ $\textup{Tors}\left(S,G_{S}\right)$ est localement connexe, puisque comme nous l'avons remarqu{\'e} dans la section 1.2, tout $G_{S}$-torseur est isomorphe au $G_{S}$-torseur trivial, localement pour la topologie {\'e}tale sur $S$.
\end{itemize}}
\end{exem}

D'ailleurs, il est int{\'e}ressant de remarquer que toute gerbe neutre est de cette forme:
\begin{lem} Soit $S$ un sch{\'e}ma. Toute $S$-gerbe neutre $\mathcal{G}$ est {\'e}quivalente {\`a} une $S$-gerbe de torseurs.
\end{lem}

\uline{\textsc{Preuve}}: soit $\mathcal{G}$ une $S$-gerbe neutre, et soit $g$ un objet de $\mathcal{G}\left(S\right)$. Alors on a une {\'e}quivalence:
	\[\xymatrix@R=5pt@C=20pt{\epsilon:&\mathcal{G}\ar[r]&\textup{Tors}\left(S,\textup{Aut}\left(g\right)\right)\\& g' \ar@{|->}[r]&\textup{Isom}\left(g,g'\right)}
\]
o{\`u} $\textup{Aut}\left(g\right)$ (\textit{resp.} $\textup{Isom}\left(g,g'\right)$) d{\'e}signe le faisceau des automorphismes de l'objet $g$ (\textit{resp.} le faisceau des isomorphismes entre $g$ et $g'$). 
\begin{flushright}
$\Box$
\end{flushright}

On d{\'e}duit imm{\'e}diatemment de ce lemme et de la d{\'e}finition de gerbe la cons{\'e}quence suivante:
\begin{coro} Toute $S$-gerbe est localement {\'e}quivalente {\`a} une gerbe de torseurs.
\end{coro}

\uline{\textsc{Preuve}}: soit $\mathcal{G}$ une $S$-gerbe. Par d{\'e}finition, il existe un recouvrement {\'e}tale $\left(S_{i}\rightarrow S\right)_{i\in I}$ tel que $\textup{Ob}\left(\mathcal{G}\left(S_{i}\right)\right)\neq\emptyset$. Soit $g_{i}$ un objet de $\mathcal{G}\left(S_{i}\right)$, pour tout $i\in I$. Alors les gerbes (sur le site {\'e}tale restreint {\`a} $S_{i}$) $\mathcal{G}_{\left|S_{i}\right.}$ et $\textup{Tors}\left(S_{i},\textup{Aut}\left(g_{i}\right)\right)$ sont {\'e}quivalentes, o{\`u} $\textup{Aut}\left(g_{i}\right)$ est le faisceau sur $S_{i}$ des automorphismes de $g_{i}$, et $\mathcal{G}_{\left|S_{i}\right.}$ est le produit fibr{\'e} $\mathcal{G}\times_{S}S_{i}$.
\begin{flushright}
$\Box$
\end{flushright}
\begin{exem}[Gerbe des alg{\`e}bres simples centrales sur un corps, etc\ldots] \textup{Soit $k$ un corps. D'apr{\`e}s l'exemple 1.4.2, le champ $\textup{Asc}\left(k,n\right)$ (\textit{resp.} $\textup{SB}\left(k,n-1\right)$) des alg{\`e}bres simples centrales sur $k$ (\textit{resp.} des vari{\'e}t{\'e}s de Severi-Brauer de dimension $n-1$) est une gerbe. En outre, les gerbes $\textup{Asc}\left(k,n\right)$ et $\textup{SB}\left(k,n-1\right)$ sont {\'e}quivalentes, puisque toutes deux sont {\'e}quivalentes {\`a} la gerbe neutre $\textup{Tors}\left(k,PGL_{n}\right)$.}

\textup{De la m{\^e}me fa\c{c}on, lorsque $S$ est un sch{\'e}ma, le champ $\textup{LBun}\left(S\right)$\label{Lbun} (\textit{resp.} $\textup{VBun}\left(n,S\right)$\label{Vbun}, \textit{resp.} $\textup{Az}\left(n,S\right)$\label{AznS}) des fibr{\'e}s en droites sur $S$ (\textit{resp.} des fibr{\'e}s vectoriels de rang $n$ sur $S$, \textit{resp.} des alg{\`e}bres d'Azumaya d'indice $n$ sur $S$) est une gerbe neutre, {\'e}quivalente {\`a} la gerbe $\textup{Tors}\left(S,\mathbb{G}_{m,S}\right)$ (\textit{resp.} $\textup{Tors}\left(S,GL_{n,S}\right)$, \textit{resp.} $\textup{Tors}\left(S,PGL_{n,S}\right)$).}
\end{exem}

Nous pr{\'e}sentons maintenant la gerbe que nous {\'e}voquerons le plus souvent au cours de ce travail.
\begin{exem}[Gerbe des mod{\`e}les d'un torseur]\index{gerbe!des mod{\`e}les!d'un torseur} \textup{Soient $k$ un corps de caract{\'e}ristique nulle, $X$ un $k$-sch{\'e}ma lisse et g{\'e}om{\'e}triquement connexe, $G$ un $k$-groupe alg{\'e}brique lin{\'e}aire, et $\bar{P}\rightarrow \bar{X}$ un $\bar{G}_{X}$-torseur. Par analogie avec les rev{\^e}tements, nous introduisons la d{\'e}finition suivante:}
\begin{defi} On dit que $\bar{P}$ est de \textbf{corps des modules $k$}\index{corps des modules} s'il repr{\'e}sente une classe de $H^{1}\left(\bar{X},\bar{G}_{X}\right)^{\Gamma}$. Il revient au m{\^e}me de dire que pour tout $\sigma\in\Gamma$, les torseurs $\bar{P}$ et $^{\sigma}\bar{P}$ sont isomorphes.
\end{defi}
\begin{defi} On conserve les hypoth{\`e}ses adopt{\'e}es ci-dessus concernant $k$, $X$ et $G$, et on suppose $\bar{P}$ de corps des modules $k$.
\begin{enumerate}[(i)]
\item $L$ {\'e}tant une extension {\'e}tale de $k$, on appelle \textbf{mod{\`e}le de $\bar{P}$}\index{mod{\`e}le d'un torseur} au-dessus de $X_{L}$ un $G_{X_{L}}$-torseur $Y_{L}\rightarrow X_{L}$ tel que les $\bar{G}_{X}$-torseurs $\bar{P}$ et $\overline{Y_{L}}=Y_{L}\times_{X_{L}}\bar{X}$ sont isomorphes.
\item On appelle \textbf{gerbe des mod{\`e}les de $\bar{P}$}, et on note $D\left(\bar{P}\right)$\label{DPbar} la $k$-gerbe dont le groupoïde fibre au-dessus d'un ouvert {\'e}tale $\left(\textup{Spec}\;L\rightarrow\textup{Spec}\;k\right)$ a pour objets les $G_{X_{L}}$-torseurs $Y_{L}\rightarrow X_{L}$ tel qu'il existe $\sigma\in\Gamma$ tel que les $\bar{G}_{X}$-torseurs $^{\sigma}\bar{P}$ et $\overline{Y_{L}}$ sont isomorphes.
	
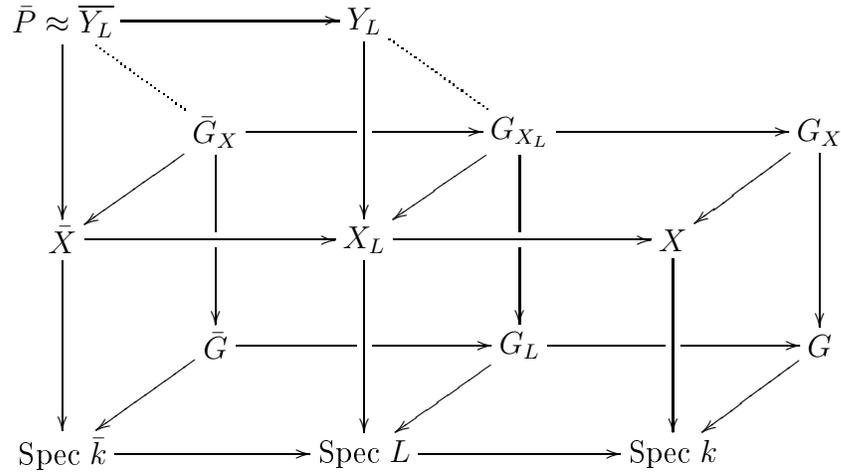
\[\xymatrix{\bar{P}\approx\overline{Y_{L}}\ar[dd]\ar[rr]\ar@{.}[dr]&&Y_{L}\ar[dd]\ar@{.}[dr]\\&\bar{G}_{X}\ar'[r][rr]\ar'[d][dd]\ar[dl]&&G_{X_{L}}\ar'[d][dd]\ar[dl]\ar[rr]&&G_{X}\ar[dl]\ar[dd]\\\bar{X}\ar[dd]\ar[rr]&&X_{L}\ar[dd]\ar[rr]&&X\ar[dd]\\&\bar{G}\ar[dl]\ar'[r][rr]&&G_{L}\ar[dl]\ar'[r][rr]&&G\ar[dl]\\\textup{Spec}\;\bar{k}\ar[rr]&&\textup{Spec}\;L\ar[rr]&&\textup{Spec}\;k}
\]

Les fl{\`e}ches du groupoïde fibre $D\left(\bar{P}\right)\left(L\right)$ sont les isomorphismes de $G_{X_{L}}$-torseurs sur $X_{L}$.
\end{enumerate}
\end{defi}
\textup{$D\left(\bar{P}\right)$ est effectivement une $k$-gerbe, car:
\begin{enumerate}[(1)]
\item c'est une cat{\'e}gorie fibr{\'e}e en groupoïdes sur $k$, avec les foncteurs de \textquotedblleft{restriction}\textquotedblright\ {\'e}vidents;
\item c'est un $k$-champ, puisque les isomorphismes de torseurs se recollent, et toute donn{\'e}e de descente sur les torseurs est effective (d'apr{\`e}s \cite{Gi2} III.1.4.1);
\item $D\left(\bar{P}\right)$ est localement connexe, puisque deux objets sont isomorphes {\`a} $\bar{P}$ (du fait que $^{\sigma}\bar{P}\approx\bar{P}$, $\forall\ \sigma\in\Gamma$);
\item $D\left(\bar{P}\right)$ est localement non-vide: en effet (d'apr{\`e}s \cite{SGA4-7}, 5.7 pour le cas ab{\'e}lien, et 5.14.(a) pour le cas g{\'e}n{\'e}ral):
	\[H^{1}\left(\bar{X},\bar{G}_{X}\right)=\lim_{\xrightarrow[L]{}}H^{1}\left(X_{L},G_{X_{L}}\right)
\]
la limite directe {\'e}tant prise sur les extensions {\'e}tales $L$ de $k$. Par suite, tout $\left[\bar{P}\right]$ est repr{\'e}sent{\'e} par une classe $\left[Y_{L_{0}}\right]\in H^{1}\left(X_{L_{0}},G_{X_{L_{0}}}\right)$, pour une extension {\'e}tale $L_{0}/k$ \textquotedblleft{suffisamment grande}\textquotedblright. Donc $D\left(\bar{P}\right)$ est localement non-vide, puisque le singleton:
	\[\left\{\textup{Spec}\;L_{0}\rightarrow\textup{Spec}\;k\right\}
\]
est un recouvrement {\'e}tale de $k$.
\end{enumerate}}

\textup{Notons, m{\^e}me si c'est une {\'e}vidence d'apr{\`e}s la d{\'e}finition, que \textbf{$D\left(\bar{P}\right)$ est une gerbe neutre si et seulement si $\bar{P}$ est d{\'e}fini sur $k$, \textit{i.e.} si et seulement si $\bar{P}$ a un mod{\`e}le sur $X$. La gerbe $D\left(\bar{P}\right)$ mesure donc l'obstruction {\`a} ce que la descente de $\bar{P}$ soit possible}.}
\begin{rem} \textup{La gerbe $D\left(\bar{P}\right)$ que l'on vient de d{\'e}finir est exactement la gerbe $D\left(c\right)$ d{\'e}crite dans \cite{Gi2} V.3.1.6.}
\end{rem}
\begin{rem} \textup{Dans le cas o{\`u} $G$ est ab{\'e}lien, et o{\`u} la condition $\bar{G}_{X}\left(\bar{X}\right)=\bar{G}\left(\bar{k}\right)$ est satisfaite, on peut associer {\`a} tout torseur $\bar{P}$ repr{\'e}sentant une classe dans $H^{1}\left(\bar{X},\bar{G}_{X}\right)^{\Gamma}$ un type\index{type d'un torseur} $\lambda_{\bar{P}}\in\textup{Hom}_{\Gamma}\left(\widehat{G},\textup{Pic}\;\bar{X}\right)$. En outre, comme on l'a rappel{\'e} dans l'introduction, on dispose des deux suites exactes: }
	\[H^{1}\left(k,G\right)\longrightarrow H^{1}\left(X,G_{X}\right)\stackrel{u}{\longrightarrow} H^{1}\left(\bar{X},\bar{G}_{X}\right)^{\Gamma}\stackrel{\delta^{1}}{\longrightarrow}H^{2}\left(k,G\right)\stackrel{}{\longrightarrow}H^{2}\left(X,G_{X}\right)
\]
\textup{et:}
	\[H^{1}\left(k,G\right)\longrightarrow H^{1}\left(X,G_{X}\right)\stackrel{\chi}{\longrightarrow} \textup{Hom}_{\Gamma}\left(\widehat{G},\textup{Pic}\;\bar{X}\right)\stackrel{\partial}{\longrightarrow}H^{2}\left(k,G\right)\stackrel{}{\longrightarrow}H^{2}\left(X,G_{X}\right)
\]

\textup{L'image de $\bar{P}$ par le morphisme $\delta^{1}$ est la classe de la gerbe $D\left(\bar{P}\right)$ (\textit{cf.} \cite{Gi2} V.3.2.1). D'autre part, la gerbe $\partial\left(\lambda_{\bar{P}}\right)$\label{Dlambda}, qui mesure l'obstruction {\`a} ce qu'il existe un $G_{X}$-torseur sur $X$ de type $\lambda_{\bar{P}}$, coïncide avec $D\left(\bar{P}\right)$ (\textit{cf.} \cite{HS} 3.7.(c)). Par cons{\'e}quent:}
\begin{pro} Les assertions suivantes sont {\'e}quivalentes:
\begin{enumerate}[(i)]
\item $\bar{P}$ est d{\'e}fini sur $k$;
\item la gerbe $D\left(\bar{P}\right)$ est neutre;
\item la gerbe $\partial\left(\lambda_{\bar{P}}\right)$ est neutre;
\item il existe un $G_{X}$-torseur sur $X$ de type $\lambda_{\bar{P}}$.
\end{enumerate}
\end{pro}
\textup{\uline{\textsc{Preuve}}: {\'e}vidente, d'apr{\`e}s la remarque pr{\'e}c{\'e}dente et les d{\'e}finitions des gerbes $D\left(\bar{P}\right)$ et $\partial\left(\lambda_{\bar{P}}\right)$.}
\begin{flushright}
$\Box$
\end{flushright}

\end{rem}
\end{exem}
\begin{exem}[Gerbe des mod{\`e}les d'un $G$-rev{\^e}tement]\index{gerbe!des mod{\`e}les!d'un rev{\^e}tement} \textup{La d{\'e}finition de la gerbe des mod{\`e}les d'un $G$-rev{\^e}tements (satisfaisant la condition corps des modules) pr{\'e}sente beaucoup de points communs avec la gerbe des mod{\`e}les d'un torseur. Nous renvoyons {\`a} la section 2 de \cite{DD1} pour une pr{\'e}sentation et une {\'e}tude d{\'e}taill{\'e}es de cette gerbe.}
\end{exem}
\begin{exem}[Gerbe des banalisations d'une alg{\`e}bre d'Azumaya]\index{gerbe!des banalisations d'une alg{\`e}bre d'Azumaya} \textup{Soit $X$ un sch{\'e}ma r{\'e}gulier et g{\'e}om{\'e}triquement irr{\'e}ductible. On suppose que $X$ est connexe, de telle sorte que les alg{\`e}bres d'Azumaya sur $X$ sont d'indice constant\footnote{\textit{I.e.} tous les $n_{i}$ sont {\'e}gaux dans l'exemple 1.2.9.} (\textit{cf.} \cite{Mi} p.143). Soit $\mathcal{A}$ une alg{\`e}bre d'Azumaya sur $X$. On appelle \textit{banalisation}\index{banalisation d'une alg{\`e}bre d'Azumaya} de $\mathcal{A}$ un couple $\left(\mathcal{E},\alpha\right)$, o{\`u} $\mathcal{E}$ est une $\mathcal{O}_{X}$-alg{\`e}bre localement libre de type fini, et $\alpha$ est un isomorphisme de $\mathcal{O}_{X}$-alg{\`e}bres:}
	\[\alpha:\textup{End}_{\mathcal{O}_{X}}\left(\mathcal{E}\right)\longrightarrow \mathcal{A}
\]
\pagebreak

\textup{On d{\'e}finit la gerbe des banalisations $\mathcal{B}\left(\mathcal{A}\right)$ de $\mathcal{A}$ ainsi: au-dessus d'un ouvert {\'e}tale $\left(X'\rightarrow X\right)$, le groupoïde fibre $\left[\mathcal{B}\left(\mathcal{A}\right)\right]\left(X'\right)$ est le groupoïde dont:}\vspace{2mm}
\textup{\begin{itemize}
\item les objets sont les banalisations de $\mathcal{A}_{\left|X'\right.}=\mathcal{A}\otimes_{\mathcal{O}_{X}}\mathcal{O}_{X'}$. Ce sont donc les couples $\left(\mathcal{E}',\alpha\right)$, $\mathcal{E}'$ {\'e}tant une $\mathcal{O}_{X'}$-alg{\`e}bre localement libre de type fini, et $\alpha$ un isomorphisme de $\mathcal{O}_{X'}$-alg{\`e}bres:
	\[\alpha:\textup{End}_{\mathcal{O}_{X'}}\left(\mathcal{E}'\right)\longrightarrow \mathcal{A_{\left|X'\right.}}
\]\vspace{2mm}
\item une fl{\`e}che entre deux banalisations $\left(\mathcal{E}_{1},\alpha_{1}\right)$ et $\left(\mathcal{E}_{2},\alpha_{2}\right)$ est un isomorphisme de $\mathcal{O}_{X'}$-alg{\`e}bres $\psi$ rendant commutatif le diagramme:
\[\xymatrix@R=40pt{\textup{End}_{\mathcal{O}_{X'}}\left(\mathcal{E}_{1}\right)\ar[rd]_{\alpha_{1}} \ar[rr]^{\psi}&& \textup{End}_{\mathcal{O}_{X'}}\left(\mathcal{E}_{2}\right)\ar[dl]^{\alpha_{2}} \\ & \mathcal{A}_{\left|X'\right.} }
\]
\end{itemize}}

\textup{Alors la th{\'e}orie de la descente assure que $\mathcal{B}\left(\mathcal{A}\right)$ est un champ, et c'est une gerbe d'apr{\`e}s l'{\'e}nonc{\'e} ci-dessous (\textit{cf.} \cite{G4} I.5.1 ou \cite{Mi} IV.4.2.1):}
\begin{pro} Soit $\mathcal{A}$ une $\mathcal{O}_{X}$-alg{\`e}bre qui est de type fini en tant que $\mathcal{O}_{X}$-module. Les assertions suivantes sont {\'e}quivalentes:
\begin{enumerate}[(i)]
\item $\mathcal{A}$ est une alg{\`e}bre d'Azumaya sur $X$;
\item il existe un recouvrement {\'e}tale $\left(X_{i}\rightarrow X\right)_{i\in I}$ tel que $\mathcal{A}_{\left|X_{i}\right.}$ soit banale, $\forall\;i\in I$; plus pr{\'e}cis{\'e}ment, il existe pour tout $i\in I$ un entier $n_{i}$ tel que: $\mathcal{A}_{\left|X_{i}\right.}\approx M_{n_{i}}\left(\mathcal{O}_{X_{i}}\right)$
\end{enumerate}
\end{pro}
\end{exem}
\begin{exem}[Gerbe des trivialisations d'un espace homog{\`e}ne]\index{gerbe!des trivialisations d'un espace homog{\`e}ne} \textup{Soient $k$ un corps de caract{\'e}ristique nulle, dont on fixe une cl{\^o}ture alg{\'e}brique $\bar{k}$, et $H$ un sous-groupe fini de $SL_{n}\left(k\right)$. Soit $V$ un espace homog{\`e}ne sous $SL_{n}$ avec isotropie $H$: on entend par l{\`a} que $V$ est une $k$-forme de $SL_{n}/H$, \textit{i.e.} on a un isomorphisme sur $\bar{k}$:}
	\[\bar{V}\approx \bar{SL_{n}}/\bar{H}
\]

\textup{On a la relation de domination de Springer \cite{Sp}:}\index{relation de domination de Springer}
	\[H^{1}\left(k,SL_{n}\right)\multimap H^{1}\left(k;SL_{n},H\right)
\]

\textup{L'espace homog{\`e}ne $V$ repr{\'e}sente une classe de $H^{1}\left(k;SL_{n},H\right)$. La gerbe $\mathcal{T}\left(V\right)$ que l'on va associer {\`a} $V$ mesure l'obstruction {\`a} ce que $\left[V\right]$ appartienne {\`a} l'image de la relation. Or, comme $H^{1}\left(k,SL_{n}\right)$ est r{\'e}duit {\`a} la classe du torseur trivial (par le th{\'e}or{\`e}me 90 de Hilbert), toute classe appartenant {\`a} l'image de la relation est triviale. D'o{\`u} l'appellation de gerbe des trivialisations de $V$. D{\'e}finissons maintenant cette gerbe: le groupoïde fibre $\left[\mathcal{T}\left(V\right)\right]\left(L\right)$ au-dessus d'un ouvert {\'e}tale $\left(\textup{Spec}\:L\rightarrow \textup{Spec}\;k\right)$ a pour objets les $SL_{n}$-torseurs $P_{L}$ (forc{\'e}ment triviaux) sur $L$ tels qu'il existe une application: $f_{L}:P_{L}\rightarrow V_{L}$, et pour fl{\`e}ches les isomorphismes de torseurs. C'est effectivement une gerbe d'apr{\`e}s \cite{Sp} 2.}
\end{exem}
\begin{exem}[Classe de Chern d'un fibr{\'e} en droites] \textup{Soient $X$ une vari{\'e}t{\'e} projective lisse complexe, $X^{an}$ la vari{\'e}t{\'e} analytique associ{\'e}e \cite{Fu}, et $\mathcal{L}$ un fibr{\'e} en droites sur $X$. Nous verrons en appendice comment interpr{\'e}ter la classe de Chern $c_{1}\left(\mathcal{L}\right)$ comme une gerbe sur le site analytique de $X$, mesurant l'obstruction {\`a} ce que $\mathcal{L}$ appartienne {\`a} l'image du morphisme:}
	\[H^{1}\left(X,\mathcal{O}_{X}\right)\longrightarrow \textup{Pic}\;X
\]
\end{exem}

Donnons maintenant deux exemples de champs qui ne sont pas des gerbes:
\begin{exem}[Le champ $\mathcal{M}_{g}$] \textup{Ce champ n'est pas une gerbe car il poss{\`e}de un ouvert\footnote{Pour la d{\'e}finition d'ouvert d'un champ alg{\'e}brique, nous renvoyons {\`a} \cite{LMB} ou \cite{Vi}.} (non-vide) de courbes ne poss{\'e}dant pas d'automorphisme non-trivial. Or nous verrons dans la section suivante que les gerbes se diff{\'e}rencient des champs par ce que les objets d'une gerbe ont tous (localement) les m{\^e}mes automorphismes.}
\end{exem}
\begin{exem}[Le champ $\pi_{\ast}\textup{Tors}\left(X,G_{X}\right)$] \textup{On consid{\`e}re une nouvelle fois un corps de caract{\'e}ristique nulle $k$, un $k$-sch{\'e}ma $\pi:X\rightarrow \textup{Spec}\;k$, et un $k$-groupe alg{\'e}brique lin{\'e}aire $G$. En g{\'e}n{\'e}ral, le champ $\pi_{\ast}\textup{Tors}\left(X,G_{X}\right)$ n'est pas une gerbe. Rappelons que ce $k$-champ a pour cat{\'e}gorie fibre au-dessus d'un ouvert {\'e}tale $\left(\textup{Spec}\;L\rightarrow\textup{Spec}\;k\right)$ la cat{\'e}gorie $Tors\left(X_{L},G_{X_{L}}\right)$ des $G_{X_{L}}$-torseurs sur $X_{L}$. C'est une cat{\'e}gorie fibr{\'e}e en groupoïdes sur $k$ munie des foncteurs de restriction {\'e}vidents, c'est un pr{\'e}champ car les isomorphismes de torseurs se recollent, et c'est un champ car les torseurs sur $X$ constituent un champ. En outre, ce champ est localement non-vide, et m{\^e}me globalement non-vide puisque le $G_{X}$-torseur trivial $X\times G_{X}$ est un objet de ce champ au-dessus de $\textup{Spec}\;k$. Mais deux objets ne sont pas n{\'e}cessairement localement isomorphes (ce champ n'est pas localement connexe): {\'e}tant donn{\'e}s $\left(\textup{Spec}\;L\rightarrow\textup{Spec}\;k\right)$ un ouvert {\'e}tale de $k$, deux $G_{X_{L}}$-torseurs $P$ et $Q$ sur $X_{L}$, il existe certes un recouvrement {\'e}tale $\left(X_{i}\rightarrow X_{L}\right)_{i\in I}$ tel que:}
	\[P_{\left|X_{i}\right.}\approx Q_{\left|X_{i}\right.},\ \forall\ i\in I.
\]

\textup{Cependant, les ouverts de ce recouvrement n'ont aucune raison de provenir d'un recouvrement {\'e}tale de $L$; plus explicitement, il n'y a aucune raison pour que les sch{\'e}mas $X_{i}$ soient de la forme $X_{L_{i}}$, o{\`u} les $L_{i}$ seraient des extensions {\'e}tales de $L$. Donc $P$ et $Q$, vus comme objets de:}
	\[\left[\pi_{\ast}\textup{Tors}\left(X,G_{X}\right)\right]\left(L\right)
\]
\textup{ne sont pas localement isomorphes.}

\textup{Cet exemple illustre donc le fait que l'image directe d'une gerbe par un morphisme de sch{\'e}mas n'est pas en g{\'e}n{\'e}ral une gerbe, alors que l'image inverse d'une gerbe est toujours une gerbe \cite{Gi2} V.1.4.2.}
\end{exem}
\end{section}
\begin{section}{Liens}
Avant de donner la d{\'e}finition de lien en g{\'e}n{\'e}ral, nous donnons juste id{\'e}e pratique de cette notion. Remarquons tout d'abord que par d{\'e}finition, deux objets $g$ et $g'$ d'une gerbe $\mathcal{G}$ sur un sch{\'e}ma $S$ sont localement isomorphes. Il s'ensuit que $g$ et $g'$ ont localement le m{\^e}me faisceau d'automorphismes, l'isomorphisme entre ces faisceaux {\'e}tant obtenu par conjugaison:
	\[\begin{array}{ccc}\textup{Aut}\left(g\right)&\longrightarrow&\textup{Aut}\left(g'\right)\\f&\longmapsto&\varphi\circ f\circ\varphi^{-1}\end{array}
\]
o{\`u} $\varphi$ est un isomorphisme entre $g$ et $g'$. L'id{\'e}e qui s'impose donc naturellement, lorsque l'on souhaite classifier les gerbes sur un sch{\'e}ma $S$, est d'associer {\`a} chaque $S$-gerbe un \textquotedblleft{faisceau en groupes}\textquotedblright\ sur $S$. C'est ce qui justifie l'introduction de la notion de lien.

\subsection*{Le champ des liens sur un sch{\'e}ma}
Soit $S$ un sch{\'e}ma. Les faisceaux de groupes sur $S$ constituent un $S$-champ not{\'e} $\left(\textup{FAGR}/S\right)$\label{FAGRS} (\textit{cf.} \cite{Gi2} II.3.4.12).\index{champ!des faisceaux de groupes} Pour tout ouvert {\'e}tale $\left(S'\rightarrow S\right)$, la cat{\'e}gorie fibre 
	\[\left(\textup{FAGR}/S\right)\left(S'\rightarrow S\right)
\]
est la cat{\'e}gorie $FAGR\left(S'\right)$\label{fagrS'} des faisceaux de groupes sur le site {\'e}tale de $S'$. Le pr{\'e}champ des liens\index{pr{\'e}champ!des liens} est construit en deux temps {\`a} partir de ce champ:
\begin{defi} On note $\left(Lien/S\right)$\label{Prechlien} la cat{\'e}gorie fibr{\'e}e sur $S$, dont la cat{\'e}gorie fibre $\left(Lien/S\right)\left(S'\right)$ au-dessus d'un ouvert {\'e}tale $\left(S'\rightarrow S\right)$ a pour objets les faisceaux de groupes sur $S'$; les morphismes entre deux objets $\mathcal{F}$ et $\mathcal{G}$ sont les sections du faisceau quotient:
	\[\textup{Int}\left(\mathcal{F}\right)\backslash{\textup{Hom}_{FAGR\left(S'\right)}\left(\mathcal{F},\mathcal{G}\right)}/\textup{Int}\left(\mathcal{G}\right)
\]
\end{defi}
\begin{pro} La $S$-cat{\'e}gorie fibr{\'e}e $\left(Lien/S\right)$ est un $S$-pr{\'e}champ, le pr{\'e}champ des liens sur $S$. On appelle \textbf{champ des liens sur $S$}\index{champ!des liens} et on note $\left(\textup{LIEN/S}\right)$\label{ChLien} le champ qui lui est associ{\'e} par le foncteur (\textit{cf.} exemple 1.4.6):
	\[\left(\textup{Pr{\'e}champ/S}\right)\longrightarrow \left(\textup{Champ/S}\right)
\]

On appelle \textbf{lien}\index{lien} sur $S$ un objet du champ $\left(\textup{LIEN/S}\right)$. Notons que l'on peut associer {\`a} tout faisceau de groupes $G$ sur $S$ un lien sur $S$, cette association {\'e}tant obtenue en gr{\^a}ce au foncteur compos{\'e}:
	\[\textup{lien}:\left(\textup{FAGR}/S\right)\longrightarrow \left(Lien/S\right)\longrightarrow\left(LIEN/S\right).
\]

Un lien sur $S$ (ou $S$-lien) est dit \textbf{repr{\'e}sentable}\index{lien!repr{\'e}sentable} s'il appartient {\`a} l'image essentielle de ce foncteur; autrement dit, un lien $\mathcal{L}$ sur $S$ est repr{\'e}sentable par un faisceau de groupes $G$ s'il existe un isomorphisme de liens sur $S$:\label{lienG}
	\[\mathcal{L}\approx \textup{lien}\;G
\]

Un lien $\mathcal{L}$ sur $S$ est dit \textbf{localement repr{\'e}sentable}\index{lien!localement repr{\'e}sentable} par un faisceau de groupes $G$ sur $S$ s'il existe un recouvrement {\'e}tale $\left(S_{i}\rightarrow S\right)$ et des isomorphismes de liens sur $S_{i}$:
	\[\mathcal{L}_{\left|S_{i}\right.}\approx \textup{lien}\left(G_{\left|S_{i}\right.}\right)
\]
Par construction m{\^e}me, tout lien sur $S$ est localement repr{\'e}sentable par un faisceau de groupes.

Enfin, un lien $\mathcal{L}$ sur $S$ est dit \textbf{r{\'e}alisable}\index{lien!r{\'e}alisable} s'il existe une $S$-gerbe $\mathcal{G}$ telle que:
	\[\textup{lien}\left(\mathcal{G}\right)\approx\mathcal{L}
\]
\end{pro}
\begin{defi} Soient $S$ un sch{\'e}ma, $G$ un faisceau de groupes sur $S$ et $\mathcal{G}$ une $S$-gerbe. On dit que $\mathcal{G}$ est \textbf{li{\'e}e par $G$}\index{gerbe!li{\'e}e} si $\textup{lien}\left(\mathcal{G}\right)$ est repr{\'e}sentable par $G$. On dit que $\mathcal{G}$ est \textbf{localement li{\'e}e par $G$}\index{gerbe!localement li{\'e}e} si $\textup{lien}\left(\mathcal{G}\right)$ est localement repr{\'e}sentable par $G$. Il suffit pour cela qu'il existe pour tout ouvert {\'e}tale $\left(S'\rightarrow S\right)$ et pour tout $g'\in \textup{Ob}\left(\mathcal{G}\left(S'\right)\right)$ des isomorphismes fonctoriels (\cite{Mi} p.144):
	\[G_{S}\left(S'\right)\longrightarrow \textup{Aut}_{\mathcal{G}\left(S'\right)}\left(g'\right)
\]
\end{defi}
\begin{exem}\textup{On se place dans la situation de l'exemple 1.5.6, et on consid{\`e}re $\bar{P}$ un $\bar{G}_{X}$-torseur de corps des modules $k$. Alors:}
\begin{lem} Si $G$ est ab{\'e}lien, alors:
\begin{enumerate}[(i)]
\item la $k$-gerbe $D\left(\bar{P}\right)$ est li{\'e}e par $\pi_{\ast}G_{X}$;
\item si de plus la condition $\bar{G}_{X}\left(\bar{X}\right)=\bar{G}\left(\bar{k}\right)$ est satisfaite, alors $D\left(\bar{P}\right)$ est li{\'e}e par $G$.
\end{enumerate}
\end{lem}

\uline{\textsc{Preuve}}: \textup{En effet, soit $\left(\textup{Spec}\;L\rightarrow \textup{Spec}\;k\right)$ un ouvert {\'e}tale tel que le groupoïde fibre $D\left(\bar{P}\right)\left(L\right)$ soit non-vide, et soit $P_{L}$ un de ses objets. La gerbe $D\left(\bar{P}\right)_{\left|\textup{Spec}\;L\right.}$ est donc neutre(!), et elle est donc {\'e}quivalente {\`a} une gerbe de torseurs \textbf{\uline{sur $L$}}. Explicitement, on a une {\'e}quivalence de $L$-gerbes \cite{Gi2} V.3.1.6.(ii):
	\[\begin{array}{rccc}\epsilon:&D\left(\bar{P}\right)_{\left|\textup{Spec}\;L\right.}&\longrightarrow&\textup{Tors}\left(L,\pi_{\ast}\textup{ad}_{G_{X_{L}}}\left(P_{L}\right)\right)\\&P'&\longmapsto&\pi_{\ast}\left(P'\wedge^{G_{X_{L}}}P_{L}^{opp}\right)\end{array}
\]
o{\`u} $P_{L}^{opp}$ d{\'e}signe le bitorseur oppos{\'e} {\`a} $P_{L}$, vu comme un $\left(\textup{ad}_{G_{X_{L}}}\left(P_{L}\right),G_{X_{L}}\right)$-bitorseur. On v{\'e}rifie ensuite que $\pi_{\ast}$ et $\textup{ad}$ commutent, ce qui est imm{\'e}diat; enfin, comme $G$ est ab{\'e}lien, le faisceau $\textup{ad}_{\pi_{\ast}G_{X_{L}}}\left(\pi_{\ast}P_{L}\right)$ n'est autre que le faisceau $\pi_{\ast}G_{X_{L}}$, d'o{\`u} le point (i).}

\textup{Pour le point (ii) on utilise le fait que la condition $\bar{G}_{X}\left(\bar{X}\right)=\bar{G}\left(\bar{k}\right)$ implique que les faisceaux $\pi_{\ast}G_{X}$ et $G$ sont isomorphes (\textit{cf.} preuve du lemme 1.2.14).}
\begin{flushright}
$\Box$
\end{flushright}
\end{exem}
\begin{exem}\textup{La gerbe des mod{\`e}les d'un $G$-rev{\^e}tement est li{\'e}e par le centre de $G$ (\textit{cf.} \cite{DD2}).}
\end{exem}
\begin{exem}\textup{La gerbe des banalisations d'une alg{\`e}bre d'Azumaya sur un sch{\'e}ma $X$ est li{\'e}e par $\mathbb{G}_{m,X}$ (\textit{cf.} \cite{Gi2} V.4.2 ou \cite{Mi} p.145).}
\end{exem}
\begin{exem}\textup{La classe de Chern d'un fibr{\'e} en droites sur une vari{\'e}t{\'e} analytique $X$, vue comme une gerbe sur le site analytique de $X^{an}$, est li{\'e}e par $\mathbb{Z}_{X^{an}}$ (\textit{cf.} appendice A).}
\end{exem}

L'{\'e}nonc{\'e} suivant est encore trait{\'e} dans le cas g{\'e}n{\'e}ral dans \cite{Gi2} (corollaire IV.1.1.7.3).
\begin{pro} Soient $S$ un sch{\'e}ma et $G$ un faisceau de groupes sur $S$. L'ensemble des classes d'isomorphie de liens localement repr{\'e}sentables par $G$ est en bijection avec l'ensemble de cohomologie $H^{1}\left(k,\textup{Out}\;G\right)$, o{\`u} $\textup{Out}\;G$ est le faisceau sur $S$ des automorphismes ext{\'e}rieurs de $G$.
\end{pro}
\begin{rem}\textup{L'ensemble $H^{1}\left(S, \textup{Out}\;G\right)$ des formes ext{\'e}rieures\index{forme!ext{\'e}rieure} est point{\'e} par la classe du $\textup{Out}\;G$-torseur trivial; ce dernier correspond au $S$-lien $\textup{lien}\;G$.}
\end{rem}
Soient $S$ un sch{\'e}ma et $G$ un faisceau de groupes sur $S$. On a toujours la suite exacte de faisceaux de groupes sur $S$:
	\[0\longrightarrow \textup{Int}\;G\longrightarrow \textup{Aut}\;G\longrightarrow \textup{Out}\;G\longrightarrow 1
\]
d'o{\`u} une suite de cohomologie (d{\'e}j{\`a} {\'e}voqu{\'e}e, voir diagramme p. 14):
	\[H^{1}\left(S, \textup{Int}\;G\right)\longrightarrow H^{1}\left(S, \textup{Aut}\;G\right)\stackrel{\lambda}{\longrightarrow}H^{1}\left(S, \textup{Out}\;G\right)
\]

L'application $\lambda$ est juste celle qui associe {\`a} une forme\footnote{Rappelons que l'on appelle $S$-forme\index{forme} d'un objet $A$ d{\'e}fini sur $S$ un objet $A'$ d{\'e}fini sur $S$, localement isomorphe {\`a} $A$ pour la topologie {\'e}tale sur $S$. Par exemple, une $k$-vari{\'e}t{\'e} de Severi-Brauer de dimension $n$ est une $k$-forme de l'espace projectif $\mathbb{P}^{n}_{k}$; une $k$-alg{\`e}bre simple centrale d'indice $n$ est une $k$-forme de l'alg{\`e}bre de matrices $M_{n}\left(k\right)$; une alg{\`e}bre d'Azumaya d'indice $n$ sur $S$ est une $S$-forme de la $\mathcal{O}_{S}$-alg{\`e}bre $M_{n}\left(\mathcal{O}_{S}\right)$, etc\ldots} (sur $S$) de $G$ le lien qu'elle repr{\'e}sente. Un lien est donc repr{\'e}sentable par $G$ s'il appartient {\`a} l'image de $\lambda$. En outre, on d{\'e}duit facilement de l'exactitude de cette suite l'{\'e}nonc{\'e}:
\begin{lem} Avec les notations adopt{\'e}es pr{\'e}c{\'e}demment, si $G'$ est une $S$-forme int{\'e}rieure de $G$, alors:
	\[\textup{lien}\;G'\approx\textup{lien}\;G
\]
\end{lem}

Une cons{\'e}quence beaucoup moins imm{\'e}diate de l'exactitude de la suite pr{\'e}c{\'e}dente est la suivante:
\begin{pro} Soient $S$ un sch{\'e}ma, $G_{S}$ un sch{\'e}ma en groupes r{\'e}ductifs sur $S$. Alors tout $S$-lien localement (pour la topologie {\'e}tale) repr{\'e}sentable par $G_{S}$ est repr{\'e}sentable par une $S$-forme de $G_{S}$.
\end{pro}

\uline{\textsc{Preuve}}: soit $\mathcal{L}$ un $S$-lien localement repr{\'e}sentable par $G_{S}$. $\mathcal{L}$ repr{\'e}sente une classe de $H^{1}\left(S,\textup{Out}\;G_{S}\right)$. Du fait que $G_{S}$ est r{\'e}ductif, la suite:
\vspace{2mm}
	\[\xymatrix{1\ar[r]&\textup{Int}\;G_{S}\ar[r]&\textup{Aut}\;G_{S}\ar[r]&\textup{\textup{Out}}\;G_{S}\ar[r]\ar@/_16pt/[l]_{\sigma}&1}
\]
est scind{\'e}e \cite{De1} p.28. Cette section induit une section de l'application:
	\[\xymatrix{H^{1}\left(S,\textup{Aut}\;G_{S}\right)\ar[r]_{\lambda}&H^{1}\left(S,\textup{Aut}\;G_{S}\right)\ar@/_20pt/[l]_{\sigma^{\left(1\right)}}}
\]\pagebreak

d'o{\`u} la conclusion.
\begin{flushright}
$\Box$
\end{flushright}
\begin{rem}\textup{Nous avons d{\'e}lib{\'e}r{\'e}ment choisi de nous restreindre {\`a} la topologie {\'e}tale, en vue de nos applications. Signalons cependant que cet {\'e}nonc{\'e} est valable dans un cadre beaucoup plus g{\'e}n{\'e}ral (\textit{cf.} \cite{Do1} V.3.2).}
\end{rem}
\end{section}
\begin{section}{Cohomologie {\`a} valeurs dans un lien}
\begin{defi} Soient $S$ un sch{\'e}ma et $\mathcal{L}$ un $S$-lien. On d{\'e}finit l'ensemble\label{H2L}
	\[H^{2}\left(S,\mathcal{L}\right)
\]
comme l'ensemble des classes d'{\'e}quivalence de $S$-gerbes de lien $\mathcal{L}$ pour la relation d'{\'e}quivalence suivante: $\mathcal{G}$ et $\mathcal{G}'$ sont dites \textbf{{\'e}quivalentes au sens de Giraud}\index{gerbes!equivalentes au sens de Giraud@{\'e}quivalentes au sens de Giraud} s'il existe une {\'e}quivalence $\epsilon:\mathcal{G}\rightarrow \mathcal{G}'$ li{\'e}e par $\textup{id}_{\mathcal{L}}$. On note $\left[\mathcal{G}\right]$ la classe d'une $S$-gerbe de lien $\mathcal{L}$.
\end{defi}
\begin{exem} \textup{Soient $S$ un sch{\'e}ma, $G_{S}$ un $S$-sch{\'e}ma en groupes, et $G'_{S}$ une $S$-forme int{\'e}rieure de $G_{S}$. Alors les gerbes $\textup{Tors}\left(S,G_{S}\right)$ et $\textup{Tors}\left(S,G'_{S}\right)$ ont {\'e}videmment m{\^e}me lien (d'apr{\`e}s le lemme 1.6.11), et ce lien n'est autre que $\textup{lien}\;G_{S}$. En outre, il existe une {\'e}quivalence
	\[\epsilon:\textup{Tors}\left(S,G_{S}\right)\longrightarrow \textup{Tors}\left(S,G'_{S}\right)
\]
mais $\epsilon$ n'est pas forc{\'e}ment li{\'e}e par l'identit{\'e}. En effet, on d{\'e}duit de la suite exacte de faisceaux:}
	\[0\longrightarrow Z\left(G_{S}\right)\longrightarrow G_{S}\longrightarrow \textup{Int}\;G_{S}\longrightarrow 0
\]
\textup{une suite longue d'ensembles point{\'e}s:}
	\[H^{1}\left(S,G_{S}\right)\stackrel{\alpha}{\longrightarrow} H^{1}\left(S,\textup{Int}\;G_{S}\right)\stackrel{\beta}{\longrightarrow} H^{2}\left(S,Z\left(G_{S}\right)\right)
\]
\textup{Comme nous le verrons bient{\^o}t, le groupe $H^{2}\left(S,Z\left(G_{S}\right)\right)$ agit simplement transitivement sur l'ensemble $H^{2}\left(S,\textup{lien}\;G_{S}\right)$; en outre, la proposition IV.3.2.6 de \cite{Gi2} assure que l'ensemble $H^{1}\left(S,\textup{Int}\;G_{S}\right)$ \textquotedblleft{agit}\textquotedblright\ transitivement (par l'interm{\'e}diaire de $\beta$) sur l'ensemble $H^{2}\left(S,\textup{lien}\;G_{S}\right)'$ des classes neutres de $S$-gerbes li{\'e}es par $G_{S}$.}

\textup{Par cons{\'e}quent, les gerbes $\textup{Tors}\left(S,G_{S}\right)$ et $\textup{Tors}\left(S,G'_{S}\right)$ ne sont {\'e}quivalentes (au sens de Giraud) que si $G'_{S}$ est une $S$-forme int{\'e}rieure de $G_{S}$ telle que:}
	\[\beta\left(\left[G'_{S}\right]\right)=0\in H^{2}\left(S,Z\left(G_{S}\right)\right)
\]

\textup{Remarquons pour conclure cet exemple que si $G'_{S}$ est une $S$-forme int{\'e}rieure de $G_{S}$, alors les gerbes $\textup{Tors}\left(S,G_{S}\right)$ et $\textup{Tors}\left(S,G'_{S}\right)$ sont {\'e}quivalentes \textit{au sens de Breen} (\textit{cf.} \cite{B3}), \textit{i.e.} repr{\'e}sentent la m{\^e}me classe dans:}
	\[H^{1}\left(S,G_{S}\rightarrow \textup{Aut}\;G_{S}\right)
\]
\end{exem}

Avant d'aller plus loin, commen\c{c}ons par quelques faits et remarques {\'e}l{\'e}mentaires concernant le $H^{2}$ {\`a} valeurs dans un lien.
\begin{fait} Soient $S$ un sch{\'e}ma et $\mathcal{L}$ un $S$-lien. L'ensemble $H^{2}\left(S,\mathcal{L}\right)$ est non-vide si et seulement si $\mathcal{L}$ est r{\'e}alisable.
\end{fait}

\uline{\textsc{Preuve}}: triviale, d'apr{\`e}s la d{\'e}finition de lien r{\'e}alisable.
\begin{flushright}
$\Box$
\end{flushright}
\begin{fait} Soient $S$ un sch{\'e}ma et $\mathcal{L}$ un $S$-lien. Soient $\mathcal{G}$ et $\mathcal{G}'$ deux $S$-gerbes {\'e}quivalentes (au sens de Giraud ou non) de lien $\mathcal{L}$. Si $\mathcal{G}$ est neutre, alors $\mathcal{G}'$ est {\'e}galement neutre. On appelle \textbf{classe neutre}\index{classe!neutre} de $H^{2}\left(S,\mathcal{L}\right)$ une classe de $S$-gerbes {\'e}quivalentes au sens de Giraud et de lien $\mathcal{L}$ dont un (donc tous, par ce qui pr{\'e}c{\`e}de) repr{\'e}sentant est neutre.
\end{fait}

\uline{\textsc{Preuve}}: si $\mathcal{G}$ est neutre, alors elle a une section au-dessus de $S$. L'{\'e}quivalence entre $\mathcal{G}$ et $\mathcal{G}'$ fournit alors une section de $\mathcal{G}'$ au-dessus de $S$. Donc $\mathcal{G}'$ est neutre {\`a} son tour.
\begin{flushright}
$\Box$
\end{flushright}
\begin{fait} Soient $S$ un sch{\'e}ma et $G_{S}$ un sch{\'e}ma en groupes sur $S$. Alors l'ensemble $H^{2}\left(S,\textup{lien}\;G_{S}\right)$ poss{\`e}de une classe privil{\'e}gi{\'e}e:
	\[\left[\textup{Tors}\left(S,G_{S}\right)\right]
\]

On l'appelle la \textbf{classe triviale}\index{classe!triviale} de $H^{2}\left(S,\textup{lien}\;G_{S}\right)$.
\end{fait}

L'ensemble $H^{2}\left(S,\textup{lien}\;G_{S}\right)$ poss{\`e}de donc toujours cette classe triviale, mais il peut aussi (d'apr{\`e}s l'exemple 1.7.2) poss{\'e}der plusieurs classes neutres, diff{\'e}rentes de la classe triviale. Explicitement, avec les notations du fait ci-dessus, soit $G'_{S}$ une $S$-forme int{\'e}rieure de $G_{S}$. Si $\beta\left(\left[G'_{S}\right]\right)\neq0$, alors la gerbe $\textup{Tors}\left(S,G'_{S}\right)$ n'est pas {\'e}quivalente (au sens de Giraud) {\`a} la gerbe $\textup{Tors}\left(S,G_{S}\right)$. Dans cette situation donc, la classe $\left[\textup{Tors}\left(S,G'_{S}\right)\right]$ est une classe neutre de $H^{2}\left(S,\textup{lien}\;G_{S}\right)$, diff{\'e}rente de la classe triviale.

Evidemment, la situation est nettement plus simple si l'on consid{\`e}re un sch{\'e}ma en groupes ab{\'e}liens:
\begin{fait} Soient $S$ un sch{\'e}ma et $G_{S}$ un $S$-sch{\'e}ma en groupes ab{\'e}liens. L'ensemble $H^{2}\left(S,\textup{lien}\;G_{S}\right)$ poss{\`e}de une unique classe neutre, qui est la classe triviale: $\left[\textup{Tors}\left(S,G_{S}\right)\right]$.
\end{fait}

\uline{\textsc{Preuve}}: puisque $G_{S}$ est ab{\'e}lien, il ne poss{\`e}de pas d'automorphisme int{\'e}rieur non-trivial. La suite exacte de cohomologie (\textit{cf.} section pr{\'e}c{\'e}dente):
	\[H^{1}\left(S, \textup{Int}\;G_{S}\right)\longrightarrow H^{1}\left(S, \textup{Aut}\;G_{S}\right)\stackrel{\lambda}{\longrightarrow}H^{1}\left(S, \textup{Out}\;G_{S}\right)
\]
se r{\'e}duit alors {\`a} l'identit{\'e}:
	\[H^{1}\left(S, \textup{Aut}\;G_{S}\right)\stackrel{\textup{id}}{\longrightarrow}H^{1}\left(S, \textup{Aut}\;G_{S}\right)
\]\pagebreak

Autrement dit, il existe un unique ({\`a} isomorphisme pr{\`e}s) $S$-lien localement repr{\'e}sentable par $G_{S}$: le lien $\textup{lien}\;G_{S}$ lui-m{\^e}me, d'o{\`u} le fait.
\begin{flushright}
$\Box$
\end{flushright}

Dans ce cas, on pourra noter $0$ cette unique classe neutre (et on pourra alors parler de \textit{classe nulle}\index{classe!nulle}). En effet:
\begin{fait} Soient $S$ un sch{\'e}ma et $G_{S}$ un $S$-sch{\'e}ma en groupes ab{\'e}liens. L'ensemble $H^{2}\left(S,\textup{lien}\;G_{S}\right)$ poss{\`e}de une loi de groupe. De plus les groupes $H^{2}\left(S,\textup{lien}\;G_{S}\right)$ et $H^{2}_{\acute{e}t}\left(S,G_{S}\right)$ sont isomorphes.
\end{fait}

\uline{\textsc{Preuve}}: c'est la proposition IV.3.5.1 de \cite{Gi2}.
\begin{flushright}
$\Box$
\end{flushright}

Enfin le th{\'e}or{\`e}me suivant est d'une importance capitale, puisqu'il permet de r{\'e}duire la cohomologie {\`a} valeurs dans un lien {\`a} celle de son centre:
\begin{fait}[Th{\'e}or{\`e}me IV.3.3.3 de \cite{Gi2}] Soient $S$ un sch{\'e}ma et $\mathcal{L}$ un $S$-lien. Alors l'ensemble $H^{2}\left(S,\mathcal{L}\right)$ est un pseudo-torseur sous $H^{2}\left(S,Z\left(\mathcal{L}\right)\right)$\footnote{\textit{I.e.} $H^{2}\left(S,\mathcal{L}\right)$ est vide ou principal homog{\`e}ne sous l'action de $H^{2}\left(S,Z\left(\mathcal{L}\right)\right)$.}, $Z\left(\mathcal{L}\right)$ d{\'e}signant le centre du lien $\mathcal{L}$.

Dans le cas particulier o{\`u} $\mathcal{L}=\textup{lien}\;G$ est repr{\'e}sentable par un faisceau de groupes $G$ sur $S$, l'ensemble $H^{2}\left(S,\textup{lien}\;G\right)$ est principal homog{\`e}ne sous l'action de $H^{2}\left(S,Z\left(G\right)\right)$ (puisque le centre de $\textup{lien}\;G$ est repr{\'e}sentable par le centre de $G$ \cite{Gi2} IV.1.5.3.(iii)).
\end{fait}

\begin{defi} Soient $S$ un sch{\'e}ma et $\mathcal{L}$ un $S$-lien. Nous dirons que l'ensemble $H^{2}\left(S,\mathcal{L}\right)$ est \textbf{inessentiel}\index{inessentiel} s'il est non-vide et s'il n'est compos{\'e} que de classes neutres.

D'apr{\`e}s le fait pr{\'e}c{\'e}dent, lorsque $G$ est un faisceau de groupes ab{\'e}liens sur $S$, l'ensemble $H^{2}\left(S,\textup{lien}\;G\right)$ est inessentiel si et seulement si il est r{\'e}duit {\`a} la classe nulle.
\end{defi}

Voici maintenant un exemple de situation particuli{\`e}rement important o{\`u} le $H^{2}$ est inessentiel, et qui est en partie une cons{\'e}quence du fait 1.7.8.
\begin{theo} Soient $k$ un corps de caract{\'e}ristique nulle et $G$ un $k$-groupe alg{\'e}brique. Alors l'ensemble $H^{2}\left(k,\textup{lien}\;G\right)$ est inessentiel dans les cas suivants:
\begin{enumerate}[(i)]
\item $G$ est semi-simple adjoint;
\item $k$ est un corps de nombres purement imaginaires et $G$ est semi-simple;
\item $k$ est un corps de nombres et $G$ est semi-simple simplement connexe.
\end{enumerate}
\end{theo}

\uline{\textsc{Preuve}}: pour les cas (ii) et (iii), c'est le th{\'e}or{\`e}me VI.3.2 de \cite{Do1}. Pour le cas o{\`u} $G$ est semi-simple adjoint, commen\c{c}ons par noter que l'ensemble $H^{2}\left(k,\textup{lien}\;G\right)$ n'est pas vide, puisqu'il poss{\`e}de au moins la classe triviale. D'apr{\`e}s le fait 1.7.8, il est donc principal homog{\`e}ne sous $H^{2}\left(k,Z\left(G\right)\right)$. Or, $G$ {\'e}tant adjoint, son centre est trivial. Donc le groupe $H^{2}\left(k,Z\left(G\right)\right)$ est nul, et l'ensemble $H^{2}\left(k,\textup{lien}\;G\right)$ est r{\'e}duit {\`a} la classe triviale. En particulier, il est inessentiel.
\begin{flushright}
$\Box$
\end{flushright}

Pour achever ce chapitre, nous donnons deux exemples d'application de ce th{\'e}or{\`e}me.

\begin{coro} Soient $k$ un corps de caract{\'e}ristique nulle, et $G$ une $k$-forme de $PGL_{n}$. Toute $k$-gerbe de lien $\textup{lien}\;G$ est neutre.
\end{coro}

\uline{\textsc{Preuve}}: c'est une cons{\'e}quence directe du th{\'e}or{\`e}me 1.7.10, puisque $PGL_{n}$ est semi-simple adjoint.
\begin{flushright}
$\Box$
\end{flushright}
\begin{coro} Si $k$ est un corps de nombres, l'application:
	\[\delta_{n}:H^{1}\left(k,PGL_{n}\right)\longrightarrow\ _{n}\textup{Br}\;k
\]
est bijective. Autrement dit, tout {\'e}l{\'e}ment d'ordre $n$ dans $\textup{Br}\;k$ est l'image par $\delta_{n}$ d'une $k$-alg{\`e}bre simple centrale d'indice $n$. En termes de gerbes, cela revient {\`a} dire que toute $k$-gerbe li{\'e}e par $\mu_{n}$ (c'est en particulier un $k$-champ de Deligne-Mumford) est la gerbe des trivialisations d'une $k$-alg{\`e}bre simple centrale d'indice $n$.
\end{coro}

\uline{\textsc{Preuve}}: de la suite exacte de faisceaux sur $k$:
	\[0\longrightarrow \mu_{n}\longrightarrow SL_{n}\longrightarrow PGL_{n}\longrightarrow 1
\]
on d{\'e}duit l'existence d'une application injective\footnote{La seule trivialit{\'e} de l'ensemble $H^{1}\left(k,SL_{n}\right)$ n'entra{\^i}ne pas l'injectivit{\'e} de $\delta^{n}$ \textit{cf.} \cite{Jah}.}:
	\[1\longrightarrow H^{1}\left(k,PGL_{n}\right)\stackrel{\delta_{n}}{\longrightarrow}H^{2}\left(k,\mu_{n}\right)=\ _{n}\textup{Br}\;k
\]

Soit $\left[\mathcal{G}\right]\in\!_{n}\textup{Br}\;k$. L'obstruction {\`a} ce que cette classe appartienne {\`a} l'image de $\delta_{n}$ est mesur{\'e}e par une gerbe dont le lien est repr{\'e}sentable par une forme de $SL_{n}$ \cite{Gi2} IV.4.2.10. Comme $SL_{n}$ est semi-simple simplement connexe, cette gerbe est neutre d'apr{\`e}s le th{\'e}or{\`e}me 1.7.10. Donc $\delta_{n}$ est bijective.
\begin{flushright}
$\Box$
\end{flushright}
\begin{rem}\textup{En fait, il est d{\'e}j{\`a} connu que cette application est bijective dans bien d'autres cas que les corps de nombres (\textit{cf.} \cite{CTGP}). Remarquons aussi qu'il existe des corps sur lequels $\delta_{n}$ n'est pas bijective; c'est le cas du corps $K_{M}$ construit par Merkurjev pour obtenir un contre-exemple {\`a} une conjecture de Kaplansky \cite{Me}. Il existe en effet sur $K_{M}$ une alg{\`e}bre simple centrale d'indice 4, mais d'esposant 2 (\textit{i.e.} dont l'image dans $\textup{Br}\;k$ est d'ordre 2). Une telle alg{\`e}bre simple centrale repr{\'e}sente donc une classe de $H^{1}\left(K_{M},PGL_{4}\right)$ dont l'image par l'isomorphisme:}
	\[\Delta:\textup{Br}_{Az}K_{M}\longrightarrow\textup{Br}\;K_{M}
\]
\textup{appartient {\`a} $_{2}\textup{Br}\;K_{M}$. Il s'ensuit que l'application:}
	\[\delta_{2}:H^{1}\left(K_{M},PGL_{2}\right)\longrightarrow H^{2}\left(K_{M},\mu_{2}\right)=\ _{2}\textup{Br}\;K_{M}
\]
\textup{n'est pas surjective.}

\textup{Par cons{\'e}quent:}
\begin{coro} Il existe une $K_{M}$-forme $SL'_{2}$ de $SL_{2}$ telle que $H^{2}\left(K_{M},\textup{lien}\;SL'_{2}\right)$ poss{\`e}de une classe non-neutre.
\end{coro}
\end{rem}
\end{section}
\end{chapter}

\begin{chapter}[Descente de torseurs: le cas ab{\'e}lien]{Descente de torseurs et points rationnels: le cas ab{\'e}lien}
\thispagestyle{empty}
Dans \cite{DD1}, D{\`e}bes et Douai montrent que l'obstruction {\`a} ce qu'un $G$-rev{\^e}tement $\bar{f}:\bar{X}\rightarrow \bar{B}$ de corps des modules $k$ soit d{\'e}fini sur $k$ est mesur{\'e}e par une gerbe (la gerbe des mod{\`e}les $\mathcal{G}\left(\bar{f}\right)$ du $G$-rev{\^e}tement) localement li{\'e}e par le centre de $G$. Si $G$ est ab{\'e}lien, ou si $Z\left(G\right)$ est un facteur direct de $G$, cette gerbe est neutre lorsque la suite exacte de groupes fondamentaux:
	\[1\longrightarrow \Pi_{\bar{k}}\left(\bar{B}^{\ast}\right)\longrightarrow \Pi_{k}\left(B^{\ast}\right) \longrightarrow \Gamma \longrightarrow 1
\]
est scind{\'e}e (o{\`u} $B^{\ast}=B-D$, $D$ {\'e}tant le lieu de ramification du $G$-rev{\^e}tement). C'est en particulier le cas lorsque la base du rev{\^e}tement poss{\`e}de un point $k$-rationnel (en dehors du lieu de ramification). L'objectif de ce chapitre est d'obtenir le m{\^e}me type d'{\'e}nonc{\'e} pour les torseurs sous un sch{\'e}ma en groupes ab{\'e}liens.

On se place donc dans la situation suivante: $k$ est un corps de caract{\'e}ristique nulle, dont on fixe une cl{\^o}ture alg{\'e}brique $\bar{k}$ et dont on note $\Gamma$ le groupe de Galois absolu; $X$ est un $k$-sch{\'e}ma g{\'e}om{\'e}triquement connexe, quasi-compact et quasi-s{\'e}par{\'e}, $\pi:X\rightarrow \textup{Spec}\;k$ est le morphisme structural, et $G$ est un $k$-groupe alg{\'e}brique lin{\'e}aire ab{\'e}lien.

Une premi{\`e}re id{\'e}e consiste {\`a} utiliser la suite spectrale de Leray\index{suite spectrale!de Leray} attach{\'e}e {\`a} cette situation:
	\[E_{2}^{p,q}=H^{p}\left(k,R^{q}\pi_{\ast}G_{X}\right)\Longrightarrow H^{p+q}\left(X,G_{X}\right)=E^{p+q}
\]

Dans ce que nous serons amen{\'e}s {\`a} consid{\'e}rer comme les \textquotedblleft{bons cas}\textquotedblright\ (\textit{i.e.} lorsque la condition $\bar{G}_{X}\left(\bar{X}\right)=\bar{G}\left(\bar{k}\right)$ est satisfaite), la suite exacte {\`a} 5 termes associ{\'e}e {\`a} la suite spectrale ci-dessus s'{\'e}crit:
	\[0\longrightarrow H^{1}\left(k,G\right)\longrightarrow H^{1}\left(X,G_{X}\right)\stackrel{u}{\longrightarrow} H^{1}\left(\bar{X},\bar{G}_{X}\right)^{\Gamma}\stackrel{\delta^{1}}{\longrightarrow}H^{2}\left(k,G\right)\stackrel{v}{\longrightarrow}H^{2}\left(X,G_{X}\right)
\]

Cette suite est {\'e}videmment parfaitement adapt{\'e}e {\`a} notre probl{\`e}me de descente, puisqu'on y lit directement l'obstruction {\`a} ce qu'un $\bar{G}_{X}$-torseur $\bar{P}$ sur $\bar{X}$ de corps des modules $k$ soit d{\'e}fini sur $k$. Plus pr{\'e}cis{\'e}ment, le morphisme $\delta^{1}$ d{\'e}fini en alg{\`e}bre homologique a une interpr{\'e}tation en termes de gerbes: c'est celui qui associe {\`a} une classe $\left[\bar{P}\right]$ la classe d'{\'e}quivalence $\left[D\left(\bar{P}\right)\right]$ de la gerbe d'un quelconque de ses repr{\'e}sentants (\textit{cf.} \cite{Gi2} V.3.1.4.1). Dire que $\bar{P}$ est d{\'e}fini sur $k$, c'est exactement dire que $\left[\bar{P}\right]$ appartient {\`a} l'image de $u$, ce qui {\'e}quivaut donc {\`a} la nullit{\'e} de la gerbe $D\left(\bar{P}\right)\in H^{2}\left(k,G\right)$. Nous verrons que l'existence d'un point $k$-rationnel sur $X$ entra{\^i}ne alors que tout $\bar{G}_{X}$ sur $\bar{X}$ de corps des modules $k$ est d{\'e}fini sur $k$, ou de fa\c{c}on {\'e}quivalente, que le morphisme $u$ est surjectif. Notons que cette approche fournit, {\`a} peu de frais, des r{\'e}sultats sur la descente des $G$-rev{\^e}tements ab{\'e}liens.

Dans la deuxi{\`e}me section, on s'int{\'e}resse {\`a} ce qui se passe en g{\'e}n{\'e}ral, si l'on n'impose plus la condition $\bar{G}_{X}\left(\bar{X}\right)=\bar{G}\left(\bar{k}\right)$. Nous montrons que le point-clef est finalement la comparaison entre les faisceaux (sur $k$) $G$ et $\pi_{\ast}G_{X}$. Pour illustrer les diff{\'e}rences entre cette situation et celle de la premi{\`e}re section, on s'int{\'e}resse en particulier {\`a} une vari{\'e}t{\'e} $X$ telle que:
	\[\mathbb{G}_{m,\bar{X}}\left(\bar{X}\right)\neq\mathbb{G}_{m,\bar{k}}\left(\bar{k}\right)
\]
Plus pr{\'e}cis{\'e}ment, on a:
	\[\bar{k}\left[X\right]^{\ast}=\bar{k}^{\ast}\oplus \mathbb{Z}
\]

Dans cette situation particuli{\`e}re, l'existence d'un point $k$-rationnel ne suffit pas {\`a} ce que tout $\mathbb{G}_{m,\bar{X}}$-torseur sur $\bar{X}$ de corps des modules $k$ soit d{\'e}fini sur $k$.

Enfin, nous nous pla\c{c}ons dans la derni{\`e}re section sur un corps de nombres, et nous utilisons un r{\'e}sultat de Skorobogatov qui assure que la descente des torseurs sous des groupes ab{\'e}liens sur des \textquotedblleft{bonnes vari{\'e}t{\'e}s}\textquotedblright\ est possible d{\`e}s qu'il existe sur $X$ des points ad{\'e}liques d'un certain type, ce qui est plus faible que de demander l'existence de points $k$-rationnels.
\begin{section}{Cons{\'e}quences de la suite spectrale de Leray}
Dans cette section, $k$ d{\'e}signe un corps de caract{\'e}ristique nulle, dont on choisit une cl{\^o}ture alg{\'e}brique $\bar{k}$; on note $\Gamma=\textup{Gal}\left(\bar{k}/k\right)$ le groupe de Galois absolu de $k$. On consid{\`e}re $G$ un $k$-groupe alg{\'e}brique ab{\'e}lien. Enfin on se donne:
	\[\pi:X\longrightarrow \textup{Spec}\;k
\]
un $k$-sch{\'e}ma g{\'e}om{\'e}triquement connexe, quasi-compact et quasi-s{\'e}par{\'e}, et on suppose satisfaite la condition:
	\[\bar{G}_{X}\left(\bar{X}\right)=\bar{G}\left(\bar{k}\right)
\]

Rappelons quelques d{\'e}finitions:
\begin{defi} Un morphisme de sch{\'e}mas $f:X\rightarrow Y$ est dit:
\begin{enumerate}[(i)]
\item \textbf{quasi-compact} si pour tout ouvert quasi-compact $U$ de $Y$, l'image r{\'e}ciproque $f^{-1}\left(U\right)$ est quasi-compacte;\index{k-sch{\'e}ma@$k$-sch{\'e}ma!quasi-compact}
\item \textbf{quasi-s{\'e}par{\'e}} \index{k-sch{\'e}ma@$k$-sch{\'e}ma!quasi-s{\'e}par{\'e}}(resp. \textbf{s{\'e}par{\'e}}) si le morphisme diagonal
	\[\Delta_{f}:X\times_{Y}X \longrightarrow Y
\]
est quasi-compact (resp. une immersion ferm{\'e}e).
\end{enumerate}

Un $k$-sch{\'e}ma est dit quasi-compact (resp. quasi-s{\'e}par{\'e}, resp. s{\'e}par{\'e}) si le morphisme structural $X\rightarrow \textup{Spec}\;k$ est quasi-compact (resp. quasi-s{\'e}par{\'e}, resp. s{\'e}par{\'e}).
\end{defi}

\begin{exem} \textup{Un morphisme s{\'e}par{\'e} est quasi-s{\'e}par{\'e}, puisqu'une immersion ferm{\'e}e est quasi-compacte; un morphisme noeth{\'e}rien est quasi-compact (\textit{cf.} \cite{EGA1} I.6.1.9); un morphisme affine est quasi-compact et s{\'e}par{\'e} (\textit{cf.} \cite{EGA1} I.9.1.3); le caract{\`e}re quasi-compact (\textit{resp.} quasi-s{\'e}par{\'e}) est stable par composition et par changement de base quelconque (\textit{cf.} \cite{EGA1} I.6.1.5, \textit{resp.} I.6.1.9);\ldots}
\end{exem}

Avec ces hypoth{\`e}ses sur $X$ et $G$, on a une suite exacte {\`a} 5 termes:
\begin{flushleft}
$\left(S1\right):\ \ 0\longrightarrow H^{1}\left(k,\pi_{\ast}G_{X}\right)\longrightarrow H^{1}\left(X,G_{X}\right)\stackrel{u}{\longrightarrow} H^{0}\left(k,R^{1}\pi_{\ast}G_{X}\right)$
\end{flushleft}
\begin{flushright}
$\stackrel{\delta^{1}}{\longrightarrow}H^{2}\left(k,\pi_{\ast}G_{X}\right)\stackrel{v}{\longrightarrow}H^{2}\left(X,G_{X}\right)$
\end{flushright}
qui est la suite exacte en basses dimensions associ{\'e}e {\`a} la suite spectrale de Leray:
	\[E_{2}^{p,q}=R^{p}\underline{\Gamma}_{\ k}\left(R^{q}\pi_{\ast}G_{X}\right)\Longrightarrow R^{p+q}\underline{\Gamma}_{\ X}\left(G_{X}\right)=E^{p+q}
\]

o{\`u}:
	\[\underline{\Gamma}_{\ k}:FAGRAB\left(k\right)\longrightarrow \mathfrak{Ab}\ \ \left(\textit{resp.}\ \underline{\Gamma}_{\ X}:FAGRAB\left(X\right)\longrightarrow \mathfrak{Ab}\right)
\]
est le foncteur qui associe {\`a} un faisceau de groupes ab{\'e}liens sur le site {\'e}tale de $k$ (\textit{resp.} de $X$) le groupe ab{\'e}lien de ses sections globales. Comme par d{\'e}finition les foncteurs d{\'e}riv{\'e}s {\`a} droite de ces foncteurs sont justement les foncteurs cohomologie, on peut r{\'e}{\'e}crire la suite spectrale ci-dessus:
	\[E_{2}^{p,q}=H^{p}\left(k,R^{q}\pi_{\ast}G_{X}\right)\Longrightarrow H^{p+q}\left(X,G_{X}\right)=E^{p+q}
\]

On peut maintenant rendre plus agr{\'e}able l'expression de la suite (S1) gr{\^a}ce au th{\'e}or{\`e}me 5.2 de \cite{SGA4-8} dont voici l'{\'e}nonc{\'e}:
\begin{theo} Soient $f:Z\rightarrow Y$ un morphisme quasi-compact et quasi-s{\'e}par{\'e} de sch{\'e}mas, $G$ un faisceau ab{\'e}lien sur $Z$, $y$ un point de $Y$, $\bar{y}$ le point g{\'e}om{\'e}trique au-dessus de $y$, relatif {\`a} une cl{\^o}ture s{\'e}parable $k\left(\bar{y}\right)$ de $k\left(y\right)$, $\bar{Y}=\textup{Spec}\left(\mathcal{O}_{Y,\bar{y}}\right)$ le sch{\'e}ma localis{\'e} strict correspondant, $\bar{Z}=Z\times_{Y}\bar{Y}$, $\bar{G}$ l'image inverse de $G$ sur $\bar{Z}$. Alors l'homomorphisme canonique:
	\[\left(R^{q}f_{\ast}G\right)_{\bar{y}}\longrightarrow H^{q}\left(\bar{Z},\bar{G}\right)
\]
est un isomorphisme, pour tout $q\geq0$.
\end{theo}
\begin{rem}\textup{Cet {\'e}nonc{\'e} reste valable (d'apr{\`e}s la remarque 5.3 de \cite{SGA4-8}) pour un faisceau de groupes $G$ sur $Z$ non-n{\'e}cessairement ab{\'e}lien, pour $q=0$ et $q=1$, en prenant pour $R^{1}f_{\ast}G$ la d{\'e}finition de \cite{Gi2} V.2.1: c'est le faisceau sur $Y$ associ{\'e} au pr{\'e}faisceau dont l'ensemble des sections $R^{1}f_{\ast}G\left(Y'\right)$ au-dessus d'un ouvert {\'e}tale $\left(Y'\rightarrow Y\right)$ est donn{\'e} par:
	\[R^{1}f_{\ast}G\left(Y'\right)=H^{1}\left(Z\times_{Y}Y',G_{\left|Z\times_{Y}Y'\right.}\right)
\]
}
\begin{coro} Sous les hypoth{\`e}ses du d{\'e}but de cette section, on a des isomorphismes:\label{R1Pi}
	\[\left(\pi_{\ast}G_{X}\right)_{\textup{Spec}\;\bar{k}}\approx \bar{G}\left(\bar{k}\right)
\]
	\[\left(R^{1}\pi_{\ast}G_{X}\right)_{\textup{Spec}\;\bar{k}}\approx H^{1}\left(\bar{X},\bar{G}_{X}\right)
\]
\end{coro}\pagebreak

\uline{\textsc{Preuve}}: \textup{il suffit d'appliquer le th{\'e}or{\`e}me ci-dessus avec $Z=X$, $Y=\textup{Spec}\;k$, $f=\pi$ et $G=G_{X}$, et d'utiliser (pour obtenir le premier isomorphisme) l'hypoth{\`e}se $\bar{G}_{X}\left(\bar{X}\right)=\bar{G}\left(\bar{k}\right)$.}
\begin{flushright}
$\Box$
\end{flushright}
\end{rem}

Par cons{\'e}quent, on peut r{\'e}{\'e}crire la suite (S1):
	\[\left(S2\right):\ 0\longrightarrow H^{1}\left(k,G\right)\longrightarrow H^{1}\left(X,G_{X}\right)\stackrel{u}{\longrightarrow} H^{1}\left(\bar{X},\bar{G}_{X}\right)^{\Gamma}\stackrel{\delta^{1}}{\longrightarrow}H^{2}\left(k,G\right)\stackrel{v}{\longrightarrow}H^{2}\left(X,G_{X}\right)
\]
\begin{rem}\textup{De fait, l'image du morphisme $v$ est incluse dans la partie transgressive $H^{2}\left(X,G_{X}\right)^{alg}$\label{H2trans}, qui est le noyau du morphisme {\'e}vident:
	\[H^{2}\left(X,G_{X}\right)\longrightarrow H^{2}\left(\bar{X},\bar{G}_{X}\right)
\]
}
\end{rem}

Ces consid{\'e}rations nous am{\`e}nent {\`a} introduire une nouvelle notation:
\begin{defi} Soit $k$ un corps. On appelle \textbf{$k$-sch{\'e}ma de type $\left(\ast\right)$}\index{k-sch{\'e}ma@$k$-sch{\'e}ma!de type a@de type $\left(\ast\right)$} un $k$-sch{\'e}ma g{\'e}om{\'e}triquement connexe, quasi-compact et quasi-s{\'e}par{\'e}.
\end{defi}

Une premi{\`e}re cons{\'e}quence de l'exactitude de la suite $\left(S2\right)$ est la suivante:
\begin{theo}[Obstruction ab{\'e}lienne {\`a} l'existence d'un point rationnel]\index{obstruction!ab{\'e}lienne} Soient $k$ un corps de caract{\'e}ristique nulle, $X$ un $k$-sch{\'e}ma de type $\left(\ast\right)$, et $G$ un $k$-groupe alg{\'e}brique ab{\'e}lien tels que $\bar{G}_{X}\left(\bar{X}\right)=\bar{G}\left(\bar{k}\right)$.

Si $X\left(k\right)\neq\emptyset$, alors tout $\bar{G}_{X}$-torseur sur $\bar{X}$ de corps des modules $k$ est d{\'e}fini sur $k$.
\end{theo}

\uline{\textsc{Preuve}}: \textup{l'existence d'un point $k$-rationnel sur $X$ entra{\^i}ne l'existence d'une r{\'e}traction du morphisme $v:H^{2}\left(k,G\right)\rightarrow H^{2}\left(X,G_{X}\right)$. Donc $v$ est injectif\footnote{Puisque dans la cat{\'e}gorie des groupes ab{\'e}liens, il est {\'e}quivalent de dire qu'un morphisme est injectif ou qu'il poss{\`e}de une r{\'e}traction.}, donc le cobord $\delta^{1}$ de la suite $\left(S2\right)$ est nul, donc le morphisme $u$ est surjectif, d'o{\`u} la conclusion.}
\begin{flushright}
$\Box$
\end{flushright}

\begin{rem} \textup{Il revient au m{\^e}me de dire, avec les notations et hypoth{\`e}ses du th{\'e}or{\`e}me, l'existence d'un point $k$-rationnel sur $X$ entra{\^i}ne l'existence d'un point $k$-rationnel de la gerbe des mod{\`e}les de tout $\bar{G}_{X}$-torseur sur $\bar{X}$ de corps des modules $k$.}
\end{rem}
Une derni{\`e}re mani{\`e}re de traduire l'{\'e}nonc{\'e} pr{\'e}c{\'e}dent est que lorsque $X$ poss{\`e}de un point rationnel, il n'existe pas d'obstruction {\`a} la descente des $\bar{G}_{X}$-torseurs de corps des modules $k$.

En outre, on peut prolonger la suite exacte $\left(S2\right)$: en effet, on peut associer {\`a} toute suite spectrale:
	\[E_{2}^{p,q}\Longrightarrow E^{p+q}
\]
une suite exacte:
	\[0\longrightarrow E_{2}^{1,0}\longrightarrow E^{1}\longrightarrow E_{2}^{0,1}\longrightarrow E_{2}^{2,0}\longrightarrow E^{2,tr}\longrightarrow E_{2}^{1,1}\longrightarrow E_{2}^{3,0}
\]
o{\`u}:
	\[E^{2,tr}=\ker\left\{E^{2}\longrightarrow E^{0,2}_{2}\right\}
\]

Sous les hypoth{\`e}ses du th{\'e}or{\`e}me 2.1.8, on obtient ainsi la suite exacte:
\begin{flushleft}
$\left(S3\right):\ 0\longrightarrow H^{1}\left(k,G\right)\longrightarrow H^{1}\left(X,G_{X}\right)\stackrel{u}{\longrightarrow} H^{1}\left(\bar{X},\bar{G}_{X}\right)^{\Gamma}\stackrel{\delta^{1}}{\longrightarrow}H^{2}\left(k,G\right)$
\end{flushleft}
\begin{flushright}
$\stackrel{v}{\longrightarrow}H^{2}\left(X,G_{X}\right)^{alg}\stackrel{w}{\longrightarrow} H^{1}\left(k,H^{1}\left(\bar{X},\bar{G}_{X}\right)\right)\stackrel{\delta^{2}}{\longrightarrow}H^{3}\left(k,G\right)$
\end{flushright}

On en d{\'e}duit donc la:
\begin{pro} Soient $k$ un corps de caract{\'e}ristique nulle, $X$ un $k$-sch{\'e}ma de type $\left(\ast\right)$, et $G$ un $k$-groupe alg{\'e}brique tels que $\bar{G}_{X}\left(\bar{X}\right)=\bar{G}\left(\bar{k}\right)$.

Si $X\left(k\right)\neq\emptyset$ alors:
\begin{enumerate}[(i)]
\item la suite:
	\[0\longrightarrow H^{1}\left(k,G\right)\longrightarrow H^{1}\left(X,G_{X}\right)\stackrel{u}{\longrightarrow}H^{1}\left(\bar{X},\bar{G}_{X}\right)^{\Gamma}\longrightarrow 0
\]
est exacte;
\item la suite:
	\[0\longrightarrow H^{2}\left(k,G\right) \stackrel{v}{\longrightarrow}H^{2}\left(X,G_{X}\right)^{alg}\stackrel{w}{\longrightarrow} H^{1}\left(k,H^{1}\left(\bar{X},\bar{G}_{X}\right)\right)\longrightarrow 0
\]
est exacte. En particulier:
	\[H^{1}\left(k,H^{1}\left(\bar{X},\bar{G}_{X}\right)\right)\approx \frac{H^{2}\left(X,G_{X}\right)^{alg}}{H^{2}\left(k,G\right)}
\]
\item le morphisme $z:H^{3}\left(k,G\right)\longrightarrow H^{3}\left(X,G_{X}\right)$ est injectif.
\end{enumerate}
\end{pro}

\uline{\textsc{Preuve}}: \textup{l'existence d'un point $k$-rationnel sur $X$ entra{\^i}ne l'existence de r{\'e}tractions des morphismes $v$ et $z$, d'o{\`u} la conclusion, en utilisant l'exactitude de la suite $\left(S3\right)$.}
\begin{flushright}
$\Box$
\end{flushright}

Nous donnons maintenant des exemples d'applications de ces propri{\'e}t{\'e}s aux groupes de Picard et de Brauer d'une $k$-vari{\'e}t{\'e}, ainsi qu'aux $G$-rev{\^e}tements ab{\'e}liens.

\begin{quote}
	$\bullet$ \uline{Application aux groupes de Picard et de Brauer}
\end{quote}

On s'int{\'e}resse donc ici au cas particulier o{\`u} $G=\mathbb{G}_{m,k}$. De mani{\`e}re {\`a} pouvoir utiliser les r{\'e}sultats pr{\'e}c{\'e}dents, on souhaite voir remplie la condition:
	\[\mathbb{G}_{m,\bar{X}}\left(\bar{X}\right)=\mathbb{G}_{m,\bar{k}}\left(\bar{k}\right)
\]
c'est-{\`a}-dire:
	\[\bar{k}\left[X\right]^{\ast}=\bar{k}^{\ast}
\]

Pour ces applications, nous consid{\'e}rons donc une $k$-vari{\'e}t{\'e} $X$ \textbf{propre}\footnote{Mais les r{\'e}sultats obtenus ici sont encore valables pour une vari{\'e}t{\'e} $X$ telle que: $\bar{k}\left[X\right]^{\ast}=\bar{k}^{\ast}$; une telle vari{\'e}t{\'e} n'est pas n{\'e}cessairement propre (on peut par exemple penser {\`a} l'espace affine $\mathbb{A}^{n}_{k}$).}. En {\'e}crivant alors la suite $\left(S3\right)$ avec $G=\mathbb{G}_{m,k}$, on obtient la suite exacte:
\begin{flushleft}
$\left(S4\right):\ \ 0\longrightarrow \textup{Pic}\;X\longrightarrow \left(\textup{Pic}\;\bar{X}\right)^{\Gamma}\longrightarrow \textup{Br}\;k\longrightarrow \textup{Br}^{alg}X$
\end{flushleft}
\begin{flushright}
$\longrightarrow H^{1}\left(k,\textup{Pic}\;\bar{X}\right)\longrightarrow H^{3}\left(k,\mathbb{G}_{m,k}\right)$
\end{flushright}
o{\`u}: $\textup{Br}^{alg}X=\ker\left\{\textup{Br}\;X\rightarrow \textup{Br}\;\bar{X}\right\}$\label{Brtr} est le groupe de Brauer \textit{transgressif}\index{groupe!de Brauer!transgressif} de $X$. De cette suite exacte, on d{\'e}duit imm{\'e}diatemment la:
\begin{pro} Si $X$ est une $k$-vari{\'e}t{\'e} propre, et si $X\left(k\right)\neq\emptyset$, alors:
\begin{enumerate}[(i)]
\item tout fibr{\'e} en droites sur $\bar{X}$ de corps des modules $k$ est d{\'e}fini sur $k$;
\item la suite:
	\[0\longrightarrow \textup{Br}\;k\longrightarrow \textup{Br}^{alg}X\longrightarrow H^{1}\left(k,\textup{Pic}\;\bar{X}\right)\longrightarrow 0
\]
est exacte;
\item le morphisme $H^{3}\left(k,\mathbb{G}_{m,k}\right)\longrightarrow H^{3}\left(X,\mathbb{G}_{m,X}\right)$ est injectif.
\end{enumerate}
\end{pro}

\uline{\textsc{Preuve}}: \textup{c'est une cons{\'e}quence imm{\'e}diate de la proposition 2.1.10.}
\begin{flushright}
$\Box$
\end{flushright}

Notons tout de suite que le $\left(i\right)$ de la proposition ci-dessus ne peut pas fournir r{\'e}ellement d'obstruction {\`a} l'existence d'un point $k$-rationnel lorsque $k$ est un corps de nombres. Plus explicitement:
\begin{pro} Soient $k$ un corps de nombres, $\mathbb{A}_{k}$ son anneau des ad{\`e}les, et soit $X$ une $k$-vari{\'e}t{\'e} propre. Si $X\left(\mathbb{A}_{k}\right)\neq\emptyset$, alors le morphisme:
	\[\textup{Pic}\;X\longrightarrow \left(\textup{Pic}\;\bar{X}\right)^{\Gamma}
\]
est un isomorphisme.
\end{pro}

\uline{\textsc{Preuve}}: soit $\bar{L}$ un fibr{\'e} en droites sur $\bar{X}$ de corps des modules $k$. La gerbe des mod{\`e}les $D\left(\bar{L}\right)$ vit dans $\textup{Br}\;k$. Mais, comme $X$ a un $k_{v}$-point pour toute place $v$ de $k$:
	\[D\left(\bar{L}\right)\otimes_{k}k_{v}=\textup{loc}_{v}\left(D\left(\bar{L}\right)\right)
\]
est neutre, pour toute place $v$ de $k$ ($\textup{loc}_{v}:\textup{Br}\;k\rightarrow\textup{Br}\;k_{v}$ {\'e}tant le morphisme {\'e}vident). En particulier, $D\left(\bar{L}\right)$ est dans le noyau de l'application $\textup{loc}$ intervenant dans la c{\'e}l{\`e}bre suite exacte:
	\[\xymatrix@C=40pt{0\ar[r]&\textup{Br}\;k\ar[r]^{\textup{loc}\ \ \ \ \ }&\displaystyle\bigoplus_{v\in\Omega_{k}}\textup{Br}\;k_{v}\ar[r]^{\ \ \ \ \ \sum\textup{inv}_{v}}&\mathbb{Q}/\mathbb{Z}\ar[r]& 0}
\]

Donc $D\left(\bar{L}\right)=0$, donc $\bar{L}$ est d{\'e}fini sur $k$.
\begin{flushright}
$\Box$
\end{flushright}

L'exactitude de la suite $\left(S4\right)$ a aussi des cons{\'e}quences sur le calcul de l'obstruction de Brauer-Manin\footnote{Pour la d{\'e}finition et la construction de cette obstruction, nous renvoyons {\`a} la section 2.3.1 du pr{\'e}sent chapitre.} d'une vari{\'e}t{\'e}. Explicitement, supposons toujours $X$ propre: l'obstruction de Brauer-Manin de $X$ associ{\'e}e {\`a} $\mathcyr{B}\left(X\right)$ (en reprenant la terminologie de \cite{CTS}, d{\'e}finition 3.1.1), que nous noterons $m_{\mathcal{H},\mathcyr{B}\left(X\right)}\left(X\right)$, est un {\'e}l{\'e}ment de 
	\[\mathcyr{B}\left(X\right)^{D}=\textup{Hom}\left(\mathcyr{B}\left(X\right),\mathbb{Q}/\mathbb{Z}\right)
\]
o{\`u} $\mathcyr{B}\left(X\right)$ est le groupe construit {\`a} partir du groupe de Brauer de la fa\c{c}on suivante: on commence par poser:\label{Bra}
	\[\textup{Br}_{a}X=\frac{\textup{Br}^{alg}X}{\textup{im}\left(\textup{Br}\;k\rightarrow \textup{Br}\;X\right)}
\]
puis on d{\'e}finit $\mathcyr{B}\left(X\right)$ (\textit{resp.} $\mathcyr{SH}^{1}\left(k,\textup{Pic}\;\bar{X}\right)$) comme le noyau de l'application de localisation:
	\[\textup{Br}_{a}X\rightarrow\displaystyle\prod_{v\in\Omega_{k}}\textup{Br}_{a}\left(X\otimes_{k}k_{v}\right)\ \left(\textit{resp. } H^{1}\left(k,\textup{Pic}\;\bar{X}\right)\rightarrow\displaystyle\prod_{v\in\Omega_{k}}H^{1}\left(k_{v},\textup{Pic}\;\bar{X}\right)\right)
\]

En utilisant l'exactitude de la suite $\left(S4\right)$ on obtient:
	\[\textup{Br}_{a}X\approx H^{1}\left(k,\textup{Pic}\;\bar{X}\right)
\]
donc en \textquotedblleft{passant aux noyaux sur toutes les places}\textquotedblright, on a:\label{BcyrX}\label{SHA1}
	\[\mathcyr{B}\left(X\right)\approx \mathcyr{SH}^{1}\left(k,\textup{Pic}\;\bar{X}\right)
\]

Le calcul de l'obstruction de Brauer-Manin de $X$ associ{\'e}e {\`a} $\mathcyr{B}\left(X\right)$ se ram{\`e}ne donc essentiellement au calcul de $\textup{Pic}\;\bar{X}$. D'o{\`u} par exemple la:
\begin{pro} Soient $k$ un corps de nombres et $X$ une $k$-vari{\'e}t{\'e} qui est:
\begin{enumerate}[(i)]
\item une vari{\'e}t{\'e} de Severi-Brauer;
\item une intersection compl{\`e}te lisse de dimension $\geq3$;
\end{enumerate}
alors:
	\[\mathcyr{B}\left(X\right)=0
\]

En particulier, l'obstruction de Brauer-Manin $m_{H,\mathcyr{B}\left(X\right)}\left(X\right)$\label{mHX} d'une telle vari{\'e}t{\'e} est nulle.
\end{pro}

\uline{\textsc{Preuve}}: elle repose sur le fait que dans les deux cas, on a $\textup{Pic}\;\bar{X}=\mathbb{Z}$.\footnote{C'est trivial pour les vari{\'e}t{\'e}s de Severi-Brauer, et nous renvoyons {\`a} la d{\'e}monstration de la proposition suivante pour le cas des intersections compl{\`e}tes.} Donc $H^{1}\left(k,\textup{Pic}\;\bar{X}\right)=0$, et \textit{a fortiori} $\mathcyr{SH}^{1}\left(k,\textup{Pic}\;\bar{X}\right)=0$. D'o{\`u} la conclusion, puisque $m_{H\mathcyr{B}\left(X\right)}\left(X\right)\in\mathcyr{B}\left(X\right)^{D}$ et $\mathcyr{B}\left(X\right)\approx \mathcyr{SH}^{1}\left(k,\textup{Pic}\;\bar{X}\right)$.
\begin{flushright}
$\Box$
\end{flushright}

Cet {\'e}nonc{\'e} n'est toutefois absolument pas surprenant dans la mesure o{\`u} le groupe de Brauer de telles vari{\'e}t{\'e}s est \textquotedblleft{trivial}\textquotedblright\ dans le sens suivant:
\begin{pro} La \textquotedblleft{conjecture de Grothendieck}\textquotedblright\ est \index{conjecture de Grothendieck sur les groupes de Brauer@Conjecture de Grothendieck sur les groupes de Brauer}\textquotedblleft{trivialement vraie}\textquotedblright\ pour les vari{\'e}t{\'e}s de Severi-Brauer et pour les vari{\'e}t{\'e}s intersections compl{\`e}tes lisses dans $\mathbb{P}^{n}$ de dimension $\geq3$ sur un corps $k$ de caract{\'e}ristique nulle. Plus pr{\'e}cis{\'e}ment, si $X$ est une telle vari{\'e}t{\'e}, alors:
	\[\textup{Br}_{Az}X\approx \textup{Br}\;k\approx \textup{Br}\;X
\]
$\textup{Br}_{Az}X$ d{\'e}signant le groupe des classes d'isomorphie d'alg{\`e}bres d'Azumaya sur $X$.
\end{pro}

\uline{\textsc{Preuve}}: si $X$ est une vari{\'e}t{\'e} de Severi-Brauer: alors il existe un entier $n$ tel que: $\bar{X}\approx\mathbb{P}^{n}_{\bar{k}}$. Par suite, $\textup{Pic}\;\bar{X}=\mathbb{Z}$, donc $H^{1}\left(k,\textup{Pic}\;\bar{X}\right)=0$, donc $\textup{Br}\;k=\textup{Br}^{alg}X$. De plus, $\textup{Br}\;\bar{X}=0$, car le groupe de Brauer d'un espace projectif sur un corps alg{\'e}briquement clos est nul \cite{G4}. Donc $\textup{Br}^{alg}X=\textup{Br}\;X$.

Si $X$ est une intersection compl{\`e}te lisse dans $\mathbb{P}^{n}$ de dimension $\geq3$: commen\c{c}ons par montrer: $\textup{Pic}\;\bar{X}=\mathbb{Z}$. D'apr{\`e}s le principe de Lefschetz\index{principe!de Lefschetz} (\textit{cf.} \cite{Harris} 15.1), $k$ {\'e}tant de caract{\'e}ristique nulle, on peut supposer $\bar{k}=\mathbb{C}$. Supposons dans un premier temps que $X$ est une hypersurface: on a le diagramme commutatif et {\`a} lignes exactes suivant:
	\[\xymatrix@C=25pt@R=25pt{\ldots \ar[r] &H^{1}\left(\bar{X},\mathcal{O}_{\bar{X}}\right) \ar[r] & H^{1}\left(\bar{X},\mathcal{O}_{\bar{X}}^{\ast}\right) \ar[r] & H^{2}\left(\bar{X},\mathbb{Z}\right) \ar[r] & H^{2}\left(\bar{X},\mathcal{O}_{\bar{X}}\right) \ar[r] & \ldots \\ \ldots \ar[r] &H^{1}\left(\mathbb{P}^{n},\mathcal{O}_{\mathbb{P}^{n}}\right) \ar[r] \ar[u] & H^{1}\left(\mathbb{P}^{n},\mathcal{O}_{\mathbb{P}^{n}}^{\ast}\right) \ar[r] \ar[u] & H^{2}\left(\mathbb{P}^{n},\mathbb{Z}\right) \ar[r] \ar[u] & H^{2}\left(\mathbb{P}^{n},\mathcal{O}_{\mathbb{P}^{n}}\right) \ar[r] \ar[u] & \ldots}
\]
les lignes {\'e}tant obtenues {\`a} partir de la suite exponentielle, et les colonnes {\`a} partir de l'inclusion: $\bar{X}\hookrightarrow\mathbb{P}^{n}$. Pour des raisons {\'e}videntes de dimension\footnote{Pour $H^{i}\left(\mathbb{P}^{n},\mathcal{O}_{\mathbb{P}^{n}}\right)$, $i=1,2$, c'est exactement le (b) du th{\'e}or{\`e}me III.5.1 de \cite{Har}; pour $H^{i}\left(\bar{X},\mathcal{O}_{\bar{X}}\right)$, $i=1,2$, c'est le (c) de l'exercice III.5.5 de \textit{loc. cit}.}, les groupes:
\begin{center}
$H^{1}\left(\bar{X},\mathcal{O}_{\bar{X}}\right)$, $H^{2}\left(\bar{X},\mathcal{O}_{\bar{X}}\right)$, $H^{1}\left(\mathbb{P}^{n},\mathcal{O}_{\mathbb{P}^{n}}\right)$ et $H^{2}\left(\mathbb{P}^{n},\mathcal{O}_{\mathbb{P}^{n}}\right)$
\end{center}
sont nuls. D'autre part, le th{\'e}or{\`e}me de la section hyperplane de Lefschetz (\textit{cf.} \cite{GH} p.156) assure que le morphisme:
	\[H^{2}\left(\mathbb{P}^{n},\mathbb{Z}\right)\longrightarrow H^{2}\left(\bar{X},\mathbb{Z}\right)
\]
est un isomorphisme. Donc:
	\[H^{1}\left(\bar{X},\mathcal{O}^{\ast}_{\bar{X}}\right)\approx H^{1}\left(\mathbb{P}^{n},\mathcal{O}^{\ast}_{\mathbb{P}^{n}}\right)
\]
c'est-{\`a}-dire:
	\[\textup{Pic}\;\bar{X}\approx\mathbb{Z}
\]

Pour obtenir cette propri{\'e}t{\'e} dans le cas o{\`u} $\bar{X}$ est une intersection compl{\`e}te lisse dans $\mathbb{P}^{n}$ (et non plus seulement une hypersurface), il suffit d'appliquer le th{\'e}or{\`e}me de la section hyperplane \textquotedblleft{suffisamment}\textquotedblright\ de fois (pr{\'e}cis{\'e}ment $n-\dim\bar{X}$ fois). D'apr{\`e}s la suite (S4), on a d{\'e}j{\`a}:
	\[\textup{Br}^{alg}X=\textup{Br}\;k
\]

La conclusion provient alors du fait que le morphisme naturel
	\[\textup{Br}_{Az}X\longrightarrow \textup{Br}\;X
\]
est toujours injectif, et de ce que $\textup{Br}\;\bar{X}=0$ dans ce cas \cite{Ma2}.
\begin{flushright}
$\Box$
\end{flushright}

\begin{quote}
	$\bullet$ \uline{Application aux $G$-rev{\^e}tements ($G$ ab{\'e}lien)}
\end{quote}

On consid{\`e}re $X$ une $k$-vari{\'e}t{\'e} projective g{\'e}om{\'e}triquement irr{\'e}ductible et $G$ un groupe fini ab{\'e}lien. Dans ce cas, $G$-rev{\^e}tements {\'e}tales et $G$-torseurs coïncident, et l'obstruction {\`a} ce qu'un $G$-rev{\^e}tement $\bar{f}:\bar{Y}\rightarrow \bar{X}$ de corps des modules $k$ soit d{\'e}fini sur $k$ est mesur{\'e} par une gerbe $\mathcal{G}\left(\bar{f}\right)$ vivant dans $H^{2}\left(k,G\right)$. Une cons{\'e}quence presque imm{\'e}diate de la suite spectrale de Leray est l'{\'e}nonc{\'e} suivant (qui est {\`a} rapprocher du th{\'e}or{\`e}me de Combes-Harbater):
\begin{theo} Avec les notations introduites plus haut, si $X$ poss{\`e}de un point $k$-rationnel en dehors du lieu de ramification de $\bar{f}$, alors le $G$-rev{\^e}tement $\bar{f}$ est d{\'e}fini sur $k$.
\end{theo}
\uline{\textsc{Preuve}}: On note $X^{\ast}=X-R$, o{\`u} $R$ est le lieu de ramification de $\bar{f}$, $\bar{X}^{\ast}=\bar{X}-\bar{R}$, et $\bar{f}^{\ast}:\bar{Y}^{\ast}\rightarrow \bar{X}^{\ast}$ la restriction de $\bar{f}$ {\`a} $\bar{X}^{\ast}$. La condition $\bar{G}\left(\bar{X}\right)=\bar{G}\left(\bar{k}\right)$ est satisfaite, car $G$ est fini, et la conclusion provient alors de ce que l'existence d'un point $k$-rationnel sur $X^{\ast}$ entra{\^i}ne l'existence d'une r{\'e}traction du morphisme $v$ dans la suite exacte:
	\[0\longrightarrow H^{1}\left(k,G\right)\longrightarrow H^{1}\left(X^{\ast},G\right)\stackrel{u}{\longrightarrow}H^{1}\left(\bar{X}^{\ast},G\right)^{\Gamma}\stackrel{\delta^{1}}{\longrightarrow}H^{2}\left(k,G\right)\stackrel{v}{\longrightarrow}H^{2}\left(X^{\ast},G\right)
\]
\begin{flushright}
$\Box$
\end{flushright}

Si l'on suppose maintenant $k=\mathbb{Q}$, et:
	\[X\left(A_{\mathbb{Q}}\right)\neq\emptyset
\]
c'est-{\`a}-dire si l'on suppose que $X$ a des points r{\'e}els et des points $p$-adiques pour tout nombre premier $p$, alors:
	\[\mathcal{G}\left(\bar{f}\right)\otimes_{\mathbb{Q}}\mathbb{R}=0
\]
et:
	\[\mathcal{G}\left(\bar{f}\right)\otimes_{\mathbb{Q}}\mathbb{Q}_{p}=0,\ \forall\;p\ premier
\]
c'est-{\`a}-dire:
	\[\mathcal{G}\left(\bar{f}\right)\in\mathcyr{SH}^{2}\left(\mathbb{Q},G\right)=0
\]

On retrouve ainsi le principe local-global de D{\`e}bes et Douai (\textit{cf.} \cite{DD2} theorem 3.8):
\begin{theo}[Principe local-global pour les $G$-rev{\^e}tements ab{\'e}liens] Avec les notations indiqu{\'e}es ci-dessus, un $G$-rev{\^e}tement $\bar{f}:\bar{Y}\rightarrow\bar{X}$ de corps des modules $Q$ est d{\'e}fini sur $\mathbb{Q}$ si et seulement si il est d{\'e}fini sur $\mathbb{R}$ et sur $\mathbb{Q}_{p}$ pour tout premier $p$. C'est en particulier le cas si $X$ poss{\`e}de des points ad{\'e}liques.
\end{theo}
\end{section}
\begin{section}{De l'importance de la condition $\bar{G}_{X}\left(\bar{X}\right)=\bar{G}\left(\bar{k}\right)$}
Nous commen\c{c}ons par un exemple (d{\^u} {\`a} J.-L. Colliot-Th{\'e}l{\`e}ne et O. Gabber) illustrant le caract{\`e}re indispensable de cette condition. Soit $k$ un corps poss{\'e}dant une extension cyclique non-triviale $L$. Il existe donc une classe $\alpha$ non-nulle dans $H^{2}\left(k,\mathbb{Z}\right)$. On consid{\`e}re le $k$-sch{\'e}ma:
	\[X=\mathbb{G}_{m}=\textup{Spec}\;k\left[T,T^{-1}\right]
\]

On fixe $\bar{k}$ une cl{\^o}ture s{\'e}parable de $k$, et on note comme d'habitude $\bar{X}=X\otimes_{k}\bar{k}$. Les fonctions inversibles sur $\bar{X}$ sont donn{\'e}es par:
	\[\bar{k}\left[X\right]^{\ast}\approx \bar{k}^{\ast}\oplus \mathbb{Z}
\]
le scindage {\'e}tant obtenu par {\'e}valuation en $T=1$:
\vspace{2mm}
	\[\xymatrix@R=5pt@C=35pt{0\ar[r] &\mathbb{Z}\ar[r]&\bar{k}\left[X\right]^{\ast}\ar[r]^{\ \ ev_{1}}&\bar{k}^{\ast}\ar@/_20pt/[l]\ar[r]&0\\&n\ar@{|->}[r]&T^{n}\\&&\mu T^{n}\ar@{|->}[r]&\mu}
\]
On en d{\'e}duit la suite exacte:
	\[H^{2}\left(k,\mathbb{Z}\right)\stackrel{i}{\longrightarrow}H^{2}\left(k,\bar{k}\left[X\right]^{\ast}\right)\longrightarrow \textup{Br}\;k
\]
Remarquons que le morphisme $i$ est injectif, car $H^{1}\left(k,\bar{k}^{\ast}\right)=0$. Notons:
	\[\beta=i\left(\alpha\right).
\]

$\beta$ est donc un {\'e}l{\'e}ment non-nul de $H^{2}\left(k,\bar{k}\left[X\right]^{\ast}\right)$. D'autre part, on a un morphisme\footnote{C'est l'edge $E^{2,0}_{2}\longrightarrow E^{2}$ de la suite spectrale de Leray:
	\[H^{p}\left(k,R^{q}\pi_{\ast}\mathbb{G}_{m,X}\right)\Longrightarrow H^{p+q}\left(X,\mathbb{G}_{m,X}\right).
\]
}:
	\[j:H^{2}\left(k,\bar{k}\left[X\right]^{\ast}\right)\longrightarrow \textup{Br}\;X
\]

Comme $X$ est affine, le monomorphisme de groupes:
	\[\Delta:\textup{Br}_{Az}X\longrightarrow \textup{Br}\;X
\]
obtenu en associant {\`a} une alg{\`e}bre d'Azumaya sur $X$ la gerbe de ses banalisations (\textit{cf.} \cite{Gi2} V.4.2) est un isomorphisme (\textit{cf.} \cite{Ga} thm. 1 p.163). Par suite il existe une un sch{\'e}ma de Severi-Brauer $Y\rightarrow X$ telle que:
	\[\Delta\left(Y\right)=j\left(\beta\right)
\]

Comme l'{\'e}valuation de $\beta$ en 1 est triviale, $Y$ poss{\`e}de un point $k$-rationnel au-dessus de $T=1$. En particulier:
	\[Y\left(k\right)\neq\emptyset
\]
De la suite spectrale de Leray et du morphisme $Y\rightarrow X$ on d{\'e}duit le diagramme commutatif et {\`a} lignes exactes suivant:
	\[\xymatrix{\textup{Pic}\;Y\ar[r]&\left(\textup{Pic}\;\bar{Y}\right)^{\Gamma}\ar[r]^{\delta_{Y}\ \ \ }&H^{2}\left(k,\bar{k}\left[Y\right]^{\ast}\right)\ar[r]^{\ \ \ \ \ u_{Y}}&\textup{Br}\;Y \\ \textup{Pic}\;X\ar[r]\ar[u]&\left(\textup{Pic}\;\bar{X}\right)^{\Gamma}\ar[r]^{\delta_{X}\ \ \ }\ar[u]&H^{2}\left(k,\bar{k}\left[X\right]^{\ast}\right)\ar[r]^{\ \ \ \ \ \ u_{X}}\ar[u]^{w}&\textup{Br}\;X\ar[u]}
\]

On a $\bar{k}\left[X\right]^{\ast}=\bar{k}\left[Y\right]^{\ast}$,\footnote{Car le morphisme $Y\rightarrow X$ est lisse avec des fibres propres et g{\'e}om{\'e}triquement int{\`e}gres. Donc une fonction inversible sur $Y$ a une restriction {\`a} la fibre g{\'e}n{\'e}rique qui est constante (je remercie D. Harari pour cet argument).} et $w$ est un isomorphisme. $\beta$ (identifi{\'e} {\`a} son image par $w$ dans $H^{2}\left(k,\bar{k}\left[Y\right]^{\ast}\right)$) est non-nul, mais $u_{Y}\left(\beta\right)=0$. Donc le morphisme $u_{Y}$ n'est pas injectif, malgr{\'e} l'existence d'un point $k$-rationnel. En particulier, le morphisme $\textup{Pic}\;Y\rightarrow\left(\textup{Pic}\;\bar{Y}\right)^{\Gamma}$ n'est pas surjectif.

\begin{pro} Pour le $k$-sch{\'e}ma $Y$ construit ci-dessus, le morphisme
	\[\textup{Pic}\;Y\rightarrow\left(\textup{Pic}\;\bar{Y}\right)^{\Gamma}
\]
n'est pas surjectif, bien que $Y$ poss{\`e}de un point $k$-rationnel.
\end{pro}
\vspace{1mm}

Nous montrons maintenant pourquoi, sur cet exemple, l'existence d'un point $k$-rationnel sur $Y$ n'entra{\^i}ne pas l'existence d'une r{\'e}traction du morphisme $v$ dans la suite exacte (c'est la suite $\left(S1\right)$ d{\'e}duite de la suite spectrale de Leray):
	\[\left(S5\right):\ \textup{Pic}\;Y\stackrel{u}{\longrightarrow} H^{0}\left(k,R^{1}\pi_{\ast}\mathbb{G}_{m,Y}\right)\stackrel{\delta^1}{\longrightarrow} H^{2}\left(k,\pi_{\ast}\mathbb{G}_{m,Y}\right)\stackrel{v}{\longrightarrow} \textup{Br}\;Y
\]
o{\`u} $\pi$ d{\'e}signe le morphisme structural $Y\rightarrow \textup{Spec}\;k$, obtenu par composition {\`a} partir des morphismes $Y\rightarrow X$ et $X\rightarrow \textup{Spec}\;k$. \newline

Commen\c{c}ons par d{\'e}crire $\delta^{1}$: soit $\mathcal{P}$ un {\'e}l{\'e}ment de $H^{0}\left(k,R^{1}\pi_{\ast}\mathbb{G}_{m,Y}\right)$; c'est une section globale du faisceau $R^{1}\pi_{\ast}\mathbb{G}_{m,Y}$, et c'est donc la donn{\'e}e:\vspace{2mm}
\begin{itemize}
\item d'une famille d'extensions {\'e}tales $\left(L_{i}/k\right)_{i\in I}$;\vspace{2mm}
\item pour tout $i\in I$, d'un $\mathbb{G}_{m,Y_{L_{i}}}$-torseur $P_{i}\rightarrow Y_{L_{i}}$;\vspace{2mm}
\item pour tout couple $\left(i,j\right)\in I\times I$, d'un isomorphisme:
	\[\varphi_{ij}:{P_{j}}_{\left|Y_{L_{ij}}\right.}\longrightarrow {P_{i}}_{\left|Y_{L_{ij}}\right.}
\]\vspace{2mm}
\end{itemize}

On peut maintenant donner une description tout-{\`a}-fait concr{\`e}te de ce qui emp{\^e}che $\mathcal{P}$ d'appartenir {\`a} l'image du morphisme $u$. Pour tout triplet $\left(i,j,k\right)\in I^{3}$, on a le diagramme suivant:
	\[\xymatrix@R=40pt{P_{k}\times_{Y_{L_{k}}}Y_{L_{ijk}} \ar[rr]^{{\varphi_{ik}}_{\left|Y_{L_{ijk}}\right.}} \ar[rd]_{{\varphi_{jk}}_{\left|Y_{L_{ijk}}\right.}} & & P_{i}\times_{Y_{L_{i}}}Y_{L_{ijk}}\\& P_{j}\times_{Y_{L_{j}}}Y_{L_{ijk}} \ar[ru]_{{\varphi_{ij}}_{\left|Y_{L_{ijk}}\right.}}}
\]
qui n'a aucune raison d'{\^e}tre commutatif. On pose alors:
	\[c_{ijk}={\varphi_{ij}}_{\left|Y_{L_{ijk}}\right.}\circ {\varphi_{jk}}_{\left|Y_{L_{ijk}}\right.}\circ \left({\varphi_{ik}}_{\left|Y_{L_{ijk}}\right.}\right)^{-1}
\]
Il est imm{\'e}diat que $\left(c_{ijk}\right)_{i,j,k\in I}$ est un 2-cocycle, qui mesure l'obstruction {\`a} ce que la famille des isomorphismes $\left(\varphi_{ij}\right)_{i,j\in I}$, qui est une donn{\'e}e de \textit{recollement} sur les $P_{i}$ soit une donn{\'e}e de \textit{descente};
	\[\xymatrix@R=40pt{P_{k}\times_{Y_{L_{k}}}Y_{L_{ijk}} \ar[rr]^{{\varphi_{ik}}_{\left|Y_{L_{ijk}}\right.}} \ar[rd]_{{\varphi_{jk}}_{\left|Y_{L_{ijk}}\right.}} & \ar@{=>}[d]^{c_{ijk}}& P_{i}\times_{Y_{L_{i}}}Y_{L_{ijk}}\\& P_{j}\times_{Y_{L_{j}}}Y_{L_{ijk}} \ar[ru]_{{\varphi_{ij}}_{\left|Y_{L_{ijk}}\right.}}}
\]
Remarquons maintenant que:
	\[c_{ijk}\in \textup{ad}_{\mathbb{G}_{m,Y_{L_{ijk}}}}\left(P_{i}\times_{Y_{L_{i}}}Y_{L_{ijk}}\right),\ \forall\ \left(i,j,k\right)\in I^{3}
\]
et du fait que $\mathbb{G}_{m}$ est ab{\'e}lien (!), on a:
	\[c_{ijk}\in \mathbb{G}_{m,Y_{L_{ijk}}}\left(Y_{L_{ijk}}\right)=\mathcal{O}_{Y_{L_{ijk}}}^{\ast}\left(Y_{L_{ijk}}\right),\ \forall\ \left(i,j,k\right)\in I^{3}
\]
Or, par d{\'e}finition du faisceau $\pi_{\ast}\mathbb{G}_{m,Y}$, ceci signifie encore que:
	\[c_{ijk}\in \pi_{\ast}\mathbb{G}_{m,Y}\left(L_{ijk}\right),\ \forall\ \left(i,j,k\right)\in I^{3}
\]
Par suite, la classe du 2-cocycle:
	\[\left[c_{ijk}\right]\in \check{H}^{2}\left(\left(L_{i}/k\right)_{i\in I},\pi_{\ast}\mathbb{G}_{m,Y}\right)\hookrightarrow \check{H}^{2}\left(k,\pi_{\ast}\mathbb{G}_{m,Y}\right)\simeq H^{2}\left(k,\pi_{\ast}\mathbb{G}_{m,Y}\right)
\]
(o{\`u} le dernier isomorphisme est d{\^u} au th{\'e}or{\`e}me III.2.17 de \cite{Mi}) est exactement l'obstruction {\`a} ce que $\mathcal{P}$ appartienne {\`a} l'image de $u$.

Passons maintenant {\`a} la description de $v$. Soit $\mathcal{G}\in H^{2}\left(k,\pi_{\ast}\mathbb{G}_{m,Y}\right)$. En utilisant une fois encore le th{\'e}or{\`e}me III.2.17 de \cite{Mi}, $\mathcal{G}$ est repr{\'e}sent{\'e} par un 2-cocycle:
	\[\left(g_{ijk}\right)_{i,j,k\in I'}\in Z^{2}\left(\left(L_{i}/k\right)_{i\in I'},\pi_{\ast}\mathbb{G}_{m,Y}\right)
\]
Par d{\'e}finition, on a donc:
	\[g_{ijk}\in \pi_{\ast}\mathbb{G}_{m,Y}\left(L_{ijk}\right),\ \forall\ \left(i,j,k\right)\in {I'}^{3}
\]
c'est-{\`a}-dire:
	\[g_{ijk}\in \mathbb{G}_{m,Y}\left(Y_{L_{ijk}}\right),\ \forall\ \left(i,j,k\right)\in {I'}^{3}
\]

On obtient donc trivialement un 2-cocycle:
	\[\left(g_{ijk}\right)_{i,j,k\in I'}\in Z^{2}\left(\left(Y_{L_{i}}/Y\right)_{i\in I'},\mathbb{G}_{m,Y}\right)
\]
ce qui fournit donc une interpr{\'e}tation particuli{\`e}rement claire du morphisme:
	\[v:H^{2}\left(k,\pi_{\ast}\mathbb{G}_{m,Y}\right)\longrightarrow H^{2}\left(Y,\mathbb{G}_{m,Y}\right)
\]\vspace{2mm}

Montrons maintenant que l'existence d'un point $k$-rationnel entra{\^i}ne l'existence d'un morphisme:
	\[r:H^{2}\left(Y,\mathbb{G}_{m,Y}\right)\longrightarrow \textup{Br}\;k
\]
Soit $\mathfrak{G}\in H^{2}\left(Y,\mathbb{G}_{m,Y}\right)$; une nouvelle application du th{\'e}or{\`e}me III.2.17 de \cite{Mi} permet de supposer que $\mathfrak{G}$ est repr{\'e}sent{\'e} par un 2-cocycle\footnote{O{\`u} $\left(Y_{\alpha}/Y\right)_{\alpha\in A}$ est un recouvrement {\'e}tale de $Y$ qui ne provient pas \textit{a priori} d'un recouvrement {\'e}tale de $k$.}:
	\[\left(\gamma_{\alpha\beta\epsilon}\right)_{\alpha,\beta,\epsilon\in A}\in Z^{2}\left(\left(Y_{\alpha}/Y\right)_{\alpha\in A},\mathbb{G}_{m,Y}\right)
\]
Pour tout $\left(\alpha,\beta,\epsilon\right)\in A^{3}$, on a\footnote{On a not{\'e}: $Y_{\alpha\beta\epsilon}=Y_{\alpha}\times_{Y}Y_{\beta}\times_{Y}Y_{\epsilon}$.}:
	\[\gamma_{\alpha\beta\epsilon}\in G_{Y}\left(Y_{\alpha\beta\epsilon}\right)
\]
Puisque $Y$ poss{\`e}de un point $k$-rationnel $y$, on peut consid{\'e}rer le produit fibr{\'e}:
	\[Y_{\alpha}\times_{Y,y}\textup{Spec}\;k
\]
C'est un $k$-sch{\'e}ma {\'e}tale, {\'e}tant obtenu par changement de base {\`a} partir du morphisme $Y_{\alpha}\rightarrow Y$, qui est {\'e}tale par hypoth{\`e}se. On note $\textup{Spec}\;K_{\alpha}$ ce sch{\'e}ma. Pour r{\'e}sumer la situation, on a donc le diagramme commutatif suivant:
	\[\xymatrix@C=60pt{&&\mathbb{G}_{m,Y}\ar[dl]\ar[dd]^{h} \\ Y_{\alpha} \ar[r]^{f_{\alpha}}&Y\ar[dd]_{\pi}\\ &&\mathbb{G}_{m,k}\ar[ld]\\\textup{Spec}\;K_{\alpha}\ar[r]_{q_{\alpha}} \ar[uu]^{y_{\alpha}}&\textup{Spec}\;k\ar@/_10pt/[uu]_{y}}
\]
Revenons maintenant {\`a} notre cocycle $\left(\gamma_{\alpha\beta\epsilon}\right)_{\alpha,\beta,\epsilon}$. Pour tout triplet $\left(\alpha,\beta,\epsilon\right)\in A^{3}$, le nouveau diagramme ci-dessous est commutatif:
	\[\xymatrix@C=60pt{&&\mathbb{G}_{m,Y}\ar[dl]\ar[dd]^{h} \\ Y_{\alpha\beta\epsilon}\ar@/^12pt/[rru]^{c_{\alpha\beta\epsilon}} \ar[r]^{f_{\alpha\beta\epsilon}}&Y\ar[dd]_{\pi}\\ &&\mathbb{G}_{m,k}\ar[ld]\\\textup{Spec}\;K_{\alpha\beta\epsilon}\ar[r]_{q_{\alpha\beta\epsilon}} \ar[uu]^{y_{\alpha\beta\epsilon}}&\textup{Spec}\;k\ar@/_10pt/[uu]_{y}}
\]
o{\`u}:
	\[K_{\alpha\beta\epsilon}=K_{\alpha}\otimes_{k}K_{\beta}\otimes_{k}K_{\epsilon}
\]
On pose:
	\[\widehat{c_{\alpha\beta\epsilon}}^{y}=h\circ c_{\alpha\beta\epsilon}\circ y_{\alpha\beta\epsilon},\ \forall\ \left(\alpha,\beta,\epsilon\right)\in A^{3}
\]
Alors:
	\[\widehat{c_{\alpha\beta\epsilon}}^{y}\in \mathbb{G}_{m,k}\left(K_{\alpha}\right),\ \forall\ \left(\alpha,\beta,\epsilon\right)\in A^{3}
\]
Ce faisant, on obtient donc un 2-cocycle:
	\[\left(\widehat{c_{\alpha\beta\epsilon}}^{y}\right)\in Z^{2}\left(\left(K_{\alpha}/k\right)_{\alpha\in A},\mathbb{G}_{m,k}\right)
\]
Il est imm{\'e}diat que l'on a ainsi d{\'e}fini un morphisme:
	\[r:H^{2}\left(Y,\mathbb{G}_{m,Y}\right)\longrightarrow \textup{Br}\;k
\]
mais du fait que $\pi_{\ast}\mathbb{G}_{m,Y}\neq \mathbb{G}_{m,k}$, ce n'est pas une r{\'e}traction du morphisme:
	\[u:H^{2}\left(k,\pi_{\ast}\mathbb{G}_{m,Y}\right)\longrightarrow H^{2}\left(Y,\mathbb{G}_{m,Y}\right)
\]
mais seulement une r{\'e}traction du morphisme:
	\[u':\textup{Br}\;k\longrightarrow H^{2}\left(Y,\mathbb{G}_{m,Y}\right)
\]
d{\'e}fini de mani{\`e}re {\'e}vidente vu ce qui pr{\'e}c{\`e}de. Pour conclure, on a le diagramme commutatif:
	\[\xymatrix@R=50pt{\textup{Pic}\;Y\ar[r]^{u\ \ \ \ \ \ \ \ }& H^{0}\left(k,R^{1}\pi_{\ast}\mathbb{G}_{m,Y}\right)\ar[r]^{\ \ \delta^{1}}& H^{2}\left(k,\pi_{\ast}\mathbb{G}_{m,Y}\right)\ar[rr]^{\ \ \ v}&& \textup{Br}\;Y\ar@/^30pt/[dl]^{r} \\ &&& \textup{Br}\;k \ar[ul]^{i}\ar[ur]_{u'}}
\]
\begin{rem} \textup{Nous avons affirm{\'e} un peu vite que $r$ est une r{\'e}traction du morphisme $u'$. Pour compl{\'e}ter la preuve de cette affirmation, il reste encore {\`a} prouver le fait suivant: si $L/k$ est une extension {\'e}tale, et en reprenant les notations habituelles:
	\[\xymatrix{Y_{L} \ar[r]^{f} \ar[d]_{\pi_{L}}& Y \ar[d]_{\pi}\\\textup{Spec}\;L \ar[r]_{q}&\textup{Spec}\;k \ar@/_15pt/[u]^{y}}
\]
alors:
	\[\textup{Spec}\;L\simeq Y_{L}\times_{f,Y,y}\textup{Spec}\;k
\]
}

\textup{Pour ce faire, on commence par remarquer que l'existence du point $k$-rationnel $y$ entra{\^i}ne l'existence d'une section:
	\[y_{L}:\textup{Spec}\;L\rightarrow Y_{L}
\]
du morphisme $\pi_{L}$. Son existence (et son unicit{\'e} d'ailleurs) est assur{\'e}e par la propri{\'e}t{\'e} universelle du produit fibr{\'e}, puisque l'on dispose du diagramme commutatif suivant:}
	
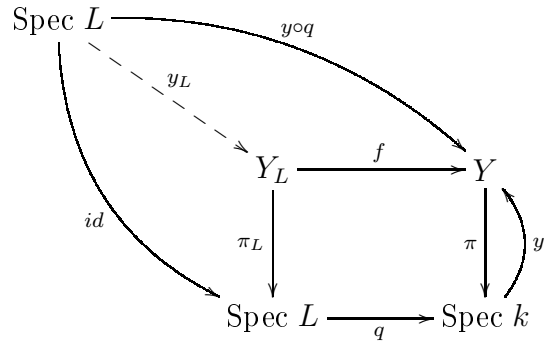
\[\xymatrix@R=40pt@C=40pt{\textup{Spec}\;L\ar@/_25pt/[ddr]_{id} \ar@/^20pt/[drr]^{y\circ q} \ar@{-->}[dr]^{y_{L}} \\&Y_{L} \ar[r]^{f} \ar[d]_{\pi_{L}}& Y \ar[d]_{\pi}\\&\textup{Spec}\;L \ar[r]_{q}&\textup{Spec}\;k \ar@/_15pt/[u]_{y}}
\]

\textup{Donnons nous maintenant un sch{\'e}ma $Z$ et deux morphismes
	\[g:Z\longrightarrow X_{L} \ \ \textup{et}\ \ h:Z\longrightarrow \textup{Spec}\;k
\]
tels que: $f\circ\;g=y\circ\;h$. On a:
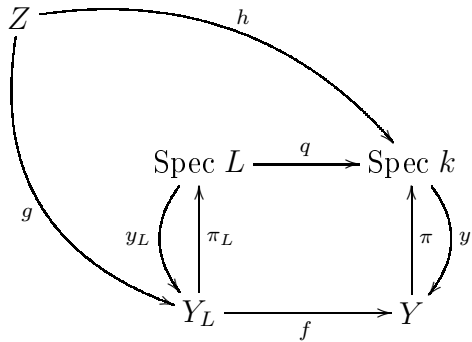
	\[\xymatrix@R=40pt@C=40pt{Z \ar@/_30pt/[ddr]_{g} \ar@/^24pt/[drr]^{h}\\&\textup{Spec}\;L \ar[r]^{q} \ar@/_15pt/[d]_{y_{L}}&\textup{Spec}\;k \ar@/^15pt/[d]^{y}\\&Y_{L} \ar[r]_{f} \ar[u]_{\pi_{L}}& Y\ar[u]_{\pi}}
\]
Il existe un unique morphisme 
	\[Y\longrightarrow \textup{Spec}\;L
\]
rendant commutatif tout le diagramme ci-dessus, le morphisme $\pi_{L}\circ\;g$. En effet, on a d'un c{\^o}t{\'e}:}
	\[f\circ\;y_{L}\circ\;\pi_{L}\circ\;g=f\circ\;g
\]
\textup{et de l'autre:}
	\[y\circ\;q\circ\;\pi_{L}\circ\;g=f\circ\;y_{L}\circ\;\pi_{L}\circ\;g=f\circ\;g
\]
\textup{Enfin:}
	\[q\circ\;\pi_{L}\circ\;g=\pi\circ\;f\circ\;g=\pi\circ\;y\circ\;h=h
\]

\textup{Par suite, le $k$-sch{\'e}ma $\textup{Spec}\;L$ est effectivement le produit fibr{\'e}:
	\[\textup{Spec}\;L\simeq Y_{L}\times_{f,Y,y}\textup{Spec}\;k
\]
}
\end{rem}
\end{section}
\begin{section}{Points ad{\'e}liques et torseurs}
Le but de cette section est d'illustrer par des applications un r{\'e}sultat de Colliot-Th{\'e}l{\`e}ne et Sansuc \cite{CTS} ({\'e}tendu par Skorobogatov \cite{Sk} des tores aux groupes multiplicatifs) assurant que l'existence de points ad{\'e}liques est suffisante pour descendre des torseurs sous des groupes ab{\'e}liens sur des vari{\'e}t{\'e}s propres. Jusqu'{\`a} la fin de ce chapitre:\vspace{2mm}
\begin{itemize}
\item $k$ est un corps de nombres, dont on note $\Omega_{k}$\label{Omegak} l'ensemble des places et $\mathbb{A}_{k}$\label{Ak} l'anneau des ad{\`e}les;\vspace{2mm}
\item $X$ est une $k$-vari{\'e}t{\'e} propre, de telle sorte que\footnote{En g{\'e}n{\'e}ral, on a seulement l'inclusion: $X\left(\mathbb{A}_{k}\right)\subset\prod_{v\in\Omega_{k}} X\left(k_{v}\right)$. Le crit{\`e}re valuatif de propret{\'e} (\textit{cf.} \cite{Har} II.4.7) fournit l'autre inclusion.}:\label{Xak}
	\[X\left(\mathbb{A}_{k}\right)=\prod_{v\in\Omega_{k}} X\left(k_{v}\right)
\]
On note toujours $\pi:X\rightarrow \textup{Spec}\;k$ le morphisme structural.
\end{itemize}

On identifie l'ensemble des points $k$-rationnels de $X$ {\`a} son image dans $\prod_{v\in\Omega_{k}} X\left(k_{v}\right)$ par le morphisme diagonal\index{morphisme!diagonal}:
	\[\begin{array}{ccc}X\left(k\right)&\longrightarrow&\displaystyle\prod_{v\in\Omega_{k}} X\left(k_{v}\right)\\x&\longmapsto&\left(x_{v}=x\circ p_{v}\right)_{v}\end{array}
\]
o{\`u}:
	\[p_{v}:\textup{Spec}\;k_{v}\longrightarrow \textup{Spec}\;k
\]
est le morphisme {\'e}vident.
	\[\xymatrix@R=50pt@C=60pt{&X\ar[d]_{\pi} \\ \textup{Spec}\;k_{v}\ar@{-->}[ur]^{x_{v}=x\circ p_{v}} \ar[r]_{p_{v}\ \ } &\textup{Spec}\;k \ar@/_20pt/[u]_{x}}
\]

On a donc {\'e}videmment l'inclusion:
	\[X\left(k\right)\subset X\left(\mathbb{A}_{k}\right)
\]

Il est donc imm{\'e}diat que $X\left(k\right)\neq\emptyset\Longrightarrow X\left(\mathbb{A}_{k}\right)\neq\emptyset$. Un probl{\`e}me difficile est d'arriver {\`a} d{\'e}terminer pour quelles vari{\'e}t{\'e}s la r{\'e}ciproque est vraie, \textit{i.e.} l'existence de points ad{\'e}liques sur $X$ entra{\^i}ne l'existence de points rationnels sur $X$; c'est ce que l'on appelle le principe de Hasse\index{principe!de Hasse}. Plus pr{\'e}cis{\'e}ment:
\begin{defi} Soit $X$ une vari{\'e}t{\'e} propre sur un corps de nombres $k$.
\begin{enumerate}[(i)]
\item Si $X\left(\mathbb{A}_{k}\right)\neq\emptyset$, on dit que $X$ a des points \textbf{partout localement}.\index{partout localement}
\item On dit que la vari{\'e}t{\'e} $X$ est un \textbf{contre-exemple au principe de Hasse} si $X$ a des points partout localement, mais ne poss{\`e}de pas de point $k$-rationnel.
\end{enumerate}
\end{defi}

\begin{quote}
	$\bullet$ \uline{Exemples de vari{\'e}t{\'e}s satisfaisant le principe de Hasse}
\end{quote}

\begin{itemize}
\item les quadriques projectives lisses \cite{Se2};\vspace{2mm}
\item les $k$-formes de $\mathbb{P}^{1}\times \mathbb{P}^{1}$ \cite{CTS};\vspace{2mm}
\item les vari{\'e}t{\'e}s de Severi-Brauer \cite{CTS};\vspace{2mm}
\item les surfaces de Del Pezzo de degr{\'e} 6 \cite{CTS};\vspace{2mm}
\item certaines intersections de quadriques \cite{CTCS};\vspace{2mm}
\item les surfaces cubiques singuli{\`e}res dans $\mathbb{P}^{3}_{k}$ (Skolem);\vspace{2mm}
\item les torseurs sous un groupe semi-simple simplement connexe (Kneser-Harder, \textit{cf.} th{\'e}or{\`e}me 4.2 de \cite{Sa}).\vspace{4mm}
\end{itemize}

\begin{quote}
	$\bullet$ \uline{Quelques contre-exemples au principe de Hasse}
\end{quote}

\begin{itemize}
\item la surface cubique de Cassels et Guy \cite{CG}: c'est la surface\index{surface!de Cassels et Guy} de $\mathbb{P}^{3}_{\mathbb{Q}}$ d'{\'e}quation:
	\[5X^{3}+9Y^{3}+10Z^{3}+12T^{3}=0
\]

Swinnerton-Dyer \cite{Sw} a {\'e}galement exhib{\'e} une surface projective lisse contre-exemple au principe de Hasse;\vspace{2mm}
\item Poonen \cite{Po} a prouv{\'e} que pour tout $t\in\mathbb{Q}$, la courbe de $\mathbb{P}^{2}$ d'{\'e}quation:
	\[5X^3+9Y^3+10Z^3+12\left(\frac{t^{12}-t^4-1}{t^{12}-t^8-1}\right)^{3}\left(X+Y+Z\right)^{3}=0
\]
est un contre-exemple au principe de Hasse;\vspace{2mm}
\item Siksek et Skorobogatov \cite{SSk} ont exhib{\'e} une courbe de Shimura contre-exemple au principe de Hasse;\vspace{2mm}
\item enfin Sarnak et Wang \cite{SW} ont construit une hypersurface lisse de $\mathbb{P}^{4}_{\mathbb{Q}}$ de degr{\'e} 1130 qui est un contre-exemple au principe de Hasse, sous r{\'e}serve que la conjecture de Lang (\textit{cf.} \cite{Lang} conjecture 1.2 p.179) tienne.
\end{itemize}
$\ $

Ce dernier contre-exemple nous int{\'e}resse particuli{\`e}rement, puisque l'obstruction de Brauer-Manin\index{obstruction!de Brauer-Manin} d'une hypersurface lisse de $\mathbb{P}^{4}_{\mathbb{Q}}$ est nulle. Comme le r{\'e}sultat de Skorobogatov que l'on veut utiliser est fortement li{\'e} {\`a} cette obstruction, nous rappelons ici sa construction.

\begin{subsection}{Construction de l'obstruction de Brauer-Manin}

Soit $X$ une $k$-vari{\'e}t{\'e} propre. On a un accouplement naturel:
	\[\left(Acc\right):\xymatrix@R=5pt{X\left(k\right)\times\textup{Br}\;X\ar[r]&\textup{Br}\;k\\\left(x,b\right)\ar@{|->}[r]&x^{\ast}b}
\]
que l'on peut d{\'e}crire de la mani{\`e}re suivante: comme d'apr{\`e}s nos hypoth{\`e}ses $\bar{k}\left[X\right]^{\ast}=\bar{k}^{\ast}$, la suite spectrale de Leray:
	\[E^{p,q}_{2}=H^{p}\left(k,R^{q}\pi_{\ast}\mathbb{G}_{m,X}\right)\Longrightarrow H^{p+q}\left(X,\mathbb{G}_{m,X}\right)=E^{p+q}
\]
fournit la suite longue de cohomologie:
	\[0\longrightarrow\textup{Pic}\;X\longrightarrow\left(\textup{Pic}\;\bar{X}\right)^{\Gamma}\longrightarrow\textup{Br}\;k\longrightarrow\textup{Br}\;X\longrightarrow H^{1}\left(k,\textup{Pic}\;\bar{X}\right)\longrightarrow0
\]

Si $x$ est un point rationnel, c'est-{\`a}-dire une section du morphisme structural:
	\[\xymatrix{X\ar[r]_{\pi\ \ \ }&\textup{Spec}\;k\ar@/_15pt/[l]_{x}}
\]
alors il fournit une r{\'e}traction not{\'e}e $x^{\ast}$ de l'edge:
	\[\xymatrix{E^{2,0}_{2}=\textup{Br}\;k\ar[r]&\textup{Br}\;X=E^{2}\ar@/^20pt/[l]^{x^{\ast}}}
\]

Le morphisme $x^{\ast}$ ainsi d{\'e}crit est celui intervenant dans la d{\'e}finition de l'accouplement $\left(Acc\right)$. 
\begin{remsub}\textup{Dans le cas o{\`u} $X$ est une vari{\'e}t{\'e} pour laquelle la \textquotedblleft{conjecture de Grothendieck}\textquotedblright\ sur les groupes de Brauer est vraie, on a une interpr{\'e}tation plus g{\'e}om{\'e}trique pour le morphisme $x^{\ast}$.}

\textup{Supposons que $X$ soit une telle vari{\'e}t{\'e}. Soit $b\in\textup{Br}\;X$. Notons $\mathcal{A}$ une alg{\`e}bre d'Azumaya sur $X$ correspondant\footnote{On entend par l{\`a} que: $\Delta\left(\left[\mathcal{A}\right]\right)=b$, o{\`u} $\Delta:\textup{Br}_{Az}X\rightarrow\textup{Br}\;X$ est le morphisme de Grothendieck.} {\`a} $b$; alors $x^{\ast}b$ n'est autre que la fibre $\mathcal{A}_{x}$ de $\mathcal{A}$ au point $x$. Evidemment, $\mathcal{A}_{x}$ est une $k$-alg{\`e}bre simple centrale, car c'est une $k\left(x\right)$-alg{\`e}bre simple centrale, et la rationnalit{\'e} de $x$ entra{\^i}ne $k\left(x\right)=k$.}
\end{remsub}

On veut maintenant rendre local l'accouplement $\left(Acc\right)$. Pour toute place $v\in\Omega_{k}$ on peut {\'e}videmment d{\'e}finir de la m{\^e}me mani{\`e}re que $\left(Acc\right)$ un accouplement:
	\[\left(Acc_v\right):\xymatrix@R=5pt{X\left(k_{v}\right)\times\textup{Br}\;X\ar[r]&\mathbb{Q}/\mathbb{Z}\\\left(x_{v},b\right)\ar@{|->}[r]&\textup{inv}_{v}\left(x_{v}^{\ast}b\right)}
\]
o{\`u} $\textup{inv}_{v}:\textup{Br}\;k_{v}\rightarrow \mathbb{Q}/\mathbb{Z}$\label{invv} est l'invariant fourni par la th{\'e}orie du corps de classes. En effet, comme $x_{v}$ est un point de $X$ {\`a} valeurs dans $\textup{Spec}\;k_{v}$, le diagramme ci-dessous est commutatif:
	\[\xymatrix@R=50pt@C=60pt{&X\ar[d]^{\pi}\\\textup{Spec}\;k_{v}\ar[ur]^{\ \ x_{v}}\ar[r]_{p_{v}}&\textup{Spec}\;k}
\]
et $x_{v}$ induit donc un morphisme:
	\[x_{v}^{\ast}:\textup{Br}\;X\longrightarrow \textup{Br}\;k_{v}
\]
\begin{remsub} \textup{Une fois encore, si $\textup{Br}_{Az}X=\textup{Br}\;X$, on a une interpr{\'e}tation g{\'e}om{\'e}trique de $x_{v}^{\ast}$. En effet, soit $b\in\textup{Br}\;X$, et soit $\mathcal{A}$ une alg{\`e}bre d'Azumaya correspondant {\`a} $b$. Dans ce cas $x_{v}^{\ast}b$\label{xvastb} n'est autre que le pullback de $\mathcal{A}$ par le morphisme $x_{v}$:}
	\[\xymatrix@R=50pt@C=30pt{&\mathcal{A}\ar[d]\\x_{v}^{\ast}b\approx\mathcal{A}\times_{X}\textup{Spec}\;k_{v}\ar@{-->}[d]\ar@{-->}[ur]&X\ar[d]^{\pi}\\\textup{Spec}\;k_{v}\ar[ur]_{x_{v}}\ar[r]_{p_{v}}&\textup{Spec}\;k}
\]
\end{remsub}

On d{\'e}finit maintenant un accouplement global {\`a} l'aide de la famille d'accouplements locaux $\left(Acc_{v}\right)$:
	\[\left\langle \bullet,\bullet\right\rangle:\xymatrix@R=5pt{X\left(\mathbb{A}_{k}\right)\times\mathcyr{B}\left(X\right)\ar[r]&\mathbb{Q}/\mathbb{Z}\\\left(\left(x_{v}\right)_{v},b\right)\ar@{|->}[r]&\displaystyle\sum_{v\in\Omega_{k}}\textup{inv}_{v}\left(x_{v}^{\ast}\tilde{b}\right)}
\]
o{\`u} $\tilde{b}$ d{\'e}signe un repr{\'e}sentant dans $\textup{Br}^{alg}X\subset\textup{Br}\;X$ de $b$. Il n'est \textit{a priori} pas {\'e}vident que $\left\langle \left(x_{v}\right)_{v},b\right\rangle$ soit bien d{\'e}fini par la formule ci-dessus. Nous rappelons bri{\`e}vement les arguments qui prouvent la coh{\'e}rence de cette d{\'e}finition.
\begin{faitsub} La somme $\displaystyle\sum_{v\in\Omega_{k}}\textup{inv}_{v}\left(x_{v}^{\ast}\tilde{b}\right)$ est finie.
\end{faitsub}

\uline{\textsc{Preuve}}: puisque $\left(x_{v}\right)_{v}\in X\left(\mathbb{A}_{k}\right)$, on a $x_{v}^{\ast}b\in \textup{Br}\;\mathcal{O}_{v}$, pour presque tout $v\in\Omega_{k}$, d'o{\`u} la conclusion puisque $\textup{Br}\;\mathcal{O}_{v}=0$ (\textit{cf.} \cite{Mi}).
\begin{flushright}
$\Box$
\end{flushright}
\begin{faitsub} La valeur de la somme $\displaystyle\sum_{v\in\Omega_{k}}\textup{inv}_{v}\left(x_{v}^{\ast}\tilde{b}\right)$ ne d{\'e}pend pas du repr{\'e}sentant $\tilde{b}\in \textup{Br}^{alg}X$ de $b$ choisi.
\end{faitsub}

\uline{\textsc{Preuve}}: soient $\tilde{b}$ et $\tilde{b'}$ deux repr{\'e}sentants de $b$. Alors il existe $c\in \textup{Br}^{cst}X$\index{groupe!de Brauer!constant} tel que\footnote{On a not{\'e}: $\textup{Br}^{cst}X=\textup{im}\left\{\textup{Br}\;k\longrightarrow\textup{Br}\;X\right\}$.}:
	\[\tilde{b}=\tilde{b'}+c
\]

Donc:
	\[\displaystyle\sum_{v\in\Omega_{k}}\textup{inv}_{v}\left(x_{v}^{\ast}\tilde{b}\right)=\displaystyle\sum_{v\in\Omega_{k}}\textup{inv}_{v}\left(x_{v}^{\ast}\tilde{b'}\right)+\displaystyle\sum_{v\in\Omega_{k}}\textup{inv}_{v}\left(x_{v}^{\ast}c\right)
\]

On d{\'e}duit de la suite exacte:
	\[\left(S7\right):\ \ \xymatrix@C=40pt{0\ar[r]&\textup{Br}\;k\ar[r]&\displaystyle\bigoplus_{v\in\Omega_{k}}\textup{Br}\;k_{v}\ar[r]^{\ \ \ \ \sum \textup{inv}_{v}}&\mathbb{Q}/\mathbb{Z}\ar[r]&0}
\]
la nullit{\'e} de la somme $\displaystyle\sum_{v\in\Omega_{k}}\textup{inv}_{v}\left(x_{v}^{\ast}c\right)$, ce qui ach{\`e}ve la preuve du fait.
\begin{flushright}
$\Box$
\end{flushright}

Les deux faits pr{\'e}c{\'e}dents nous assurent donc que l'accouplement $\left\langle \bullet,\bullet\right\rangle$ est bien d{\'e}fini. Nous {\'e}non\c{c}ons maintenant sa propri{\'e}t{\'e} fondamentale:
\begin{prosub} Avec les notations introduites pr{\'e}c{\'e}demment:
\begin{enumerate}[(i)]
\item La valeur de l'accouplement:
	\[\left\langle \left(x_{v}\right)_{v},b\right\rangle
\]
ne d{\'e}pend pas du point ad{\'e}lique de $X$ choisi.
\item Si $\left(x_{v}\right)_{v}\in X\left(\mathbb{A}_{k}\right)$ appartient {\`a} l'image du morphisme diagonal:
	\[X\left(k\right)\longrightarrow\prod_{v\in\Omega_{k}}X\left(k_{v}\right)
\]
alors: $\left\langle \left(x_{v}\right)_{v},b\right\rangle=0$.
\end{enumerate}
\end{prosub}

\uline{\textsc{Preuve}}: soient $\left(x_{v}\right)_{v}\in X\left(\mathbb{A}_{k}\right)$ et $b\in\mathcyr{B}\left(X\right)$. Choisissons un repr{\'e}sentant $\tilde{b}\in\textup{Br}^{alg}X$ de $b$. Pour toute place $v$ de $k$, on a:
	\[\textup{loc}_{v}\left(b\right)\in \textup{Br}^{cst}X_{v}
\]
ce qui prouve (i).


Soient maintenant $b\in\mathcyr{B}\left(X\right)$ et $\tilde{b}$ un repr{\'e}sentant de $b$. On a:
	\[x_{v}^{\ast}\tilde{b}=\left(x\circ p_{v}\right)^{\ast}\tilde{b}=p_{v}^{\ast}\left(x^{\ast}\tilde{b}\right),\ \forall\;v\in\Omega_{k},
\]
et $x^{\ast}\tilde{b}\in\textup{Br}\;k$.
\pagebreak

En utilisant une nouvelle fois l'exactitude de la suite $\left(S7\right)$, on en d{\'e}duit finalement que:
	\[\left\langle \left(x_{v}\right)_{v},b\right\rangle=0
\]
\begin{flushright}
$\Box$
\end{flushright}

En conclusion, on a donc d{\'e}fini un morphisme de groupes:
	\[m_{\mathcal{H,\mathcyr{B}\left(X\right)}}\left(X\right):\xymatrix@R=5pt{\mathcyr{B}\left(X\right)\ar[r]&\mathbb{Q}/\mathbb{Z}\\b\ar@{|->}[r]&\displaystyle\sum_{v\in\Omega_{k}}\textup{inv}_{v}\left(x_{v}^{\ast}\tilde{b}\right)}
\]
o{\`u} $\left(x_{v}\right)_{v}\in X\left(\mathbb{A}_{k}\right)$ est quelconque, et $\tilde{b}$ est un repr{\'e}sentant de $b$ dans $\textup{Br}^{alg}X$. De la proposition pr{\'e}c{\'e}dente, on d{\'e}duit alors le:
\begin{theosub}$\ $

\begin{center}
Si $X\left(k\right)\neq\emptyset$, alors $m_{\mathcal{H}}\left(X\right)=0\in \mathcyr{B}\left(X\right)^{D}$.
\end{center}
\end{theosub}

Notons enfin pour terminer ces rappels que (voir p. \pageref{BcyrX}):
	\[\left(\clubsuit\right)\ \ \ \ \mathcyr{B}\left(X\right)^{D}\approx\mathcyr{SH}^{1}\left(k,\textup{Pic}\;\bar{X}\right)^{D}
\]
puisque l'on d{\'e}duit de la suite exacte:
	\[\textup{Br}\;k\longrightarrow \textup{Br}^{alg}X\longrightarrow H^{1}\left(k,\textup{Pic}\;\bar{X}\right)\longrightarrow 0
\]
l'isomorphisme:
	\[\textup{Br}_{a}X\approx H^{1}\left(k,\textup{Pic}\;\bar{X}\right)
\]
puis l'isomorphisme $\left(\clubsuit\right)$, en \textquotedblleft{passant aux noyaux sur toutes les places}\textquotedblright\ et en dualisant.
\end{subsection}
\begin{subsection}{Exemples de calculs de $m_{\mathcal{H,\mathcyr{B}\left(X\right)}}\left(X\right)$}
\begin{exemsub}[Vari{\'e}t{\'e}s de Severi-Brauer et intersections compl{\`e}tes lisses]
\textup{D'apr{\`e}s la proposition 2.1.13, on a le:}
\begin{lemsub} Si $V$ est une $k$-vari{\'e}t{\'e} de Severi-Brauer ou une intersection compl{\`e}te lisse de dimension $\geq3$, alors $m_{\mathcal{H,\mathcyr{B}\left(X\right)}}\left(V\right)=0$.
\end{lemsub}
\end{exemsub}
\begin{exemsub}[Surfaces de Del Pezzo] \textup{Rappelons que l'on appelle \textit{surface de Del Pezzo}\index{surface!de Del Pezzo} une surface projective lisse dont le diviseur anticanonique est ample. Le plan projectif est un exemple de surface de Del Pezzo\footnote{Puisque $K_{\mathbb{P}^{2}}=-3\mathbb{P}^{1}$, d'apr{\`e}s \cite{Har} II.8.20.1.}. Un exemple moins trivial est fourni par la surface cubique non-singuli{\`e}re de $\mathbb{P}^{3}$. Cette derni{\`e}re est obtenue en {\'e}clatant $\mathbb{P}^{2}$ en 6 points en position g{\'e}n{\'e}rale (non co-coniques et 3 {\`a} 3 non-align{\'e}s). Donc son groupe de Picard est $\mathbb{Z}^{7}$ (un premier exemplaire de $\mathbb{Z}$ correspond {\`a} $\mathbb{P}^{2}$, et les 6 autres sont fournis par les diviseurs exceptionnels correspondant aux {\'e}clatements successifs). D'une mani{\`e}re plus g{\'e}n{\'e}rale, une surface de Del Pezzo est toujours obtenue en {\'e}clatant $\mathbb{P}^{2}$ en un certain nombre de points (\textit{cf.} \cite{Ma2}), et son groupe de Picard est donc toujours de la forme $\mathbb{Z}^{N}$. Si $\Gamma$ op{\`e}re trivialement sur $\mathbb{Z}^{N}$ (c'est alors un $\Gamma$-module de permutation), on a (\textit{cf.} \cite{CTS}):}
	\[\mathcyr{SH}^{1}\left(k,\mathbb{Z}^{N}\right)=0
\]
d'o{\`u} la nullit{\'e} de $m_{\mathcal{H},\mathcyr{B}\left(X\right)}\left(X\right)$ dans ce cas.
\end{exemsub}
\begin{exemsub}[Espaces homog{\`e}nes sous $SL_{n}$ avec isotropie finie]
\textup{Soit $\bar{H}$ un sous-groupe fini (non-n{\'e}cessairement ab{\'e}lien) de $SL_{n}\left(\bar{k}\right)$.\footnote{En particulier, $\bar{H}$ est un $\bar{k}$-groupe alg{\'e}brique; il ne provient pas n{\'e}cessairement d'un groupe d{\'e}fini sur $k$.} On consid{\`e}re alors le $\bar{k}$-espace homog{\`e}ne:}
	\[\bar{V}=SL_{n}\left(\bar{k}\right)/\bar{H}
\]

\textup{On choisit $X$ une $k$-forme\footnote{Il est un peu abusif d'utiliser parler de \textquotedblleft{$k$-forme}\textquotedblright\ ici; nous voulons juste dire que $X$ est une $k$-vari{\'e}t{\'e} telle que: $X\otimes_{k}\bar{k}\approx\bar{V}$.} de $\bar{V}$. D'apr{\`e}s le corollaire 4.6 de \cite{FI}, on a:}
	\[\textup{Pic}\;\bar{X}=\widehat{\bar{H}}
\]

\textup{Il s'ensuit que l'obstruction de Brauer-Manin d'une telle vari{\'e}t{\'e} $X$ n'est en g{\'e}n{\'e}ral pas nulle.}
\begin{lemsub} Soient $n$ un entier, et $\bar{H}$ un $\bar{k}$-groupe alg{\'e}brique qui est un sous-groupe de $SL_{n}\left(\bar{k}\right)$. Soit $X$ une $k$-forme de $SL_{n}\left(\bar{k}\right)/\bar{H}$. Alors:
	\[m_{\mathcal{H}}\left(X\right)\in \mathcyr{SH}^{1}\left(k,\widehat{\bar{H}}\right)
\]
\end{lemsub}

\textup{Nous terminons ces exemples en faisant le lien entre cette remarque et la dualit{\'e} de Tate-Poitou. On suppose maintenant $\bar{H}$ \textit{ab{\'e}lien} fini, et on note toujours $X$ une forme de $SL_{n}\left(\bar{k}\right)/\bar{H}$. On peut associer {\`a} $X$ la gerbe de ses trivialisations, c'est-{\`a}-dire la gerbe $\mathcal{G}_{X}$ qui mesure l'obstruction {\`a} ce que $X$ soit domin{\'e} par un $SL_{n}$-torseur sur $k$ (n{\'e}cessairement trivial). Dans la terminologie de Springer \cite{Sp}, la donn{\'e}e de $X$ correspond {\`a} la donn{\'e}e d'un $1$-cocycle {\`a} valeurs dans $SL_{n}\left(\bar{k}\right)/\bar{H}$, et $\mathcal{G}_{X}$ est alors la classe du $2$-cocycle {\`a} valeurs dans $\bar{H}$ mesurant ce qui emp{\^e}che de le relever en un $1$-cocycle {\`a} valeurs dans $SL_{n}\left(\bar{k}\right)$.}

\textup{Supposons que $X$ ait un point partout localement. Alors\footnote{Pour {\^e}tre compl{\`e}tement rigoureux, il e{\^u}t fallu dire: $\mathcal{G}_{X}$ repr{\'e}sente une classe dans $\mathcyr{SH}^{2}\left(k,\bar{H}\right)$\label{SH2ab}. Cet abus est justifi{\'e} par le fait que si un repr{\'e}sentant dans une classe d'{\'e}quivalence est une gerbe neutre, alors tous les repr{\'e}sentants de cette classe sont des gerbes neutres.}:}
	\[\mathcal{G}_{X}\in\mathcyr{SH}^{2}\left(k,\bar{H}\right)
\]

\textup{Moyennant la d{\'e}finition de la cohomologie {\`a} valeurs dans un topos localement annel{\'e} (\textit{cf.} \cite{SGA4-5}) on peut d{\'e}finir l'obstruction de Brauer-Manin\footnote{Nous renvoyons au chapitre IV pour la d{\'e}finition de l'obstruction de Brauer-Manin d'une gerbe.} de la gerbe $\mathcal{G}_{X}$, not{\'e}e $m_{\mathcal{H}}\left(\mathcal{G}_{X}\right)$ \cite{DEZ}. Celle-ci a le bon go{\^u}t de satisfaire la propri{\'e}t{\'e} suivante:}
\begin{prosub} Avec les notations introduites ci-dessus:
\begin{enumerate}[(i)]
\item $m_{\mathcal{H}}\left(\mathcal{G}_{X}\right)=m_{\mathcal{H}}\left(X\right)$;
\item la dualit{\'e} de Tate-Poitou:\index{dualit{\'e} de Tate-Poitou!}
	\[\mathcyr{SH}^{2}\left(k,\bar{H}\right)\times\mathcyr{SH}^{1}\left(k,\widehat{\bar{H}}\right)\longrightarrow \mathbb{Q}/\mathbb{Z}
\]
est explicitement r{\'e}alis{\'e}e gr{\^a}ce {\`a} l'obstruction de Brauer-Manin des gerbes, dans le sens o{\`u}:\label{mHGer}
	\[m_{\mathcal{H}}:\xymatrix@R=5pt{\mathcyr{SH}^{2}\left(k,\bar{H}\right)\ar[r]&\mathcyr{SH}^{1}\left(k,\widehat{\bar{H}}\right)^{D}\\\mathcal{G}\ar@{|->}[r]&m_{\mathcal{H}}\left(\mathcal{G}\right)}
\]
est un isomorphisme de groupes.
\end{enumerate}
\end{prosub}

\uline{\textsc{Preuve}}: \textup{c'est une cons{\'e}quence de la proposition 3.2 de \cite{DEZ}.}
\begin{flushright}
$\Box$
\end{flushright}
\end{exemsub}
\end{subsection}
\begin{subsection}{Brauer-Manin orthogonalit{\'e} et descente}
Apr{\`e}s ces quelques exemples et avant d'{\'e}noncer le th{\'e}or{\`e}me de Skorobogatov, il nous faut introduire un peu de terminologie:
\begin{defisub} Un point ad{\'e}lique $\left(x_{v}\right)_{v}\in X\left(\mathbb{A}_{k}\right)$ est dit \textbf{Brauer-Manin orthogonal {\`a}} $b\in \textup{Br}\;X$ si:
	\[\sum_{v\in\Omega_{k}}\textup{inv}_{v}\left(x_{v}^{\ast}b\right)=0
\]
cette somme {\'e}tant finie d'apr{\`e}s le fait 2.3.4. Soit $B$ une partie de $\textup{Br}\;X$; le point ad{\'e}lique $\left(x_{v}\right)_{v}$ est dit \textbf{Brauer-Manin orthogonal {\`a}} $B$ s'il est Brauer-Manin orthogonal {\`a} tout $b\in B$. On note:\label{XAkB}
	\[X\left(\mathbb{A}_{k}\right)^{B}=\left\{\left(x_{v}\right)_{v}\in X\left(\mathbb{A}_{k}\right)\left|\ \sum_{v\in\Omega_{k}}\textup{inv}_{v}\left(x_{v}^{\ast}b\right)=0,\;\forall\;b\in B\right.\right\}
\]
l'ensemble des points ad{\'e}liques Brauer-Manin orthogonaux {\`a} $B$.
\end{defisub}
\begin{exemsub}\textup{D'apr{\`e}s les exemples pr{\'e}c{\'e}dents tout point ad{\'e}lique sur une vari{\'e}t{\'e} de Severi-Brauer (resp. une intersection compl{\`e}te lisse dans $\mathbb{P}^{n}$ de dimension $\geq3$, resp. une surface de Del Pezzo) $X$ est Brauer-Manin orthogonal {\`a} $\textup{Br}\;X$, \textit{i.e}:}
	\[X\left(\mathbb{A}_{k}\right)^{\textup{Br}\;X}=X\left(\mathbb{A}_{k}\right)
\]
\end{exemsub}

Soit $B$ une partie de $\textup{Br}\;X$. On a la suite d'inclusions:
\begin{equation}
	X\left(k\right)\subset X\left(\mathbb{A}_{k}\right)^{\textup{Br}\;X}\subset X\left(\mathbb{A}_{k}\right)^{B}\subset X\left(\mathbb{A}_{k}\right)
\end{equation}
o{\`u} la premi{\`e}re inclusion d{\'e}coule de la proposition 2.3.1.5(ii), et les deux autres sont triviales. Nous allons justement introduire un nouveau maillon dans cette cha{\^i}ne. Soit $M$ un $k$-groupe alg{\'e}brique ab{\'e}lien, tel que $\bar{M}_{X}\left(\bar{X}\right)=\bar{M}\left(\bar{k}\right)$ (\textit{e.g.} $M$ fini). De la suite spectrale \textquotedblleft{des Ext}\textquotedblright\ \cite{CTS}:\index{suite spectrale!des Ext}
	\[Ext^{p}_{\left(\textup{Spec}\;k\right)_{\acute{e}t}}\left(M,R^{q}\pi_{\ast}\mathbb{G}_{m,X}\right)\Longrightarrow Ext_{X_{\acute{e}t}}^{p+q}\left(M_{X};\mathbb{G}_{m,X}\right)
\]
on d{\'e}duit la suite {\`a} 5 termes:
	\[H^{1}\left(k,M\right)\longrightarrow H^{1}\left(X,M_{X}\right)\stackrel{type}{\longrightarrow} \textup{Hom}_{\Gamma}\left(\widehat{M\left(\bar{k}\right)},\textup{Pic}\;\bar{X}\right)\stackrel{\partial}{\longrightarrow}H^{2}\left(k,M\right)\longrightarrow H^{2}\left(X,M_{X}\right)
\]
\begin{defisub} Soit $Y\rightarrow X$ un $M_{X}$-torseur. On appelle \textbf{type de $Y$}\index{type d'un torseur} l'image de $\left[Y\right]$ par le morphisme du m{\^e}me nom dans $\textup{Hom}_{\Gamma}\left(\widehat{M\left(\bar{k}\right)},\textup{Pic}\;\bar{X}\right)$. Si $\widehat{M\left(\bar{k}\right)}=\textup{Pic}\;\bar{X}$, on appelle \textbf{torseur universel sur $X$}\index{torseur!universel} un $M_{X}$-torseur dont le type est l'identit{\'e} de $\textup{Pic}\;\bar{X}$.
\end{defisub}
\begin{remsub}\textup{L'exactitude de la suite ci-dessus a la cons{\'e}quence suivante: pour tout $\lambda\in \textup{Hom}_{\Gamma}\left(\widehat{M\left(\bar{k}\right)},\textup{Pic}\;\bar{X}\right)$}:
	\[\partial\left(\lambda\right)=0\Leftrightarrow \textup{il existe un $M_{X}$-torseur sur $X$ de type $\lambda$}.
\]
\end{remsub}

La suite {\`a} 5 termes que l'on vient d'{\'e}voquer n'est {\'e}videmment pas sans rapport avec celle d{\'e}duite de la suite spectrale de Leray:
	\[H^{1}\left(k,M\right)\longrightarrow H^{1}\left(X,M_{X}\right)\longrightarrow H^{1}\left(\bar{X},\bar{M}_{X}\right)^{\Gamma}\stackrel{D}{\longrightarrow}H^{2}\left(k,M\right)\longrightarrow H^{2}\left(X,M_{X}\right)
\]

Nous renvoyons {\`a} l'appendice B de \cite{HS} pour une comparaison d{\'e}taill{\'e}e de ces deux suites spectrales. Les groupes $H^{1}\left(\bar{X},\bar{M}_{X}\right)^{\Gamma}$ et $\textup{Hom}_{\Gamma}\left(\widehat{M\left(\bar{k}\right)},\textup{Pic}\;\bar{X}\right)$ sont isomorphes (\textit{cf.} \cite{CTS},\cite{HS}), et on note:
	\[\tau:H^{1}\left(\bar{X},\bar{M}_{X}\right)^{\Gamma}\longrightarrow \textup{Hom}_{\Gamma}\left(\widehat{M\left(\bar{k}\right)},\textup{Pic}\;\bar{X}\right)
\]
cet isomorphisme. Ceci nous am{\`e}ne {\`a} introduire une nouvelle d{\'e}finition:
\begin{defisub} Soit $\bar{Y}\rightarrow\bar{X}$ un $\bar{M}_{X}$-torseur sur $\bar{X}$ de corps des modules $k$. On appelle \textbf{type de $\bar{Y}$} l'image de $\bar{Y}$ par le morphisme $\tau$.
\end{defisub}

On d{\'e}duit des remarques pr{\'e}c{\'e}dentes le:
\begin{lemsub} Soit $\bar{Y}\rightarrow\bar{X}$ un $\bar{M}_{X}$-torseur de corps des modules $k$. Les assertions suivantes sont {\'e}quivalentes:
\begin{enumerate}[(i)]
\item $\bar{Y}\rightarrow\bar{X}$ est d{\'e}fini sur $k$;
\item la gerbe des mod{\`e}les $D\left(\bar{Y}\right)$ est neutre;
\item la gerbe $\partial\left(\tau\left(\bar{Y}\right)\right)$ est neutre;
\item il existe un $M_{X}$-torseur sur $X$ de type $\tau\left(\bar{Y}\right)$.
\end{enumerate}

En outre, ces conditions sont {\'e}videmment satisfaites lorsque $X\left(k\right)\neq\emptyset$.
\end{lemsub}

Soit $\lambda\in \textup{Hom}_{\Gamma}\left(\widehat{M\left(\bar{k}\right)},\textup{Pic}\;\bar{X}\right)$. On peut lui associer un morphisme de groupes:
	\[\lambda_{\ast}:H^{1}\left(k,\widehat{\bar{M}}\right)\longrightarrow H^{1}\left(k,\textup{Pic}\;\bar{X}\right)
\]
l'image d'un $\widehat{\bar{M}}$-torseur $Y\rightarrow \textup{Spec}\;k$ par $\lambda_{\ast}$ {\'e}tant donn{\'e}e par extension du groupe structural {\`a} l'aide de $\lambda$: explicitement, le torseur $\lambda_{\ast}Y$ n'est autre que le $\textup{Pic}\;\bar{X}$-torseur:
	\[\xymatrix{Y\wedge^{\lambda}\textup{Pic}\;\bar{X}\ar[d]\\\textup{Spec}\;k}
\]
obtenu {\`a} l'aide du produit contract{\'e} d{\'e}fini \textit{via} le morphisme $\lambda$ \cite{Gi2}. Notons maintenant $r$ le morphisme naturel\footnote{C'est le morphisme: $\ker\left\{E^{2}\rightarrow E^{2,0}_{2}\right\}\longrightarrow E^{1,1}_{2}$ d{\'e}duit de la suite spectrale de Leray.}:
	\[r:\textup{Br}^{alg}X\longrightarrow H^{1}\left(k,\textup{Pic}\;\bar{X}\right)
\]

\begin{defisub} Avec les notations adopt{\'e}es ci-dessus, on pose:\label{Brlambda}
	\[\textup{Br}_{\lambda}X=r^{-1}\left[\lambda_{\ast}H^{1}\left(k,\widehat{\bar{M}}\right)\right]
\]
\end{defisub}
On a {\'e}videmment la cha{\^i}ne d'inclusions:

\begin{equation}
	X\left(k\right)\subset X\left(\mathbb{A}_{k}\right)^{\textup{Br}\;X}\subset X\left(\mathbb{A}_{k}\right)^{\textup{Br}_{\lambda}X}\subset X\left(\mathbb{A}_{k}\right)
\end{equation}

Le point fondamental est alors le r{\'e}sultat de Skorobogatov:
\begin{theosub}[Theorem 3, \cite{Sk}] Supposons que $\bar{M}$ soit un $\Gamma$-module de type fini en tant que groupe ab{\'e}lien. Si $X\left(\mathbb{A}_{k}\right)^{\textup{Br}_{\lambda}X}\neq\emptyset$, alors les conditions {\'e}quivalentes du lemme 2.3.3.6 sont satisfaites.
\end{theosub}

Autrement dit:
\begin{corosub} L'existence d'un point ad{\'e}lique sur $X$ Brauer-Manin orthogonal {\`a} $\textup{Br}_{\lambda}X$ entra{\^i}ne l'existence d'un mod{\`e}le pour tout $\bar{M}_{X}$-torseur sur $\bar{X}$ de corps des modules $k$ et de type $\lambda$.
\end{corosub}

\begin{exemsub}[Application aux $G$-rev{\^e}tements sur des espaces homog{\`e}nes]
\textup{On fixe $H$ un $k$-groupe alg{\'e}brique ab{\'e}lien fini, et on choisit un entier $n$ de telle sorte que $H$ se r{\'e}alise comme un sous-groupe de $SL_{n}\left(k\right)$ (on peut par exemple prendre $n=\left|H\right|$, mais un tel $n$ n'est {\'e}videmment pas unique). On choisit $X$ une $k$-forme de $SL_{n}\left(\bar{k}\right)/\bar{H}$. On a:}
	\[\textup{Pic}\;\bar{X}=\widehat{\bar{H}}
\]

\textup{Soit $\bar{Y}\rightarrow\bar{X}$ un $\widehat{\bar{H}}$-torseur sur $\bar{X}$ de corps des modules $k$. Du fait que $\bar{H}$ est ab{\'e}lien fini, c'est aussi un $\widehat{\bar{H}}$-rev{\^e}tement de $\bar{X}$. On note $\tau\left(\bar{Y}\right)$ le type du rev{\^e}tement $\bar{Y}$. On d{\'e}duit de ce qui pr{\'e}c{\`e}de le:}

\begin{lemsub} S'il existe sur $X$ un point ad{\'e}lique Brauer-Manin orthogonal {\`a} $\textup{Br}_{\tau\left(\bar{Y}\right)}X$, alors le $\widehat{\bar{H}}$-rev{\^e}tement $\bar{Y}\rightarrow\bar{X}$ est d{\'e}fini sur $k$.
\end{lemsub}

\begin{remsub} \textup{Une fois encore, on se rend compte que l'{\'e}nonc{\'e} de Skorobogatov fournit donc une condition suffisante beaucoup plus faible que l'existence d'un point rationnel pour qu'un rev{\^e}tement ab{\'e}lien soit d{\'e}fini sur son corps des modules. Il convient malgr{\'e} tout de temp{\'e}rer ce r{\'e}sultat par la remarque suivante: le th{\'e}or{\`e}me 2.3.3.8 ne donne aucune information quant aux rev{\^e}tements sur des vari{\'e}t{\'e}s de Severi-Brauer, etc\ldots\ D'ailleurs, plus g{\'e}n{\'e}ralement, cet {\'e}nonc{\'e} est \textquotedblleft{vide}\textquotedblright\ pour des $k$-vari{\'e}t{\'e}s telles que $\textup{Pic}\;\bar{X}$ soit sans torsion (\textit{e.g.} des vari{\'e}t{\'e}s rationnelles \cite{CTS}). En effet, puisque $M$ est fini, il existe des types int{\'e}ressants (non-triviaux) dans}
	\[\textup{Hom}_{\Gamma}\left(\widehat{\bar{M}},\textup{Pic}\;\bar{X}\right)
\]
\textup{si $\textup{Pic}\;\bar{X}$ poss{\`e}de de la torsion.}
\end{remsub}
\begin{remsub}\textup{Evidemment, on aimerait beaucoup pouvoir {\'e}tendre ce th{\'e}or{\`e}me au cas non-ab{\'e}lien. Les obstacles sont de deux sortes. Premi{\`e}rement, il faut faire une croix sur les suites spectrales. La cons{\'e}quence la plus f{\^a}cheuse de cette disparition est que l'on n'a plus aucune raison d'avoir un isomorphisme:}
	\[H^{1}\left(\bar{X},\bar{M}_{X}\right)^{\Gamma}\approx \textup{Hom}_{\Gamma}\left(\widehat{M\left(\bar{k}\right)},\textup{Pic}\;\bar{X}\right)
\]

\textup{Adieu donc notre belle correspondance entre types et rev{\^e}tements de corps des modules $k$!}

\textup{La seconde cons{\'e}quence, non moins f{\^a}cheuse, est que la dualit{\'e} de Tate-Poitou ne tient plus. Plus pr{\'e}cis{\'e}ment, lorsque $H$ n'est plus ab{\'e}lien, l'obstruction de Brauer-Manin (des gerbes) permet juste de d{\'e}finir une application \cite{DEZ}:}
	\[m_{\mathcal{H}}:\mathcyr{SH}^{2}\left(k,\textup{lien}\;H\right)\longrightarrow \mathcyr{SH}^{1}\left(k,\widehat{\bar{H}}\right)^{D}
\]
\textup{mais cette application n'a plus aucune raison d'{\^e}tre un isomorphisme. Or le th{\'e}or{\`e}me de Skorobogatov utilise de mani{\`e}re essentielle la dualit{\'e} de Tate-Poitou\ldots}
\end{remsub}
\end{exemsub}
\end{subsection}
\end{section}
\end{chapter}

\begin{chapter}[Descente de torseurs: le cas non-ab{\'e}lien]{Descente de torseurs et points rationnels: le cas non-ab{\'e}lien}
\thispagestyle{empty}
L'objectif de ce chapitre est de comprendre ce qui emp{\^e}che de g{\'e}n{\'e}raliser d'une mani{\`e}re franchement satisfaisante au cas non-ab{\'e}lien les r{\'e}sultats du chapitre pr{\'e}c{\'e}dent. On consid{\`e}re toujours un corps $k$ de caract{\'e}ristique nulle, $X$ un $k$-sch{\'e}ma et un $k$-groupe alg{\'e}brique $G$ non-n{\'e}cessairement ab{\'e}lien; il n'est donc plus question d'utiliser des suites spectrales. Le plan que nous suivrons dans ce chapitre est le suivant. Dans un premier temps, nous donnerons une interpr{\'e}tation purement topologique de la suite d'ensembles point{\'e}s:
	\[0\longrightarrow H^{1}\left(k,\pi_{\ast}G_{X}\right)\longrightarrow H^{1}\left(X,G_{X}\right)\stackrel{u}{\longrightarrow} H^{0}\left(k,R^{1}\pi_{\ast}\bar{G}_{X}\right)
\]
et nous d{\'e}terminerons les conditions sur $G$ permettant de \textquotedblleft{prolonger}\textquotedblright\ cette suite au cas non-ab{\'e}lien. Dans un deuxi{\`e}me temps, nous donnerons une preuve du r{\'e}sultat principal inspir{\'e}e de la preuve du th{\'e}or{\`e}me de Combes-Harbater pour les rev{\^e}tements (pour laquelle je tiens {\`a} remercier Michel Emsalem). Nous nous int{\'e}resserons ensuite au cas des groupes non-ab{\'e}liens finis, pour lesquels (assez ironiquement d'ailleurs) tout fonctionne {\`a} merveille. Enfin, nous t{\^a}cherons de passer en revue les obstructions {\`a} ce que l'on puisse dans le cas g{\'e}n{\'e}ral calculer \textquotedblleft{explicitement}\textquotedblright\ le lien de la gerbe des torseurs d'un $\bar{G}_{X}$-torseur $\bar{P}\rightarrow\bar{X}$ de corps des modules $k$, et en particulier l'obstruction {\`a} ce que ce lien soit repr{\'e}sentable (ou juste localement repr{\'e}sentable) par un $k$-groupe alg{\'e}brique, m{\^e}me lorsque le groupe $G$ dont on part est le plus sympathique possible (par exemple semi-simple simplement connexe). La clef de ce chapitre est de toute fa\c{c}on li{\'e}e au lemme fondamental (lemme \ref{lemf}), dont nous rappelons ici l'{\'e}nonc{\'e}:
\vspace{4mm}
\newline
\textbf{Lemme fondamental.} \textit{Soient $S$ un sch{\'e}ma, $G_{S}$ un sch{\'e}ma en groupes sur $S$, et $P$ un $G_{S}$-torseur sur $S$. Alors $\textup{ad}_{G_{S}}\left(P\right)$ est une $S$-forme int{\'e}rieure de $G_{S}$; autrement dit:}
	\[\textup{ad}_{G_{S}}\left(P\right)\textit{ repr{\'e}sente une classe de } H^{1}\left(S,\textup{Int}\;G_{S}\right).
\]

\textit{En particulier, si $G_{S}$ est ab{\'e}lien, alors:}
	\[\textup{ad}_{G_{S}}\left(P\right)\approx G_{S}
\]
\pagebreak

\begin{section}{Une interpr{\'e}tation topologique}
Avant de d{\'e}crire compl{\`e}tement la suite
	\[0\longrightarrow H^{1}\left(k,\pi_{\ast}G_{X}\right)\longrightarrow H^{1}\left(X,G_{X}\right)\stackrel{u}{\longrightarrow} H^{0}\left(k,R^{1}\pi_{\ast}\bar{G}_{X}\right)
\]
il nous semble indispensable ici, au vu des consid{\'e}rations qui vont suivre, de faire quelques rappels sur les diff{\'e}rents ensembles point{\'e}s intervenant dans cette suite.

Soient $Y$ un $k$-sch{\'e}ma et $G$ un $k$-groupe alg{\'e}brique (non-n{\'e}cessairement ab{\'e}lien). Ainsi on a une bijection (d'apr{\`e}s le corollaire 4.7 p.123 de \cite{Mi}):
	\[\textup{EHP}\left(G_{Y}/Y\right)\longrightarrow H^{1}\left(Y,G_{Y}\right)
\]
Rappelons bri{\`e}vement comment l'on obtient cette bijection.\vspace{2mm}

\begin{flushleft}
$\ \bullet$ \uline{1-cocycle associ{\'e} {\`a} un $G_{Y}$-torseur sur $Y$}
\end{flushleft}\vspace{2mm}

Soit $P$ un $G_{Y}$-torseur sur $Y$. Alors il existe un recouvrement {\'e}tale $\left(Y_{\alpha}/Y\right)_{\alpha\in A}$ tel que $P_{\left|Y_{\alpha}\right.}$ soit trivial, pour tout $\alpha\in A$. Le choix d'une section
	\[p_{\alpha}\in P\left(Y_{\alpha}\right)
\]
pour tout $\alpha$ fournit un isomorphisme de torseurs:
	\[\eta^{P}_{\alpha}:P_{\left|Y_{\alpha}\right.}\longrightarrow G_{Y_{\alpha},d},\ \forall \;\alpha\in A
\]
d{\'e}fini par:
	\[\eta^{P}_{\alpha}\left(p_{\alpha}\right)=e,\ \forall \;\alpha\in A
\]

Soit maintenant $\left(\alpha,\beta\right)\in A\times A$. La simple transitivit{\'e} de l'action de $G_{Y}\left(Y_{\alpha\beta}\right)$ (o{\`u} on a {\'e}videmment not{\'e} $Y_{\alpha\beta}=Y_{\alpha}\times_{Y}Y_{\beta}$) entra{\^i}ne l'existence et l'unicit{\'e} d'un $g_{\alpha\beta}^{P}\in G_{Y}\left(Y_{\alpha\beta}\right)$ tel que:
	\[{p_{\beta}}_{\left|Y_{\alpha\beta}\right.}={p_{\beta}}_{\left|Y_{\alpha\beta}\right.}.\;g_{\alpha\beta}
\]
cette {\'e}galit{\'e} ayant lieu dans $P\left(Y_{\alpha\beta}\right)$
Il est imm{\'e}diat\footnote{En utilisant encore la simple transitivit{\'e} de l'action de $G_{Y}$.} que la famille $\left(g_{\alpha\beta}\right)_{\alpha,\beta}\in A$ est un 1-cocycle {\`a} valeurs dans $G_{Y}$. 

La classe de ce 1-cocycle ne d{\'e}pend pas du choix des sections locales de $P$; en effet, soit $\left(p'_{\alpha}\right)_{\alpha}$ une autre famille de sections, avec:
	\[p_{\alpha}\in P\left(Y_{\alpha}\right),\ \forall\;\alpha\in A
\]
Alors:
	\[\forall\;\alpha\in A,\ \exists!\;h_{\alpha}^{P}\in G_{Y}\left(Y_{\alpha}\right)\ t.q:\ p'_{\alpha}=p_{\alpha}.\;h^{P}_{\alpha}\ \in P\left(Y_{\alpha}\right)
\]
\pagebreak

Au-dessus de $Y_{\alpha\beta}$, on dispose des relations:
	\[\left\{\begin{array}{c}{p'_{\alpha}}_{\left|Y_{\alpha\beta}\right.}={p_{\alpha}}_{\left|Y_{\alpha\beta}\right.}.\;{h^{P}_{\alpha}}_{\left|Y_{\alpha\beta}\right.}\\\ \\{p'_{\beta}}_{\left|Y_{\alpha\beta}\right.}={p_{\beta}}_{\left|Y_{\alpha\beta}\right.}.\;{h^{P}_{\beta}}_{\left|Y_{\alpha\beta}\right.}\\\ \\{p_{\beta}}_{\left|Y_{\alpha\beta}\right.}={p_{\alpha}}_{\left|Y_{\alpha\beta}\right.}.\;g^{P}_{\alpha\beta}\\\ \\{p'_{\beta}}_{\left|Y_{\alpha\beta}\right.}={p'_{\alpha}}_{\left|Y_{\alpha\beta}\right.}.\;{g'}^{P}_{\alpha\beta}\end{array}\right.
\]

D'une part, on d{\'e}duit de la premi{\`e}re et de la derni{\`e}re relation la nouvelle condition:
	\[{p'_{\beta}}_{\left|Y_{\alpha\beta}\right.}={p_{\alpha}}_{\left|Y_{\alpha\beta}\right.}.\;h^{P}_{\alpha}.\;{g'}^{P}_{\alpha\beta}
\]
et on obtient d'autre part {\`a} l'aide des deux autres relations:
	\[{p'_{\beta}}_{\left|Y_{\alpha\beta}\right.}={p_{\alpha}}_{\left|Y_{\alpha\beta}\right.}.\;{g}^{P}_{\alpha\beta}.\;h^{P}_{\beta}
\]
Par suite:
	\[{g'}^{P}_{\alpha\beta}=\left(h^{P}_{\alpha}\right)^{-1}.\;{g}^{P}_{\alpha\beta}.\;h^{P}_{\beta},\ \forall\;\left(\alpha,\beta\right)\in A\times A
\]
ce qui prouve justement que les 1-cocycles $\left(g_{\alpha\beta}\right)_{\alpha,\beta\in A}$ et $\left(g'_{\alpha\beta}\right)_{\alpha,\beta\in A}$ sont cohomologues.
\begin{rem}\textup{Explicitement, il revient au m{\^e}me de dire que le $G_{Y}$-torseur $P$ est obtenu en recollant les torseurs triviaux $G_{Y_{\alpha},d}$ {\`a} l'aide des $g_{\alpha\beta}$. Plus pr{\'e}cis{\'e}ment, le choix d'une famille de sections locales $\left(p_{\alpha}\right)_{\alpha\in A}$ pour $P$ entra{\^i}ne l'existence d'une famille d'automorphismes $\varphi_{\alpha\beta}$ des torseurs triviaux $G_{Y_{\alpha\beta},d}$, ce que l'on illustre sur par le diagramme suivant:}
	
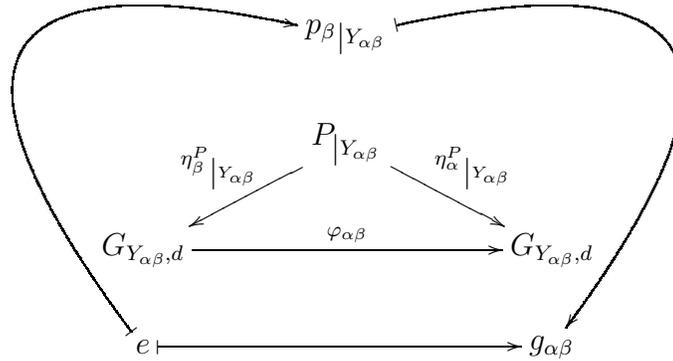
\[\xymatrix@R=20pt@C=40pt{&{p_{\beta}}_{\left|Y_{\alpha\beta}\right.} \ar@{|->}@/^100pt/[dddr] \\ & P_{\left|Y_{\alpha\beta}\right.} \ar[dl]_{{\eta^{P}_{\beta}}_{{\left|Y_{\alpha\beta}\right.}}} \ar[dr]^{{\eta^{P}_{\alpha}}_{{\left|Y_{\alpha\beta}\right.}}} \\G_{Y_{\alpha\beta},d} \ar[rr]^{\varphi_{\alpha\beta}} && G_{Y_{\alpha\beta},d} \\e\ar@{|->}[rr] \ar@{|->}@/^100pt/[uuur] &&g_{\alpha\beta}}
\]
\end{rem}
Pour achever ces rappels sur les torseurs, montrons que deux torseurs $P$ et $P'$ isomorphes donnent lieu {\`a} des 1-cocycles cohomologues. On peut supposer qu'il existe un recouvrement {\'e}tale $\left(Y_{\alpha}/Y\right)_{\alpha\in A}$ de $Y$ trivialisant $P$ et $P'$ (il suffit de prendre \textquotedblleft{l'intersection}\textquotedblright\ d'un recouvrement {\'e}tale trivialisant $P$ et d'un recouvrement {\'e}tale trivialisant $P'$). On note alors $\left(g^{P}_{\alpha\beta}\right)_{\alpha,\beta\in A}$ (\textit{resp.} $\left(g^{P'}_{\alpha\beta}\right)_{\alpha,\beta\in A}$) le 1-cocycle associ{\'e} comme pr{\'e}c{\'e}demment {\`a} $P$ (\textit{resp.} {\`a} $P'$) via le choix d'une famille de sections locales $\left(p_{\alpha}\right)_{\alpha\in A}$ (\textit{resp.} $\left(p'_{\alpha}\right)_{\alpha\in A}$). On suppose donc qu'il existe un isomorphisme de $G_{Y}$-torseurs sur $Y$:
	\[f:P'\longrightarrow P
\]
En utilisant la simple transitivit{\'e} de l'action de $G_{Y}\left(Y_{\alpha}\right)$ sur $P\left(Y_{\alpha}\right)$ pour tout $\alpha\in A$, on obtient:
	\[\forall\;\alpha\in A,\;\exists!\;h^{PP'}_{\alpha}\in G_{Y}\left(Y_{\alpha}\right)\ t.q:\ f_{\left|Y_{\alpha}\right.}\left(p'_{\alpha}\right)=p_{\alpha}.\;h^{PP'}_{\alpha}
\]

Pour tout couple $\left(\alpha,\beta\right)\in A\times A$, on a encore, au-dessus de $Y_{\alpha\beta}$:
	\[\left\{\begin{array}{c}f_{\left|Y_{\alpha\beta}\right.}\left({p'_{\alpha}}_{\left|Y_{\alpha\beta}\right.}\right)={p_{\alpha}}_{\left|Y_{\alpha\beta}\right.}.\;{h^{PP'}_{\alpha}}_{\left|Y_{\alpha\beta}\right.}\\\ \\f_{\left|Y_{\alpha\beta}\right.}\left({p'_{\beta}}_{\left|Y_{\alpha\beta}\right.}\right)={p_{\beta}}_{\left|Y_{\alpha\beta}\right.}.\;{h^{PP'}_{\beta}}_{\left|Y_{\alpha\beta}\right.}\\\ \\{p_{\beta}}_{\left|Y_{\alpha\beta}\right.}={p_{\alpha}}_{\left|Y_{\alpha\beta}\right.}.\;g^{P}_{\alpha\beta}\\\ \\{p'_{\beta}}_{\left|Y_{\alpha\beta}\right.}={p'_{\alpha}}_{\left|Y_{\alpha\beta}\right.}.\;g^{P'}_{\alpha\beta}\end{array}\right.
\]
De la deuxi{\`e}me et de la troisi{\`e}me relation, on d{\'e}duit:
	\[f_{\left|Y_{\alpha\beta}\right.}\left({p'_{\beta}}_{\left|Y_{\alpha\beta}\right.}\right)={p_{\alpha}}_{\left|Y_{\alpha\beta}\right.}.\;g_{\alpha\beta}^{P}.\;{h^{PP'}_{\beta}}_{\left|Y_{\alpha\beta}\right.}
\]
tandis que les deux extr{\^e}mes donnent:
	\[f_{\left|Y_{\alpha\beta}\right.}\left({p'_{\beta}}_{\left|Y_{\alpha\beta}\right.}\right)={p_{\alpha}}_{\left|Y_{\alpha\beta}\right.}.\;{h^{PP'}_{\alpha}}_{\left|Y_{\alpha\beta}\right.}.\;g_{\alpha\beta}^{P'}
\]
La simple transitivit{\'e} de l'action de $G_{Y}\left(Y_{\alpha\beta}\right)$ sur $P\left(Y_{\alpha\beta}\right)$ pour tout couple $\left(\alpha,\beta\right)\in A\times A$ permet alors de conclure:
	\[g_{\alpha\beta}^{P'}=\left({h^{PP'}_{\alpha}}_{\left|Y_{\alpha\beta}\right.}\right)^{-1}.\;g_{\alpha\beta}^{P}.\;{h^{PP'}_{\beta}}_{\left|Y_{\alpha\beta}\right.},\ \forall \left(\alpha,\beta\right)\in A\times A
\]

Ce qui ach{\`e}ve de prouver que l'ensemble des classes d'isomorphie de $G_{Y}$-torseurs sur $Y$ est en bijection avec les classes de cohomologie de 1-cocycles {\`a} valeurs dans $G_{Y}$.\vspace{2mm}

\begin{flushleft}
$\ \bullet$ \uline{Description de l'ensemble $H^{0}\left(k,R^{1}\pi_{\ast}G_{X}\right)$}
\end{flushleft}\vspace{2mm}

C'est l'ensemble des sections globales du $k$-faisceau $R^{1}\pi_{\ast}G_{X}$, qui est le faisceau associ{\'e} au pr{\'e}faisceau:
	\[\left(\textup{Spec}\;L\rightarrow \textup{Spec}\;k\right)\rightsquigarrow H^{1}\left(X_{L},G_{X_{L}}\right)
\]
$L$ d{\'e}signant dans la formule ci-dessus une extension {\'e}tale de $k$. Une section globale est donc une classe:
	\[\left[\left(L_{i}/k\right)_{i\in I},\left(P_{i}\right)_{i\in I},\left(\varphi_{ij}\right)_{i,j\in I}\right]
\]
o{\`u}:
\begin{enumerate}[(i)]
\item $L_{i}$ est une extension {\'e}tale de $k$, $\forall\;i\in I$;
\item pour tout $i\in I$, $P_{i}$ repr{\'e}sente une classe de $H^{1}\left(X_{L_{i}},G_{X_{L_{i}}}\right)$;
\item pour tout $\left(i,j\right)\in I\times I$:
	\[\varphi_{ij}:P_{j}\times_{X_{L_{j}}}X_{L_{ij}}\longrightarrow P_{i}\times_{X_{L_{i}}}X_{L_{ij}}
\]
est un isomorphisme de $G_{X_{L_{ij}}}$-torseurs sur $X_{L_{ij}}$, o{\`u} l'on a not{\'e} $L_{ij}=L_{i}\otimes_{k}L_{j}$.
\end{enumerate}
Deux telles classes
	\[\left[\left(L_{i}/k\right)_{i\in I},\left(P_{i}\right)_{i\in I},\left(\varphi_{ij}\right)_{i,j\in I}\right]\ \ \textup{et}\ \ \left[\left(L_{i'}/k\right)_{i'\in I},\left(P'_{i'}\right)_{i'\in I'},\left(\varphi_{i'j'}\right)_{i',j'\in I'}\right]
\]
sont {\'e}quivalentes si et seulement si il existe un recouvrement {\'e}tale\footnote{Il e{\^u}t fallu {\'e}crire: le recouvrement{\'e}tale $\left(\textup{Spec}\;L_{\alpha}\rightarrow\textup{Spec}\;k\right)_{\alpha\in A}$.} $\left(L_{\alpha}/k\right)_{\alpha\in A}$ raffinant l'intersection des deux pr{\'e}c{\'e}dents, et des isomorphismes:
	\[\psi_{\alpha}:P'_{i}\times_{X_{L_{i}}}X_{L_{\alpha}}\longrightarrow P_{i}\times_{X_{L_{i}}}X_{L_{\alpha}},\ \forall\;\alpha\in A
\]
compatibles avec les isomorphismes $\varphi_{ij}$ et $\varphi_{i'j'}$ dans le seul sens raisonnable que l'on puisse donner {\`a} cette formulation.

Nous passons maintenant {\`a} la description de la suite exacte d'ensembles point{\'e}s tant attendue.
\vspace{2mm}

\begin{flushleft}
$\ \bullet$ \uline{Description de l'application $a:H^{1}\left(k,\pi_{\ast}G_{X}\right)\longrightarrow H^{1}\left(X,G_{X}\right)$}
\end{flushleft}\vspace{2mm}

Soit $T$ un repr{\'e}sentant d'une classe de $H^{1}\left(k,\pi_{\ast}G_{X}\right)$; c'est donc la donn{\'e}e d'une famille d'extensions {\'e}tales $\left(L_{i}/k\right)_{i\in I}$ et d'un 1-cocycle:
	\[\left(\tau_{ij}\right)_{i,j}\in Z^{1}\left(\left(L_{i}/k\right),\pi_{\ast}G_{X}\right)
\]
Par cons{\'e}quent, pour tout $\left(i,j\right)\in I\times I$, on a:
	\[\tau_{ij}\in \pi_{\ast}G_{X}\left(L_{ij}\right)
\]
Or, par d{\'e}finition de l'image directe d'un faisceau:
	\[\pi_{\ast}G_{X}\left(L_{ij}\right)=G_{X}\left(\pi^{-1}\left(L_{ij}\right)\right)=G_{X}\left(X_{L_{ij}}\right)
\]
Notons maintenant:
	\[\widetilde{\tau_{ij}}\in G_{X}\left(X_{L_{ij}}\right)
\]
l'{\'e}l{\'e}ment $\tau_{ij}$, vu comme une section du faisceau $G_{X}$. Il est alors {\'e}vident que l'on obtient un 1-cocycle:
	\[\left(\widetilde{\tau_{ij}}\right)\in Z^{1}\left(\left(X_{L_{i}}/X\right),G_{X}\right)
\]

On peut alors d{\'e}finir $a\left(\left[T\right]\right)$ comme la classe du $G_{X}$-torseur sur $X$ correspondant {\`a} ce 1-cocycle. Pour cela, il suffit de v{\'e}rifier l'ind{\'e}pendance du repr{\'e}sentant choisi pour la classe de $T$: soit donc $\left(\tau'_{ij}\right)_{i,j\in I}$ un 1-cocycle {\`a} valeurs dans $\pi_{\ast}G_{X}$ cohomologue au 1-cocycle $\left(\tau_{ij}\right)_{i,j}$. Il existe donc une 1-cocha{\^i}ne:
	\[\left(h_{i}\right)_{i\in I}\in C^{1}\left(\left(L_{i}/k\right)_{i\in I},\pi_{\ast}G_{X}\right)
\]
telle que:
	\[\tau'_{ij}=\left({h_{j}}_{\left|L_{ij}\right.}\right)^{-1}.\;\tau_{ij}.\;{h_{i}}_{\left|L_{ij}\right.}\;\in \left[\pi_{\ast}G_{X}\right]\left(L_{ij}\right),\ \forall\ \left(i,j\right)\in I\times I
\]
Or cette relation implique {\'e}videmment la suivante:
	\[\widetilde{\tau'_{ij}}=\left(\widetilde{{h_{j}}}_{\left|X_{L_{ij}}\right.}\right)^{-1}.\;\widetilde{\tau_{ij}}.\;{\widetilde{h_{i}}}_{\left|X_{L_{ij}}\right.}\;\in \left[\pi_{\ast}G_{X}\right]\left(L_{ij}\right),\ \forall\ \left(i,j\right)\in I\times I
\]
o{\`u} $\widetilde{h_{i}}$ d{\'e}signe l'{\'e}l{\'e}ment $h_{i}$, qui appartient \textit{a priori} {\`a}:
	\[\left[\pi_{\ast}G_{X}\right]\left(L_{i}\right)
\]
vu comme un {\'e}l{\'e}ment de:
	\[G_{X}\left(X_{L_{i}}\right)
\]
Il s'ensuit que les 1-cocyles $\left(\widetilde{\tau_{ij}}\right)_{i,j\in I}$ et $\left(\widetilde{\tau'_{ij}}\right)_{i,j\in I}$ sont cohomologues, ce qui prouve la coh{\'e}rence de la d{\'e}finition de l'application $a$. A de tr{\`e}s l{\'e}g{\`e}res modifications pr{\`e}s, ce sont les m{\^e}mes arguments qui permettent de prouver la trivialit{\'e} du noyau (au sens des ensembles point{\'e}s\footnote{\textit{I.e.} la pr{\'e}image de la classe privil{\'e}gi{\'e}e de $H^{1}\left(X,G_{X}\right)$ par $a$ est la classe privil{\'e}gi{\'e}e de $H^{1}\left(k,\pi_{\ast}G_{X}\right)$.} bien s{\^u}r) de l'application $a$. En outre, une cons{\'e}quence imm{\'e}diate des pr{\'e}c{\'e}dents calculs est la suivante:
\begin{lem} L'image de l'application $a:H^{1}\left(k,\pi_{\ast}G_{X}\right)\longrightarrow H^{1}\left(X,G_{X}\right)$ s'identifie {\`a} l'ensemble des classes d'isomorphie de $G_{X}$-torseurs sur $X$ trivialis{\'e}s par un recouvrement {\'e}tale de $X$ \textquotedblleft{provenant de $k$}\textquotedblright. Plus pr{\'e}cis{\'e}ment, l'image de l'application $a$ est constitu{\'e}e des classes d'isomorphie de $G_{X}$-torseurs $P\rightarrow X$ pour lesquels il existe une famille d'extensions {\'e}tales $\left(L_{i}/k\right)_{i\in I}$ telle que:
	\[P\times_{X}X_{L_{i}}\simeq G_{X_{L_{i}},d},\ \forall\;i\in I
\]
\end{lem}
\begin{rem} \textup{Dans le contexte des chapitres pr{\'e}c{\'e}dents, l'image de $a$ est donc constitu{\'e}e des classes d'isomorphie de $G_{X}$-torseurs $P\rightarrow X$ tels que $\bar{P}\simeq \bar{G}_{X,d}$.}
\end{rem}
\begin{rem}\textup{En reprenant les notations du paragraphe V.3.3.1 de \cite{Gi2}, on peut encore exprimer le fait pr{\'e}c{\'e}dent en disant que l'image de $a$ est la cat{\'e}gorie:
	\[Tors\left(X,G_{X}\right)^{\left(\textup{Spec}\;k\right)_{\acute{e}t}}
\]
cette derni{\`e}re {\'e}tant une \uline{cat{\'e}gorie} {\'e}quivalente {\`a} la gerbe
	\[\textup{Tors}\left(k,\pi_{\ast}G_{X}\right)
\]
}
\end{rem}
\pagebreak

\begin{flushleft}
$\ \bullet$ \uline{Description de l'application $u:H^{1}\left(X,G_{X}\right)\longrightarrow H^{0}\left(k,R^{1}\pi_{\ast}G_{X}\right)$}
\end{flushleft}\vspace{2mm}

On obtient {\'e}videmment cette application en associant {\`a} la classe $\left[P\right]\in H^{1}\left(X,G_{X}\right)$ la classe:
	\[\left[\left(k/k\right),\left(P\right),\left(\textup{id}_{P}\right)\right]
\]
Il est maintenant imm{\'e}diat, vues les pr{\'e}c{\'e}dentes descriptions de l'application $a$ d'une part, et des sections globales du faisceau $R^{1}\pi_{\ast}G_{X}$ d'autre part, que le noyau (au sens non-ab{\'e}lien) est effectivement constitu{\'e} des classes d'isomorphie de $G_{X}$-torseurs trivialis{\'e}s par un recouvrement {\'e}tale de $X$ \textquotedblleft{provenant de $k$}\textquotedblright, ce qui ach{\`e}ve de prouver l'exactitude de la suite d'ensembles point{\'e}s:
	\[0\longrightarrow H^{1}\left(k,\pi_{\ast}G_{X}\right)\stackrel{a}{\longrightarrow} H^{1}\left(X,G_{X}\right)\stackrel{u}{\longrightarrow} H^{0}\left(k,R^{1}\pi_{\ast}G_{X}\right)
\]
\vspace{2mm}

\begin{flushleft}
$\ \bullet$ \uline{Obstruction non-ab{\'e}lienne}
\end{flushleft}\vspace{2mm}

On rentre maintenant dans le vif du sujet avec l'{\'e}tude de la surjectivit{\'e} de l'application $u$. Consid{\'e}rons une classe:
	\[\left[\left(L_{i}/k\right)_{i\in I},\left(P_{i}\right)_{i\in I},\left(\varphi_{ij}\right)_{i,j\in I}\right]
\]
dans $H^{0}\left(k,R^{1}\pi_{\ast}G_{X}\right)$. En choisissant un repr{\'e}sentant de cette classe, on dispose donc d'une famille de torseurs $\left(P_{i}\right)_{i\in I}$, et d'isomorphismes\footnote{De $G_{X_{L_{ij}}}$-torseurs sur $X_{L_{ij}}$.}:
	\[\varphi_{ij}:P_{j}\times_{X_{L_{j}}}X_{L_{ij}}\longrightarrow P_{j}\times_{X_{L_{i}}}X_{L_{ij}}
\]
En d'autres termes, la famille $\left(\varphi_{ij}\right)_{i,j\in I}$ constitue une donn{\'e}e de recollement pour les torseurs $P_{i}$, et on veut mesurer ce qui l'emp{\^e}che d'{\^e}tre une donn{\'e}e de descente.

Pour tout triplet $\left(i,j,k\right)\in I^{3}$, on a le diagramme suivant:
	\[\xymatrix@R=40pt{P_{k}\times_{X_{L_{k}}}X_{L_{ijk}} \ar[rr]^{{\varphi_{ik}}_{\left|X_{L_{ijk}}\right.}} \ar[rd]_{{\varphi_{jk}}_{\left|X_{L_{ijk}}\right.}} & \ar@{=>}[d]^{c_{ijk}}& P_{i}\times_{X_{L_{i}}}X_{L_{ijk}}\\& P_{j}\times_{X_{L_{j}}}X_{L_{ijk}} \ar[ru]_{{\varphi_{ij}}_{\left|X_{L_{ijk}}\right.}}}
\]
o{\`u} l'on a not{\'e}:
	\[c_{ijk}={\varphi_{ij}}_{\left|X_{L_{ijk}}\right.}\circ {\varphi_{jk}}_{\left|X_{L_{ijk}}\right.}\circ \left({\varphi_{ik}}_{\left|X_{L_{ijk}}\right.}\right)^{-1}
\]
En particulier:
	\[c_{ijk}\in\textup{ad}_{G_{X_{L_{ijk}}}}\left(P_{i}\times_{X_{L_{i}}}X_{L_{ijk}}\right),\ \forall\;\left(i,j,k\right)\in I^{3}
\]
c'est-{\`a}-dire encore:
	\[c_{ijk}\in\left[\pi_{\ast}\textup{ad}_{G_{X}}\left(P_{i}\right)\right]\left(L_{ijk}\right),\ \forall\;\left(i,j,k\right)\in I^{3}
\]
\begin{rem}\textup{Dans le cas ab{\'e}lien, on a {\'e}videmment:
	\[c_{ijk}\in\left[\pi_{\ast}G_{X}\right]\left(L_{ijk}\right),\ \forall\;\left(i,j,k\right)\in I^{3}
\]
et les $c_{ijk}$ constituent donc un 2-cocycle {\`a} valeurs dans $\pi_{\ast}G_{X}$; plus pr{\'e}cis{\'e}ment, la classe de ce 2-cocycle est justement l'image de la classe
	\[\left[\left(L_{i}/k\right)_{i\in I},\left(P_{i}\right)_{i\in I},\left(\varphi_{ij}\right)_{i,j\in I}\right]
\]
par le cobord
	\[H^{0}\left(k,R^{1}\pi_{\ast}G_{X}\right)\stackrel{\delta^{1}}{\longrightarrow} H^{2}\left(k,\pi_{\ast}G_{X}\right)
\]
Dans le cas non-ab{\'e}lien, la famille $\left(c_{ijk}\right)_{i,j,k\in I}$ pourrait encore {\^e}tre vue comme un 2-cocycle, mais {\`a} valeurs dans un \textit{syst{\`e}me de coefficients} \textit{cf.} \cite{Do1}, pr{\'e}cis{\'e}ment celui d{\'e}fini par la famille de faisceaux locaux $\left(\pi_{\ast}\textup{ad}_{G_{X}}\left(P_{i}\right)\right)_{i\in I}$.}
\end{rem}
Il s'agit maintenant de comprendre pourquoi (dans le cas ab{\'e}lien comme dans le cas g{\'e}n{\'e}ral) ce 2-cocycle devient trivial, c'est-{\`a}-dire de comprendre pourquoi les $P_{i}$ se recollent effectivement, \textquotedblleft{une fois que l'on passe {\`a} $X$}\textquotedblright. Notons tout d'abord que l'on peut, comme dans la description de l'application $a$, associer au 2-cocycle constitu{\'e} par les
	\[c_{ijk}\in \left[\pi_{\ast}\textup{ad}_{G_{X}}\left(P_{i}\right)\right]\left(L_{ijk}\right)
\]
un 2-cocycle sur $X$:
	\[\widetilde{c_{ijk}}\in \left[\textup{ad}_{G_{X}}\left(P_{i}\right)\right]\left(X_{L_{ijk}}\right)
\]
(dans le cas ab{\'e}lien, on a $\widetilde{c_{ijk}}\in G_{X}\left(X_{L_{ijk}}\right)$). 
Or ce dernier peut {\^e}tre rendu trivial, puisque l'on dispose d'assez d'ouverts sur $X$ pour pouvoir rendre tous les torseurs $P_{i}$ triviaux.
\begin{rem}\textup{Une autre mani{\`e}re de dire les choses est que la gerbe \uline{sur $X$} engendr{\'e}e par les torseurs $P_{i}$ est la gerbe des $G_{X}$-torseurs sur $X$, du fait que l'on peut localiser suffisamment de mani{\`e}re {\`a} rendre tous les $P_{i}$ triviaux, ce qui n'est {\'e}videmment pas le cas si l'on se restreint aux ouverts {\'e}tales de $X$ \textquotedblleft{provenant de $k$}\textquotedblright.}
\end{rem}

Il reste encore {\`a} voir que l'existence d'un point $k$-rationnel sur $X$ entra{\^i}ne l'existence d'une r{\'e}traction de l'application:
	\[\left[\left(c_{ijk}\right)_{i,j,k}\right]\longmapsto \left[\left(\widetilde{c_{ijk}}\right)\right]
\]
On note maintenant
	\[x: \textup{Spec}\;k\longrightarrow X
\]
un point $k$-rationnel de $X$. Dans le cas ab{\'e}lien, la r{\'e}traction est obtenue en envoyant la classe de $\left(\widetilde{c_{ijk}}\right)$ sur la classe du 2-cocycle:
	\[\left(q\circ\widetilde{c_{ijk}}\circ x_{ijk}\right)
\]
o{\`u}:
	\[x_{ijk}:\textup{Spec}\;L_{ijk}\longrightarrow X_{ijk}
\]
est induit par $x$, et o{\`u}
	\[q:G_{X}\longrightarrow G
\]
est le morphisme {\'e}vident (puisque $G_{X}$ est le produit fibr{\'e} $G\times_{\textup{Spec}\;k}X$). Evidemment, le 2-cocycle $\left(q\circ\widetilde{c_{ijk}}\circ x_{ijk}\right)_{i,j,k\in I}$ ainsi obtenu n'est qu'un cocycle {\`a} valeurs dans $G$, et ne correspond au cocycle initial $\left(c_{ijk}\right)_{i,j,k\in I}$ que sous r{\'e}serve que la condition:
	\[\bar{G}_{X}\left(\bar{X}\right)=\bar{G}\left(\bar{k}\right)
\]
soit satisfaite.

Dans le cas non-ab{\'e}lien, on applique exactement le m{\^e}me raisonnement, pour aboutir {\`a} la m{\^e}me conclusion, {\`a} la diff{\'e}rence (de taille!) pr{\`e}s que l'on doit imposer la condition beaucoup plus contraignante suivante, si l'on souhaite pouvoir descendre \uline{tous} les $G_{\bar{X}}$-torseurs sur $\bar{X}$ de corps des modules $k$:\vspace{2mm}

\uline{\textbf{Condition} $\left(\star\right)$}: \textit{Pour toute $\bar{X}$-forme int{\'e}rieure $G'$ de $\bar{G}_{X}$, on a:}
	\[H^{0}\left(\bar{X},G'\right)=H^{0}\left(\bar{k},G'_{\bar{x}}\right)
\]\vspace{2mm}

Une justification plus rigoureuse de ce fait est apport{\'e}e dans la section suivante.
\end{section}
\begin{section}[Obstruction non-ab{\'e}lienne]{Obstruction non-ab{\'e}lienne {\`a} l'existence d'un point $k$-rationnel}
Soit $\bar{P}\rightarrow \bar{X}$ un $\bar{G}_{X}$-torseur de corps des modules $k$. On note $G'=\textup{ad}_{G_{\bar{X}}}\left(\bar{P}\right)$ le faisceau sur $\bar{X}$ de ses automorphismes. D'apr{\`e}s le lemme \ref{lemf}, c'est une forme int{\'e}rieure sur $\bar{X}$ de $\bar{G}_{X}$. On note encore:
	\[\bar{x}:\textup{Spec}\;\bar{k}\longrightarrow \bar{X}
\]
le morphisme induit par $x$, et $G'_{\bar{x}}$ la fibre en ce point de $G'$. On dispose d'une fl{\`e}che naturelle:
	\[\varphi_{\bar{x}}:H^{0}\left(\bar{X},G'\right)\longrightarrow H^{0}\left(\textup{Spec}\;\bar{k},G'_{\bar{x}}\right)
\]
\begin{rem}\textup{Lorsque $G$ est ab{\'e}lien, demander que la condition $\bar{G}_{X}\left(\bar{X}\right)=\bar{G}\left(\bar{k}\right)$ soit satisfaite {\'e}quivaut encore {\`a} demander que l'application ci-dessus soit bijective (et cette hypoth{\`e}se n'est pas trop contraignante, du fait qu'il n'y a pas de $\bar{X}$-formes int{\'e}rieures de $\bar{G}_{X}$ autres que $\bar{G}_{X}$ lui-m{\^e}me).}
\end{rem}
\begin{theo} Soient $X$ un $k$-sch{\'e}ma, et $G$ un $k$-groupe lin{\'e}aire, et $\bar{P}\rightarrow \bar{X}$ un $G_{\bar{X}}$-torseur de corps des modules $k$. On suppose que $X$ poss{\`e}de un point $k$-rationnel $x$, et on suppose satisfaite la condition suivante:

\textbf{Condition $\left(\star_{\bar{P}}\right)$:} En notant $G'=\textup{ad}_{G_{\bar{X}}}\left(\bar{P}\right)$ et $\bar{x}$ un point g{\'e}om{\'e}trique associ{\'e} au point $k$-rationnel $x$, on a:
	\[H^{0}\left(\bar{X},G'\right)=H^{0}\left(\bar{k},G'_{\bar{x}}\right)
\]

Alors $\bar{P}\rightarrow\bar{X}$ est d{\'e}fini sur $k$.
\end{theo}

\uline{\scshape{Preuve}}. On note encore $G'=\textup{ad}_{G_{\bar{X}}}\left(\bar{P}\right)$ le faisceau sur $\bar{X}$ des automorphismes de $\bar{P}$. Soit $\bar{p}$ un point dans la fibre $\bar{f}^{-1}\left(\bar{x}\right)$. Puisque $\bar{P}$ est de corps des modules $k$, il existe un isomorphisme:
	\[\varphi_{\sigma}:\bar{P}\longrightarrow \ ^{\sigma}\bar{P},\ \forall\ \sigma\in \Gamma
\]

Par ailleurs $\bar{P}_{\bar{x}}$ est un $\left(\bar{G}_{X}\right)_{\bar{x}}$-torseur, de groupe d'automorphismes $G'_{\bar{x}}$, ce dernier agissant simplement transitivement {\`a} droite sur $P_{\bar{x}}$. Comme $x$ est un point $k$-rationnel, tous les isomorphismes $\varphi_{\sigma}$ respectent la fibre au-dessus de $\bar{x}$, pour tout $\sigma\in \Gamma$. De plus, puisque la condition $\left(\star\right)$ est satisfaite, on peut composer tout $\varphi_{\sigma}$ {\`a} droite par un {\'e}l{\'e}ment de $G'\left(\bar{X}\right)$ de telle sorte que l'on ait:
	\[\varphi_{\sigma}\left(\bar{p}\right)=\ ^{\sigma}\bar{p}
\]

En utilisant toujours la condition $\left(\star\right)$, cette propri{\'e}t{\'e} d{\'e}termine uniquement $\varphi_{\sigma}$. Il s'ensuit imm{\'e}diatemment que pour tout couple $\left(\tau,\sigma\right)$ d'{\'e}l{\'e}ments de $\Gamma$, le diagramme ci-dessous commute:
	\[\xymatrix{\bar{P}\ar[rr]^{\varphi_{\tau\sigma}}\ar[dr]_{\varphi_{\tau}} && ^{\tau\sigma}\bar{P} \\&^{\tau}\bar{P}\ar[ur]_{^{\tau}\varphi_{\sigma}} }
\]
Donc $\bar{f}:\bar{P}\rightarrow \bar{X}$ est d{\'e}fini sur $k$.
\begin{flushright}
$\Box$
\end{flushright}
\begin{rem}\textup{Evidemment, la condition $\left(\star\right)$ est tr{\`e}s lourde. On aurait envie qu'elle soit v{\'e}rifi{\'e}e par exemple pour les sch{\'e}mas en groupes lin{\'e}aires au-dessus de vari{\'e}t{\'e}s projectives. Ce n'est absolument pas le cas, comme le montre l'exemple suivant: on consid{\`e}re la droite projective $X=\mathbb{P}^{1}_{\overline{\mathbb{Q}}}$, munie du recouvrement par les ouverts standards $\left\{U_{0},U_{\infty}\right\}$, o{\`u}:
	\[U_{0}=\left\{\left(x:y\right)/y\neq0\right\}\ \ et\ \ U_{\infty}=\left\{\left(x:y\right)/x\neq0\right\}
\]
Nous noterons:
	\[U_{0\infty}=U_{0}\cap U_{\infty}
\]
On construit maintenant une forme int{\'e}rieure $G'$ du faisceau $GL_{2,\mathbb{P}^{1}}$ en recollant les faisceaux $GL_{2,U_{0}}$ et $GL_{2,U_{\infty}}$ au-dessus de $U_{0\infty}$ {\`a} l'aide de la conjugaison:
	\[\xymatrix{\chi:&GL_{2,U_{0,\infty}}=\left(GL_{2,U_{0}}\right)_{\left|U_{0\infty}\right.}\ar[r]& GL_{2,U_{0,\infty}}=\left(GL_{2,U_{\infty}}\right)_{\left|U_{0\infty}\right.}\\&A\ar@{|->}[r]&M_{0\infty}.\;A.\;M_{0\infty}^{-1}}
\]
o{\`u}:
	\[M_{0\infty}=\left(\begin{array}{cc}x^{-1}&0\\0&x\end{array}\right)
\]
Alors la section:
	\[A_{0}=\left(\begin{array}{cc}1&x\\0&1\end{array}\right)
\]
de $GL_{2,U_{0}}$ fournit une section globale du faisceau $G'$
	\[\left(\textup{puisque:}\ \ M_{0\infty}.\;A_{0}.\;M_{0\infty}^{-1}=\left(\begin{array}{cc}1&x^{-1}\\0&1\end{array}\right)\right)
\]
non-constante. On ne peut donc pas esp{\'e}rer que $G'$ provienne d'un $\mathbb{Q}$-groupe alg{\'e}brique ici.}
\end{rem}
\end{section}
\begin{section}{Le cas fini}
\begin{theo} Soient $k$ un corps de caract{\'e}ristique nulle, $X$ un $k$-sch{\'e}ma g{\'e}om{\'e}triquement connexe, et $G$ un $k$-groupe fini. Si $X\left(k\right)\neq\emptyset$, alors tout $\bar{G}_{X}$-torseur sur $\bar{X}$ de corps des modules $k$ est d{\'e}fini sur $k$.
\end{theo}

\uline{\textsc{Preuve}}. Notons $x:\textup{Spec}\;k\rightarrow X$ un point $k$-rationnel de $X$, et soit $\bar{P}\rightarrow\bar{X}$ un $\bar{G}_{X}$-torseur de corps des modules $k$. On a alors un morphisme naturel de $k$-champs\footnote{En effet, $D\left(\bar{P}\right)$ est une gerbe, mais $\pi_{\ast}\textup{Tors}\left(X,G_{X}\right)$ n'est qu'un champ (\textit{cf.} exemple 1.5.18). Par ailleurs, ce morphisme fait de $D\left(\bar{P}\right)$ une \textit{sous-gerbe maximale} du champ $\pi_{\ast}\textup{Tors}\left(X,G_{X}\right)$ (\textit{cf.} \cite{Gi1} p.113).}:
	\[\epsilon:D\left(\bar{P}\right)\longrightarrow \pi_{\ast}\textup{Tors}\left(X,G_{X}\right)
\]

La composante au-dessus d'un ouvert {\'e}tale $\left(\textup{Spec}\;L\rightarrow\textup{Spec}\;k\right)$ du morphisme en question est le foncteur {\'e}vident:
	\[\xymatrix@R=5pt@C=20pt{\epsilon_{L}:&D\left(\bar{P}\right)\left(L\right)\ar[r]&Tors\left(X_{L},G_{X_{L}}\right)\\& P_{L} \ar@{|->}[r]&P_{L}}
\]

En utilisant la formule d'adjonction pour les champs\index{formule d'adjonction} (\textit{cf.} \cite{Gi1} I.5.2.1):
	\[\underline{\textup{Cart}}_{X_{\acute{e}t}}\left(\pi^{\ast}D\left(\bar{P}\right),\textup{Tors}\left(X,G_{X}\right)\right)\approx \underline{\textup{Cart}}_{\left(\textup{Spec}\;k\right)_{\acute{e}t}}\left(D\left(\bar{P}\right),\pi_{\ast}\textup{Tors}\left(X,G_{X}\right)\right)
\]
on obtient l'existence d'un morphisme de gerbes (puisque l'image inverse d'une gerbe est une gerbe, \cite{Gi2} V.1.4.2):
	\[\phi:\pi^{\ast}D\left(\bar{P}\right)\longrightarrow \textup{Tors}\left(X,G_{X}\right)
\]

On veut prouver que $\phi$ est une {\'e}quivalence. Notons $\mathcal{L}_{\bar{P}}$ le lien de la gerbe $D\left(\bar{P}\right)$. Comme la gerbe image inverse par $\pi$ de $D\left(\bar{P}\right)$ est li{\'e}e par l'image inverse du lien de $D\left(\bar{P}\right)$ (\textit{cf.} \cite{Gi2} V.1.4.2), le morphisme $\phi$ est li{\'e} par un morphisme:
	\[\Lambda:\pi^{\ast}\mathcal{L}_{\bar{P}}\longrightarrow \textup{lien}\;G_{X}
\]

D'autre part, il existe un recouvrement {\'e}tale $\left(\textup{Spec}\;L_{i}\rightarrow\textup{Spec}\;k\right)_{i\in I}$, et des mod{\`e}les locaux pour $\bar{P}$ au-dessus de ces ouverts. Plus pr{\'e}cis{\'e}ment, pour tout $i\in I$, il existe un $G_{X_{i}}$-torseur\footnote{Pour all{\'e}ger les notations, on a not{\'e} $X_{i}$ le sch{\'e}ma $X\otimes_{k}L_{i}$.} $P_{i}$ sur $X_{i}$ tel que:
	\[P_{i}\times_{X_{i}}\bar{X}\approx \bar{P}
\]

Au-dessus d'un ouvert $\left(X_{i}\rightarrow X\right)$, le morphisme $\Lambda$ est repr{\'e}sent{\'e} par le morphisme naturel de faisceaux de groupes sur le site {\'e}tale de $X_{i}$ (\textit{cf.} \cite{Gi1} IV.3.5.4):
	\[\lambda_{i}:\pi^{\ast}\pi_{\ast}ad_{G_{X_{i}}}\left(P_{i}\right)\longrightarrow ad_{G_{X_{i}}}\left(P_{i}\right)
\]
Du fait que ces $ad_{G_{X_{i}}}\left(P_{i}\right)$ sont finis, le morphisme d'adjonction ci-dessus est un isomorphisme (d'apr{\`e}s la remarque 3.1.(d) de \cite{Mi}) et il s'ensuit que dans ce cas le morphisme 
	\[\phi:\pi^{\ast}D\left(\bar{P}\right)\longrightarrow \textup{Tors}\left(X,G_{X}\right)
\]
est un isomorphisme de gerbes. Donc la gerbe $\pi^{\ast}D\left(\bar{P}\right)$ est neutre. Par suite, la gerbe $x^{\ast}\pi^{\ast}D\left(\bar{P}\right)$ est neutre. Comme $\pi\circ x=\textup{id}_{\textup{Spec}\;k}$, on en d{\'e}duit que la gerbe $D\left(\bar{P}\right)$, qui est {\'e}quivalente {\`a} la gerbe $x^{\ast}\pi^{\ast}D\left(\bar{P}\right)$, est elle aussi neutre, ce qui ach{\`e}ve la preuve.
\begin{flushright}
$\Box$
\end{flushright}
\end{section}
\begin{section}{Remarques}
\begin{subsection}{Probl{\`e}mes de repr{\'e}sentabilit{\'e}}
Il ressort des deux chapitres pr{\'e}c{\'e}dents qu'il existe des obstacles de taille {\`a} ce que l'on puisse g{\'e}n{\'e}raliser au cas non-ab{\'e}lien la suite {\`a} 5 termes d{\'e}duite de la suite spectrale de Leray. Le premier, et non des moindres, est (pour citer \cite{Gi2}) \textquotedblleft{qu'il n'existe pas de notion raisonnable pour l'image directe d'un lien}\textquotedblright.

On ne peut en g{\'e}n{\'e}ral pas dire grand chose du lien de la gerbe $D\left(\bar{P}\right)$, m{\^e}me lorsque le $k$-groupe $G$ initial jouit des meilleures propri{\'e}t{\'e}s dont l'on puisse r{\^e}ver. Il serait par exemple tentant de croire que si $G$ est un $k$-groupe r{\'e}ductif (\textit{e.g.} $GL_{n}\left(k\right)$) alors le lien de la gerbe des mod{\`e}les d'un $\bar{G}_{X}$-torseur de corps des modules $k$ est localement repr{\'e}sentable par une $k$-forme de $G$. Il n'en est rien en g{\'e}n{\'e}ral, et voici (me semble-t-il) une explication possible: m{\^e}me si l'on conna{\^i}t parfaitement la structure des $k$-groupes r{\'e}ductifs (et plus g{\'e}n{\'e}ralement des sch{\'e}mas en groupes r{\'e}ductifs, gr{\^a}ce {\`a} \cite{De1}), on ne peut pas arriver {\`a} avoir des informations exploitables sur $\textup{lien}\;D\left(\bar{P}\right)$, attendu que ce $k$-lien provient dans un certain sens de formes int{\'e}rieures sur $X$ du sch{\'e}ma en groupes $G_{X}$. Grossi{\`e}rement, le passage de $X$ {\`a} $k$ brouille toutes les pistes. 
\end{subsection}
\begin{subsection}{La condition corps des modules et le champ $\pi_{\ast}\textup{Tors}\left(X,G_{X}\right)$}
Consid{\'e}rons de nouveau un corps $k$ de caract{\'e}ristique nulle, $\pi:X\rightarrow\textup{Spec}\;k$ un $k$-sch{\'e}ma, un $k$-groupe alg{\'e}brique $G$ et un $\bar{G}_{X}$-torseur $\bar{P}\rightarrow \bar{X}$, que l'on ne suppose pas pour l'instant de corps des modules $k$, autrement dit: $\left[\bar{P}\right]\in H^{1}\left(\bar{X},\bar{G}_{X}\right)$.\vspace{2mm}
\begin{itemize}
\item \uline{Si $\bar{P}$ est de corps des modules $k$}: on peut lui associer sa gerbe des mod{\`e}les $D\left(\bar{P}\right)$. C'est une $k$-gerbe. On dispose de plus d'un monomorphisme de $k$-champs {\'e}vident (que nous avions d{\'e}j{\`a} {\'e}voqu{\'e} dans la preuve du lemme 3.1.2):
	\[\epsilon:D\left(\bar{P}\right)\longrightarrow \pi_{\ast}\textup{Tors}\left(X,G_{X}\right)
\]
dont la composante au-dessus d'un ouvert {\'e}tale $\left(\textup{Spec}\;L\rightarrow\textup{Spec}\;k\right)$ est le foncteur:
	\[\xymatrix@R=5pt@C=20pt{\epsilon_{L}:&D\left(\bar{P}\right)\left(L\right)\ar[r]&Tors\left(X_{L},G_{X_{L}}\right)\\& P_{L} \ar@{|->}[r]&P_{L}}
\]
dont la d{\'e}finition est coh{\'e}rente, puisqu'un objet de cette gerbe au-dessus de $\textup{Spec}\;L$ est en particulier un $G_{X_{L}}$-torseur sur $X_{L}$. Dans un certain sens (que nous n'allons pas tarder {\`a} pr{\'e}ciser) le $k$-champ $\pi_{\ast}\textup{Tors}\left(X,G_{X}\right)$ \textquotedblleft{contient}\textquotedblright\ les gerbes des mod{\`e}les de tous les $\bar{G}_{X}$-torseurs sur $\bar{X}$ de corps des modules $k$;\vspace{3mm}
\item \uline{Si $\bar{P}$ n'est pas de corps des modules $k$}: on peut toujours d{\'e}finir la $k$-cat{\'e}gorie fibr{\'e}e en groupoïdes $D\left(\bar{P}\right)$ en prenant comme cat{\'e}gorie fibre au-dessus d'un ouvert {\'e}tale $\left(\textup{Spec}\;L\rightarrow\textup{Spec}\;k\right)$ le groupoïde dont:\vspace{2mm}
\begin{itemize}
\item \textbf{les objets} sont toujours les $G_{X_{L}}$-torseurs $P_{L}$ sur $X_{L}$ tels qu'il existe un automorphisme $\sigma\in\Gamma$ et un isomorphisme de $\bar{G}_{X}$-torseurs sur $\bar{X}$:
	\[^{\sigma}\bar{P}\stackrel{\approx}{\longrightarrow} P_{L}\otimes_{L}\bar{k}
\]
\item \textbf{les fl{\`e}ches} sont les isomorphismes de $G_{X_{L}}$-torseurs sur $X_{L}$.\vspace{3mm}
\end{itemize}
$D\left(\bar{P}\right)$ est toujours un $k$-champ. Il est m{\^e}me localement non-vide, puisque comme (\textit{cf.} \cite{SGA4-7} 5.7 et 5.14):
	\[H^{1}\left(\bar{X},\bar{G}_{X}\right)=\lim_{\xrightarrow[L]{}}H^{1}\left(X_{L},G_{X_{L}}\right)
\]
la limite directe {\'e}tant prise sur les extensions {\'e}tales de $k$, il existe une extension {\'e}tale $L'$ et un {\'e}l{\'e}ment $\left[P_{L'}\right]$ de $H^{1}\left(X_{L'},G_{X_{L'}}\right)$ tels que: $P_{L'}\otimes_{L'}\bar{k}\approx \bar{P}$. 

Mais ce champ $D\left(\bar{P}\right)$ n'est pas une $k$-gerbe, car deux objets ne sont pas localement isomorphes, du fait que pour $\sigma\in\Gamma$, les torseurs $\bar{P}$ et $^{\sigma}\bar{P}$ ne sont pas en g{\'e}n{\'e}ral isomorphes. Cependant, il est possible d'associer {\`a} $\bar{P}$ une gerbe au-dessus de son corps des modules $k_{\bar{P}}$, ce dernier {\'e}tant l'intersection de toutes les extensions {\'e}tales $L$ de $k$ telles que:\footnote{Pour donner une interpr{\'e}tation topologique du corps des modules, disons que l'ouvert {\'e}tale $\left(\textup{Spec}\;k_{\bar{P}}\rightarrow\textup{Spec}\;k\right)$ est le plus grand ouvert {\'e}tale sur lequel on peut associer une donn{\'e}e de recollement {\`a} $\bar{P}$.}
	\[^{\tau}\bar{P}\approx \bar{P},\ \forall\;\tau\in\textup{Gal}\left(\bar{k}/L\right).
\]
\end{itemize}
\vspace{2mm}

On peut traduire de deux mani{\`e}res les remarques ci-dessus:\vspace{2mm}
\begin{itemize}
\item \uline{dans le langage de Giraud} (\textit{cf.} \cite{Gi1} et \cite{Gi2}), la premi{\`e}re remarque signifie que la gerbe $D\left(\bar{P}\right)$ est une section globale du faisceau des sous-gerbes maximales de $\pi_{\ast}\textup{Tors}\left(X,G_{X}\right)$, et la seconde signifie que la gerbe $D\left(\bar{P}\right)$ est une section au-dessus de l'ouvert {\'e}tale $\left(\textup{Spec}\;k_{\bar{P}}\rightarrow\textup{Spec}\;k\right)$ de ce faisceau. Ces remarques sont contenues dans l'{\'e}nonc{\'e} suivant:
\begin{propri}[Lemme V.3.1.5 de \cite{Gi2}] Soient $f:E'\rightarrow E$ un morphisme de sites et $A$ un faisceau de groupes sur $E'$. On a un isomorphisme canonique d'ensembles sur $E$ entre $R^{1}f_{\ast}A$ et le faisceau des sous-gerbes maximales du champ $f_{\ast}\textup{Tors}\left(E',A\right)$.
\end{propri}

En appliquant ce lemme {\`a} la pr{\'e}sente situation, \textit{i.e.} en prenant pour $E$ (\textit{resp.} $E'$) le site {\'e}tale de $k$ (\textit{resp.} de $X$), pour $f$ le morphisme de sites induit par le morphisme structural $\pi$ et pour $A$ le faisceau $G_{X}$, on obtient ainsi:
\begin{lem} On a un isomorphisme canonique d'ensembles sur $k$ entre les classes d'isomorphie de $\bar{G}_{X}$-torseurs sur $\bar{X}$ et le faisceau des sous-gerbes maximales du champ $\pi_{\ast}\textup{Tors}\left(X,G_{X}\right)$.
\end{lem}

En effet, on peut associer {\`a} un $G_{X}$-torseur sur $\bar{X}$ une $k_{\bar{P}}$-gerbe, et celle-ci est une sous-gerbe maximale du $k_{\bar{P}}$-champ $\pi_{\ast}\textup{Tors}\left(X,G_{X}\right)\otimes_{k}k_{\bar{P}}$, \textit{i.e.} une section au-dessus de $k_{\bar{P}}$ du faisceau des sous-gerbes maximales de $\pi_{\ast}\textup{Tors}\left(X,G_{X}\right)$;\vspace{3mm}
\item \uline{dans le langage plus moderne de Laumon et Moret-Bailly} (\textit{cf.} \cite{LMB}), on traduit la premi{\`e}re remarque en disant que $\bar{P}$ (\textit{resp.} $\left[\bar{P}\right]$) est un point $k$-rationnel\footnote{Au sens de la d{\'e}finition 5.2 de \cite{LMB}.}\index{point!k-rationnel@$k$-rationnel d'un champ} du champ $\pi_{\ast}\textup{Tors}\left(X,G_{X}\right)$ (\textit{resp.} une section globale du faisceau {\'e}tale grossier \footnote{Nous renvoyons {\`a} \cite{LMB} 3.19 ou {\`a} la preuve du lemme 3.3.2.2 pour la d{\'e}finition de faisceau grossier attach{\'e} {\`a} un champ.} $R^{1}\pi_{\ast}G_{X}$ associ{\'e} {\`a} ce champ), et la seconde en disant que $\bar{P}$ (\textit{resp.} $\left[\bar{P}\right]$) est un point {\`a} valeurs dans $\textup{Spec}\;k_{\bar{P}}$ du champ $\pi_{\ast}\textup{Tors}\left(X,G_{X}\right)$ (\textit{resp.} une section au-dessus de $k_{\bar{P}}$ du faisceau grossier $R^{1}\pi_{\ast}G_{X}$). Dans les deux cas, la gerbe $D\left(\bar{P}\right)$ est la gerbe r{\'e}siduelle en le point $\bar{P}$ du champ $\pi_{\ast}\textup{Tors}\left(X,G_{X}\right)$. L'analogue du lemme V.3.1.5 de \cite{Gi2} est fourni par l'{\'e}nonc{\'e} ci-dessous:
\begin{propri}[Corollaire 11.4 de \cite{LMB}] Soient $S$ un sch{\'e}ma, $\mathfrak{X}$ un $S$-champ alg{\'e}brique noeth{\'e}rien, et $\xi$ un point du $S$-champ alg{\'e}brique $\mathfrak{X}$. L'application qui {\`a} $\xi$ associe $\mathcal{G}_{\bar{\xi}}$ (la gerbe r{\'e}siduelle\index{gerbe!r{\'e}siduelle} du champ $\mathfrak{X}$ au point $\xi$) d{\'e}finit une bijection de l'ensemble des classes d'isomorphie de couples
	\[\left(\mathcal{G}\stackrel{A}{\longrightarrow} \textup{Spec}\;K,\mathcal{G}\stackrel{I}{\longrightarrow}\mathfrak{X}\right)
\]
o{\`u} $\mathcal{G}$ est une gerbe fppf sur un $S$-corps $K$ et o{\`u} $\mathcal{G}\stackrel{I}{\longrightarrow}\mathfrak{X}$ est un monomorphisme de $S$-champs.
\end{propri}

Nous utiliserons cette propri{\'e}t{\'e} (o{\`u}, avec nos hypoth{\`e}ses, on peut remplacer \textquotedblleft{fppf}\textquotedblright\ par \textquotedblleft{{\'e}tale}\textquotedblright) avec $S=\textup{Spec}\;k$ et $\mathfrak{X}=\pi_{\ast}\textup{Tors}\left(X,G_{X}\right)$. A tout $\bar{G}_{X}$-torseur sur $\bar{X}$, on peut associer sa gerbe des mod{\`e}les $D\left(\bar{P}\right)$; comme on l'a vu plus haut, c'est une $k_{\bar{P}}$-gerbe ($k_{\bar{P}}$ {\'e}tant le corps des modules de $\bar{P}$), que l'on peut toujours \textquotedblleft{plonger}\textquotedblright\ naturellement dans le $k$-champ $\pi_{\ast}\textup{Tors}\left(X,G_{X}\right)$ (essentiellement parce qu'un mod{\`e}le de $\bar{P}$ au-dessus d'une extension {\'e}tale $L$ est un $G_{X_{L}}$-torseur sur $X_{L}$, donc en particulier un objet du groupoïde fibre $\pi_{\ast}\textup{Tors}\left(X,G_{X}\right)\left(L\right)$). Alors le morphisme $A$ du corollaire 11.4 de \cite{LMB} n'est autre que le morphisme structural:
	\[D\left(\bar{P}\right)\longrightarrow \textup{Spec}\;k_{\bar{P}}
\]
de la $k_{\bar{P}}$-gerbe des mod{\`e}les, et le monomorphisme $I$ est le morphisme naturel:
	\[D\left(\bar{P}\right)\longrightarrow \pi_{\ast}\textup{Tors}\left(X,G_{X}\right)
\]
\end{itemize}
\vspace{2mm}
Achevons ces consid{\'e}rations par un {\'e}nonc{\'e} (trivial) faisant le lien entre les deux points de vue que l'on vient de comparer:
\begin{lemsub} Avec les notations adopt{\'e}es tout au long de cette section, le faisceau $R^{1}\pi_{\ast}G_{X}$ est le faisceau grossier (sur le site {\'e}tale de $k$) attach{\'e} au $k$-champ $\pi_{\ast}\textup{Tors}\left(X,G_{X}\right)$.
\end{lemsub}

\uline{\textsc{Preuve}}: le faisceau grossier attach{\'e} au $k$-champ $\pi_{\ast}\textup{Tors}\left(X,G_{X}\right)$ est le faisceau associ{\'e} au pr{\'e}faisceau:
	\[\left(\textup{Spec}\;L\rightarrow\textup{Spec}\;k\right)\rightsquigarrow \left\{\textup{classes d'isomorphie d'objets de }\pi_{\ast}\textup{Tors}\left(X,G_{X}\right)\left(L\right)\right\}
\]
soit encore:
	\[\left(\textup{Spec}\;L\rightarrow\textup{Spec}\;k\right)\rightsquigarrow \left\{\textup{classes d'isomorphie d'objets de }Tors\left(X_{L},G_{X_{L}}\right)\right\}
\]

C'est donc le faisceau associ{\'e} au pr{\'e}faisceau:
	\[\left(\textup{Spec}\;L\rightarrow\textup{Spec}\;k\right)\rightsquigarrow H^{1}\left(X_{L},G_{X_{L}}\right)
	\]

Ce qui est exactement la d{\'e}finition du faisceau $R^{1}\pi_{\ast}G_{X}$.
\begin{flushright}
$\Box$
\end{flushright}
\newpage

\thispagestyle{empty}
\end{subsection}
\end{section}
\end{chapter}

\begin{chapter}[Obstruction de Brauer-Manin des gerbes]{Obstruction de Brauer-Manin des gerbes}
\thispagestyle{empty}
Cette section est le fruit d'un travail en commun avec Jean-Claude Douai et Michel Emsalem. A peu de choses pr{\`e}s, le contenu de cette section est celui de notre texte, disponible sur le serveur \textquotedblleft{arxiv}\textquotedblright, {\`a} l'adresse suivante:
\begin{center}
arxiv:math.AG/0303231
\end{center}

Nous montrons ici comment associer {\`a} une gerbe d{\'e}finie sur un corps de nombres une obstruction de Brauer-Manin mesurant, comme dans le cas des vari{\'e}t{\'e}s, le d{\'e}faut d'existence d'une section globale. La motivation est la suivante: soient $k$ un corps de nombres, et $V$ un $k$-espace homog{\`e}ne sous $SL_{n}$ avec isotropie $H$. On peut lui associer la gerbe $\mathcal{G}\left(V\right)=\mathcal{G}$ de ses trivialisations, \textit{i.e.} la gerbe qui mesure l'obstruction {\`a} ce que $V$ soit domin{\'e} par un $SL_{n}$-torseur sur $k$, autrement dit qu'il soit trivial. D'autre part, si $V$ a des points partout localement, on peut lui associer son obstruction de Brauer-Manin $m_{\mathcal{H}}\left(V\right)$. Si cette derni{\`e}re est non-nulle, $V$ n'a pas de point $k$-rationnel.

Or plusieurs espaces homog{\`e}nes peuvent avoir la m{\^e}me gerbe $\mathcal{G}$ des trivialisations. Dans un souci d'{\'e}conomie, on veut d{\'e}finir une obstruction de Brauer-Manin de $\mathcal{G}$, dont la non-nullit{\'e} est une obstruction {\`a} ce que $\mathcal{G}$ soit neutre, et par cons{\'e}quent une obstruction {\`a} ce que chacun des espaces homog{\`e}nes ayant $\mathcal{G}$ pour gerbe des trivialisations ait un point $k$-rationnel.

Une fois de plus, ce raisonnement est inspir{\'e} de la th{\'e}orie des rev{\^e}tements. Plus explicitement, si $k$ est un corps de nombres, $\bar{f}$ un $\bar{k}$-rev{\^e}tement de corps des modules $k$, et $\mathcal{G}\left(\bar{f}\right)$ la $k$-gerbe associ{\'e}e {\`a} $\bar{f}$ (\textit{i.e.} la gerbe des mod{\`e}les de $\bar{f}$,). Dans \cite{DDM}, D{\`e}bes, Douai et Moret-Bailly ont introduit la notion de vari{\'e}t{\'e} de descente associ{\'e}e {\`a} $\bar{f}$. Si $V$ est une telle vari{\'e}t{\'e}, la gerbe $\mathcal{G}\left(\bar{f}\right)$ est alors isomorphe au champ quotient $\left[V/GL_{n}\right]$, pour un $n$ idoine. En fait, il existe une infinit{\'e} possible de telles vari{\'e}t{\'e}s de descente\index{vari{\'e}t{\'e}!de descente} $V$ correspondant {\`a} une infinit{\'e} de choix possibles pour l'entier $n$. Soit maintenant $K$ une extension de $k$: tout $K$-point de $V$ d{\'e}finit un $K$-point de $\mathcal{G}$, et r{\'e}ciproquement tout $K$-point de $\mathcal{G}$ se rel{\`e}ve en un $K$-point de $V$. Il s'ensuit que si $V$ et $V'$ sont deux vari{\'e}t{\'e}s de descente correspondant {\`a} la m{\^e}me $k$-gerbe $\mathcal{G}$, alors:
	\[V\left(K\right)\neq\emptyset\Leftrightarrow V'\left(K\right)\neq\emptyset
\]

Forts de ces observations, on veut comparer les invariants de $V$ et $V'$; ils ne d{\'e}pendent que de $\mathcal{G}$. En particulier:
	\[\textup{Br}_{a}V\approx \textup{Br}_{a}V'\approx \textup{Br}_{a}\mathcal{G},
\]
et
	\[\textup{Pic}\:V\approx \textup{Pic}\:V'\approx \textup{Pic}\:\mathcal{G}.
\]

Ceci nous am{\`e}ne {\`a} calculer l'invariant de Brauer-Manin $m_{\EuScript{H}}\left(V\right)$ de $V$, {\`a} introduire l'invariant de Brauer-Manin $m_{\EuScript{H}}\left(\mathcal{G}\right)$\index{invariant!de Brauer-Manin d'une gerbe} de la gerbe $\mathcal{G}$, puis {\`a} prouver que $m_{\EuScript{H}}\left(V\right)=m_{\EuScript{H}}\left(\mathcal{G}\right)$, l'int{\'e}r{\^e}t de cette {\'e}galit{\'e} {\'e}tant sa validit{\'e} pour toute vari{\'e}t{\'e} de descente $V$ correspondant {\`a} $\mathcal{G}$ (plus loin, nous dirons que $V$ est une \textit{pr{\'e}sentation}\index{pr{\'e}sentation} de $\mathcal{G}$).

Tout ce qui pr{\'e}c{\`e}de s'{\'e}tend aux $k$-gerbes quelconques localement li{\'e}es par un groupe fini (pour des raisons {\'e}videntes, de telles gerbes seront appel{\'e}es des \textit{gerbes de Deligne-Mumford}\index{gerbe!de Deligne-Mumford}). L'application $m_{\EuScript{H}}$ qui {\`a} une classe de $k$-gerbes $\left[\mathcal{G}\right]$ associe l'invariant de Brauer-Manin d'un de ses repr{\'e}sentants peut alors {\^e}tre vue comme une g{\'e}n{\'e}ralisation de la dualit{\'e} de Tate-Poitou\index{dualit{\'e} de Tate-Poitou} dans le cas ab{\'e}lien (nous renvoyons au th{\'e}or{\`e}me 4.4.1 pour un {\'e}nonc{\'e} pr{\'e}cis); cet invariant vit dans le groupe de Tate-Shafarevich:
	\[\mathcyr{SH}^{1}\left(k,\widehat{\bar{H}}\right)^{D}
\]
o{\`u} $\bar{H}$ est le groupe d'automorphismes d'un objet de $\mathcal{G}\left(\textup{Spec}\:\bar{k}\right)$.
\begin{section}{Rappels}
\begin{subsection}[Calcul de $\textup{Br}_{a}$ pour un espace homog{\`e}ne sous $SL_{n}$]{Calcul de $\textup{Br}_{a}V$ dans le cas o{\`u} $V$ est un espace homog{\`e}ne de $SL_{n}$ avec isotropie $H$}
$\ $
\newline

Dans tout ce chapitre, $k$ d{\'e}signe un corps de nombres et $\Omega_{k}$ l'ensemble de ses places. Soit $V$ une $k$-vari{\'e}t{\'e}\footnote{Par $k$-vari{\'e}t{\'e}, on entend ici $k$-sch{\'e}ma s{\'e}par{\'e} de type fini.} alg{\'e}brique lisse, g{\'e}om{\'e}triquement irr{\'e}ductible. De la suite spectrale:
	\[H^{p}\left(k,H^{q}_{\acute{e}t}\left(\bar{V},\mathbb{G}_{m}\right)\right)\Longrightarrow H^{p+q}_{\acute{e}t}\left(V,\mathbb{G}_{m}\right)
\]
on d{\'e}duit la suite exacte longue:
\begin{equation}
	\xymatrix{0 \ar[r] & H^{1}\left(k,\bar{k}\left[V\right]^{\ast}\right) \ar[r] & \textup{Pic}\:V \ar[r] & \textup{Pic}\:\bar{V}^{\textup{Gal}\left(\bar{k}/k\right)} \ar[r] &H^{2}\left(k,\bar{k}\left[V\right]^{\ast}\right) \ar`dr_l[lll]`^dr[lll] [dll]  & \\ &&\textup{Br}_{1}V \ar[r]&H^{1}\left(k,\textup{Pic}\:\bar{V}\right) \ar[r] & H^{3}\left(k,\bar{k}\left[V\right]^{\ast}\right)}
\end{equation}
Posons:
	\[U\left(\bar{V}\right)=\frac{\bar{k}\left[V\right]^{\ast}}{\bar{k}^{\ast}}
\]
La suite exacte $\left(4.1\right)$ fournit une nouvelle suite exacte:
\begin{equation}
	\textup{Pic}\:\bar{V}^{\textup{Gal}\left(\bar{k}/k\right)}\rightarrow H^{2}\left(k,U\left(\bar{V}\right)\right) \rightarrow \textup{Br}_{a}V \rightarrow H^{1}\left(k,\textup{Pic}\:\bar{V}\right) \rightarrow H^{3}\left(k,U\left(\bar{V}\right)\right)
\end{equation}

Supposons que $V$ est un $k$-espace homog{\`e}ne sous un $k$-groupe alg{\'e}brique semi-simple simplement connexe $\widetilde{G}$ (\textit{e.g.} $SL_{n}$) avec isotropie un groupe fini, \textit{i.e}: il existe un $\bar{k}$-groupe fini $\bar{H}$ tel que:
	\[\bar{V}=\widetilde{G}\left(\bar{k}\right)/\bar{H}
\]
Nous avons alors la suite exacte:
	\[0\longrightarrow U\left(\bar{V}\right) \longrightarrow U\left(\widetilde{G}\left(\bar{k}\right)\right)
\]
provenant de la fibration $\widetilde{G}\left(\bar{k}\right)\rightarrow \bar{V}$. Or on sait (d'apr{\`e}s le lemme 6.5 (iii) de \cite{Sa}) que:

\[U\left(\widetilde{G}\left(\bar{k}\right)\right)=\widehat{\widetilde{G}\left(\bar{k}\right)}=0
\]
La suite exacte $\left(4.2\right)$ se r{\'e}duit alors {\`a} l'isomorphisme:
	\[\textup{Br}_{a}V\stackrel{\sim}{\longrightarrow}H^{1}\left(k,\textup{Pic}\:\bar{V}\right)
\]
\begin{remsub} \textup{Notons au passage que cet isomorphisme est valable plus g{\'e}n{\'e}ralement lorsque $V$ est une vari{\'e}t{\'e} alg{\'e}brique propre (\textit{e.g.} projective) d{\'e}finie sur un corps de nombres $k$. Car dans cette situation, d'une part le groupe $H^3\left(k,\mathbb{G}_{m}\right)$ est nul, et d'autre part $\bar{k}\left[V\right]^{\ast}$ se r{\'e}duit {\'e}videmment aux constantes.}

\textup{Par exemple, le groupe $\textup{Br}_{a}V$ (donc \textit{a fortiori} $\mathcyr{B}\left(V\right)$) est nul lorsque $V$ est une $k$-vari{\'e}t{\'e} de Severi-Brauer, ou une $k$-vari{\'e}t{\'e} projective lisse qui est une intersection compl{\`e}te de dimension $\geq3$ (\textit{cf.} la proposition 2.1.13).}
\end{remsub}
\end{subsection}
\begin{subsection}{Exemples}
\begin{enumerate}[(i)]
\item Si $H=0$, alors $\bar{V}\approx\widetilde{G}\left(\bar{k}\right)$, et $\textup{Pic}\:\bar{V}=0$ (car $\textup{Pic}\:\widetilde{G}\left(\bar{k}\right)=0$ par le corollaire 4.5 de \cite{FI}).
\item Si $H=\mu$ est un $k$-sous-groupe central de $\widetilde{G}$, alors $V=G=\widetilde{G}/\mu$ est semi-simple, et $\textup{Pic}\:\widetilde{G}\left(\bar{k}\right)=\widehat{\mu\left(\bar{k}\right)}$ (par le corollaire 4.6 de \cite{FI}), d'o{\`u}:

\[\textup{Br}_{a}V=\textup{Br}_{a}G=H^{1}\left(\bar{k}/k,\widehat{\mu\left(\bar{k}\right)}\right)=H^{1}\left(k,\widehat{\mu}\right)
\]
$\textup{Pic}\:G$ et $\textup{Br}_{a}G$ sont justiciables de la philosophie de Kottwitz: ce sont des invariants des groupes semi-simples qui sont nuls lorsque $G=\widetilde{G}$ est simplement connexe. Ils peuvent donc s'exprimer en fonction du centre $Z\left(^{L}G\right)$ du dual de Langlands de $G$ \cite{Ko}.

Cette remarque vaut encore pour:
	\[\mathcyr{B}\left(G\right)=\ker\left\{\textup{Br}_{a}G\rightarrow{\prod_{v\in\Omega_{k}}\textup{Br}_{a}G_{v}}\right\}
\]
\item Prenons pour $V$ un $k$-tore $T$. Alors (\textit{cf.} le lemme 6.9 de \cite{Sa}):
	\[\textup{Pic}\:\bar{T}=H^{1}\left(k,\widehat{T}\right)\ \ et\ \ \textup{Br}_{a}T=H^{2}\left(k,\widehat{T}\right)
\]
En outre:
	\[\mathcyr{B}\left(T\right)\approx\mathcyr{SH}^{2}\left(k,\widehat{T}\right)\approx\mathcyr{SH}^{1}\left(k,\widehat{T}\right)^{D}
\]
le deuxi{\`e}me isomorphisme {\'e}tant fourni par la dualit{\'e} de Kottwitz \cite{Ko} qui {\'e}tend aux tores celle de Tate-Poitou.
\item Consid{\'e}rons maintenant un $k$-espace homog{\`e}ne de $SL_{n}$ avec isotropie un groupe fini; on suppose donc qu'il existe un groupe fini $\bar{H}$ tel que:
	\[\bar{V}\approx SL_{n,\bar{k}}/\bar{H}
\]
Alors $\textup{Pic}\:\bar{V}\approx\widehat{\bar{H}}$ (\textit{cf.} \cite{BK2}). On dispose en effet de la $\bar{k}$-fibration:
	\[SL_{n,\bar{k}}\longrightarrow SL_{n,\bar{k}}/\bar{H}
\]
{\`a} laquelle est attach{\'e}e la suite spectrale de Hochschild-Serre:
	\[E^{p,q}_{2}=H^{p}\left(\bar{H},H^{q}\left(SL_{n},\mathbb{G}_{m}\right)\right)\Longrightarrow H^{p+q}_{\acute{e}t}\left(\bar{V},\mathbb{G}_{m}\right)=E^{p+q}
\]
Dans cette derni{\`e}re, le terme $E_{2}^{0,1}$ est nul\footnote{En effet,
$\textup{Pic}\:SL_{n}=H^{1}\left(SL_{n,\bar{k}},\mathbb{G}_{m}\right)=\textup{Hom}\left(\Pi_{1}\left(SL_{n}\right),\mathbb{G}_{m}\right)=0$ puisque $SL_{n}$ est simplement connexe.}, donc:

\[H^{1}_{\acute{e}t}\left(\bar{V},\mathbb{G}_{m}\right)=H^{1}\left(\bar{H},\mathbb{G}_{m}\right)=\textup{Hom}\left(\bar{H},\mathbb{G}_{m}\right)=\widehat{\bar{H}}
\]
D'o{\`u} la:
\end{enumerate}
\begin{prosub} Soit $V$ un $k$-espace homog{\`e}ne d'un groupe semi-simple simplement connexe $\widetilde{G}$ avec isotropie un groupe fini $\bar{H}$. Alors d'apr{\`e}s \cite{BK2}:
	\[\textup{Br}_{a}V\approx H^{1}\left(k,\widehat{\bar{H}}\right)
\]

En outre, si $k$ est un corps de nombres, et si on suppose que $V$ a des points partout localement (\textit{i.e.} si $V_{v}=V\otimes_{k}k_{v}$ a un $k_{v}$-point, pour toute place $v$ de $k$), alors:
	\[\mathcyr{B}\left(V\right)=\ker\left\{\textup{Br}_{a}V\rightarrow{\prod_{v\in\Omega_{k}}\textup{Br}_{a}V_{v}}\right\}\approx \mathcyr{SH}^{1}\left(k,\widehat{\bar{H}}\right)
\]
\end{prosub}
\begin{corosub} Sous les hypoth{\`e}ses et notations de la proposition pr{\'e}c{\'e}dente, si $\bar{H}$ est sans caract{\`e}re, alors $\textup{Br}_{a}V=\mathcyr{B}\left(V\right)=0$.
\end{corosub}
\begin{exemsub} \textup{Le corollaire s'applique donc si $\bar{H}=SL\left(2,\mathbb{F}_{p}\right)$ avec $p\neq2,3$, ou encore si $\bar{H}=\EuScript{A}_{n}$ avec $n\geq5$. Plus g{\'e}n{\'e}ralement, il faut et il suffit que $\bar{H}$ soit {\'e}gal {\`a} son groupe d{\'e}riv{\'e}.}
\end{exemsub}
\end{subsection}
\end{section}
\begin{section}[Interpr{\'e}tation en termes de champ quotient]{Interpr{\'e}tation comme champ quotient des $k$-gerbes localement li{\'e}es par un groupe alg{\'e}brique fini}
On s'int{\'e}resse donc dans cette section aux $k$-gerbes qui sont des champs de Deligne-Mumford \cite{LMB}. Rappelons d'abord la proposition 5.1 de \cite{DDM}:
\begin{pro} Soient $k$ un corps, et $\mathcal{G}$ une $k$-gerbe (pour la topologie {\'e}tale) qui est un champ de Deligne-Mumford. Alors:
\begin{enumerate}
\item Il existe une $k$-alg{\`e}bre $L$ avec action {\`a} gauche d'un groupe fini $\Gamma$ admettant $k$ comme anneau des invariants telle que $\mathcal{G}$ soit isomorphe au champ quotient $\left[\textup{Spec}\:L/\Gamma\right]$;
\item Il existe un $k$-sch{\'e}ma affine $V$, un entier $n\geq0$, une action {\`a} droite de $GL_{n,k}$ sur $V$ et un $1$-morphisme $\pi:V\rightarrow\mathcal{G}$ avec les propri{\'e}t{\'e}s suivantes:
\begin{enumerate}[(i)]
\item $\pi$ induit un isomorphisme du champ quotient $\left[V\right/GL_{n,k}]$ vers $\mathcal{G}$;
\item $V$ est lisse et g{\'e}om{\'e}triquement irr{\'e}ductible;
\item l'action de $GL_{n,k}$ sur $V$ est transitive et {\`a} stabilisateurs finis;
\item pour chaque extension $K$ de $k$, chaque objet de $\mathcal{G}\left(K\right)$ se rel{\`e}ve en un point de $V\left(K\right)$ via $\pi$.
\end{enumerate}
En particulier, {\`a} cause de (iii) et (iv), si $K$ est une extension de $k$ telle que $\mathcal{G}\left(K\right)\neq\emptyset$,\footnote{Ce qui signifie que l'ensemble d'objets $Ob\left(\mathcal{G}\left(K\right)\right)$ de la cat{\'e}gorie fibre de $\mathcal{G}$ au-dessus de $\textup{Spec}\:K$ est non-vide; par la suite, nous fairons syst{\'e}matiquement cet abus de langage.} la $K$-vari{\'e}t{\'e} $V\otimes_{k}K$ est isomorphe au quotient de $GL_{n,K}$ par un groupe fini.
\end{enumerate}
\end{pro}
\begin{rems}$\ $

\begin{enumerate}[(a)]
\item \textup{Dans la remarque 5.2(b) de \cite{DDM}, il est montr{\'e} que l'on peut en fait prendre pour $k$-alg{\`e}bre $L$ une extension galoisienne finie de $k$, auquel cas $\Gamma$ est l'ensemble des couples $\left(\sigma,\varphi\right)$ o{\`u} $\sigma\in\:\textup{Gal}\left(L/k\right)$ et $\varphi:\sigma x\stackrel{\sim}{\rightarrow}x$ est un isomorphisme dans la cat{\'e}gorie (en fait, le groupoïde) $\mathcal{G}\left(L\right)$. Il y a une structure de groupe sur $\Gamma$ pour laquelle la projection naturelle $\Gamma\rightarrow \textup{Gal}\left(L/k\right)$ est un morphisme surjectif, et le noyau $H=H\left(L\right)$ est le stabilisateur fini dont l'existence est donn{\'e}e par le (iii) de la proposition 4.2.1.}

\textup{A la gerbe $\mathcal{G}$ est aussi associ{\'e}e une extension (cf. \cite{Gi2}, \cite{Sp}) 
	\[\mathcal{G}\leftrightsquigarrow\left(\EuScript{E}\right):1\rightarrow{H}\rightarrow\Gamma\rightarrow{\textup{Gal}\left(L/k\right)}\rightarrow{1}
\]
d{\'e}finissant une action ext{\'e}rieure $\mathcal{L}_{H}$ de $\textup{Gal}\left(L/k\right)$ sur $H=H\left(L\right)$, et une classe de 2-cohomologie not{\'e}e $\left[\mathcal{G}\right]$ dans $H^{2}\left(L/k,\mathcal{L}_{H}\right)\hookrightarrow H^{2}\left(k,\mathcal{L}_{H}\right)$.}
\item \textup{On obtient les m{\^e}mes conclusions en rempla\c{c}ant dans la proposition 4.2.1 $GL_{n}$ par $SL_{n}$ (\textit{cf.} remarque 5.2(c) de \cite{DDM}).}
\end{enumerate}
\end{rems}
$\ $
\vspace{1mm}
\newline
\begin{center}
\uline{Une construction fondamentale}
\end{center}
$\ $
\newline
Partons de l'extension $\left(\EuScript{E}\right)$ de la remarque pr{\'e}c{\'e}dente:
	\[\left(\EuScript{E}\right):1\rightarrow{H}\rightarrow\Gamma\rightarrow{\textup{Gal}\left(L/k\right)}\rightarrow{1}
\]
$\Gamma$ est un groupe fini; on peut donc le plonger dans $SL_{n}$ pour un certain $n$, ce qui conduit au diagramme suivant:
	\[\left(D\right):\xymatrix{1\ar[r] & H \ar[r] \ar@{=}[d] & \Gamma \ar[r] \ar[d] & \textup{Gal}\left(L/k\right) \ar[r] \ar@{-->}[d] & 1\\ 1 \ar[r] & H \ar[r] & SL_{n,\bar{k}} \ar@{-->}[r] & SL_{n,\bar{k}}/H \ar[r] & 1}
\]
$SL_{n,\bar{k}}/H$ n'est pas un groupe, puisque $H$ n'est pas n{\'e}cessairement normal dans $SL_{n,\bar{k}}$. C'est seulement un $k$-espace homog{\`e}ne (toujours au sens de Springer \cite{Sp}), d'o{\`u} la pr{\'e}sence des pointill{\'e}s dans le diagrammme pr{\'e}c{\'e}dent. La fl{\`e}che verticale
	\[\xymatrix{\textup{Gal}\left(L/k\right) \ar@{-->}[d]\\ SL_{n,\bar{k}}/H }
\]
donne lieu {\`a} un $1$-cocycle dans $Z^{1}\left(L/k;SL_{n},H\right)$, qui repr{\'e}sente pr{\'e}cis{\'e}ment la classe du $k$-espace homog{\`e}ne $V$ du (2) de la proposition 4.2.1. La $k$-gerbe $\mathcal{G}\approx\left[V/SL_{n}\right]$ (associ{\'e}e {\`a} $\left(\EuScript{E}\right)$) s'interpr{\`e}te alors comme la gerbe des rel{\`e}vements du $k$-espace homog{\`e}ne $V$ {\`a} $SL_{n}$. En d'autres termes, $\left[\mathcal{G}\right]$ est l'image de $V$ par le cobord (\textit{cf.} \cite{Sp}, \cite{Do1})
	\[\delta^{1}:Z^{1}\left(L/k;SL_{n},H\right)\longrightarrow H^{2}\left(k,\mathcal{L}_{H}\right)
\]
Dans la suite, nous appellerons \textbf{pr{\'e}sentation de $\mathcal{G}=\left[V/SL_{n}\right]$}\index{pr{\'e}sentation} un couple $\left(V,\pi\right)$ comme dans la proposition 4.2.1.
\begin{rem} \textup{La proposition 4.2.1 traduit en particulier le fait que les deux groupes de Brauer d'un sch{\'e}ma (cohomologique et \textquotedblleft{Azumaya}\textquotedblright) coïncident lorsque ce sch{\'e}ma est un corps. En effet, si $\mathcal{G}\in H^{2}\left(k,\mathbb{G}_{m}\right)$, alors il existe un $n$ tel que $\mathcal{G}\in H^{2}\left(k,\mu_{n}\right)$, puisque le groupe de Brauer d'un corps, et plus g{\'e}n{\'e}ralement d'un sch{\'e}ma r{\'e}gulier \cite{G4}, est de torsion. Par cons{\'e}quent, $\mathcal{G}$ est un champ de Deligne-Mumford, et il existe un entier $n'$ et un $k$-espace homog{\`e}ne $V$ de $SL_{n'}$ avec isotropie $\mu_{n}$ tels que }
	\[\mathcal{G}\approx\left[V/SL_{n'}\right]
\]

\textup{Remarquons que l'entier $n'$ peut {\^e}tre diff{\'e}rent de $n$, et c'est en particulier le cas pour le corps $k_{M}$ construit par Merkurjev dans son article sur la conjecture de Kaplansky \cite{Me}: il existe un {\'e}l{\'e}ment de $H^{2}\left(k_{M},\mu_{2}\right)$ qui ne provient pas d'un {\'e}l{\'e}ment de $H^{1}\left(k_{M},PGL_{2}\right)=H^{1}\left(k_{M};SL_{2},\mu_{2}\right)$ (mais qui est atteint par un {\'e}l{\'e}ment de $H^{1}\left(k_{M},PGL_{4}\right)$).}

\textup{Cependant, l'espace homog{\`e}ne $V$ n'est autre que la vari{\'e}t{\'e} de Severi-Brauer pr{\'e}-image de $\mathcal{G}$ par le morphisme naturel
	\[\textup{Br}_{Az}k\longrightarrow H^{2}\left(k,\mathbb{G}_{m}\right)
\]
et on retrouve ainsi le point de vue de \cite{EHKV}.}
\end{rem}
\end{section}
\begin{section}[Invariant de Brauer-Manin d'une $k$-gerbe]{Invariant de Brauer-Manin d'une $k$-gerbe localement li{\'e}e par un groupe fini}
Les $k$-champs alg{\'e}briques (en particulier les $k$-gerbes qui sont des champs de Deligne-Mumford) sont des g{\'e}n{\'e}ralisations de la notion de sch{\'e}ma. Par suite, il est tout-{\`a}-fait naturel de d{\'e}finir l'obstruction de Brauer-Manin d'une $k$-gerbe de mani{\`e}re analogue {\`a} celle d'un $k$-sch{\'e}ma.

Soit donc $\mathcal{G}$ une $k$-gerbe, qui est un champ de Deligne-Mumford. Le site {\'e}tale de $\mathcal{G}$
est d{\'e}fini au chapitre 12 de \cite{LMB}. Le groupe de Brauer \index{groupe!de Brauer!d'une gerbe}
$\textup{Br}\:\mathcal{G}=H^{2}_{\acute{e}t}\left(\mathcal{G},\mathbb{G}_{m}\right)$\label{Brger} est d{\'e}fini dans \cite{SGA4-5}
(o{\`u}, d'une mani{\`e}re plus g{\'e}n{\'e}rale, la cohomologie d'un topos localement annel{\'e} est
d{\'e}finie). Il existe une suite spectrale (cf. \cite{SGA4-5}, prop. 5.3):
	\[E_{2}^{p,q}=H^{p}\left(k,H^{q}_{\acute{e}t}\left(\mathcal{G},\mathbb{G}_{m}\right)\right)\Longrightarrow H^{p+q}_{\acute{e}t}\left(\mathcal{G},\mathbb{G}_{m}\right)=E^{p+q}
\]
qui permet de d{\'e}finir $\textup{Br}^{cst}\mathcal{G}$, $\textup{Br}^{alg}\mathcal{G}$ et $\textup{Br}_{a}\mathcal{G}$ comme pour les vari{\'e}t{\'e}s. Plus pr{\'e}cis{\'e}ment:
	\[\left\{\begin{array}{l}\textup{Br}^{cst}\mathcal{G}=\textup{im}\;\left\{E_{2}^{2,0}=\textup{Br}\;k\longrightarrow \textup{Br}\;\mathcal{G}=E^{2}\right\}\\\textup{Br}^{alg}\mathcal{G}=\ker\;\left\{E^{2}=\textup{Br}\mathcal{G}\longrightarrow\textup{Br}\overline{\mathcal{G}}=E^{0,2}_{2}\right\}\\\textup{Br}_{a}\mathcal{G}=\displaystyle\frac{\textup{Br}^{alg}\mathcal{G}}{\textup{Br}^{cst}\mathcal{G}}\end{array}\right.
\]

Soit $\left(V,\pi\right)$ une pr{\'e}sentation de $\mathcal{G}$ (d'apr{\`e}s nos hypoth{\`e}ses, on a donc $\mathcal{G}\approx\left[V/SL_{n}\right]$). Posons:
	\[U\left(\bar{\mathcal{G}}\right)=\frac{\bar{k}\left[\mathcal{G}\right]^{\ast}}{\bar{k}^{\ast}}
\]
On a:
	\[U\left(\bar{\mathcal{G}}\right)\subset U\left(\bar{V}\right) \subset U\left(SL_{n,\bar{k}}\right)=\widehat{SL_{n,\bar{k}}}=0
\]
L'analogue de la suite exacte $\left(2\right)$ associ{\'e}e {\`a} la suite spectrale pr{\'e}c{\'e}dente implique alors:
	\[\begin{array}{rl}\textup{Br}_{a}\mathcal{G}&\approx H^{1}\left(k,\textup{Pic}\:\bar{\mathcal{G}}\right)\\ &\approx H^{1}\left(k,\textup{Hom}\left(\Pi_{1}\left(\bar{\mathcal{G}}\right),\mathbb{G}_{m}\right)\right)\\&\approx H^{1}\left(k,\textup{Hom}\left(\bar{H},\mathbb{G}_{m}\right)\right)\\\end{array}
\]
car il est bien connu que $\Pi_{1}\left(\bar{\mathcal{G}}\right)=\bar{H}$ (\textit{cf.} \cite{No}), et finalement:
\begin{equation}
	\textup{Br}_{a}\mathcal{G}\approx H^{1}\left(k,\widehat{\bar{H}}\right)
\end{equation}

Si $K$ est un corps qui est une extension quelconque de $k$, nous pouvons d{\'e}finir un accouplement:
	\[\begin{array}{ccl}\mathcal{G}\left(K\right)\times \textup{Br}\:\mathcal{G}& \rightarrow & \textup{Br}\:K\\ \left(x,b\right) & \longmapsto & b\left(x\right)\end{array}
\]
o{\`u} l'image $b\left(x\right)$ peut {\^e}tre interpr{\'e}t{\'e}e de diff{\'e}rentes fa\c{c}ons ($\mathcal{B}$ d{\'e}signe ci-dessous un repr{\'e}sentant de $b$):
\begin{enumerate}[(i)]
\item $b\left(x\right)$ est la gerbe r{\'e}siduelle de $\mathcal{B}$ au point $x$ du champ alg{\'e}brique $\mathcal{G}=\left[V/SL_{n}\right]$;
\item comme dans la section 2.3, $x$ peut {\^e}tre vu comme une section au dessus de $\textup{Spec}\:K$ du (1-)morphisme structural $\mathcal{G}\rightarrow{\textup{Spec}\:k}$; autrement dit, c'est un (1-)morphisme rendant commutatif le diagramme (de morphismes de champs) suivant:
	\[\xymatrix{&\mathcal{G} \ar[d]\\\textup{Spec}\:K \ar[r] \ar[ur]^{x}&\textup{Spec}\:k}
\]
$\mathcal{B}$ {\'e}tant une gerbe sur $\mathcal{G}$, on peut consid{\'e}rer la gerbe image inverse $x^{\ast}\mathcal{B}$ de $\mathcal{B}$ par le morphisme $x$; la gerbe $x^{\ast}\mathcal{B}$ ainsi obtenue correspond exactement {\`a} $b\left(x\right)$; elle est obtenue par pull-back {\`a} partir de $x$ et de $\mathcal{B}$:
	\[\xymatrix{&\mathcal{B} \ar@{-}[d]\\x^{\ast}\mathcal{B} \ar@{-}[d] \ar@{--}[ur]^{x^{\ast}}&\mathcal{G} \ar[d]\\\textup{Spec}\:K \ar[r] \ar[ur]^{x}&\textup{Spec}\:k}
\]
\end{enumerate}

De plus, si $H$ ab{\'e}lien, toute $K$-section de la gerbe $\mathcal{G}$ (\textit{i.e.} tout objet de $\mathcal{G}\left(K\right)$) est un $H$-torseur sur $\textup{Spec}\:K$ ($\mathcal{G}$ est par d{\'e}finition localement {\'e}quivalente {\`a} la gerbe $Tors\:H$; l'existence d'une $K$-section implique que $\mathcal{G}_{\left|K\right.}$ est {\'e}quivalente {\`a} la gerbe des $H$-torseurs sur $K$). Si on suppose que $K$ est un corps local, on obtient alors l'{\'e}nonc{\'e} suivant:
\begin{pro}[Cas local] Soient $K$ un corps local, $\mathcal{G}$ une $K$-gerbe li{\'e}e par un groupe ab{\'e}lien fini $H$ et $\left(V,\pi\right)$ une pr{\'e}sentation de $\mathcal{G}$. Le diagramme suivant est commutatif:
	\[\xymatrix@C=2pt@R=25pt{\left(Acc.1\right)&V\left(K\right)\ar@{->>}[d]&\times&\textup{Br}_{a}V\ar[dr] \\\left(Acc.2\right)&\mathcal{G}\left(K\right)\ar@{->>}[d]&\times&\textup{Br}_{a}\mathcal{G}\ar @<1ex> [u]^{\approx}\ar[r]&\mathbb{Q}/\mathbb{Z} \\\left(Acc.3\right)& H^{1}\left(K,\bar{H}\right)&\times&H^{1}\left(K,\widehat{\bar{H}}\right) \ar @<1ex> [u]^{\approx} \ar[ur]}
\]
o{\`u}:
\begin{itemize}
\item l'accouplement $\left(Acc.1\right):V\left(K\right)\times \textup{Br}_{a}V\longrightarrow\mathbb{Q}/\mathbb{Z}$ est d{\'e}fini comme suit: {\`a} un point $x$ de $V\left(K\right)$, et {\`a} une classe $b$ dans $\textup{Br}_{a}V$, on associe:
	\[\left\langle x,b\right\rangle=\left[s_{x}\left(b\right)\right]_{x}
\]
$s_{x}$ d{\'e}signant la section\footnote{En effet, l'existence d'un point $K$-rationnel entra{\^i}ne que la suite:
	\[0\longrightarrow \textup{Br}\;K\longrightarrow \textup{Br}_{1}V\longrightarrow \textup{Br}_{a}V\longrightarrow 0
\]
est scind{\'e}e.} induite par $x$ de la projection canonique $p$ (\textit{cf.} \cite{BK1}):
	\[\xymatrix{\textup{Br}_{1}X \ar[r]_{p}&\textup{Br}_{a}X \ar@/_12pt/[l]_{s_{x}}}
\]
\item l'accouplement $\left(Acc.2\right):\mathcal{G}\left(K\right)\times \textup{Br}_{a}\mathcal{G}\longrightarrow\mathbb{Q}/\mathbb{Z}$ est d{\'e}fini de la m{\^e}me mani{\`e}re que $\left(Acc.1\right)$;
\item l'accouplement $\left(Acc.3\right):H^{1}\left(K,\bar{H}\right)\times H^{1}\left(K,\widehat{\bar{H}}\right)\longrightarrow\mathbb{Q}/\mathbb{Z}$ est l'accouplement de Tate pour les corps locaux.
\end{itemize}
\end{pro}
$\ $
\newline

$k$ {\'e}tant un corps de nombres, nous pouvons d{\'e}finir pour toute place $v$ de $k$ l'accouplement:
	\[\begin{array}{ccl} \mathcal{G}\left(k_{v}\right)\times \textup{Br}_{1}\:\mathcal{G}& \rightarrow & \mathbb{Q}/\mathbb{Z}\\ \left(x,b\right) & \longmapsto & \textup{inv}_{v}\left(b\left(x\right)\right)\end{array}
\]
o{\`u} comme d'habitude $\textup{inv}_{v}$ est l'invariant donn{\'e} par la th{\'e}orie du corps de classes, et $b\left(x\right)$ est la classe dans $\textup{Br}\:k_{v}$ de $\mathcal{B}_{x}$, o{\`u} $\mathcal{B}$ est un repr{\'e}sentant de $b$. Supposons maintenant que $\mathcal{G}\left(k_{v}\right)$ soit non-vide pour toute place $v$ de $k$, et restreignons nous au sous-groupe $\mathcyr{B}\left(\mathcal{G}\right)$ de $\textup{Br}_{a}\mathcal{G}$ d{\'e}fini par:
	\[\mathcyr{B}\left(\mathcal{G}\right)=\ker\left\{\textup{Br}_{a}\mathcal{G}\longrightarrow \prod_{v\in\Omega_{k}}{\textup{Br}_{a}\:\mathcal{G}_{\left|k_{v}\right.}}\right\}
\]
Nous d{\'e}finissons ainsi un accouplement:
	\[\begin{array}{rccl}\left\langle.\:,\:.\right\rangle: & $\(\displaystyle\prod_{v\in\Omega_{k}}{\mathcal{G}\left(k_{v}\right)}\)$ \times \mathcyr{B}\left(\mathcal{G}\right)& \longrightarrow & \mathbb{Q}/\mathbb{Z}\\ & \left(\left(x_{v}\right)_{v},b\right) & \longmapsto & $\(\displaystyle\sum_{v\in\Omega_{k}}{\textup{inv}_{v}\left(\widetilde{b}\left(x_{v}\right)\right)}\)$\end{array}
\]
o{\`u} $\widetilde{b}$ est un relev{\'e} de $b$ dans $\textup{Br}_{1}\mathcal{G}$. Par analogie avec la d{\'e}finition usuelle de cet accouplement (\textit{cf.} section 2.3), $\left\langle
\left(x_{v}\right)_{v},b\right\rangle$ ne d{\'e}pend pas de $\left(x_{v}\right)_{v}$, et $\left\langle
\left(x_{v}\right)_{v},b\right\rangle\neq0$ est une obstruction {\`a} l'existence d'une section $k$-rationnelle
de $\textup{Spec}\:k$ (\textit{i.e.} d'un objet de la cat{\'e}gorie fibre $\mathcal{G}\left(k\right)$). Nous obtenons de cette fa\c{c}on un {\'e}l{\'e}ment bien d{\'e}fini:
	\[m_{\mathcal{H}}\left(\mathcal{G}\right)\in \mathcyr{B}\left(\mathcal{G}\right)^{D}=\textup{Hom}\left(\mathcyr{B}\left(\mathcal{G}\right),\mathbb{Q}/\mathbb{Z}\right)=\textup{Hom}\left(\mathcyr{SH}^{1}\left(k,\widehat{\bar{H}}\right),\mathbb{Q}/\mathbb{Z}\right)
\]
puisque: $\mathcyr{B}\left(\mathcal{G}\right)\approx\mathcyr{SH}^{1}\left(k,\widehat{\bar{H}}\right)$ (c'est une cons{\'e}quence imm{\'e}diate de l'isomorphisme $\left(4.3\right)$).
\begin{pro}[Cas global] Soient $k$ un corps de nombres et $\mathcal{G}$ une $k$-gerbe. Pour toute pr{\'e}sentation $\left(V,\pi\right)$ de $\mathcal{G}$:
	\[\mathcyr{B}\left(V\right)\stackrel{\sim}{\longleftarrow}\mathcyr{B}\left(\mathcal{G}\right)
\]
et $m_{\mathcal{H}}\left(\mathcal{G}\right)$ est {\'e}gale {\`a} l'image de $m_{\mathcal{H}}\left(V\right)$ par l'isomorphisme:
	\[\mathcyr{B}\left(V\right)^{D}\stackrel{\sim}{\longrightarrow}\mathcyr{B}\left(\mathcal{G}\right)^{D}
\]
\end{pro}
$\ $
\vspace{1mm}
\newline

L'isomorphisme $\mathcyr{B}\left(V\right)\stackrel{\sim}{\longleftarrow}\mathcyr{B}\left(\mathcal{G}\right)$ r{\'e}sulte de l'isomorphisme 
	\[\mathcyr{SH}^{1}\left(k,\textup{Pic}\:\bar{V}\right)\stackrel{\sim}{\longrightarrow}\mathcyr{SH}^{1}\left(k,\textup{Pic}\:\bar{\mathcal{G}}\right)
\]
ce dernier {\'e}tant induit par les isomorphismes compos{\'e}s
	\[H^{1}\left(k,\textup{Pic}\:\bar{V}\right)\approx H^{1}\left(k,\widehat{\bar{H}}\right)\approx H^{1}\left(k,\textup{Pic}\:\bar{\mathcal{G}}\right)
\]
Nous avons le diagramme commutatif:
	\[\xymatrix@C=2pt@R=10pt{\prod{V\left(k_{v}\right)}\ar @<-1ex> @{->>}[dd]\ \  \times & \mathcyr{B}\left(V\right) \ar[dr] \\&&& \mathbb{Q}/\mathbb{Z}\\\prod{\mathcal{G}\left(k_{v}\right)}\ \  \times & \mathcyr{B}\left(\mathcal{G}\right) \ar@<1ex>[uu]^{\approx} \ar[ur]}
\]
dans lequel la surjectivit{\'e} de la fl{\`e}che de gauche provient de la proposition 4.2.1(iv). On a vu que le
calcul de $m_{\mathcal{H}}\left(\mathcal{G}\right)$ ne d{\'e}pendait pas de la famille $\left(x_{v}\right)_{v}$
choisie dans \(\displaystyle\prod_{v\in\Omega_{k}}{\mathcal{G}\left(k_{v}\right)}\). De la m{\^e}me mani{\`e}re, on
sait que le calcul de $m_{\mathcal{H}}\left(V\right)$ ne d{\'e}pend pas non plus de la famille
$\left(y_{v}\right)_{v}$ choisie dans \(\displaystyle\prod_{v\in\Omega_{k}}{V\left(k_{v}\right)}\). Pour calculer
$m_{\mathcal{H}}\left(V\right)$, on peut donc prendre pour $\left(y_{v}\right)_{v}$ n'importe quel rel{\`e}vement de la famille $\left(x_{v}\right)_{v}$. On en d{\'e}duit que $m_{\mathcal{H}}\left(\mathcal{G}\right)$ n'est autre que l'application compos{\'e}e:
	\[\mathcyr{B}\left(\mathcal{G}\right)\stackrel{\sim}{\longrightarrow}\mathcyr{B}\left(V\right)\xrightarrow{m_{\mathcal{H}}\left(V\right)}\mathbb{Q}/\mathbb{Z}
\]
o{\`u} $m_{\mathcal{H}}\left(V\right)\in \mathcyr{B}\left(V\right)^{D}$ est donn{\'e}e par:
	\[b\longmapsto \left\langle \left(y_{v}\right)_{v},b\right\rangle
\]

Dans la suite, nous verrons l'{\'e}l{\'e}ment $m_{\mathcal{H}}\left(V\right)$ comme un {\'e}l{\'e}ment de
$\mathcyr{SH}^{1}\left(k,\widehat{\bar{H}}\right)^{D}$.
\end{section}
\begin{section}[$1/2$-th{\'e}or{\`e}me de Tate-Poitou non-ab{\'e}lien]{$1/2$-th{\'e}or{\`e}me de Tate-Poitou pour les groupes non-ab{\'e}liens}
\begin{theo} Soient $k$ un corps de nombres, $H$ un $k$-groupe fini, $\mathcal{L}_{H}$ un $k$-lien localement repr{\'e}sentable par $H$. L'application
	\[\begin{array}{rccc}m_{\mathcal{H}}: & \mathcyr{SH}^{2}\left(k,\mathcal{L}_{H}\right) & \longrightarrow & \mathcyr{SH}^{1}\left(k,\widehat{\bar{H}}\right)^{D}\\ & \left[\mathcal{G}\right] & \longmapsto & m_{\mathcal{H}}\left(\mathcal{G}\right)\end{array}
\]
o{\`u} $\mathcyr{SH}^{2}\left(k,\mathcal{L}_{H}\right)$ d{\'e}signe l'ensemble des classes d'{\'e}quivalence de gerbes localement li{\'e}es par $H$ admettant partout localement une section (\textit{i.e.} qui sont partout localement neutres), se factorise par
	\[\xymatrix{\mathcyr{SH}^{2}\left(k,\mathcal{L}_{H}\right) \ar[r]^{ab} \ar[dr]_{m_{\mathcal{H}}} & \mathcyr{SH}^{2}\left(k,\frac{\bar{H}}{\left[\bar{H},\bar{H}\right]}\right) \ar[d]_{\approx}\\ & \mathcyr{SH}^{1}\left(k,\widehat{\frac{\bar{H}}{\left[\bar{H},\bar{H}\right]}}\right)^{D}}
\]
o{\`u} $ab$ est l'application d'ab{\'e}lianisation naturelle, et l'isomorphisme vertical est fourni par la dualit{\'e} de Tate-Poitou.
\end{theo}
\pagebreak

\begin{rem}
\item \textup{On peut {\'e}tendre le th{\'e}or{\`e}me 4.4.1 au cas o{\`u} $H$ est un $k$-groupe lin{\'e}aire,
\textit{i.e.} au cas o{\`u} les $k$-gerbes consid{\'e}r{\'e}es ne sont plus de Deligne-Mumford. Ceci peut de
faire en rempla\c{c}ant dans la construction fondamentale (de la section 2) $SL_{n}$ par $GL_{n}$. Le th{\'e}or{\`e}me 4.1 peut alors {\^e}tre compl{\'e}t{\'e} par les deux r{\'e}sultats suivants:}
\begin{enumerate}
\item \textup{Si $H$ est un $k$-tore $T$, $m_{\mathcal{H}}$ prend ses valeurs dans $\mathcyr{SH}^{1}\left(k,X^{\ast}\left(T\right)\right)^{D}$:
	\[m_{\mathcal{H}}:\mathcyr{SH}^{2}\left(k,T\right)\longrightarrow\mathcyr{SH}^{1}\left(k,X^{\ast}\left(T\right)\right)^{D}
\]
et coïncide avec l'isomorphisme donn{\'e} par la dualit{\'e} de Kottwitz pour les tores \cite{Ko}.}
\item \textup{Si $H$ est un $k$-groupe semi-simple, alors $\mathcyr{SH}^{1}\left(k,\widehat{\bar{H}}\right)=0$ et 
	\[m_{\mathcal{H}}:\mathcyr{SH}^{2}\left(k,\mathcal{L}_{H}\right)\longrightarrow\mathcyr{SH}^{1}\left(k,\widehat{\bar{H}}\right)^{D}=0
\]
est l'application nulle. On sait d'apr{\`e}s \cite{Bo2} dans le cas semi-simple (resp. d'apr{\`e}s \cite{Do1} dans le cas semi-simple simplement connexe) que toutes les classes de
	\[\mathcyr{SH}^{2}\left(k,\mathcal{L}_{H}\right)\ \left(resp.\textup{ de }H^{2}\left(k,\mathcal{L}_{H}\right)\right)
\]
sont neutres. Ainsi l'obstruction de Brauer-Manin est la seule dans le cas semi-simple. Compte tenu de la remarque $\left(1\right)$ pr{\'e}c{\'e}dente, on en d{\'e}duit que le m{\^e}me r{\'e}sultat vaut dans le cas des groupes r{\'e}ductifs connexes, puis dans le cas des groupes connexes (\cite{Bo1}).}
\end{enumerate}
\end{rem}
\newpage
\thispagestyle{empty}
\end{section}
\end{chapter}
\newpage
\thispagestyle{empty}
$\ $

\newpage
\addcontentsline{toc}{chapter}{Index}
\markboth{INDEX}{INDEX}
\thispagestyle{empty}
\printindex
\newpage
\thispagestyle{empty}
\begin{flushleft}
\thispagestyle{empty}
\end{flushleft}
\addcontentsline{toc}{chapter}{Glossaire des notations}
\markboth{GLOSSAIRE DES NOTATIONS}{GLOSSAIRE DES NOTATIONS}
\huge 
\begin{flushleft}
\textbf{Glossaire des notations}
\end{flushleft}
\begin{flushleft}
\thispagestyle{empty}
\end{flushleft}
\normalsize
$\ $\newline\vspace{1mm}
$\textup{Br}_{Az}X$,  groupe de Brauer-Azumaya de $X$, \pageref{Braz}\newline\vspace{2mm}
$\textup{Br}\;X$,  groupe de Brauer cohomologique de $X$, \pageref{BrX}\newline\vspace{2mm}
$S_{\acute{e}t}$,  site {\'e}tale de $S$, \pageref{Setsite}\newline\vspace{2mm}
$\mathbb{G}_{a,S}$,  groupe additif de $S$, \pageref{Gadd}\newline\vspace{2mm}
$\mathbb{G}_{m,S}$,  groupe multiplicatif de $X$, \pageref{Gmult}\newline\vspace{2mm}
$\mu_{n,S}$,  faisceau des racines $\textup{n}^{\grave{e}mes}$ de l'unit{\'e} sur $S$, \pageref{mun}\newline\vspace{2mm}
$\widetilde{S}_{\acute{e}t}$,  topos {\'e}tale de $S$, \pageref{Settopos}\newline\vspace{2mm}
$\textup{ad}_{G_{S}}\left(P\right)$,  faisceau des automorphismes du $G_{S}$-torseur $P$, \pageref{adG}\newline\vspace{2mm}
$Tors\left(S,G_{S}\right)$,  cat{\'e}gorie des $G_{S}$-torseurs sur le site {\'e}tale de $S$, \pageref{Torscat}\newline\vspace{2mm}
$Bitors\left(S;H_{S},G_{S}\right)$,  cat{\'e}gorie des $\left(H_{S},G_{S}\right)$-bitorseurs sur le site {\'e}tale de $S$, \pageref{BitorsSHG}\newline\vspace{2mm}
$\textup{Tors}\left(k,G\right)$,  gerbe des $G$-torseurs sur le site {\'e}tale de $k$, \pageref{ChTorskG}\newline\vspace{2mm}
$\textup{Asc}\left(k,n\right)$,  gerbe des alg{\`e}bres simples centrales d'indice $n$ sur le site {\'e}tale de $k$, \pageref{Asckn}\newline\vspace{2mm}
$\textup{SB}\left(k,n-1\right)$,  gerbe des $k$-vari{\'e}t{\'e}s de Severi-Brauer de dimension $n-1$, \pageref{SBkn-1}\newline\vspace{2mm}
$\textup{Tors}\left(S,G_{S}\right)$,  gerbe des $G_{S}$-torseurs sur le site {\'e}tale de $S$, \pageref{TorsSG}\newline\vspace{2mm}
$\textup{LBun}\left(S\right)$,  gerbe des fibr{\'e}s en droites sur le site {\'e}tale de $S$, \pageref{Lbun}\newline\vspace{2mm}
$\textup{VBun}\left(n,S\right)$,  gerbe des fibr{\'e}s vectoriels de rang $n$ sur le site {\'e}tale de $S$, \pageref{Vbun}\newline\vspace{2mm}
$\textup{Az}\left(n,S\right)$,  gerbe des alg{\`e}bres d'Azumaya d'indice $n$ sur le site {\'e}tale de $S$, \pageref{AznS}\newline\vspace{2mm}
$D\left(\bar{P}\right)$,  $k$-gerbe des mod{\`e}les de $\bar{P}$, \pageref{DPbar}\newline\vspace{2mm}
$\partial\left(\lambda_{\bar{P}}\right)$,  gerbe associ{\'e}e au type $\lambda_{\bar{P}}$, \pageref{Dlambda}\newline\vspace{2mm}
$\left(\textup{FAGR}/S\right)$,  champ des faisceaux de groupes sur le site {\'e}tale de $S$, \pageref{FAGRS}\newline\vspace{2mm}
$FAGR\left(S'\right)$,  cat{\'e}gorie des faisceaux de groupes sur le site {\'e}tale de $S'$, \pageref{fagrS'}\newline\vspace{2mm}
$\left(Lien/S\right)$,  pr{\'e}champ des liens sur le site {\'e}tale de $S$, \pageref{Prechlien}\newline\vspace{2mm}
$\left(\textup{LIEN/S}\right)$,  champ des liens sur le site {\'e}tale de $S$, \pageref{ChLien}\newline\vspace{2mm}
$\textup{lien}\;G$,  lien repr{\'e}sent{\'e} par $G$, \pageref{lienG}\newline\vspace{2mm}
$H^{2}\left(S,\mathcal{L}\right)$,  ensemble des classes d'{\'e}quivalence de $S$-gerbes de lien $\mathcal{L}$, \pageref{H2L}\newline\vspace{2mm}
$R^{1}\pi_{\ast}G_{X}$,  premier foncteur d{\'e}riv{\'e} {\`a} droite du faisceau $\pi_{\ast}G_{X}$, \pageref{R1Pi}\newline\vspace{2mm}
$H^{2}\left(X,G_{X}\right)^{tr}$,  partie transgressive de $H^{2}\left(X,G_{X}\right)$, \pageref{H2trans}\newline\vspace{2mm}
$\textup{Br}^{alg}X$,  groupe de Brauer transgressif de $X$, \pageref{Brtr}\newline\vspace{2mm}
$\textup{Br}_{a}X$,  groupe de Brauer transgressif de $X$ modulo les constantes, \pageref{Bra}\newline\vspace{2mm}
$\mathcyr{B}\left(X\right)$,  noyau de $\textup{Br}_{a}X\rightarrow\displaystyle\prod_{v\in\Omega_{k}}\textup{Br}_{a}\left(X\otimes_{k}k_{v}\right)$, \pageref{BcyrX}\newline\vspace{2mm}
$\mathcyr{SH}^{1}\left(k,\textup{Pic}\;\bar{X}\right)$,  noyau de $H^{1}\left(k,\textup{Pic}\;\bar{X}\right)\rightarrow\displaystyle\prod_{v\in\Omega_{k}}H^{1}\left(k_{v},\textup{Pic}\;\bar{X}\right)$, \pageref{SHA1}\newline\vspace{2mm}
$m_{H}\left(X\right)$,  obstruction de Brauer-Manin de $X$, \pageref{mHX}\newline\vspace{2mm}
$\Omega_{k}$,  ensemble des places du corps de nombres $k$, \pageref{Omegak}\newline\vspace{2mm}
$\mathbb{A}_{k}$,  anneau des ad{\`e}les du corps de nombres $k$, \pageref{Ak}\newline\vspace{2mm}
$X\left(\mathbb{A}_{k}\right)$,  ensemble des points ad{\'e}liques de $X$, \pageref{Xak}\newline\vspace{2mm}
$\textup{inv}_{v}$,  invariant local, \pageref{invv}\newline\vspace{2mm}
$x_{v}^{\ast}b$, image inverse d'une alg{\`e}bre d'Azumaya par un $k_{v}$-point, \pageref{xvastb}\newline\vspace{2mm}
$\mathcyr{SH}^{2}\left(k,\bar{H}\right)$,  noyau de $H^{2}\left(k,\bar{H}\right)\rightarrow\displaystyle\prod_{v\in\Omega_{k}}H^{2}\left(k_{v},\bar{H}\right)$, \pageref{SH2ab}\newline\vspace{2mm}
$m_{\mathcal{H}}\left(\mathcal{G}\right)$,  obstruction de Brauer-Manin de la gerbe $\mathcal{G}$, \pageref{mHGer}\newline\vspace{2mm}
$X\left(\mathbb{A}_{k}\right)^{B}$,  points ad{\'e}liques de $X$ Brauer-Manin orthogonaux {\`a} $B$, \pageref{XAkB}\newline\vspace{2mm}
$\textup{Br}_{\lambda}X$,  partie de $\textup{Br}\;X$ associ{\'e}e au type $\lambda$, \pageref{Brlambda}\newline\vspace{2mm}
$\textup{Br}\:\mathcal{G}$,  groupe de Brauer de la gerbe $\mathcal{G}$, \pageref{Brger}\newline\vspace{2mm}
\newpage
\thispagestyle{empty}
$\ $
\newpage
\thispagestyle{empty}
$\ $
\newpage
\thispagestyle{empty}
\begin{center}
\Large{\textbf{\scshape{R{\'e}sum{\'e}}}}
\end{center}
\vspace{2mm}

Soient $k$ un corps de caract{\'e}ristique nulle et $G$ un $k$-groupe alg{\'e}brique lin{\'e}aire. Il est bien connu que si $G$ est ab{\'e}lien, les torseurs sous $G_{X}$ sur un $k$-sch{\'e}ma $\pi:X\rightarrow \textup{Spec}\;k$ fournissent une obstruction {\`a} l'existence de points $k$-rationnels sur $X$, puisque la suite spectrale de Leray:
	\[R^{p}\:\uline{\Gamma}_{\ k}\left(R^{q}\pi_{\ast}G_{X}\right)\Longrightarrow R^{p+q}\uline{\Gamma}_{\ X}\left(G_{X}\right)
\]
donne dans les bons cas (\textit{e.g.} $X$ propre) une suite exacte de groupes:
	\[0\longrightarrow H^{1}\left(k,G\right)\longrightarrow H^{1}\left(X,G_{X}\right)\stackrel{u}{\longrightarrow} H^{1}\left(\bar{X},\bar{G}_{X}\right)^{\textup{Gal}\left(\bar{k}/k\right)}\stackrel{\delta^{1}}{\longrightarrow}H^{2}\left(k,G\right)\stackrel{}{\longrightarrow}H^{2}\left(X,G_{X}\right)
\]
sur laquelle on peut directement lire l'obstruction {\`a} ce qu'un $\bar{G}_{X}$-torseur $\bar{P}\rightarrow\bar{X}$ de corps des modules $k$ soit d{\'e}fini sur $k$, \textit{i.e.} qu'il provienne par extension des scalaires {\`a} la cl{\^o}ture alg{\'e}brique $\bar{k}$ de $k$ d'un $G_{X}$-torseur $P\rightarrow X$. Le point crucial est que cette obstruction est mesur{\'e}e par une gerbe, qui est neutre lorsque $X$ poss{\`e}de un point $k$-rationnel. On essaye ici d'{\'e}tendre ce r{\'e}sultat au cas non-commutatif, et on en d{\'e}duit (sous certaines conditions) des obstructions cohomologiques non-ab{\'e}liennes {\`a} l'existence de points $k$-rationnels sur $X$, et des r{\'e}sultats sur la descente des torseurs.\vspace{3mm}\newline
\normalsize{\textbf{Mots-cl{\'e}s.} Points rationnels, (bi-)torseurs, champs, gerbes, cohomologie non-ab{\'e}lienne, obstruction de Brauer-Manin.}
\vspace{4mm}
\begin{center}
\Large{\textbf{\scshape{Abstract}}}
\end{center}
\vspace{2mm}

Let $k$ be a field of characteristic $0$ and $G$ a linear algebraic $k$-group. When $G$ is abelian, it is well known that torsors under $G_{X}$ over a $k$-scheme $\pi:X\rightarrow \textup{Spec}\;k$ provide an obstruction to the existence of $k$-rational points on $X$, since Leray spectral sequence:
	\[R^{p}\:\uline{\Gamma}_{\ k}\left(R^{q}\pi_{\ast}G_{X}\right)\Longrightarrow R^{p+q}\uline{\Gamma}_{\ X}\left(G_{X}\right)
\]
gives rise (when $X$ is \textquotedblleft{nice}\textquotedblright, \textit{e.g.} $X$ smooth and proper) to an exact sequence of groups:
	\[0\longrightarrow H^{1}\left(k,G\right)\longrightarrow H^{1}\left(X,G_{X}\right)\stackrel{u}{\longrightarrow} H^{1}\left(\bar{X},\bar{G}_{X}\right)^{\textup{Gal}\left(\bar{k}/k\right)}\stackrel{\delta^{1}}{\longrightarrow}H^{2}\left(k,G\right)\stackrel{}{\longrightarrow}H^{2}\left(X,G_{X}\right)
\]
This sequence gives an obstruction for a $\bar{G}_{X}$-torsor $\bar{P}\rightarrow\bar{X}$ with field of moduli $k$ to be defined over $k$, \textit{i.e.} to be obtained by extension of scalars to the algebraic closure $\bar{k}$ of $k$ from a $G_{X}$-torsor $P\rightarrow X$. This obstruction is measured by a gerbe, which is neutral if $X$ possesses a $k$-rational point. We try to extend this result to the non-commutative case, and in some cases, we deduce non-abelian cohomological obstruction to the existence of $k$-rational points on $X$, and results about descent of torsors. \vspace{3mm}\newline
\normalsize{\textbf{Keywords.} Rational points, (bi-)torsors, stacks, gerbes, non-abelian cohomology, Brauer-Manin obstruction.}
\end{document}